\documentclass[a4paper,12pt,openright,twoside]{book}
\pdfoutput=1
\usepackage[latin1]{inputenc}
\usepackage{amssymb,amsmath,amsthm}
\allowdisplaybreaks
\usepackage{amsfonts}
\usepackage{indentfirst}
\usepackage[top=30mm, bottom=30mm, left=35mm, right=25mm]{geometry}
\usepackage{graphicx}
\usepackage{float}
\usepackage{bbm}
\usepackage[cal=boondoxo]{mathalfa}
\usepackage{xfrac}
\usepackage{upgreek}
\usepackage{verbatim}
\usepackage{tikz}
\usetikzlibrary{matrix,arrows}
\usetikzlibrary{decorations.pathmorphing}
\usetikzlibrary{positioning}
\usepackage{tikz-cd}
\usetikzlibrary{calc}
\usepackage[labelsep=space,labelfont={bf}]{caption}

\usepackage{perpage}
\MakePerPage{footnote}
\usepackage{url}
\usepackage{hyperref}
\usepackage{fancyhdr}
\usepackage[nottoc,notlot,notlof]{tocbibind}
\usepackage{imakeidx}
\makeindex[intoc]
\newtheorem{thm}{Theorem}[section]
\newtheorem{prp}[thm]{Proposition}
\newtheorem{cor}[thm]{Corollary}

\theoremstyle{definition}
\newtheorem{dfn}[thm]{Definition}
\newtheorem{rem}[thm]{Remark}
\newtheorem{ntn}[thm]{Notation}
\newenvironment{Abs}{\begin{center}
		\huge\textbf{Abstract} \\[1.5cm] \begin{minipage}{.75\textwidth}\normalsize}{\end{minipage}	\end{center}}
\renewcommand{\qedsymbol}{$\blacksquare$}
\newenvironment{mypr}{\noindent \textbf{Proof.}}{\hfill\qedsymbol}
\newcommand{\mcg}[1]{\mathcal{#1}}
\newcommand{\ag}[2]{\left\langle #1,#2\right\rangle}
\newcommand{\agg}[3]{\left\langle #1,#2,#3\right\rangle}
\newcommand{\xc}{\xrightarrow{~\cong~}}
\newcommand{\cs}{CSMC}
\newcommand{\cat}[1]{\textup{\textbf{#1}}}
\newcommand{\prm}[1]{#1^\prime}
\newcommand{\pprm}[1]{#1^{\prime\prime}}
\newcommand{\ppprm}[1]{#1^{\prime\prime\prime}}
\newcommand{\tu}[1]{\textup{#1}}
\newcommand{\op}[1]{#1^{\tu{op}}}
\newcommand{\eqd}{\stackrel{\text{\fontsize{5}{1}\selectfont def}}{=}}
\newcommand{\nat}{\stackrel{\bullet}{\longrightarrow}}
\newcommand{\msf}[1]{\mathsf{#1}}
\DeclareMathOperator*{\Lm}{Lim}
\DeclareMathOperator*{\Clm}{Colim}
\DeclareMathOperator*{\dju}{\stackrel{\bullet}{\bigcup}}
\newcommand{\chu}{\cat{Chu}}
\newcommand{\sfr}[2]{\sfrac{#1}{#2}}
\newcommand{\bs}{\backslash}
\newcommand{\dlong}{\raisebox{-1pt}{$\def\arraystretch{0.5}\begin{array}{c} \longrightarrow\\ \longrightarrow \end{array}$}} 
\newcommand{\bb}[1]{\mathbb{#1}}
\newcommand{\bm}[1]{\mathbbm{#1}}
\newcommand{\eb}[1]{\mathbf{#1}}
\newcommand{\dlg}[3]{{#1}^{#2}_{#3}}
\newcommand{\fk}[1]{\mathfrak{#1}}
\newcommand{\Gr}{\mathrm{Gr}}
\newcommand{\ob}{\mathrm{Obj}}
\newcommand{\mor}{\mathrm{Mor}}
\newcommand{\dom}{\mathrm{Dom}}
\newcommand{\cod}{\mathrm{Cod}}
\newcommand{\id}{\mathrm{id}}
\newcommand{\rvr}[1]{\reflectbox{\tu{R}}_{#1}}
\newcommand{\pb}[3]{#1\times_{#2}#3}
\newcommand{\uo}[1]{#1_{\mathrm{obj}}}
\newcommand{\um}[1]{#1_{\mathrm{mor}}}
\newcommand{\bt}{\bowtie}
\newcommand{\T}{\mathrm{T}}
\newcommand{\uda}{\updownarrow}
\newcommand{\lra}{\leftrightarrow}
\newcommand{\dlc}{\bb{DLC}}
\newcommand{\xR}[1]{\stackrel{~\bm{#1}~}{\Longrightarrow}}
\newcommand{\wt}[1]{\widetilde{#1}}
\newcommand{\ol}[1]{\overline{#1}}
\newcommand{\flu}{\cat{Fluid}_\chu}
\newcommand{\iw}[1]{\textbf{#1}\index{#1}}
\newcommand{\tpar}[2]{\vspace*{#1}{\par \large #2}}

\begin{document}
	\frontmatter

	\begin{titlepage}
		\begin{center}
			\vspace*{0.3cm}{\Huge \textbf{The Double Category of Paired Dialgebras on the Chu Category}}
			\tpar{0.3cm}{(English Translation)}
			\rule{0pt}{1cm}
			\tpar{0.3cm}{~\\~\\}
			\tpar{0.3cm}{Aydin Manzouri}
			\tpar{0.3cm}{aydin.manzouri@gmail.com}
			\tpar{0.3cm}{May 2017}
			\tpar{3cm}{A Thesis Submitted in Partial Fulfillment of}
			\tpar{0.3cm}{the Requirements for the Degree of}
			\tpar{0.3cm}{Master of Science in Pure Mathematics}
			\tpar{1cm}{Faculty of Mathematical Sciences}
			\tpar{0.3cm}{Shahid Beheshti University}
			\tpar{0.3cm}{Tehran, Iran}
			\tpar{2cm}{\textbf{Supervisor:}}
			\tpar{0.3cm}{Professor Mojgan Mahmoudi}
			\tpar{1cm}{\textbf{Advisors:}}
			\tpar{0.3cm}{Professor Mohammad-Mehdi Ebrahimi}
			\tpar{0.3cm}{Professor Amir Daneshgar, Sharif University of Technology}
			
		\end{center}
	\end{titlepage}
	\rule{0pt}{5cm}
	\newpage
	\begin{Abs}
		\par \indent The scientific and practical needs of the twenty-first century lead humankind to convergence of the specialized and diverse branches of science and technology. This convergence reveals the need for new mathematical theories capable of providing common languages and frameworks to be utilized by professionals from different fileds in solving interdisciplinary and challenging problems.
		\par \indent The present thesis is done in the same direction. Here, we develop a new formalism with the central idea of ``unification of various mathematical branches''. For this purpose, we utilize three major tools from today's mathematics, each of which possessing a unifying nature itself: category theory and especially the theory of ``double cateogries'', the theory of ``universal dialgebra'', and the ``Chu construction''. With the aid of these tools, we define and study a double category that subsumes a significant portion of the formalisms usual within the body of mathematics and theoretical computer science. We show that this double category possesses the properties of ``horizontal self-duality'' and ``vertical self-duality''. Also, we perform a primary investigation about existence of binary horizontal products and coproducts in this category. Finally, we give some suggestions for future work.
	\end{Abs}

	\setcounter{tocdepth}{3}
	\tableofcontents

	\chapter{Preface}
	Convergence of the various scientific and technological fields, towards solving challenging multi-dimensional problems, is a need that is felt in the twenty-first century far more than anytime in the past. Today, we are the heirs of an enormous amount of scientific and technological achievements, and those continue to accumulate at accelerating rates. Within this roaring ocean of information, the degrees of specialization of the study areas have reached such levels that, frequently, effective transference of findings from one professional to another becomes a serious challenge itself.
	\par In addition, we are the heirs of huge amounts of human and environmental problems. Today, few people in the scientific community ever doubt the role of man in global warming. Non-degradable plastic and rubber waste, various toxic chemicals, petrochemicals, and heavy metals that humans have been producing and releasing into nature since the beginning of the twentieth century have polluted lands and waters; those substances interfere with ecological cycles, they cause death of organisms and cause undesirable changes in food chains, and the result of all these returns to us humans in the form of environmental disasters and health threats. On the other hand, the emergence of drug-resistant pathogens due to overconsumption of antibiotics has made the medical sciences face new challenges in fighting infectious diseases. Yet another set of alarming issues is the crises of water and food shortages and hard-to-treat epidemics and pandemics throughout the world, especially in the impoverished areas. Along with all these problems, one should add the potential dangers of the new technologies. Today, nanotechnology, biotechnology, and information technology are already forming new industrial revolutions throughout the world while many of their side-effects are unknown to us.
	\par For solving difficult multi-dimensional problems in such an extremely complex era, effective cooperation of professionals from different disciplines is a clear necessity. But any effective cooperation requires using a ``potent common language'' by all involved; therefore, mathematics manifests its value as the common language of sciences and technologies once again. Although various mathematical theories developed up to now have offered brilliant services to natural sciences and engineering fields, currently we need a new mathematics which is even more capable of unifying the infrastructures as well as the bodies of different branches of science and technology.
	\par The crucial point to be considered is that \textit{it is impossible for mathematics to succeed in such unification when it is internally fragmented (as we witness today); rather, the process of convergence and unification has to start within mathematics itself. By ``unification of mathematics'' we do \textbf{not} mean a theory principally for the foundations of mathematics (such as axiomatic set theory), but one which is able to encompass the essence of different theories within the \textbf{body} of mathematics and which provides the working mathematician with the necessary toolbox for effective cross-disciplinary communication of \textbf{ideas}.} 
	\par The present thesis is done in the same direction. For this purpose, we utilize three major tools from today's mathematics, each of which possessing a unifying nature itself: category theory and especially the theory of ``double cateogries'', the theory of ``universal dialgebra'', and the ``Chu construction''.
	\par Eilenberg and Mac Lane officially founded the theory of categories in their historic paper in 1945 \cite{Awod}. There were various generalizations of the theory later, of which we point to Ehresmann's theory of \textit{double categories} (see \cite{Fiore-Thoma} and the references within). The notion of double category can be seen as a ``two-dimensional generalization'' of the notion of category. The formalism developed in the current thesis makes use of the language of double categories.
	\par From a categorical viewpoint, algebraic theories can have their corresponding ``dual theories''. These duals are referred to as \textit{coalgebras}. Whereas the theory of universal algebra \cite{Grat, Burr} recognizes and studies the common patterns among specific algebraic theories such as the theories of groups, rings, linear algebra, etc., the theory of \textit{universal coalgebra} does the same with specific coalgebraic theories such as automata, transition systems, Petri nets, and event systems \cite{Rutt,Jaco}. 
	\par Again, categorically one can think of \textit{dialgebras} \cite{Poll} as common generalizations of algebras and coalgebras. Also, the theory of \textit{universal dialgebra} \cite{Vout} has been developed as a common generalization of universal algebra and universal coalgebra. Universal dialgebra recognizes and studies the common patterns among specific dialgebraic theories.
	\par From a different perspective, tracing the property of ``self-duality'' in categories (see Definition \ref{selfdual}) leads us to a family of categories called the \textit{Chu construction} \cite{Barr-hist}.  Existence of the self-duality property is essential for the purposes of the present work (see the discussion given in Section \ref{funda-motivs}). It is well-known that $ \cat{Set} $, the category of sets and functions, lacks this property; therefore, in order to achieve a unified formalism for solving interdisciplinary problems, one has to develop the theory of universal dialgebra based on some other category which, along with other desirable properties, possesses self-duality. Consequently, the Chu construction enters the scene as an appropriate substitute for $ \cat{Set} $.
	\par In this thesis, after introducing the above tools, we define and study a formalism which subsumes a significant portion of the formalisms usual within the body of mathematics and theoretical computer science. The general structure of the thesis is as follows. In Chapter 1 the preliminaries are given. Next, the Chu construction is introduced in Chapter 2. The contents of that chapter are mainly from \cite{Barr-hist,Bif,Pratt-chu-sp,Pratt-Gamut}. Next, universal dialgebra is given in Chapter 3, with the materials mainly inspired by \cite{Poll,Vout}.
	\par Then, the main formalism of the thesis is given in Chapter 4. In that chapter, after introducing double categories, we introduce the $ \dlc $ construction (see Section \ref{the-main-forma}). We show that this double category possesses the property of ``Klein-invariance'' (see Definition \ref{Klein-inv}), which is the conjunction of the two properties of horizontal and vertical self-duality. Also, we do a primary investigation concerning existence of binary horizontal products and coproducts in $ \dlc $. Finally, a few suggestions for future work will be given in Chapter 5.
	

	\mainmatter

	\chapter{Introduction}
	This chapter discusses the general ideas and motivations behind the present work and provides the preliminaries for the next chapters.\\
	\par \noindent \textbf{Thesis organization.} In Section \ref{funda-motivs} we speak of the fundamental motivations; in Section \ref{quick-rev} we quickly review the essentials of basic category theory; and in Section \ref{internal-cats-section}, internal categories are introduced.
	\par Chapter 2 introduces the Chu construction and discusses a number of its properties. Chapter 3 studies universal dialgebra on a given base category $ \mcg{C} $ in general and on the Chu construction in particular. Chapter 4 begins with introducing the formalism of double categories and proceeds towards the main formalism of the present work. Finally, Chapter 5 gives the conclusions and directions for future work.
	
	\section{The fundamental motivations}
	\label{funda-motivs}
	Today, in the twenty-first century, we live at a time of ever-accelerating changes, of which scientific and technological advancements constitute a significant portion. Virtually every aspect of human life is constantly flooded with scientific discoveries and technological innovations. This situation was forecast with relatively high accuracy in \cite{Roco} in 2002:
	\vspace{2mm}
	\par \textit{``We stand at the threshold of a New Renaissance in science and technology, based on a comprehensive understanding of the structure and behavior of matter from the nanoscale up to the most complex system yet discovered, the human brain. A coherent science and engineering approach based on the unity of nature and its holistic investigation will lead to technological convergence and a more efficient societal structure...''}
	\vspace{2mm}
	\par The article then explains this ``technological convergence'' as follows [ibid]:
	\vspace{2mm}
	\par \textit{``The phrase `convergent technologies' refers to the synergistic combination of four major `NBIC' (Nano-Bio-Info-Cogno) provinces of science and technology, each of which is currently progressing at a rapid rate: (a) nanoscience and nanotechnology; (b) biotechnology and biomedicine, including genetic engineering; (c) information technology, including advanced computing and communications; (d) cognitive science, including cognitive neuroscience.} 
	\par \textit{Accelerated scientific and social progress can be achieved by combining research methods and results across these provinces in duos, trios, and the full quartet.''}
	\vspace{2mm}
	\par On the other hand, there is a crucial fact concerning the nature of scientific research itself [ibid]:
	\vspace{2mm}
	\par \textit{``The sciences have reached a watershed at which they must combine if they are to continue to advance. The New Renaissance must be based on a holistic view of science and technology that envisions new technical possibilities and focuses on people.''}
	\vspace{2mm}
	\par Finally, the article lays emphasis on the key roles of mathematics, computer science, and system approach in the New Renaissance [ibid]:
	\vspace{2mm}
	\par \textit{``Developments in system approach, mathematics, and computation in conjunction with NBIC allow us for the first time to understand the natural world, and social events and humanity as closely coupled complex, hierarchical systems. Applied both to particular research problems and to the over-all organization of the research enterprise, this complex system approach provides holistic awareness of opportunities for integration, in order to obtain maximum synergism along main directions of progress.''}
	\vspace{2mm}	
	\par From the above, it follows that development of mathematics in the twenty-first century shall be heavily influenced by the requirements of intensely-converging sciences and technologies. This shall result in more unified mathematical theories, capable of solving more diverse, more complex, more hierarchical, and more challenging problems arising in the NBIC hybrid. Needless to say, this historical process also requires the pure and applied parts of mathematics to unify and function as a whole.
	\par Now, the present work has been inspired by ideas and motivations inclusive of what described above. The author has been trying to develop a new and unified theoretical framework in which those aspects of pure mathematics that are of (potential or actual) practical importance may be embodied. The result is a theory that interweaves three threads of existing mathematical theories, each of which having a unifying nature itself: category theory, universal dialgebra, and the Chu construction.\\
	
	\par \noindent \textbf{Category theory.} The twentieth century witnessed a number of triumphs in solidification of mathematical thinking, one of which was the formulation of various axiomatic set theories as candidate foundations of mathematics. Moreover, Eilenberg and Mac Lane officially formulated category theory in their historic paper (``General theory of natural equivalences'') in 1945 \cite{Awod}. Contrary to set theory which places emphasis on sets and membership, category theory prioritizes \textit{arrows}, which are \textit{relationships between} different mathematical objects. This way, category theory sheds light on even the farthest corners of mathematical realm and unifies structures so distant from each other that would seem totally unrelated otherwise. It is also worth mentioning that in addition to the set-theoretic foundations, there are suggestions for category-theoretic foundations of mathematics (e.g. see \cite{nLab-found} and the references therein). Therefore, as foundational systems for mathematics, category theory possesses much more unifying power than set theory.
	\par Now, ordinary categories themselves can be seen as kinds of ``one-dimensional structures'', in which objects play the role of ``points'' and arrows that of ``one-dimensional entities between points''. These can be generalized to higher dimensions in various ways. The present work utilizes a definite kind of generalization known as ``double categories'' \cite{Fiore-Thoma,GranPar-Lims}. These will be introduced in Chapter 4.\\
	
	\par \noindent \textbf{Universal dialgebra.} From the viewpoint of mathematical logic, groups, rings, fields, lattices, and other algebraic structures can be viewed as \textit{models} of particular theories. Whereas abstract algebra studies those models, \textbf{universal algebra} deals with algebraic or equational \textit{theories} in general \cite{Grat,nLab-uni-alg}. Gr\"{a}tzer gives the following notes on the latter subject in his book \cite{Grat}:
	\vspace{2mm}
	\par \textit{``In A. N. Whitehead's book on Universal Algebra, published in 1898, the term universal algebra had very much the same meaning that it has today....}
	\par \textit{Thus universal algebra is the study of finitary operations on a set, and the purpose of research is to find and develop the properties which such diverse algebras as rings, fields, Boolean algebras, lattices, and groups may have in common.''}
	\vspace{2mm}
	\par Now, with the aid of category theory, the above picture can be expanded in a number of ways. First of all, the collection of finitary operations on a given set may be replaced with the notion of ``$ F $-algebra'', where $ F $ is an endofunctor (see \ref{endofunctor-def}) on the category of sets . The advantage of this approach is that instead of an indexed set of fundamental operations, we have only one mapping \cite{Dene}. Next, $ F $-algebras can be categorically ``dualized'' (see \ref{dual-statement}) to yield $ F $-coalgebras. $ F $-coalgebras provide the (categorical) basis for the theory of \textbf{universal coalgebra}. Quoting Rutten \cite{Rutt}:
	\vspace{2mm}
	\par \textit{``These observations, then, have led to the development in the present paper of a general theory of coalgebras called `universal coalgebra', along the lines of universal algebra. Universal algebra ... deals with the features common to the many well-known examples of algebras such as groups, rings, etc. The central concepts are $ \Sigma $-algebra, homomorphism of $ \Sigma $-algebras, and congruence relation. The corresponding notions ... on the coalgebra side are: coalgebra, homomorphism of coalgebras, and bisimulation equivalence. These notions constitute the basic ingredients of our theory of universal coalgebra. (More generally, the notion of substitutive relation corresponds to that of bisimulation relation; hence congruences, which are substitutive equivalence relations, correspond to bisimulation equivalences.) Adding to this the above-mentioned observation that various dynamical systems (automata, transition systems, and many others as we shall see) can be represented as coalgebras, makes us speak of universal coalgebra as `a theory of systems'. We shall go even as far as, at least for the context of the present paper, to consider coalgebra and system as synonyms.''}
	\vspace{2mm}
	\par Also, Jacobs \cite{Jaco} writes:
	\vspace{2mm}
	\par \textit{``Mathematics is about the formal structures underlying counting, measuring, transforming etc.... In more recent decades also `dynamical' features have become a subject of research. The emergence of computers has contributed to this development. Typically, dynamics involves a `state of affairs', which can possibly be observed and modified....}
	\par \textit{Both mathematicians and computer scientists have introduced various formal structures to capture the essence of state-based dynamics, such as automata (in various forms), transition systems, Petri nets and event systems. The area of coalgebra has emerged within theoretical computer science with a unifying claim. It aims to be the mathematics of computational dynamics. It combines notions and ideas from the mathematical theory of dynamical systems and from the theory of state-based computation. The area of coalgebra is still in its infancy but promises a perspective on uniting, say, the theory of differential equations with automata and process theory and with biological and quantum computing, by providing an appropriate semantical basis with associated logic. The theory of coalgebras may be seen as one of the original contributions stemming from the area of theoretical computer science.''}
	\vspace{2mm}
	\par Thirdly, the assumption that universal algebra and universal coalgebra merely deal with set-based operations and dynamics may be dropped. In other words, one can study algebras and coalgebras induced by arbitrary endofunctors $ F:\mcg{C}\longrightarrow\mcg{C} $, where $ \mcg{C} $ may be any category, not necessarily the category of sets.
	\par Finally, universal algebra and universal coalgebra themselves can be unified. So far, the situation is roughly as follows: universal algebra unifies specific algebraic structures, universal coalgebra does the same with specific state-based dynamical systems, and these two are categorical duals to each other. Now, unification of the two theories results in the more recent theory of \textbf{universal dialgebra}. The notion of \textit{dialgebra} was investigated in \cite{Poll} as a common generalization of algebras and coalgebras. Dialgebras subsume both features of algebraic structuring and coalgebraic dynamics; also, there are many interesting examples of dialgebras that are neither algebra nor coalgebra.
	\par However, it was not until 2010 that the first paper on ``universal dialgebra'' as a new theory was published. In that year, Voutsadakis \cite{Vout} gave a systematic treatment of the subject for the first time, and showed how the fundamental results in both universal algebra and universal coalgebra can be viewed as special cases of those in universal dialgebra. There, the notion of $ \ag{F}{G} $-dialgebra is introduced as the categorical formalism of dialgebras, where $ F,G $ are arbitrary endofunctors on the category of sets. Although Voutsadakis's work deals exclusively with set-based dialgebras, as he himself states: \textit{``... many of the results (in fact most of them, if adequately translated) will be easily seen to hold in arbitrary categories''}. This is what we will do to some extent in Chapter 3.
	\vspace{2mm}
	\par An interesting fact is that ``being a dialgebra'' is a \textit{self-dual} statement (see \ref{self-dual-property-def}); that is, the dual of a dialgebra is again a dialgebra. This can be contrasted with the above-mentioned fact that algebras and coalgebras are duals to each other. In other words, whereas the process of dualization interchanges algebraic and coalgebraic \textit{concepts}, the same process does not affect dialgebras \textit{conceptually}.
	\par However, this does NOT mean that dialgebras are perfectly immune to dualization. Indeed, there is a subtle fact concerning the dualization process. For example, when $ \mcg{C}=\cat{Set} $ (the category of sets and functions), quoting Gumm \cite{Gumm-UnCoLo}:
	\par \textit{``Universal coalgebra is dual to universal algebra over the dual of the category of sets. As the category $ \cat{Set} $ is not self-dual, universal algebra cannot simply be translated to deliver a corresponding theory of coalgebras. It can, however, serve as a formidable source of inspiration.''}
	\vspace{2mm}
	\par The same issue affects dialgebras, too. Again, this will be discussed in Chapter 3.
	\par A remedy for the above problem is to replace $ \cat{Set} $ with a \textit{self-dual category} (see \ref{selfdual}). Consequently, in search for a self-dual substitute for $ \cat{Set} $, the author has found the third essential ingredient of the current thesis.\\
	
	\par \noindent \textbf{The Chu construction.} Thu Chu construction is exactly what the author has found appropriate for the purposes of the present work. Quoting Barr \cite{Barr-hist}:
	\vspace{2mm}
	\par \textit{``In 1975, I began a sabbatical leave at the ETH in Z\"{u}rich, with the idea of studying duality in categories in some depth. By this, I meant not such things as the duality between Boolean algebras and Stone spaces, nor between compact and discrete abelian groups, but rather self-dual categories such as complete semi-lattices, finite abelian groups, and locally compact abelian groups. Moreover, I was interested in the possibilities of having a category that was not only self-dual but one that had an internal hom and for which the duality was implemented as the internal hom into a `dualizing object'....}
	\par \textit{The desired properties were what I subsequently called $ * $-autonomy....}
	\par \textit{By the end of the year, I had in fact produced a moderate number of examples of $ * $-autonomous categories. One of them was a full subcategory of topological abelian groups that included all the locally compact abelian (LCA) groups in such a way that the duality restricted to them was the well-known duality of LCA groups....}
	\par \textit{Another example was a full subcategory of the category of locally convex topological vector spaces....}
	\par \textit{Thus I ended up with a category whose objects were pairs $ \msf{E}=(E,\prm{E}) $ of vector spaces equipped with a pairing $E\otimes \prm{E}\longrightarrow \bb{C}$} [with $\bb{C}$ being the set of complex numbers]. \textit{A map from $ \msf{E} $ to $ \msf{F} $ is a pair of linear maps $ (f,\prm{f}) $ in which $ f:E\longrightarrow F $ and $ \prm{f}:\prm{F}\longrightarrow\prm{E} $ (note the direction reversal) such that $ \langle fv,w\rangle =\langle v,\prm{f}w\rangle$ whenever $ v\in E $ and $ w\in\prm{F} $. There is no topology assumed and no continuity on the linear maps....''}
	\vspace{2mm}
	\par He then continues to describe how he had formulated a new collection of $ * $-autonomous categories; and finally [ibid]:
	\vspace{2mm}
	\par \textit{``It seemed clear that this gave a $*$-autonomous category, but there were a number of unpleasant details to be verified. Since my student, Po-Hsiang Chu needed a master's project, so I asked him to verify them, which he did.... I now had expanded from six to infinity the repertory of known $ * $-autonomous categories.''}
	\vspace{2mm}
	\par Eventually, the new collection of $ * $-autonomous categories gets named after P. Chu's work in his thesis in 1979. The details of the construction will be given in Chapter 2 (see Section \ref{Chu-V}). The Chu construction, especially the \textit{set-based} Chu construction, has many interesting features (see Chapter 2), among which is self-duality. Those features altogether make the (set-based) Chu construction an excellent replacement for $ Set $, the category of sets and functions, for serving as a base category for a theory of universal dialgebra.\\
	
	\par \noindent \textbf{The present work.} Although developing a theory of universal dialgebra on a self-dual basis has been one of the initial motivations for the author, the present work goes beyond it. Indeed, the formalism that we will be developing in Chapter 4 subsumes and interrelates \textit{all possible theories of universal dialgebra on the set-based Chu construction}. More precisely, assuming a set-based Chu category $ \chu $, we will develop a double categorical structure in which every horizontal ``cross-section'' corresponds to a ``counter-current pairing'' of \textit{two} theories of universal dialgebra on $ \chu $, while every vertical ``cross-section'' corresponds to interrelations between different dialgebraic theories for a fixed \textit{pair} of objects of $ \chu $. Moreover, we will show that the new double category is \textit{Klein-invariant} (see Definition \ref{Klein-inv}). This way, a two-dimensional unified framework emerges that subsumes a significant portion of what pure and applied mathematicians and computer scientists deal with in everyday work. The author hopes that the results of this work will prove to be beneficial for future scientific and technological advancements.
	
	\section{A quick review of basic category theory}
	\label{quick-rev}	
	In this section, we take a brief look at the basic definitions and constructions of category theory. We assume the reader's familiarity with some kind of axiomatic set theory (e.g., refer to \cite{Monk}). For the purpose of the present work, we need the concepts of ``class'' and ``conglomerate'' besides that of ``set'' and ``function''. The material given in this section is mainly borrowed from \cite{Adam,Borc1,CWM,Lart}.
	\rule{0pt}{0.1cm}
	
	\subsubsection*{Classes}
	Intuitively speaking, classes\index{class} consist of ``large collections of sets''. In particular, the following are required.
	
	\begin{enumerate}
		\item The members of each class are sets.
		
		\item For every ``property'' $P$, one can form the class of all sets with prroperty $P$.
		\par Hence there is the largest class $\mcg{U}$ of all sets, called the \textbf{universe}. Classes are precisely the subcollections of $\mcg{U}$. Thus, given classes $A$ and $B$, one can form such classes as \(A\cup B,A\cap B,\) and $A\times B$. Because of this, there is no problem in defining functions between classes, equivalnce relations, etc.
		
		\item If $X_1,X_2,...,X_n$ are classes, then so is the $n$-tuple $(X_1,X_2,...,X_n)$.
		
		\item Every set is a class.
		
		\item There is no surjection from a set to a proper class.
		\par This means that sets have ``fewer'' elements than proper classs.
	\end{enumerate}
	\rule{0pt}{0.1cm}
	
	\subsubsection*{Conglomerates}
	The concept of ``conglomerate''\index{conglomerate} has been created to deal with ``collections of classes''. The following are required.
	
	\begin{enumerate}
		\item Every class is a conglomerate.
		
		\item For every property $P$, one can form the conglomerate of all classes with property $P$.
		
		\item Conglomerates are closed under analogues of the usual set-theoretic constructions (pairing, union, intersection, product, etc.)
		\par Thus, we can form the conglomerate $\prm{\mcg{U}}$ of all classes, as well as functions between conglomerates.
	\end{enumerate}
	
	\begin{rem}
		\label{proper-classes-conglos}
		Classes that are not sets are called \textbf{proper classes}\index{proper class}. Similarly, conglomerates that are not classes are called \textbf{proper conglomerates}\index{proper conglomerate}. For example, $\mcg{U}$ is a proper class, while $\prm{\mcg{U}}$ is a proper conglomerate.
	\end{rem}
	\rule{0pt}{0.1cm}
	
	\subsubsection*{Categories}
	There are various ways to define a category\index{category}. We state two of them since each will be useful in some parts of the current work.
	\begin{dfn}
		\label{def-cateogry-first}
		\textbf{(The first definition of categories.)} A \textbf{category} is a quadruple $\mcg{C}=\left\langle \mcg{O},\hom,\circ,\id \right\rangle $ consisting of:
		\begin{itemize}
			\item a \textit{class} $\mcg{O}$, whose elements are called \textbf{objects}\index{object} of the category or $\mcg{C}$-\textbf{objects}; this class is also denoted by $\uo{\mcg{C}}$; we prefer this latter notation, and we will be using it in the sequel;
			
			\item for every pair $X,Y$ of objects, a \textit{set} $\hom(X,Y)$, whose elements are called\\ \textbf{morphisms}\index{morphism} or \textbf{arrows}\index{arrow} or $\mcg{C}$-\textbf{morphisms} or $\mcg{C}$-\textbf{arrows} from $X$ to $Y$; this set is also denoted by $\hom_\mcg{C}(X,Y)$ or by $\mcg{C}(X,Y)$; we prefer the last notation, and we will be using it in the sequel;
			
			\item for each triple $X,Y,Z$ of objects, a function \[\circ:\mcg{C}(X,Y)\times \mcg{C}(Y,Z)\longrightarrow \mcg{C}(X,Z),\]
			called the \iw{composition}, which assigns, to each pair of morphisms $\langle f,g\rangle$, the\\ composite morphism $g\circ f$ or just $gf$; each of these two notations has its own advantages and therefore, we will be using both in the sequel;
			
			\item for each object $X$, an element $\id_X$ or $1_X$ of $\mcg{C}(X,X)$, called the \iw{identity morphism} or \iw{identity arrow} on $X$; we prefer the notation ``$1_X$'', and we will use it most of the time.
		\end{itemize}
		These data are subject to the following axioms.
		\begin{enumerate}
			\item \iw{Associativity Axiom.} Given morphisms \(f\in \mcg{C}(W,X),~ g\in \mcg{C}(X,Y),~\)\\ $h \in \mcg{C}(Y,Z)$, the following equality holds: \[h\circ(g\circ f)=(h\circ g)\circ f.\] Hence we can safely discard parentheses and write $hgf$.
			
			\item \iw{Identity Axiom.} Given a morphism $f\in \mcg{C}(X,Y)$, the following equalities hold: \[1_Y\circ f=f=f\circ 1_X.\]
			
			\item \iw{Disjointness Axiom.} The sets $\mcg{C}(X,Y)$ are pairwise disjoint.
		\end{enumerate}
		The union of all sets $\mcg{C}(X,Y)$ is denoted by $\um{\mcg{C}}$ and is called the \textbf{class of morphisms} or the \textbf{class of arrows} of $\mcg{C}$, or equivalently, the \textbf{class of} $\mcg{C}$-\textbf{morphisms} or the \textbf{class of} $\mcg{C}$-\textbf{arrows}.
	\end{dfn}
	
	\rule{0pt}{0.1cm}
	\begin{rem}
		\label{def-cateogry-rems}
		\begin{enumerate}
			\item A morphism $f\in\mcg{C}(A,B)$ is represented by the notation\\ $f:A\longrightarrow B$. The object $A$ is called the \iw{domain} or the \iw{source} of $f$, and the object $B$ is called the \iw{codomain} or the \iw{target} of $f$. We may write \[\dom(f)=A,~\cod(f)=B.\] Two arrows $f,g$ with $\dom(g)=\cod(f)$ are called \iw{composable}.
			
			\item We will often write $X\in\mcg{C}$ instead of $X\in\ob(\mcg{C})$, and \(f\in\mcg{C}\) instead of $f\in\mor(\mcg{C})$.
			
			\item Also, in the situation of a diagram like that of \ref{comm-diag}, we say that the diagram is \textbf{commutative} (or the diagram \textbf{commutes})\index{commutative diagram} if, taking any two vertices, the composite of the arrows along any path from the first vertex to the second is equal to the composite along any other path from the first to the second; e.g., when $g\circ f=k\circ h$ in this diagram.
			\begin{figure}[h]
				\begin{center}
					\begin{tikzpicture}[commutative diagrams/every diagram]
					
					\node (N1) at (1.5,1.5cm) {$B$};
					
					\node (N2) at (-1.5,1.5cm) {$A$};
					
					\node (N3) at (-1.5,-1.5cm) {$C$};
					
					\node (N4) at (1.5,-1.5cm) {$D$};	
					
					\path[commutative diagrams/.cd, every arrow, every label]
					(N2) edge node{$f$} (N1)
					(N2) edge node{$h$} (N3)
					(N3) edge node{$k$} (N4)
					(N1) edge node{$g$} (N4);
					\end{tikzpicture}
				\end{center}
				\caption{} \label{comm-diag}
			\end{figure}
			An analogous terminology holds for diagrams of arbitrary shape.
			
			\item $1_X$ is the only arrow from $X$ to $X$ that plays the role of an identity for the composition. Indeed, if $\prm{1}_X\in\mcg{C}(X,X)$ is another such morphism, then
			\[1_X=1_X\circ \prm{1}_X=\prm{1}_X.\]
			
			\item Axiom (3) guarantees that each $\mcg{C}$-morphism has a \textit{unique} domain and a \textit{unique} codomain.
		\end{enumerate}
	\end{rem}
	
	Now we state the second definition:
	
	\begin{dfn}
		\label{def-cateogry-second}
		\textbf{(The second definition of categories.)} Alternatively, a category is a sextuple $\mcg{C}=\langle \uo{\mcg{C}},\um{\mcg{C}},s,t,\id,\circ\rangle$, where:
		\begin{itemize}
			\item $\uo{\mcg{C}},\um{\mcg{C}}$ are the classes of objects and morphisms of $\mcg{C}$, respectively;
			
			\item $s,t:\um{\mcg{C}}\dlong\uo{\mcg{C}}$ are two functions assigning, to every morphism $f$, its source and target, respectively;
			
			\item \(\circ:\um{\mcg{C}}\times\um{\mcg{C}}\longrightarrow\um{\mcg{C}}\) is a \textit{partial function} called the composition, which assigns, to any pair of morphisms $f,g$ with $t(f)=s(g)$, their composite morphism $g\circ f=gf$;
			
			\item $\id:\uo{\mcg{C}}\longrightarrow\um{\mcg{C}}$ is a function which assigns to each object $X$ a morphism $\id_X$ or $1_X$.
		\end{itemize}
		These data are subject to the following axioms:
		\begin{enumerate}
			\item Source and target are respected by composition:
			\[s(g\circ f)=s(f)~~~~\text{and}~~~~t(g\circ f)=t(g).\]
			
			\item Source and target are respected by identities:
			\[s(1_X)=X=t(1_X).\]
			
			\item Composition is associative; that is, whenever \(t(f)=s(g)\) and $t(g)=s(h)$ we have
			\[(h\circ g)\circ f=h\circ(g\circ f).\]
			
			\item Composition satisfies the \textbf{left and right unit laws}; if $s(f)=X$ and $t(f)=Y$, then
			\[1_Y\circ f=f=f\circ 1_X.\]
		\end{enumerate}
	\end{dfn}
	
	\begin{rem}
		\label{usefulness-of-each-dfn-cats}
		Each of Definitions \ref{def-cateogry-first} and \ref{def-cateogry-second} has its own advantages when working in different theories. For example:
		\begin{itemize}
			\item Definition \ref{def-cateogry-first} works best when introducing the notion of \textit{``monoidal category''}, a topic we will be introducing in Section \ref{monoidal-cats-section}. Monoidal categories provide the foundations for Chapter 2. It is also interesting to note that the (ordinary) study of \textit{enriched categories} utilizes the language of monoidal categories. There, too, Defintion \ref{def-cateogry-first} generalizes quite nicely to the notion of ``enriched category'' (e.g., see Chapter 6 of \cite{Borc2}).
			
			\item In contrast, Definition \ref{def-cateogry-second} generalizes elegantly and straightforwardly to the notion of \textit{internal category}, to be introduced in Section \ref{internal-cats-section}. Internal categories are one of the mainstream approaches to the development of the theory of \textit{double categories}. Double categories will be introduced in Chapter 4 and they provide us with the necessary language to be used in the main formalism of this thesis.
		\end{itemize}
	\end{rem}
	
	\begin{ntn}
		\label{notation-for-cats}
		We will be using capital letters $A,B,C,X,Y,Z,...$ for objects and small letters $f,g,h,k,...$ for morphisms in a category. We will be using calligraphic letters $\mcg{A,B,C,D,...}$ to denote categories in general, and boldface upright letters for named categories (such as \cat{Set}, \cat{Grph}, etc.). More specific notations for different situations will come in the sequel.
	\end{ntn}
	\rule{0pt}{0.5cm}
	Next, let us take a look at ``size issues'':
	
	\begin{dfn}
		\label{large-vs-small-cat}
		A category $\mcg{C}$ is also called a \iw{large category}. On the other hand, $\mcg{C}$ is called a \iw{small category} whenever $\uo{\mcg{C}}$ is a set.
	\end{dfn}
	
	\begin{rem}
		\label{ob-small-c-small}
		\begin{enumerate}
			\item When $\uo{\mcg{C}}$ is a set, then $\um{\mcg{C}}$ must be a set, so that the small category $\mcg{C}$ must \textit{also} be a set.
			
			\item Obviously, every small category is at the same time a large category, but \textit{not} the converse. Large categories that are not small are called \textbf{properly large categories}.
		\end{enumerate}
	\end{rem}
	
	\begin{dfn}
		\label{very-large-cat}
		A \iw{very large category} $\mcg{C}$ is defined in the same way as (the first definition of) a category except that instead of requiring $\uo{\mcg{C}}$ to be a class and $\mcg{C}(X,Y)$ to be a set for all $X,Y$, we require them all to be conglomerates.
	\end{dfn}
	
	\begin{rem}
		\label{very-large-cat-rems}
		\begin{enumerate}
			\item Since conglomerates are closed under union, we find out that $\um{\mcg{C}}$ must also be a conglomerate for the very large category $\mcg{C}$.
			
			\item Again, clear from the definition is the fact that every (large or small) category is at the same time a very large category, but the converse statement is \textit{false} in general. Very large categories that are not large are called \textbf{properly very large categories}.
			
		\end{enumerate}
	\end{rem}
	
	\begin{prp}
		\label{sets-form-Set}
		Sets and functions between them constitute a (properly large)\\ category, which is denoted by \cat{Set}.
	\end{prp}
	
	\rule{0pt}{0.1cm}
	Other examples of categories are:
	\begin{itemize}
		\item the empty category $\eb{0}$ with no objects and no arrows;
		
		\item the category $\eb{1}$ with one object and one (identity) arrow;
		
		\item the category $\eb{2}$ with two objects, two identity arrows, and one non-identity arrow:
		\[\bullet\xrightarrow{~~~~~~~~~~}\bullet\] (the identities are not shown);
		
		\item every set might be seen as a \textit{discrete category}, i.e., a category in which all morphisms are identities;
		
		\item \textit{monoids} are categories with only one object, and the arrows are all from that single object to itself;
		
		\item topological spaces and continuous mappings constitute a category denoted by \cat{Top};
		
		\item groups and group homomorphisms make a category \cat{Grp};
		
		\item $\cat{Vect}_\Bbbk$, the category of vector spaces and linear transformations over a field $\Bbbk$;
		
		\item $ \cat{Aut} $, the category of all deterministic Moore automata and simulations between them;
		
		\item \cat{Grph}, the category of graphs and graph homomorphisms;
		
		\item $ \cat{Mon} $, the category of monoids and homomorphisms between them;
		
		\item $ \cat{Pos} $, the category of posets and monotone mappings; and
		
		\item $ \cat{Rel} $, the category of sets and binary relations between them.
	\end{itemize}
	
	\rule{0pt}{0.1cm}
	
	\subsubsection*{The duality principle}
	Quoting Mac Lane \cite{CWM}, ``categorical duality is the process `Reverse all arrows'\mbox{\tiny ~}''. Below, we explain this more.
	
	\begin{dfn}
		\label{dual-statement}
		Let $\mcg{S}$ be a statement about (large/very large) categories. The \textbf{dual}\index{dual statement} of $\mcg{S}$ is formed by making the following replacements throughout in $\mcg{S}$:
		\begin{itemize}
			\item ``domain'' by ``codomain'' and vice versa;
			
			\item ``$h=gf$'' by ``$h=fg$'' for any arrows $f,g,h$;
			
			\item arrows and composites are reversed; and
			
			\item logic (and, or, not, then,...) is unchanged.
		\end{itemize}
		The result is Table \ref{dual-statements-table}.
		\begin{table}[h]
			\begin{center}
				\begin{tabular}{c|c}
					Statement $\mcg{S}$~ & ~Dual Statement $\mcg{S}^d$\\[1mm] \hline
					\rule{0pt}{0pt} & \rule{0pt}{0pt} \\
					\(f:X\longrightarrow Y\) & \(f:Y\longrightarrow X\)\\[1mm]
					$A=\dom(f)$ & $A=\cod(f)$\\[1mm]
					\(i=1_A\) & \(i=1_A\)\\[1mm]
					$h=gf$ & $h=fg$\\[1mm]
					$f$ is monic & $f$ is epic\\[1mm]
					$g$ is an isomorphism & $g$ is an isomorphism\\[1mm]
					$P$ is initial & $P$ is terminal\\[1mm]
					etc. & etc.\\
				\end{tabular}
			\end{center}
			\caption{} \label{dual-statements-table}
		\end{table}
	\end{dfn}
	
	\begin{rem}
		\label{dual-of-dual-statement}
		Note that the dual of the dual is the original statement: $(\mcg{S}^d)^d=\mcg{S}$. If a statement involves a diagram, the dual statement involves that diagram with all arrows reversed.
	\end{rem}
	
	The dual of each of the axioms for a category is also an axiom. Hence in any proof of a theorem about an arbitrary category from the axioms, replacing each statement by its dual gives a valid proof (of the dual conclusion). This is the \textit{duality principle}:
	
	\begin{prp}
		\label{the-duality-principle}
		\tu{\iw{The duality principle.}} If a statement $\mcg{S}$ about a (large/very large) category is a consequence of the axioms, so is the dual statement $\mcg{S}^d$.
	\end{prp}
	
	\par The duality principle provides a very useful tool in complicated situations. Having proved a theorem, we immediately gain a dual theorem by applying the duality principle. No proof of the dual theorem is needed to be given.
	
	\begin{dfn}
		\label{self-dual-property-def}
		A statement (or property) $\mcg{S}$ is called \textbf{self-dual}\index{self-dual statement} if $\mcg{S}^d=\mcg{S}$.
	\end{dfn}
	
	For example, ``being an identity morphism'' is a self-dual property.
	
	\begin{dfn}
		\label{opposite-of-cats}
		For any (large/very large) category $\mcg{C}$, the \textbf{dual}\index{dual category} (or \textbf{opposite}\index{opposite category}) \textbf{(large/very large) category of} $\mcg{C}$ is the (large/very large) category $\op{\mcg{C}}$ with the \textit{same} objects as of $\mcg{C}$ but with the arrows of $\mcg{C}$ reversed; that is, \[\op{\mcg{C}}(X,Y)=\mcg{C}(Y,X)\]
		for every $X,Y$, and \[g\op{\circ}f=f\circ g\]
		for every $f,g$.
	\end{dfn}
	
	It is easy to see that $\op{\mcg{C}}$ is indeed a (large/very large) category.
	
	\begin{prp}
		\label{C-op-op}
		For every (large/very large) category $\mcg{C}$, $\op{(\op{\mcg{C}})}=\mcg{C}$.
	\end{prp}
	
	Because of the way dual categories are defined, every statement $\mcg{S}_{\op{\mcg{C}}}(X)$ concerning an object $X$ in the category $\op{\mcg{C}}$ can be translated into a logically equivalent statement $\mcg{S}^d_{\mcg{C}}(X)$ concerning the object $X$ in the cateogry $\mcg{C}$. Obviously, in general, $\mcg{S}^d_{\mcg{C}}(X)$ is \textit{not} equivalent to $\mcg{S}_{\mcg{C}}(X)$.
	\rule{0pt}{0.1cm}
	
	\subsubsection*{Functors}
	Up to now, all of our arguments have been about objects and arrows within a single category. Now we define a ``morphism'' from one category to another.
	\begin{dfn}
		\label{cov-functor-def}
		A \textbf{(covariant) functor}\index{covariant functor}\index{functor} $F$ from a category $\mcg{C}$ to a category $\mcg{D}$, denoted by\\ \(F:\mcg{C\longrightarrow D}\), consists of the following:
		\begin{itemize}
			\item a function \[F_0:\uo{\mcg{C}}\longrightarrow\uo{\mcg{C}};\] often, by some abuse of notation, the image of $X\in\uo{\mcg{C}}$ is written $F(X)$ or just $FX$;
			
			\item for every pair $X,Y$ of objects of $\mcg{C}$, a function
			\[F_1:\mcg{C}(X,Y)\longrightarrow\mcg{D}(F_0(X),F_0(Y));\]
			again by abuse of notation, the image of $f\in\mcg{C}(X,Y)$ is often written $F(f)$ or just $Ff$.
		\end{itemize}
		These data are subject to the following axioms:
		\begin{enumerate}
			\item \textbf{Preservation of composition.} For every pair $f\in\mcg{C}(X,Y),~g\in\mcg{C}(Y,Z)$,
			\[F(g\circ f)=F(g)\circ F(f),\]
			where the left hand side composition is in $\mcg{C}$ whereas the right hand side composition is in $\mcg{D}$.
			
			\item \textbf{Preservation of identities.} For every $X\in\mcg{C}$, \[F(1_X)=1_{F(X)}.\]
		\end{enumerate}
	\end{dfn}
	
	\begin{dfn}
		\label{functor-compo}
		Given two functors $F=\ag{F_0}{F_1}:\mcg{A}\longrightarrow \mcg{B}$ and\\ $G=\ag{G_0}{G_1}:\mcg{B}\longrightarrow\mcg{C}$, their \iw{composition} is defined as
		\[G\circ F\eqd\ag{G_0\circ F_0}{G_1\circ F_1}:\mcg{A\longrightarrow C}.\]
	\end{dfn}
	
	This composition is obviously \textit{associative}. The \iw{identity functor} \[1_\mcg{C}\eqd\ag{1_{\uo{\mcg{C}}}}{1_{\um{\mcg{C}}}}:\mcg{C}\longrightarrow\mcg{C}\]
	is clearly an identity for the above composition law. Therefore:
	
	\begin{prp}
		\label{small-cats-form-Cat}
		Small categories and functors between them form a (properly large) category, which is denoted by \cat{Cat}.
	\end{prp}
	
	\begin{prp}
		\label{large-cats-form-CAT}
		Large categories and functors between them form a (properly very large) category, which is denoted by \cat{CAT}.
	\end{prp}
	
	\begin{ntn}
		\label{notation-for-functors}
		We sometimes use the simplified notations $FA$ and $Ff$ rather than $F(A)$ and $F(f)$. Also, we sometimes denote the action on \textit{both} objects and morphisms by \[F(A\xrightarrow{~f~}B)~=~FA\xrightarrow{~Ff~}FB.\]
	\end{ntn}
	
	Besides identity functors introduced above, there are some other important examples of functors, including:
	\begin{itemize}
		\item the \iw{constant functor}: for any categories $\mcg{A,B}$ and any $\mcg{B}$-object $B$, there is a functor \(\tu{K}_B:\mcg{A\longrightarrow B}\) with value $B$, defined by 
		\[\tu{K}_B(A\xrightarrow{~f~}\prm{A})~=~B\xrightarrow{~1_B~}B;\]
		
		\item the \iw{forgetful functor} or the \iw{underlying set functor}: let $\mcg{C}$ be a \textit{construct} over \cat{Set}, that is, a category of structured sets (such as groups, rings, etc.) and structure-preserving maps (group homomorphisms, ring homomorphisms, etc.); then there is a functor $U:\mcg{C}\longrightarrow\cat{Set}$, where in each case $U(A)$ is the underlying set of $A\in\mcg{C}$, and $U(f)$ is the underlying function of $f\in\mcg{C}$.
	\end{itemize}
	
	Let us now investigate some properties of functors.
	
	\begin{prp}
		\label{functors-presv-isos}
		All functors $F:\mcg{C\longrightarrow D}$ \iw{preserve} isomorphisms; i.e., whenever $C\xrightarrow{~f~}\prm{C}$ is a $\mcg{C}$-isomorphism, $Ff$ is a $\mcg{D}$-isomorphism.
	\end{prp}
	
	\begin{dfn}
		\label{iso-functor}
		\begin{enumerate}
			\item A functor $F:\mcg{C\longrightarrow D}$ is called an \textbf{isomorphism (functor)}\index{isomorphism functor} provided that $F$ is an isomorphism arrow in either \cat{Cat} or \cat{CAT}. This means that there exists another functor $G:\mcg{B\longrightarrow A}$ satisfying 
			\[GF=1_\mcg{C}~~~~\text{and}~~~~FG=1_\mcg{D}.\]
			We may denote the situation by $F:\mcg{C\cong D}$.
			
			\item The categories $\mcg{C,D}$ are said to be \iw{isomorphic} provided that they are isomorphic objects in \cat{Cat} or in \cat{CAT}. In other words, if there exists an isomorphism functor $F:\mcg{C\cong D}$.
		\end{enumerate}
	\end{dfn}
	
	\begin{dfn}
		\label{full-faithful-full-embed}
		Let $F=\ag{F_0}{F_1}:\mcg{C\longrightarrow D}$ be a functor.
		\begin{enumerate}
			\item $F$ is called \iw{faithful} provided that all the restrictions
			\[F_1\left|_{\mcg{C}(X,\prm{X})}\right.:\mcg{C}(X,{\prm{X}})\longrightarrow\mcg{D}(F_0(X),F_0(\prm{X})) \]
			are injective functions.
			
			\item $F$ is called \iw{full} provided that all the above restrictions are surjective functions.
			
			\item $F$ is called \iw{injective on objects} provided that $F_0$ is an injective function.
			
			\item $F$ is called a \iw{embedding} provided that $F_1$ is injective; i.e., $F$ is \iw{injective on morphisms}.
			
			\item $F$ is called a \iw{full embedding} provided that $F$ is full and an embedding.
		\end{enumerate}
	\end{dfn}
	
	\begin{prp}
		\label{iffs-for-full-embed}
		A functor is:
		\begin{enumerate}
			\item an embedding if and only if it is faithful and injective on objects;
			
			\item a full embedding if and only if it is full, faithful, and injective on objects; 
			
			\item an isomorphism if and only if it is full, faithful, and bijective on objects.
		\end{enumerate}
	\end{prp}
	
	\vspace{6pt}
	Now, we introduce the notion of ``self-duality'' for categories. 
	\begin{dfn}
		\label{selfdual}
		A category $\mcg{C}$ is called \textbf{self-dual}\index{self-dual category} whenever $\mcg{C}$ is isomorphic to $\op{\mcg{C}}$.
	\end{dfn}
	For example, $ \cat{Rel} $ is self-dual but $ \cat{Set} $ and $ \cat{Pos} $ are not.
	
	\begin{rem}
		Note that this concept may be defined slightly differently elsewhere (using \textit{equivalence of categories} instead of isomorphism), but we have chosen Definition \ref{selfdual} so that it fits the purposes of the present work.
	\end{rem}
	
	Next, we turn to ``contravariant'' functors.
	\begin{dfn}
		\label{ctv-functor-def}
		A \iw{contravariant functor} $F:\mcg{C\longrightarrow D}$ is defined in the same way as a covariant functor except that the axiom of preservation of composition changes here to:
		\[F(g\circ f)=F(f)\circ F(g).\]
	\end{dfn}
	
	In the present work we will be dealing with both covariant and contravariant functors. By default, functors will be assumed to be covariant unless stated otherwise.\\
	
	\par It turns out that every contravariant functor may be translated into (two) covariant analogues. The process is described below.
	
	\begin{dfn}
		\label{the-reverser}
		For every category $\mcg{C}$, the \iw{reverser functor} $\rvr{\mcg{C}}:\mcg{C}\longrightarrow \op{\mcg{C}}$ is defined as 
		\[\rvr{\mcg{C}}(X\xrightarrow{~f~}Y)~\eqd~X\xleftarrow{~\op{f}~}Y.\]
	\end{dfn}
	
	It is clear that for every category $\mcg{C}$, \(\rvr{\mcg{C}}\) is a (contravariant) isomorphism functor:
	\[\rvr{\op{\mcg{C}}}\circ\rvr{\mcg{C}}=1_\mcg{C}~~~~ \text{and}~~~~\rvr{\mcg{C}}\circ\rvr{\op{\mcg{C}}}=1_{\op{\mcg{C}}}.\]
	
	Now we can easily change contravariant functors into their covariant analogues, using the reversers:
	
	\begin{prp}
		\label{ctv-to-cov-by-rvr}
		Every contravariant functor $F:\mcg{C\longrightarrow D}$ can be regarded as either of the following covariant functors:
		\[F\circ\rvr{\op{\mcg{C}}}:\op{\mcg{C}}\longrightarrow\mcg{D},\] or
		\[\rvr{\mcg{D}}\circ F:\mcg{C}\longrightarrow \op{\mcg{D}}.\]
	\end{prp}
	
	In the light of Proposition \ref{ctv-to-cov-by-rvr}, we frequently substitute contravariant functors $\mcg{C\longrightarrow D}$ with either $\op{\mcg{C}}\longrightarrow \mcg{D}$ or \(\mcg{C}\longrightarrow\op{\mcg{D}}\). However, there are some exceptions to this rule. For an important example, see Definition \ref{contra-tilde-def}.
	
	\begin{ntn}
		\label{F-op-A-op-B-op}
		Let $F:\mcg{A\longrightarrow B}$ be a (covariant) functor. We denote the composite functor
		\[\op{\mcg{A}}\xrightarrow{~~\rvr{\op{\mcg{A}}}~~}\mcg{A}\xrightarrow{~~F~~}\mcg{B}\xrightarrow{~~\rvr{\mcg{B}}~~}\op{\mcg{B}}\]
		by $\op{F}:\op{\mcg{A}}\longrightarrow\op{\mcg{B}}$. Notice that $\op{F}$ is a \textit{covariant} functor. In the literature, it is called the \textbf{dual}\index{dual functor} (or \textbf{opposite}\index{opposite functor}) functor to $F$.
	\end{ntn}
	
	\begin{rem}
		\label{dual-statements-involving-functors}
		Obviously, $\op{(\op{F})}=F$. To form the dual of a categorical statement that involves functors, we make the same statement, but with each category and each functor replaced by its dual. Then, we translate this back into a statement about the original categories and functors.
	\end{rem}
	
	\begin{prp}
		\label{each-functor-selfdual}
		\begin{enumerate}
			\item Each of the following properties of functors is self-dual: ``isomorphism'', ``embedding'', ``faithful'', ``full'', ``isomorphism-dense'', and ``equivalence''.
			
			\item The notion of ``self-dual category'' introduced in Definition \ref{selfdual} is self-dual in the sense of Definition \ref{dual-statement}.
		\end{enumerate}			
	\end{prp}
	\rule{0pt}{0.1cm}
	
	\subsubsection*{Subcategories}
	\begin{dfn}
		\label{subcat-def}
		\begin{enumerate}
			\item A category $\mcg{C}$ is said to be a \iw{subcategory} of a (large/very large) category $\mcg{D}$ provided that the following conditions are satisfied:
			\begin{itemize}
				\item[(a)] $\uo{\mcg{C}}\subseteq\uo{\mcg{D}}$;
				
				\item[(b)] for each $C,\prm{C}\in\uo{\mcg{C}}$, we have \(\mcg{C}(C,\prm{C})\subseteq\mcg{D}(C,\prm{C})\);
				
				\item[(c)] for each $\mcg{C}$-object $C$, the $\mcg{D}$-identity on $C$ is the \(\mcg{C}\)-identity on $C$;
				
				\item[(d)] the composition law in $\mcg{C}$ is the restriction of the composition law in $\mcg{D}$ to the morphisms of $\mcg{C}$.
			\end{itemize}
			
			\item $\mcg{C}$ is called a \iw{full subcategory} of $\mcg{D}$ if, in addition to the above, for each\\ $C,\prm{C}\in\uo{\mcg{C}}$, we have \(\mcg{C}(C,\prm{C})=\mcg{D}(C,\prm{C})\).
		\end{enumerate}
	\end{dfn}
	
	As an important example, \cat{Cat} is a full subcategory of \cat{CAT}. More examples will come in the sequel.
	
	\begin{rem}
		\label{about-subcats}
		\begin{enumerate}
			\item From the above definition it follows that a full subcategory of a category $\mcg{D}$ can be specified by merely specifying its object class within $\mcg{D}$.
			
			\item The conditions (a), (b), and (c) of part (1) of the above definition do \textit{not} imply (c).
			
			\item If $F:\mcg{A\longrightarrow B}$ is a full functor or is injective on objects, then the image of $\mcg{A}$ under $F$ (i.e. the category formed by $F_0(\ob(\mcg{A}))$ and $F_1(\mor(\mcg{A}))$) is a subcategory of \(\mcg{B}\). However, for arbitrary functors $G:\mcg{C\longrightarrow D}$, the image of $\mcg{C}$ under $F$ need not be a subcategory of $\mcg{D}$.
		\end{enumerate}
	\end{rem}
	
	For every subcategory $\mcg{C}$ of a (large/very large) category $\mcg{D}$ there is an associated \iw{inclusion functor} $E:\mcg{C\hookrightarrow D}$. Each such inclusion is (a) an embedding, and (b) a full functor if and only if $\mcg{C}$ is a full subcategory of $\mcg{D}$. As the next proposition shows, inclusions of subcateogries are (up to isomorphism) precisely the embedding functors.
	
	\begin{prp}
		\label{subcat-embed-faith}
		A functor $F:\mcg{A\longrightarrow B}$ is a (full) embedding if and only if there exists a (full) subcategory $\mcg{C}$ of $\mcg{B}$ with inclusion functor $E:\mcg{C\longrightarrow B}$ and an isomorphism $G:\mcg{A\longrightarrow C}$ with $F=EG$.
	\end{prp}
	\rule{0pt}{0.1cm}
	
	\subsubsection*{Reflective and coreflective subcategories}
	\begin{dfn}
		\label{reflectiv-subcat-def}
		Let $\mcg{C}$ be a subcategory of $\mcg{D}$, and let $D$ be a $\mcg{D}$ object.
		\begin{enumerate}
			\item An $\mcg{C}$\textbf{-reflection}\index{reflection} (or $\mcg{C}$\textbf{-reflection arrow}) for $D$ is a $\mcg{D}$-morphism \(D\xrightarrow{~r~}C\) from $D$ to a $\mcg{C}$-object $C$ with the following universal property:\\
			for any $\mcg{D}$-morphism $D\xrightarrow{~f~}\prm{C}$ into some $\mcg{C}$-object $\prm{C}$, there exists a \textit{unique} $\mcg{C}$-morphism $\prm{f}:C\longrightarrow\prm{C}$ such that \[f=\prm{f}r.\]
			
			\item $\mcg{C}$ is called a \iw{reflective subcategory} $\mcg{D}$ provided that each $\mcg{D}$-object has a $\mcg{C}$-reflection.
		\end{enumerate}
	\end{dfn} 
	
	\begin{prp}
		\label{reflections-are-unique}
		Reflections are essentially unique, i.e.,
		\begin{enumerate}
			\item if $r_1:D\longrightarrow C_1$ and $r_2:D\longrightarrow C_2$ are $\mcg{C}$-reflections for $D$, then there exists a $\mcg{C}$-isomorphism $k:C_1\longrightarrow C_2$ such that \[r_2=kr_1;\]
			
			\item if $r_1:D\longrightarrow C_1$ is a $\mcg{C}$-reflection for $D$ and $k:C_1\longrightarrow C_2$ is a $\mcg{C}$-isomorphism, then $kr_1:D\longrightarrow C_2$ is a $\mcg{C}$-reflection for $D$.
		\end{enumerate}
	\end{prp}
	
	Reflections yield the following important proposition:
	
	\begin{prp}
		\label{reflection-functor}
		Let $\mcg{C}$ be a reflective subcategory of $\mcg{D}$, and for each $\mcg{D}$-object $D$ let $r_D:D\longrightarrow C_D$ be a $\mcg{C}$-reflection arrow. Then there exists an unique functor $R:\mcg{D\longrightarrow C}$ such that the following conditions are met:
		\begin{enumerate}
			\item \(\forall D\in\mcg{D},~ R(D)=C_D\);
			\item for each $\mcg{D}$-morphism $f:D\longrightarrow\prm{D}$, the diagram
			\begin{center}
				\mbox{
					\begin{tikzpicture}[commutative diagrams/every diagram]
					
					\node (N1) at (1.5,1.5cm) {$R(D)$};
					
					\node (N2) at (-1.5,1.5cm) {$D$};
					
					\node (N3) at (-1.5,-1.5cm) {$\prm{D}$};
					
					\node (N4) at (1.5,-1.5cm) {$R(\prm{D})$};	
					
					\path[commutative diagrams/.cd, every arrow, every label]
					(N2) edge node{$r_D$} (N1)
					(N2) edge node{$f$} (N3)
					(N3) edge node{$r_{\prm{D}}$} (N4)
					(N1) edge node{$R(f)$} (N4);
					\end{tikzpicture}	
				}
			\end{center}
			commutes.
		\end{enumerate}
	\end{prp}
	
	\begin{dfn}
		\label{reflector}
		A functor $R:\mcg{D\longrightarrow C}$ constructed according to the above proposition is called a \iw{reflector} for $\mcg{C}$.
	\end{dfn}
	
	The dual of the concept \textit{reflective subcateogry} is \textit{coreflective subcategory}. That is: 
	\begin{dfn}
		\label{coreflective}
		$\mcg{A}$ is a \iw{coreflective subcategory} of $\mcg{B}$ whenever $\op{\mcg{A}}$ is a reflective subcategory of $\op{\mcg{B}}$.
	\end{dfn}
	
	Thus, dualizing the descriptions given in \ref{reflectiv-subcat-def}, \ref{reflections-are-unique}, \ref{reflection-functor}, and \ref{reflector} gives a taste of coreflections, coreflective subcategories, and coreflector functors.
	\rule{0pt}{0.1cm}
	
	\subsubsection*{Natural transformations}
	\begin{dfn}
		\label{NT-def}
		Consider two functors $F,G:\mcg{C\dlong D}$. A \iw{natural transformation} $\mu:F\xrightarrow{~~\bullet~~}G$ from $F$ to $G$ is a class of $\mcg{D}$-morphisms $\langle \mu_X:FX\longrightarrow GX\rangle_{X\in\mcg{C}}$, indexed by the $\mcg{C}$-objects, and such that for every morphism $f:X\longrightarrow Y$ Diagram \ref{NT-def-diag} commutes.
		\begin{figure}[H]
			\begin{center}
				\begin{tikzpicture}[commutative diagrams/every diagram]
				
				\node (N1) at (1.5,1.5cm) {$GX$};
				
				\node (N2) at (-1.5,1.5cm) {$FX$};
				
				\node (N3) at (-1.5,-1.5cm) {$FY$};
				
				\node (N4) at (1.5,-1.5cm) {$GY$};	
				
				\path[commutative diagrams/.cd, every arrow, every label]
				(N2) edge node{$\mu_X$} (N1)
				(N2) edge node{$Ff$} (N3)
				(N3) edge node{$\mu_Y$} (N4)
				(N1) edge node{$Gf$} (N4);
				\end{tikzpicture}
			\end{center}
			\caption{} \label{NT-def-diag}
		\end{figure}
	\end{dfn}
	
	\begin{prp}
		\label{functor-cat-def}
		Let $\mcg{A}$ be a small category and $\mcg{B}$ a category. The functors from $\mcg{A}$ to $\mcg{B}$ and the natural transformarions between them constitute a category. That category is small if $\mcg{B}$ is small.
	\end{prp}
	
	The above-defined category is called the \iw{functor category} from $\mcg{A}$ to $\mcg{B}$ and is denoted by $\mcg{B^A}$.
	
	\begin{dfn}
		\label{nat-iso-def}
		Let $F,G:\mcg{A\longrightarrow B}$ be functors.
		\begin{enumerate}
			\item A natural transformation $\mu:F\xrightarrow{~\bullet~}G$ whose components $\mu_A$ are isomorphisms is called a \iw{natural isomorphism}.
			
			\item $F$ and $G$ are said to be \textbf{naturally isomorphic} (denoted by $F\cong G$) provided that there exists a natural isomorphism from $F$ to $G$.
		\end{enumerate}
	\end{dfn}
	As a trivial example, an \textit{identity natural transformation} from a functor to itself is automatically a natural isomorphism.
	
	\begin{prp}
		\label{ref-coref-nat-iso}
		\begin{enumerate}
			\item If $\mcg{A}$ is a reflective subcategory of $\mcg{B}$, then any two reflectors for $\mcg{A}$ are naturally isomorphic.
			
			\item A dual statement holds for coreflective subcategories.
		\end{enumerate}
	\end{prp}
	
	\begin{prp}
		\label{eqv-iff-isos}
		A functor $\mcg{A}\xrightarrow{~F~}\mcg{B}$ is an equivalence if and only if there exists a functor $\mcg{B}\xrightarrow{~G~}\mcg{A}$ such that \(1_\mcg{A}\cong GF\) and \(FG\cong 1_\mcg{B}\).
	\end{prp}
	From the above result it is clear that every isomorphism of two categories is automatically an equivalence between them.
	\rule{0pt}{0.1cm}
	
	\subsubsection*{Bifunctors and endofunctors}
	\begin{dfn}
		\label{bifunctor-def}
		A \iw{bifunctor} or \textbf{functor of two variables} is a functor \(F:\mcg{C_1\times C_2\longrightarrow D}\) whose domain is the product of two categories $\mcg{C_1,C_2}$. In such case, then, $F$ is called a \textit{bifunctor from} $\mcg{C_1}$ \textit{and} $\mcg{C_2}$ \textit{to} $\mcg{D}$, and the situation is also denoted by \[F(-,-):\mcg{C_1\times C_2\longrightarrow D},\] where the first placeholder ``$(-)$'' accepts objects and morphisms from $\mcg{C_1}$, while the second placeholder does the same with $\mcg{C_2}$.
	\end{dfn}
	
	\begin{rem}
		\label{fixation-of-one-var}
		Given a bifunctor $F$ as above, fixation of one of its variables yields another (ordinary) functor. That is, assuming a \textit{fixed} $\mcg{C_1}$-object $C_1$, we can derive another functor
		\[\prm{F}\eqd F(C_1,-):\mcg{C_2\longrightarrow D},\] 
		\[\prm{F}(A\xrightarrow{~f~}B)~\eqd~F(C_1,A)\xrightarrow{~F(1_{C_1},f)~}F(C_1,B).\]
		Similarly, fixation of any $\mcg{C_2}$-object $C_2$ yields another functor
		\[F''\eqd F(-,C_2):\mcg{C_1\longrightarrow D}.\]
	\end{rem}
	
	~\\
	\par Famous examples of bifunctors include the following:
	\begin{itemize}
		\item The \textit{product} bifunctor
		\[(-)\times(-):\mcg{C\times C\longrightarrow C},\]
		taking any $\ag{A}{B}$ to $A\times B$, and any \(\ag{h}{k}\) to $h\times k$.
		
		\item The \textit{coproduct} bifunctor 
		\[(-)\sqcup(-):\mcg{C\times C\longrightarrow C},\]
		taking any $\ag{A}{B}$ to $A\sqcup B$, and any \(\ag{h}{k}\) to $h\sqcup k$.
		
		\item The \iw{hom functor}
		\[\hom(-,-):\op{\mcg{C}}\times\mcg{C}\longrightarrow\cat{Set}\]
		on any small category $\mcg{C}$; this bifunctor takes any pair of objects $\ag{A}{B}$ to their ``hom-set'' $\mcg{C}(A,B)$, and its action on morphisms is defined by formulas
		\[\hom(1_A,f)(g)\eqd f\circ g,~\forall g\in\mcg{C},\]
		and
		\[\hom(f,1_A)(\op{g})\eqd g\circ f,~\forall \op{g}\in\op{\mcg{C}}.\]
		The restriction $\hom(A,-):\mcg{C}\longrightarrow\cat{Set}$ is called the \textbf{covariant hom functor} for every fixed $A\in\op{\mcg{C}}$, while the restriction \(\hom(-,A):\op{\mcg{C}}\longrightarrow\cat{Set}\) is called the \textbf{contravariant hom functor} for every fixed $A\in\mcg{C}$. Specifically, when $\mcg{C}=\cat{Set}$, $\hom(A,B)$ is the \textit{exponential} $B^A$, i.e., the set of maps from $A$ to $B$.
	\end{itemize}
	
	Furthermore, we will introduce two other famous bifunctors in Chapter 2: the \textit{tensor product} and the \textit{internal hom}. These two are fundamental to the theory of monoidal categories, and, in particular, $*$-autonomous categories.
	
	\begin{dfn}
		\label{endofunctor-def}
		Let $\mcg{C}$ be a category. An \iw{endofunctor} on $\mcg{C}$ is a functor
		\[F:\mcg{C\longrightarrow C}\]
		from $\mcg{C}$ to itself.
	\end{dfn}
	
	Endofunctors constitute the basis for the topic of \textit{universal dialgebra} which we will introduce in Chapter 3.
	\rule{0pt}{0.1cm}
	
	\subsubsection*{Cartesian closed categories}
	\begin{dfn}
		\label{CCC-def}
		A category $\mcg{C}$ is called \iw{cartesian closed} if it has finite products, and for each $\mcg{C}$-object $A$ the functor $(-)\times A:\mcg{C\longrightarrow C}$ is left adjoint (to some other functor).
	\end{dfn}
	
	\begin{ntn}
		\label{ntn-for-ccc}
		The right adjoint for the above functor is denoted by\\ $(-)^A:\mcg{C\longrightarrow C}$, and is called the \iw{exponential}. We call the objects $B^A$ \textbf{power objects}. For every $B\in\mcg{C}$, the member $B^A\times A\longrightarrow B$ of the counit of the adjunction is denoted by $\tu{ev}_B$ and is called the \iw{evaluation} (at $B$).
	\end{ntn}
	
	There is a natural isomorphism between the morphisms from a binary product\\ $B\times A\xrightarrow{~f~}C$ and the morphisms to an exponential object $B\xrightarrow{~g~}C^A$. There is a famous terminolgy for this bijective correspondence in the literature:
	
	\begin{ntn}
		\label{curry-decurry}
		Let $B,C$ be any two $\mcg{C}$-objects. The bijection
		\[\mcg{C}(B\times A,C)\xrightarrow{~\cong~}\mcg{C}(B,C^A)\]	
		is called \iw{currying}, while its inverse is called \iw{decurrying} (or \iw{uncurrying}). Beware, however, that the terms ``currying'' and ``decurrying'' are used in the even broader context of \textit{monoidal categories}. We will be dealing with this issue in Chapter 2.
	\end{ntn}
	\rule{0pt}{0.1cm}
	
	\subsubsection*{Binary products and pullbacks in \cat{CAT}}
	Because of later usage in the study of (large) double cateogries, products and pullbacks in \cat{CAT} deserve special attention. In this last part of the first section of the current chapter, we take a look at the precise formulations of products and pullbacks in \cat{CAT}, the properly very large category of all large categories and functors between them.
	
	\begin{rem}
		\label{since-Cat-subcat-of-CAT}
		Since \cat{Cat} is a full subcategory of \cat{CAT}, the formulations to be given here apply to \cat{Cat} as well.
	\end{rem}
	
	\cat{CAT} has binary products as well as pullbacks. The precise formulation of these two constructions goes as follows.\\
	
	\par \textbf{Binary products.} For large categories $\mcg{A,B}$, their product is a triple \mbox{$\agg{\mcg{A\times B}}{\pi_1}{\pi_2}$}, where 
	\begin{itemize}
		\item $\mcg{A\times B}$ is the large category consisting of: 
		\begin{itemize}
			\item objects: all pairs \(\ag{X}{Y}\) with $X\in\ob(\mcg{A})$ and $Y\in\ob(\mcg{B})$,
			
			\item morphisms: all pairs $\ag{f}{g}:\ag{X}{Y}\longrightarrow\ag{\prm{X}}{\prm{Y}}$, with \mbox{$X\xrightarrow{~f~}\prm{X}$} in $\mcg{A}$, and \mbox{$Y\xrightarrow{~g~}\prm{Y}$} in $\mcg{B}$,
			
			\item composition: performed componentwise:
			\[\ag{h}{k}\circ\ag{f}{g}=\ag{hf}{kg},\]
			and
			
			\item identities: pairs $\ag{1_X}{1_Y}$ of an $\mcg{A}$-identity $1_X$ together with a $\mcg{B}$-identity $1_Y$;
		\end{itemize}
		
		\item $\pi_1$ and $\pi_2$ are the first and second projection functors, respectively.
	\end{itemize}
	
	The above product possesses a universal property, as described below.\\
	\textit{The universality of binary product.} For every diagram 
	\[\mcg{A}\xleftarrow{~F~}\mcg{C}\xrightarrow{~G~}\mcg{B}\]
	in \cat{CAT}, there exists a \textit{unique} functor $(F,G):\mcg{C\longrightarrow A\times B}$ making Diagram \ref{CAT-prod-diag} commute.
	
	\begin{figure}[h]
		\begin{center}
			\begin{tikzpicture}[commutative diagrams/every diagram]
			
			\node (N1) at (0,2.5cm) {$\mcg{C}$};
			
			\node (N2) at (-2.5,0cm) {$\mcg{A}$};
			
			\node (N3) at (0,0cm) {$\mcg{A\times B}$};
			
			\node (N4) at (2.5,0cm) {$\mcg{B}$};	
			
			\path[commutative diagrams/.cd, every arrow, every label]
			(N1) edge node[swap]{$F$} (N2)
			(N1) edge[dashed] node[near end]{$(F,G)$} (N3)
			(N1) edge node{$G$} (N4)
			(N3) edge node[swap]{$\pi_2$} (N4)
			(N3) edge node{$\pi_1$} (N2);
			\end{tikzpicture}	
		\end{center}
		\caption{} \label{CAT-prod-diag}
	\end{figure}
	
	The functor $(F,G)$ takes every $\mcg{C}$-object $C$ to the pair $(F,G) (C)=\ag{FC}{GC}$; also, it sends every $\mcg{C}$-arrow $C\xrightarrow{~f~}\prm{C}$ to 
	\[(F,G)(f)\eqd\ag{Ff}{Gf}:\ag{FC}{GC}\longrightarrow\ag{F\prm{C}}{G\prm{C}}.\]
	
	\par \textbf{Pullbacks.} Given any diagram 
	\[\mcg{A}\xrightarrow{~H~}\mcg{C}\xleftarrow{~K~}\mcg{B}\] 
	in \cat{CAT}, its pullback is the triple
	\[\agg{\mcg{A\times_\mcg{C}B}}{\pi_1}{\pi_2},\]
	where
	\begin{itemize}
		\item $\mcg{A\times_\mcg{C}B}$ is the large category consisting of:
		\begin{itemize}
			\item objects: all pairs \(\ag{X}{Y}\) with $X\in\ob(\mcg{A})$ and $Y\in\ob(\mcg{B})$ such that\\ \(HX=KY\),
			
			\item morphisms: all pairs $\ag{f}{g}:\ag{X}{Y}\longrightarrow\ag{\prm{X}}{\prm{Y}}$, with \mbox{$X\xrightarrow{~f~}\prm{X}$} in $\mcg{A}$, and \mbox{$Y\xrightarrow{~g~}\prm{Y}$} in $\mcg{B}$, such that \(Hf=Kg\),
			
			\item composition: performed componentwise; for consecutive pairs of\\ \mbox{\((\mcg{A\times_C B})\)}-morphisms \mbox{$\ag{X_1}{Y_1}\xrightarrow{~\ag{f}{g}~}\ag{X_2}{Y_2}$} and \mbox{$\ag{X_2}{Y_2}\xrightarrow{~\ag{h}{k}~}\ag{X_3}{Y_3}$} we have
			\[H(hf)=H(h)H(f)=K(k)K(g)=K(kg);\] hence we have
			\[\ag{h}{k}\circ\ag{f}{g}=\ag{hf}{kg}\]
			as a valid composition for $\mcg{A\times_C B}$;
			
			\item identities: pairs $\ag{1_X}{1_Y}$ of an $\mcg{A}$-identity $1_X$ together with a $\mcg{B}$-identity $1_Y$; again, we have $F(1_X)=G(1_Y)$;
		\end{itemize}
		
		\item $\pi_1$ and $\pi_2$ are the first and second projection functors, respectively.		
	\end{itemize}
	
	The above pullback possesses a universal property, as described below.\\
	
	\noindent \textit{The universality of pullback.} For every commutative square\\
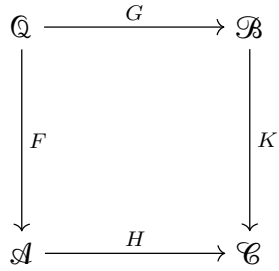
	\begin{center}
		\mbox{
			\begin{tikzpicture}[commutative diagrams/every diagram]
			
			\node (N1) at (1.5,1.5cm) {$\mcg{B}$};
			
			\node (N2) at (-1.5,1.5cm) {$\mcg{Q}$};
			
			\node (N3) at (-1.5,-1.5cm) {$\mcg{A}$};
			
			\node (N4) at (1.5,-1.5cm) {$\mcg{C}$};	
			
			\path[commutative diagrams/.cd, every arrow, every label]
			(N2) edge node{$G$} (N1)
			(N2) edge node{$F$} (N3)
			(N3) edge node{$H$} (N4)
			(N1) edge node{$K$} (N4);
			\end{tikzpicture}
		}
	\end{center}
	in \cat{CAT}, there exists a \textit{unique} functor $(F,G):\mcg{Q\longrightarrow A\times B}$ making Diagram \ref{CAT-pullb-diag} commute. The functor $(F,G)$ takes every $\mcg{Q}$-object $Q$ to the pair $(F,G) (Q)=\ag{FQ}{GQ}$; also, it sends every $\mcg{Q}$-arrow $Q\xrightarrow{~f~}\prm{Q}$ to 
	\[(F,G)(f)\eqd\ag{Ff}{Gf}:\ag{FQ}{GQ}\longrightarrow\ag{F\prm{Q}}{G\prm{Q}}.\]
	
	\begin{figure}[h]
		\begin{center}
			\begin{tikzpicture}[->,>=angle 90,commutative diagrams/every diagram]
			
			\node (N1) at (2,1cm) {$\mcg{B}$};
			
			\node (N2) at (-1,1cm) {$\mcg{A\times_C B}$};
			
			\node (N3) at (-1,-1cm) {$\mcg{A}$};
			
			\node (N4) at (2,-1cm) {$\mcg{C}$};
			
			\node (N5) at (-3,3cm) {$\mcg{Q}$};
			
			\path[commutative diagrams/.cd, every arrow, every label]
			(N2) edge node{$\pi_2$} (N1)
			(N1) edge node{$K$} (N4)
			(N3) edge node{$H$} (N4)
			(N2) edge node[swap]{$\pi_1$} (N3)
			(N5) edge[dashed] node{$(F,G)$} (N2);
			
			\draw[bend left=40] (N5) to node [above] {$G$} (N1);
			
			\draw[bend right=40] (N5) to  node [below] {$F$} (N3);
			\end{tikzpicture}
		\end{center}
		\caption{} \label{CAT-pullb-diag}
	\end{figure}
	
	\rule{0pt}{1cm}

	\section{Internal categories}
	\label{internal-cats-section}
	In this section, we intend to generalize the very notion of ``category'' to that of \textit{internal category}. Recall from Remark \ref{usefulness-of-each-dfn-cats} that, of the two equivalent definitions of categories we introduced, the second definition (i.e., Definition \ref{def-cateogry-second}) is preferred when studying internal categories. According to that definition, a category $\mcg{C}$ has a class $\uo{\mcg{C}}$ of objects and a class $\um{\mcg{C}}$ of morphisms, together with two mappings
	\[s,t:\um{\mcg{C}}\dlong\uo{\mcg{C}},\]
	which map an arrow, respectively, to its source and target.\\
	
	\par Now suppose that we are given a category $\mcg{E}$. By ``defining an internal category in $\mcg{E}$'' we mean a process of ``emulating'' the formalism of Definition \ref{def-cateogry-second} by using objects and morphisms taken from the category $\mcg{E}$. The result of this process is, then, a structure $\mcg{D}$ that behaves much like an ordinary category. In order to achieve this goal, the category $\mcg{E}$ must possess certain objects and morphisms; in particular, besides other things, it must have an ``object of objects'', an ``object of morphisms'', two parallel ``source'' and ``target'' morphisms from the former object to the latter, and finally, pullbacks to be built on the source and target morphisms. 
	\par Since we will be frequently working with pullbacks, especially with arrows \textit{between} pullbacks, it is necessary to be able to view the notion of pullback as a \textit{bifunctor}. But this requires that we introduce the notion of \textit{slice category} \cite{Lart} firstly.
	
	\begin{dfn}
		\label{slice-cat-def}
		Let $\mcg{C}$ be any (large/very large) category, and let $A$ be a $\mcg{C}$-object. The \textbf{(large/very large) slice category}\index{slice category} $\mcg{C}/A$ is defined as the category with the following data:
		\begin{itemize}
			\item \textbf{Objects.} The objects of $\mcg{C}/A$ are $\mcg{C}$-morphisms with codomain $A$.
			
			\item \textbf{Morphisms.} Given $\mcg{C}/A$-objects \(B\xrightarrow{~f~}A\) and \(C\xrightarrow{~g~}A\), a $\mcg{C}/A$-morphism is a commutative triangle of $\mcg{C}$:
			\begin{center}
				\mbox{
					\begin{tikzpicture}[commutative diagrams/every diagram]
					
					\node (N1) at (0,-1cm) {$A$};
					
					\node (N2) at (2,1cm) {$C$};
					
					\node (N3) at (-2,1cm) {$B$};
					
					\path[commutative diagrams/.cd, every arrow, every label]
					(N2) edge node{$g$} (N1)
					(N3) edge node{$h$} (N2)
					(N3) edge node[swap]{$f$} (N1);
					\end{tikzpicture}	
				}
			\end{center}
			We write $h:f\longrightarrow g$.
			
			\item \textbf{Composition.} Composition is defined by pasting the triangles together, so the composition $k\circ h$ from $\mcg{C}$ is also the composite in $\mcg{C}/A$. 
			
			\item \textbf{Identities.} For any $B\xrightarrow{~f~}A$, the identity arrow $1_B:B\longrightarrow B$ in $\mcg{C}$ is also the identity arrow $1_f:f\longrightarrow f$ in $\mcg{C}/A$.
		\end{itemize}
	\end{dfn}
	
	Now we turn to pullbacks.
	
	\begin{prp}
		\label{slice-prof-iff-pullb}
		Let $\mcg{E}$ be a (large/very large) category, and let $A$ be an $\mcg{E}$-object. Then, for any diagram \[B\xrightarrow{~f~}A\xleftarrow{~g~}C\]
		in $\mcg{E}$, the triple $\agg{P}{\pi_1}{\pi_2}$ is a pullback in $\mcg{E}$ if and only if $\agg{P}{\pi_1}{\pi_2}$ is the (binary) product of $f,g$ in the slice category $\mcg{E}/A$.
	\end{prp}
	
	Proposition \ref{slice-prof-iff-pullb} tells us that a pullback diagram in a category is essentially the same thing as the binary product of the corresponding arrows in the corresponding slice category. Therefore, we obtain the following useful result.
	
	\begin{prp}
		\label{pullb-as-bif}
		Suppose a (large/very large) category $\mcg{E}$ with pullbacks. Given any $\mcg{E}$-object $A$, there exists a bifunctor 
		\[\pb{(-)}{A}{(-)}:(\mcg{E}/A)\times(\mcg{E}/A)\longrightarrow \mcg{E}/A\]
		taking any pair of $\mcg{E}/A$-objects $\ag{B\xrightarrow{~f~}A}{C\xrightarrow{~g~}A}$ to the $\mcg{E}/A$-object \(\pb{B}{A}{C}\xrightarrow{~\delta~}A\), where $\delta$ is the diagonal of the pullback diagram $\agg{\pb{B}{A}{C}}{\pi_1}{\pi_2}$ (see Diagram \ref{pullb-as-bif-diag-1}).\\
		\begin{figure}[h]
			\begin{center}
				\begin{tikzpicture}[commutative diagrams/every diagram]
				
				\node (N1) at (1.5,1.5cm) {$C$};
				
				\node (N2) at (-1.5,1.5cm) {$\pb{B}{A}{C}$};
				
				\node (N3) at (-1.5,-1.5cm) {$B$};
				
				\node (N4) at (1.5,-1.5cm) {$A$};
				
				\path[commutative diagrams/.cd, every arrow, every label]
				(N2) edge node{$\pi_2$} (N1)
				(N2) edge node[swap]{$\pi_1$} (N3)
				(N3) edge node[swap]{$f$} (N4)
				(N1) edge node[swap]{$g$} (N4)
				(N2) edge[dashed] node[swap]{$\delta$} (N4);
				\end{tikzpicture}
			\end{center}	
			\caption{} \label{pullb-as-bif-diag-1}
		\end{figure}
		\par Moreover, the above functor sends any pair of $\mcg{E}/A$-morphisms \(\ag{B\xrightarrow{~h}\prm{B}}{C\xrightarrow{~k~}\prm{C}}\) to the $\mcg{E}/A$ morphism \(\pb{h}{A}{k}\), that is, the unique arrow making Diagram \ref{pullb-as-bif-diag-2} commute.
		\begin{figure}[h]
			\begin{center}
				\begin{tikzpicture}[commutative diagrams/every diagram]
				
				\node (N1) at (2,2cm) {$C$};
				
				\node (N2) at (-1,3cm) {$\pb{B}{A}{C}$};
				
				\node (N3) at (-4,2cm) {$B$};
				
				\node (N4) at (2,-1cm) {$\pb{\prm{B}}{A}{\prm{C}}$};
				
				\node (N5) at (5,-2cm) {$\prm{C}$};
				
				\node (N6) at (-1,-2cm) {$\prm{B}$};
				
				\node (N7) at (-1.5,1cm) {$A$};
				
				\node (N8) at (1.5,-3cm) {$A$};
				
				\path[commutative diagrams/.cd, every arrow, every label]
				(N7) edge node{$1_A$} (N8)
				(N2) edge node{$\pi_2$} (N1)
				(N2) edge node[swap]{$\pi_1$} (N3)
				(N3) edge[thick] node[swap]{$h$} (N6)
				(N1) edge[thick] node[swap]{$k$} (N5)
				(N4) edge node[swap]{$\prm{\pi}_2$} (N5)
				(N4) edge[-,line width=6pt,draw=white] node{} (N6)
				(N4) edge node[near end]{$\prm{\pi}_1$} (N6)
				(N1) edge node[swap,near start]{$g$} (N7)
				(N3) edge node{$f$} (N7)
				(N5) edge node{$\prm{g}$} (N8)
				(N6) edge node[swap]{$\prm{f}$} (N8)
				(N2) edge node[swap]{$\delta$} (N7)
				(N4) edge node{$\prm{\delta}$} (N8)
				(N2) edge[-,line width=6pt,draw=white] node{} (N4)
				(N2) edge[very thick] node[near end]{$\pb{h}{A}{k}$} (N4);	
				\end{tikzpicture}
			\end{center}
			\caption{} \label{pullb-as-bif-diag-2}
		\end{figure}
	\end{prp}
	
	Now that we are equipped with the above results, we can give the main definition of this section. The material below is mainly from Chapter 8 of \cite{Borc1}.
	
	\begin{dfn}
		\label{internal-cat-def}
		Let $\mcg{E}$ be a (large/very large) category with pullbacks. By an \iw{internal category} in $\mcg{E}$ we mean a sextuple \(\mcg{D}=\langle A,B,s,t,i,\bowtie\rangle\) consisting of:
		\begin{enumerate}
			\item an object \(A\in\uo{\mcg{E}}\), called the \iw{object of objects},
			
			\item an object $B\in\uo{\mcg{E}}$, called the \iw{object of arrows},
			
			\item two morphisms $s,t:B\dlong A$ in $\um{\mcg{E}}$, called, respectively, \iw{source} and \iw{target},
			
			\item an arrow $i:A\longrightarrow B$ in $\um{\mcg{E}}$, called \iw{identity},
			
			\item an arrow $\bowtie:B\times_A B\longrightarrow B$ in $\um{\mcg{E}}$, called \iw{composition}, where the pullback \(\agg{B\times_A B}{\pi_1}{\pi_2}\) is that of $t,s$ (Diagram \ref{pullb-s-t-def}).
			\begin{figure}[h]
				\begin{center}
					\begin{tikzpicture}[commutative diagrams/every diagram]
					
					\node (N1) at (2,1cm) {$B$};
					
					\node (N2) at (-1,1cm) {$B\times_A B$};
					
					\node (N3) at (-1,-1cm) {$B$};
					
					\node (N4) at (2,-1cm) {$A$};
					
					\path[commutative diagrams/.cd, every arrow, every label]
					(N2) edge node{$\pi_2$} (N1)
					(N1) edge node{$s$} (N4)
					(N3) edge node{$t$} (N4)
					(N2) edge node[swap]{$\pi_1$} (N3);
					\end{tikzpicture}
				\end{center}
				\caption{} \label{pullb-s-t-def}
			\end{figure}
		\end{enumerate}
		These data must satisfy the following axioms:
		\begin{itemize}
			\item[\textbf{A1}]~ $si=1_A=ti$ ~\textit{(existence of ``identity arrows'')};
			
			\item[\textbf{A2}]~ $t\pi_2=t~\circ\bowtie~$ and $~s\pi_1=s~\circ\bowtie$ ~\textit{(match-up of the target and source of composable pairs)};
			
			\item[\textbf{A3}]~ $\bowtie\circ~(is,1_B)=1_B=\text{\tiny~}\bowtie\circ~(1_B,it)$ ~\textit{(identity property for identity arrows, see Remark \ref{existence-of-uniques} below)}, where the notation $(f,g)$ denotes the unique arrow to the pullback;
			
			\item[\textbf{A4}]~ $\bowtie\circ~(1_B\times_A\bowtie)=\text{\tiny~}\bowtie\circ~(\bowtie\times_A 1_B)$ ~\textit{(associativity, see Remark \ref{pullb-as-bifunctor} below)}.
		\end{itemize}
	\end{dfn}
	
	\begin{rem}
		\label{existence-of-uniques}
		For Axiom \textbf{A3}, the equations
		\[tis=s1_B~~~~\text{and}~~~~t1_B=sit\]
		imply existence of the arrows  $(is,1_B)$ and $(1_B,it)$, respectively.
	\end{rem}
	
	\begin{rem}
		\label{pullb-as-bifunctor}
		For Axiom \textbf{A4}, we have two isomorphic objects of ``composable pairs'', namely $\pb{B}{A}{(\pb{B}{A}{B})}$ and $\pb{(\pb{B}{A}{B})}{A}{B}$. The former is seen as the pullback of $t,s\pi_1$, while the latter is seen as the pullback of $t\pi_2,s$. Viewing the pullback $\pb{(-)}{A}{(-)}$ as a bifunctor (Proposition \ref{pullb-as-bif}), we find the following arrows in $\mcg{E}$:
		\[\pb{1_B}{A}{\bowtie}:\pb{B}{A}{(\pb{B}{A}{B})}\longrightarrow\pb{B}{A}{B},\]
		\[\pb{\bowtie}{A}{1_B}:\pb{(\pb{B}{A}{B})}{A}{B}\longrightarrow\pb{B}{A}{B}.\]
	\end{rem}
	\rule{0pt}{2mm}	
	
	\section{Monoidal Categories}
	\label{monoidal-cats-section}
	In Section \ref{quick-rev} we gave two equivalent definitions for categories. Then, we noted in Remark \ref{usefulness-of-each-dfn-cats} that each of those definitions is used in some parts of the current work. Here we want to use the first one (Definition \ref{def-cateogry-first}) as a `pattern' and enter the topic of monoidal categories based on that. The material of this section is mainly from \cite{Borc2}.
	\begin{dfn}
		\label{mon-cat}
		A \iw{monoidal category} $\mathcal{V}$ consists of:
		\begin{enumerate}
			\item a category $\mcg{V}$;
			
			\item a bifunctor $(-)\boxtimes(-):\mcg{V}\times\mcg{V}\longrightarrow\mcg{V}$, called the \iw{tensor product}. We write \(A\boxtimes B\) for the image under $\boxtimes$ of the pair $\ag{A}{B}$;
			
			\item an object $I\in\mcg{V}$, called the \iw{unit}; 
			
			\item for every triple $A,B,C$ of objects an \iw{assoicator} isomorphism 
			\[a_{ABC}:(A\boxtimes B)\boxtimes C\xc A\boxtimes(B\boxtimes C)\]
			natural in $A,B,C$;
			
			\item for every object $A$, a \iw{left unitor} isomorphism \[l_A:I\boxtimes A\xc A\] natural in $A$;
			
			\item for every object $A$, a \iw{right unitor} isomorphism \[r_A:A\boxtimes I\xc A\] natural in $A$.			
		\end{enumerate}
		
		These data must satisfy the following requirements:
		
		\begin{enumerate}
			
			\item Diagram \ref{asso-coh}~is commutative for every quadruple of objects $A,B,C,D$ (\iw{associativity coherence});
			\begin{figure}[h]
				\begin{center}
					\begin{tikzpicture}[commutative diagrams/every diagram]
					
					\node (N1) at (2.5,2cm) {$(A\boxtimes B)\boxtimes(C\boxtimes D)$};
					
					\node (N2) at (-2.5,2cm) {$((A\boxtimes B)\boxtimes C)\boxtimes D$};
					
					\node (N3) at (-2.5,0cm) {$(A\boxtimes(B\boxtimes C))\boxtimes D$};
					
					\node (N4) at (-2.5,-2cm) {$A\boxtimes((B\boxtimes C)\boxtimes D)$};
					
					\node (N5) at (2.5,-2cm) {$A\boxtimes(B\boxtimes(C\boxtimes D))$};
					
					\path[commutative diagrams/.cd, every arrow, every label]
					(N1) edge node[swap]{$a_{A,B,C\boxtimes D}$} (N5)
					(N2) edge node[swap]{$a_{A\boxtimes B,C,D}$} (N1)
					(N2) edge node{$a_{ABC}\boxtimes 1$} (N3)
					(N3) edge node{$a_{A,B\boxtimes C,D}$} (N4)
					(N4) edge node{$1\boxtimes a_{BCD}$} (N5);
					\end{tikzpicture}
				\end{center}
				\caption{} \label{asso-coh}
			\end{figure}
			
			\item Diagram \ref{unit-coh}~is commutative for every pair $A,B$ (\iw{unit coherence}).
			\begin{figure}[h]
				\begin{center}
					\begin{tikzpicture}[commutative diagrams/every diagram]
					
					\node (N1) at (2,1cm) {$(A\boxtimes I)\boxtimes B$};
					
					\node (N2) at (-2,0cm) {$A\boxtimes B$};
					
					\node (N3) at (2,-1cm) {$A\boxtimes(I\boxtimes B)$};
					
					\path[commutative diagrams/.cd, every arrow, every label]
					(N1) edge node[swap,near start]{$r_A\boxtimes 1$} (N2)
					(N3) edge node[near start]{$1\boxtimes l_B$} (N2)
					(N1) edge node{$a_{AIB}$} (N3);
					\end{tikzpicture}
				\end{center}
				\caption{} \label{unit-coh}
			\end{figure}
		\end{enumerate}
	\end{dfn}
	
	\begin{rem}
		\label{sqtensor-ntn}
		Usually, the tensor product is denoted by ``$ \otimes $'' in the literature. However, as we will be dealing with \textit{two kinds} of tensors in the sequel (one for a base monoidal category $ \mcg{V} $ and one for the Chu construction \textit{over} $ \mcg{V} $), for the sake of clarity of expressions, we prefer to preserve ``$ \otimes $'' exclusively for the tensor product of the Chu construction and use ``$ \boxtimes $'' for monoidal categories in general.
	\end{rem}
	
	\begin{dfn}
		\label{symm-mon-cat}
		With the notation as above, a monoidal category is \textbf{symmetric}\index{symmetric monoidal category} when, moreover, an isomorphism \[s_{AB}:A\boxtimes B\xc B\boxtimes A\] is given for every pair $A,B$ of objects, natural in $A,B$. These isomorphisms must be such that:
		
		\begin{enumerate}
			
			\item Diagram \ref{asso-coh-symm}~is commutative for every triple $A,B,C$ of objects (\iw{associativity coherence with symmetry});
			\begin{figure}
				\begin{center}
					\begin{tikzpicture}[commutative diagrams/every diagram]
					
					\node (N1) at (2.5,2cm) {$(B\boxtimes A)\boxtimes C$};
					
					\node (N2) at (-2.5,2cm) {$(A\boxtimes B)\boxtimes C$};
					
					\node (N3) at (-2.5,0cm) {$A\boxtimes(B\boxtimes C)$};
					
					\node (N4) at (-2.5,-2cm) {$(B\boxtimes C)\boxtimes A$};
					
					\node (N5) at (2.5,-2cm) {$B\boxtimes(C\boxtimes A)$};
					
					\node (N6) at (2.5,0cm) {$B\boxtimes(A\boxtimes C)$};
					
					\path[commutative diagrams/.cd, every arrow, every label]
					(N1) edge node[swap]{$a_{BAC}$} (N6)
					(N6) edge node[swap]{$1\boxtimes s_{AC}$} (N5)
					(N2) edge node[swap]{$s_{AB}\boxtimes 1$} (N1)
					(N2) edge node{$a_{ABC}$} (N3)
					(N3) edge node{$s_{A,B\boxtimes C}$} (N4)
					(N4) edge node{$a_{BCA}$} (N5);
					\end{tikzpicture}
				\end{center}
				\caption{} \label{asso-coh-symm}
			\end{figure}
			
			\item Diagram \ref{unit-coh-symm}~is commutative for every object $A$ (\iw{unit coherence with symmetry});
			\begin{figure}[h]
				\begin{center}
					\begin{tikzpicture}[commutative diagrams/every diagram]
					
					\node (N1) at (2,1cm) {$A\boxtimes I$};
					
					\node (N2) at (-1,0cm) {$A$};
					
					\node (N3) at (2,-1cm) {$I\boxtimes A$};
					
					\path[commutative diagrams/.cd, every arrow, every label]
					(N1) edge node[swap]{$r_A$} (N2)
					(N3) edge node{$l_A$} (N2)
					(N1) edge node{$s_{AI}$} (N3);
					\end{tikzpicture}
				\end{center}
				\caption{} \label{unit-coh-symm}
			\end{figure}
			
			\item Diagram \ref{symm-axiom}~is commutative for every pair $A,B$ (\textit{the symmetry axiom}).
			\begin{figure}[h]
				\begin{center}
					\begin{tikzpicture}[commutative diagrams/every diagram]
					
					\node (N1) at (0,2cm) {$B\boxtimes A$};
					
					\node (N2) at (-2,0cm) {$A\boxtimes B$};
					
					\node (N3) at (2,0cm) {$A\boxtimes B$};
					
					\path[commutative diagrams/.cd, every arrow, every label]
					(N2) edge node{$s_{AB}$} (N1)
					(N1) edge node{$s_{BA}$} (N3);
					
					\path (N2) edge [double,thick,double distance=3pt] (N3);
					\end{tikzpicture}
				\end{center}
				\caption{} \label{symm-axiom}
			\end{figure}
		\end{enumerate}
	\end{dfn}
	
	\begin{dfn}
		\label{bicl}
		With the notation as above, a monoidal category $\mcg{V}$ is called \iw{biclosed} when for each object $B\in \mcg{V}$, both functors \[(-)\boxtimes B:\mcg{V}\longrightarrow\mcg{V}\text{~~~~,~~~~}B\boxtimes (-):\mcg{V}\longrightarrow\mcg{V}\] have right adjoints.
	\end{dfn}
	
	\begin{dfn}
		\label{csmc}
		A biclosed symmetric monoidal category is called a \iw{closed symmetric monoidal category (\cs)}.
	\end{dfn}
	
	\par Since in a symmetric monoidal category both functors $(-)\boxtimes B$ and $B\boxtimes (-)$ are naturally isomorphic, thus one obviously has:
	
	\begin{prp}
		\label{eqv-smc}
		The following are equivalent for a symmetric monoidal category $\mcg{V}$:
		\begin{itemize}
			\item[\textup{(a)}] $\mcg{V}$ is a \cs;
			\item[\textup{(b)}] for each object $B\in\mcg{V}$, the functor $(-)\boxtimes B:\mcg{V}\longrightarrow\mcg{V}$ has a right adjoint;
			\item[\textup{(c)}] for each object $B\in\mcg{V}$, the functor $B\boxtimes (-):\mcg{V}\longrightarrow\mcg{V}$ has a right adjoint.
		\end{itemize}
	\end{prp}
	
	\par The following is an important example of a \cs.
	
	\begin{dfn}
		\label{ccc}
		A category $\mcg{V}$ is \iw{cartesian closed} when it admits all finite products and, for every object $B\in\mcg{V}$, the functor $(-)\times B:\mcg{V}\longrightarrow\mcg{V}$ has a right adjoint, generally written $(-)^B:\mcg{V}\longrightarrow\mcg{V}$.
	\end{dfn}
	
	\begin{prp}
		\label{ccc-csmc}
		Every cartesian closed category is a \cs, with the cartesian product as the tensor product.
	\end{prp}
	
	\begin{mypr}
		\label{ccc-csmc-pr}
		Existence of the required natural isomorphisms and the various coherence conditions follows immediately from the universal property of the product. The tensor unit is just the terminal object.
	\end{mypr}
	
	\begin{ntn}
		\label{dash-bullet}
		When $\mcg{V}$ is a \cs, we write
		\begin{align*}
		B\multimap(-): &\text{~}\mcg{V}\longrightarrow\mcg{V}\\
		&\text{~}C\mapsto B\multimap C
		\end{align*}
		for the right adjoint to the functor $(-)\boxtimes B$. In particular, the isomorphisms \[\mcg{V}(B,B)\cong\mcg{V}(I\boxtimes B,B)\cong\mcg{V}(I,B\multimap B)\] yield a ``unit'' morphism \(u_B:I\longrightarrow B\multimap B\), corresponding with the identity on $B$. In an analogous way the isomorphisms \[\mcg{V}(C,C)\cong\mcg{V}(C\boxtimes I,C)\cong\mcg{V}(C,I\multimap C)\] yield a morphism \[i_C:C\longrightarrow I\multimap C\] corresponding with the identity on $C$. It is also useful to consider the (counit) ``evaluation morphisms'' \[\textup{ev}_{AB}:(A\multimap B)\boxtimes A\longrightarrow B\] corresponding by adjunction with the identity on $A\multimap B$, and the ``composition morphisms'' \[c_{ABC}:(A\multimap B)\boxtimes(B\multimap C)\longrightarrow A\multimap C\] corresponding by adjunction with the composite of Diagram \ref{ev_AB}:
		\begin{figure}[h]
			\begin{center}
				\begin{tikzpicture}[commutative diagrams/every diagram]
				
				\node (N1) at (0,4cm) {$(A\multimap B)\boxtimes(B\multimap C)\boxtimes A$};
				
				\node (N2) at (0,2cm) {$(A\multimap B)\boxtimes A\boxtimes(B\multimap C)$};
				
				\node (N3) at (0,0cm) {$B\boxtimes (B\multimap C)$};
				
				\node (N4) at (0,-2cm) {$(B\multimap C)\boxtimes B$};
				
				\node (N5) at (0,-4cm) {$C$};
				
				\path[commutative diagrams/.cd, every arrow, every label]
				(N1) edge node{$\cong$} (N2)
				(N2) edge node{$\textup{ev}_{AB}\boxtimes 1$} (N3)
				(N3) edge node{$\cong$} (N4)
				(N4) edge node{$\textup{ev}_{BC}$} (N5);
				
				\end{tikzpicture}
			\end{center}
			\caption{} \label{ev_AB}
		\end{figure}
	\end{ntn}
	
	\begin{prp}
		\label{int-hom}
		On a \cs $\mcg{V}$ we get a bifunctor 
		\begin{align*}
		(-)\multimap(-): &\text{~}\mcg{V}^\textup{op}\times\mcg{V}\longrightarrow\mcg{V}\\
		&\text{~}\ag{A}{B}\mapsto A\multimap B
		\end{align*}
		called the \iw{internal hom} bifunctor on $\mcg{V}$, whose composite with the hom functor\\ \mbox{$\mcg{V}(I,-):\mcg{V}\longrightarrow\cat{Set}$} is just
		\begin{align*}
		\mcg{V}^\textup{op}\times\mcg{V}&\longrightarrow\cat{Set}\\
		\ag{A}{B}&\mapsto\mcg{V}(A,B).
		\end{align*}
	\end{prp}
	
	\begin{mypr}
		Given $f:A\longrightarrow \prm{A}$ in $\mcg{V}$, the arrow $f\multimap B:(\prm{A}\multimap B)\longrightarrow (A\multimap B)$ corresponds by adjunction with the composite \[(\prm{A}\multimap B)\boxtimes A\xrightarrow{1\boxtimes f}(\prm{A}\multimap B)\boxtimes\prm{A}\xrightarrow{\textup{ev}_{\prm{A}B}}B.\]
		It is routine to check the bi-functoriality of $(-)\multimap(-)$.
		\par On the other hand, the isomorphisms \[\mcg{V}(I,B\multimap C)\cong\mcg{V}(I\boxtimes B,C)\cong\mcg{V}(B,C)\] prove the second assertion.
	\end{mypr}
	
	\begin{prp}
		\label{uB-iC}
		In a closed symmetric monoidal category $\mcg{V}$:
		\begin{enumerate}
			\item the morphisms $u_B:I\longrightarrow B\multimap B$ are natural in $B$;
			
			\item the morphisms $i_C:C\longrightarrow I\multimap C$ are isomorphisms;
			
			\item the morphisms $i_C$ are natural in $C$;
			
			\item for all $A,B,C,D\in \mcg{V}$, Diagrams \ref{uB-iC-1}~and \ref{uB-iC-2}~commute.
		\end{enumerate}
		\begin{figure}[h]
			\begin{center}
				\begin{tikzpicture}[commutative diagrams/every diagram]
				
				\node (N1) at (3.5,2cm) {$(A \multimap C)\boxtimes(C\multimap D)$};
				
				\node (N2) at (-3.5,2cm) {$((A\multimap B)\boxtimes (B \multimap C))\boxtimes (C\multimap D)$};
				
				\node (N3) at (-3.5,0cm) {$(A\multimap B)\boxtimes((B\multimap C)\boxtimes(C\multimap D))$};
				
				\node (N4) at (-3.5,-2cm) {$(A\multimap B)\boxtimes(B\multimap D)$};
				
				\node (N5) at (3.5,-2cm) {$A\multimap D$};
				
				\path[commutative diagrams/.cd, every arrow, every label]
				(N1) edge node[swap]{$c_{ACD}$} (N5)
				(N2) edge node[swap]{$c_{ABC}\boxtimes 1$} (N1)
				(N2) edge node{$a_{(A\multimap B)(B\multimap C)(C\multimap D)}$} (N3)
				(N3) edge node{$1\boxtimes c_{BCD}$} (N4)
				(N4) edge node{$c_{ABD}$} (N5);
				\end{tikzpicture}
			\end{center}
			\caption{} \label{uB-iC-1}
		\end{figure}
		
		\begin{figure}[h]
			\begin{center}
				\begin{tikzpicture}[commutative diagrams/every diagram]
				
				\node (N1) at (0,0cm) {$A\multimap B$};
				
				\node (N2) at (4,0cm) {$(A\multimap B)\boxtimes(B\multimap B)$};
				
				\node (N3) at (4,2cm) {$(A\multimap B)\boxtimes I$};
				
				\node (N4) at (-4,2cm) {$I\boxtimes(A\multimap B)$};
				
				\node (N5) at (-4,0cm) {$(A\multimap A)\boxtimes(A\multimap B)$};
				
				\path[commutative diagrams/.cd, every arrow, every label]
				(N2) edge node{$c_{ABB}$} (N1)
				(N3) edge node{$1\boxtimes u_B$} (N2)
				(N3) edge node[swap]{$r_{A\multimap B}$} (N1)
				(N4) edge node{$l_{A\multimap B}$} (N1)
				(N4) edge node[swap]{$u_A\boxtimes 1$} (N5)
				(N5) edge node[swap]{$c_{AAB}$} (N1);
				\end{tikzpicture}
			\end{center}
			\caption{} \label{uB-iC-2}
		\end{figure}
	\end{prp}
	
	\begin{mypr}
		The inverse of the morphism is the composite \[I\multimap C\xrightarrow{r_{I\multimap C}^{-1}}(I\multimap C)\boxtimes I\xrightarrow{\tu{ev}_{IC}}C.\] The rest of the proof is routine computations.
	\end{mypr}
	
	\par There are also two other important concepts when dealing with monoidal categories: namely ``monoidal functor'' and ``monoidal transformation'' \cite{Kelly,nLab-monf}.
	
	\begin{dfn}
		\label{lax-mon-functor}
		Let $\agg{\mcg{V}}{\boxtimes_\mcg{V}}{I_\mcg{V}}$ and $\agg{\mcg{W}}{\boxtimes_\mcg{W}}{I_\mcg{W}}$ be monoidal categories. A \iw{lax monoidal functor} (or \iw{weak monoidal functor}) between them is:
		\begin{enumerate}
			\item a functor \[F:\mcg{V}\longrightarrow\mcg{W};\]
			\item a morphism \[\varepsilon:I_\mcg{W}\longrightarrow F(I_\mcg{V});\]
			\item a natural transformation \[\mu_{XY}:F(X)\boxtimes_\mcg{W}F(Y)\nat F(X\boxtimes_\mcg{V}Y)\] for all $X,Y\in \mcg{V}$, satisfying the following conditions:
		\end{enumerate}
		\begin{itemize}
			\item[(a)] (\iw{Associativity}) For all objects $X,Y,Z\in \mcg{V}$ Diagram \ref{lax-asso}~commutes.\\
			\begin{figure}[h]
				\begin{center}
					\begin{tikzpicture}[commutative diagrams/every diagram]
					
					\node (N1) at (4,2cm) {$F(X)\boxtimes_\mcg{W}(F(Y)\boxtimes_\mcg{W}F(Z))$};
					
					\node (N2) at (-4,2cm) {$(F(X)\boxtimes_\mcg{W}F(Y))\boxtimes_\mcg{W}F(Z)$};
					
					\node (N3) at (-4,0cm) {$F(X\boxtimes_\mcg{V}Y)\boxtimes_\mcg{W}F(Z)$};
					
					\node (N4) at (-4,-2cm) {$F((X\boxtimes_\mcg{V}Y)\boxtimes_\mcg{V}Z)$};
					
					\node (N5) at (4,-2cm) {$F(X\boxtimes_\mcg{V}(Y\boxtimes_\mcg{V}Z))$};
					
					\node (N6) at (4,0cm) {$F(X)\boxtimes_\mcg{W}F(Y\boxtimes_\mcg{V}Z)$};
					
					\path[commutative diagrams/.cd, every arrow, every label]
					(N1) edge node[swap]{$1\boxtimes_\mcg{W}\mu_{YZ}$} (N6)
					(N2) edge node{$a_{F(X),F(Y),F(Z)}^\mcg{W}$} (N1)
					(N2) edge node{$\mu_{XY}\boxtimes_\mcg{W}1$} (N3)
					(N3) edge node{$\mu_{X\boxtimes_\mcg{V}Y,Z}$} (N4)
					(N4) edge node{$F(a^\mcg{V}_{XYZ})$} (N5)
					(N6) edge node[swap]{$\mu_{X,Y\boxtimes_\mcg{V}Z}$} (N5);
					\end{tikzpicture}
				\end{center}
				\caption{} \label{lax-asso}
			\end{figure}
			\mbox{In this diagram, $a^\mcg{V}$ and $a^\mcg{W}$ denote the associators of the respective monoidal categories.}
			
			\item[(b)] (\iw{Unitality}) For all $X\in\mcg{V}$, Diagrams \ref{left-unitality} and \ref{right-unitality} commute, where $l^\mcg{V},l^\mcg{W},r^\mcg{V},r^\mcg{W}$ denote the left and right unitors of the two monoidal categories, respectively.
			\begin{figure}[H]
				\begin{center}
					\begin{tikzpicture}[commutative diagrams/every diagram]
					
					\node (N1) at (2.2,2cm) {$F(I_\mcg{V})\boxtimes_\mcg{W}F(X)$};
					
					\node (N2) at (-2.2,2cm) {$I_\mcg{W}\boxtimes_\mcg{W}F(X)$};
					
					\node (N3) at (-2.2,-1cm) {$F(X)$};
					
					\node (N4) at (2.2,-1cm) {$F(I_\mcg{V}\boxtimes_\mcg{V}X)$};
					
					\path[commutative diagrams/.cd, every arrow, every label]
					(N2) edge node{$\varepsilon\boxtimes_\mcg{W}1$} (N1)
					(N2) edge node[swap]{$l_{F(X)}^\mcg{W}$} (N3)
					(N4) edge node[swap]{$F(l_{X}^\mcg{V})$} (N3)
					(N1) edge node{$\mu_{I_\mcg{V},X}$} (N4);
					\end{tikzpicture}
				\end{center}
				\caption{} \label{left-unitality}
			\end{figure}
			
			\begin{figure}[H]
				\begin{center}
					\begin{tikzpicture}[commutative diagrams/every diagram]
					
					\node (N1) at (2.2,2cm) {$F(X)\boxtimes_\mcg{W}F(I_\mcg{V})$};
					
					\node (N2) at (-2.2,2cm) {$F(X)\boxtimes_\mcg{W}I_\mcg{W}$};
					
					\node (N3) at (-2.2,-1cm) {$F(X)$};
					
					\node (N4) at (2.2,-1cm) {$F(X\boxtimes_\mcg{V}I_\mcg{V})$};
					
					\path[commutative diagrams/.cd, every arrow, every label]
					(N2) edge node{$1\boxtimes_\mcg{W}\varepsilon$} (N1)
					(N2) edge node[swap]{$r_{F(X)}^\mcg{W}$} (N3)
					(N4) edge node[swap]{$F(r_{X}^\mcg{V})$} (N3)
					(N1) edge node{$\mu_{X,I_\mcg{V}}$} (N4);
					\end{tikzpicture}
				\end{center}
				\caption{} \label{right-unitality}
			\end{figure}
		\end{itemize}
	\end{dfn}
	
	\begin{dfn}
		\label{strong-mon-functor}
		If $\varepsilon$ and all $\mu_{XY}$ in Definition \ref{lax-mon-functor} are isomorphisms, then $F$ is called a \iw{monoidal functor}.
	\end{dfn}
	
	\begin{dfn}
		\label{mon-trans}
		A \iw{monoidal transformation} between monoidal functors is a natural transformation that respects the extra structure in an obvious way.
	\end{dfn}
	This way, we arrive at the end of Chapter 1. In the next chapter we will introduce the Chu construction.

	\chapter{The Chu Construction}
	\label{The-Chu-Construction}
	In Section \ref{funda-motivs} we pointed to the historical roots of the Chu construction. Now, we introduce the precise formalism of the Chu construction in this chapter. 
	\par We start with $ * $-autonomous categories; next, we will introduce the general form of the Chu construction. We will show that for every \cs~ such as $ \mcg{V} $, categories $ \mcg{V},\op{\mcg{V}} $ are the coreflective and reflective subcategories of $ \chu(\mcg{V},\Gamma) $, respectively (see Proposition \ref{coreflective-reflective}). Additionally, we will prove that if $ \mcg{V} $ is bicomplete then so is $ \chu(\mcg{V},\Gamma) $ for any $ \Gamma\in\mcg{V} $ (Proposition \ref{V-bicomp}). 
	\par Next, we will focus on the Chu construction on sets (Section \ref{Chu-Set}) and will study some of its important properties. Particularly, we will point to extensional, separable, and biextensional Chu spaces (Subsection \ref{eChu-sChu-bChu}); we will show that $ \chu(\cat{Set},\Gamma) $ is balanced (Proposition \ref{chu-balanced}); we will introduce the multiplicative as well as additive connectives of linear logic by the Chu construction (Subsection \ref{LL-Chu}); we will have a look at ``realizations'' in $ \chu $ (Subsection \ref{realize}); and finally, we will study endofunctors on $ \chu $ (Subsection \ref{endofunc}) and will show that from every endofunctor $ F:\cat{Set}\longrightarrow\cat{Set} $, one can obtain an endofunctor $ \hat{F}:\chu\longrightarrow\chu $ (at least in two ways). 
	
	\section{$*$-Autonomous categories}
	\label{star-auto-cats}
	\par Now we are at the position to introduce the notion of ``$*$-autonomous category'' \cite{Barr-LL,nLab-star}. There are two equivalent definitions for it; we choose one and prove the other as an equivalent property for such categories.
	
	\begin{dfn}
		\label{star-auto}
		Let $\mcg{V}$ be a \cs. $\mcg{V}$ is said to be \iw{$*$-autonomous} (read ``star-autonomous'') if it has a \iw{dualizing object}: an object $\perp$ such that the canonical morphism \[d:A\longrightarrow(A\multimap\perp)\multimap\perp,\] which is the transpose of the evaluation map \[\tu{ev}_{A,\perp}:(A\multimap\perp)\boxtimes A\longrightarrow\perp,\] is an isomorphism for all $A$.
	\end{dfn}
	
	\begin{prp}
		\label{dual-fun}
		A closed symmetric monoidal category $\mcg{V}$ is $*$-autonomous if and only if it is equipped with a full and faithful functor \[(-)^*:\op{\mcg{V}}\longrightarrow\mcg{V}\] such that there is a natural isomorphism \[\mcg{V}(A\boxtimes B,C^*)\cong\mcg{V}(A,(B\boxtimes C)^*).\] The functor $(-)^*$ is called the \iw{dualization functor}, and $A^*$ is said to be the \textbf{dual}\index{dual object} of $A$ for any object $A$.
	\end{prp}
	
	\begin{mypr}
		($\Longrightarrow$)  Define the dualization functor as the internal hom to the dualizing object: \[(-)^*\eqd(-)\multimap\perp.\] Then the morphism $d_A$ is natural in $A$, so that there is a natural isomorphism \[d:1_{\mcg{V}}\nat(-)^{**}.\] We also have 
		\begin{align*}
		\mcg{V}(A\boxtimes B,C^*)&=\mcg{V}(A\boxtimes B,C\multimap\perp)\\
		&\cong\mcg{V}((A\boxtimes B)\boxtimes C,\perp)\\
		&\cong\mcg{V}(A\boxtimes(B\boxtimes C),\perp)\\
		&\cong\mcg{V}(A,(B\boxtimes C)\multimap\perp)\\
		&=\mcg{V}(A,(B\boxtimes C)^*).
		\end{align*}
		This yields the desired result.\\
		\par ($\Longleftarrow$)  Conversely, define the object $\perp$ as the dual of the tensor unit: \[\perp\eqd I^*.\] 
	\end{mypr}
	
	\rule{0pt}{0.5cm}
	
	From Proposition \ref{dual-fun}~it follows that:
	
	\begin{cor}
		\label{star-selfdual}
		Every $*$-autonomous category is self-dual.
	\end{cor}
	
	\begin{prp}
		\label{tensor-vs-bullet}
		The following hold for every pair of objects $A,B$ in a $*$-autonomous category $\mcg{V}$:
		\[A\multimap B\cong(A\boxtimes B^*)^*\text{~~~~~and~~~~~}A\boxtimes B\cong(A\multimap B^*)^*.\]
	\end{prp}
	
	\begin{mypr}
		We know that there is an adjunction \[(-)\boxtimes A\dashv A\multimap(-)\] for every $A\in\mcg{V}$. Thus, to prove the left hand side isomorphism, it suffices to show that the functor \[G\eqd(A\boxtimes(-)^*)^*:\mcg{V}\longrightarrow\mcg{V}\] is right adjoint to $F\eqd(-)\boxtimes A:\mcg{V}\longrightarrow\mcg{V}$. By Proposition \ref{dual-fun}~we find the natural isomorphisms \[\mcg{V}(FB,C)\cong\mcg{V}(B\otimes A,C)\cong\mcg{V}(B\boxtimes A,C^{**})\cong\mcg{V}(B,(A\boxtimes C^*)^*)\cong\mcg{V}(B,GC)\] for all $A,B,C$. Therefore, $F\dashv G$, and by uniqueness, \[(A\boxtimes (-)^*)^*\cong A\multimap(-).\] The other isomorphism is dual to the first one.
	\end{mypr}
	
	\begin{rem}
		\label{compact-closed}
		For a $*$-autonomous category $\mcg{V}$, whenever there is a natural isomorphism \[(A\boxtimes B)^*\cong A^*\boxtimes B^*,\] the category $\mcg{V}$ is called a \iw{compact closed category}. Putting $A=I,B=\perp$ in the above isomorphisms we find out that in compact closed categories we always have $\perp\cong I$, i.e., the tensor unit is the same thing as the dualizing object. A standard example is $\cat{FinVect}_\Bbbk$, the category of finite-dimensional vector spaces over some field $\Bbbk$ together with the usual tensor product, in which the field $\Bbbk$ itself plays the roles of the tensor unit and the dualizing object simultaneously. This type of categories will be of little interest in the present work.
	\end{rem}
	
	\begin{dfn}
		\label{lin-dist-cat}
		Conversely, if a $*$-autonomous category $\mcg{V}$ is \textit{not} compact closed, then the tensor product induces another bifunctor
		\begin{align*}
		(-)\odot(-):&\mcg{V}\times\mcg{V}\longrightarrow\mcg{V}\\
		&A\odot B\eqd(A^*\boxtimes B^*)^*
		\end{align*}
		called the \iw{par operation}, making $\mcg{V}$ into a \iw{linearly distributive category}. From the definition of par it is clear that $\odot$ is a \textit{symmetric} operation: $A\odot B\cong B\odot A$. For more information on linearly distributive categories, the reader is referred to \cite{nLab-ldc} and the references therein. In the present thesis, it is this type of $*$-autonomous categories that we mainly work with; more precisely, we will be working with the Chu construction (see Sections \ref{Chu-V}~and \ref{Chu-Set}).
	\end{dfn}
	
	\par Note that the par operation is often (but not always) denoted by some special symbol in the literature, namely an upside-down ampersand ``\rotatebox{180}{\raisebox{-7pt}{\&}}'', which is borrowed from linear logic. However, this symbol is difficult to work with--both in handwriting and in typesetting--and so we do not use it here.
	
	\begin{rem}
		It is noteworthy that in a $*$-autonomous category, seen as a linearly distributive category, we have two different monoidal structures $\ag{\boxtimes}{I}$ and $\ag{\odot}{\perp}$. 
	\end{rem}
	
	\par Another important aspect of $*$-autonomous categories is their \textit{internal logic}. The internal logic of $*$-autonomous categories is the multiplicative fragment of classical linear logic; conversely, $*$-autonomous categories are the categorical semantics of classical linear logic (Proposition 1.5 in \cite{Seely}). In the sequel, we will have brief looks at the basic operations of classical linear logic based on an important family of $*$-autonomous categories: the categories of the form $\cat{Chu}(\cat{Set},\Gamma)$ (see \ref{LL-Chu}). For more information on topics of internal logic and categorical semantics, see \cite{CatSemLL}.\\
	
	\rule{0pt}{0.5cm}	
	
	\section{The Chu construction}
	\label{Chu-V}
	The Chu construction is a general method for constructing a $*$-autonomous category from a given \cs~ \cite{Barr-the,Barr-hist,nLab-chu}. It is named after Po-Hsiang Chu, a student of Michael Barr, who introduced the construction in his master's thesis at McGill University. The formalism has been extensively studied by Pratt and others (\cite{Pratt-FullComp,Pratt-chu-sp,Pratt-SemBr,Pratt-rlz,Pratt-AutoQ,Pratt-Gamut,Pratt-chu-guide}) for its potential applications in theoretical computer science.
	
	\begin{dfn}
		\label{Chu(V,Gamma)}
		Let $\mcg{V}$ be a \cs~with pullbacks, and let $\Gamma$ be an object of $\mcg{V}$. We define a category $\cat{Chu}(\mcg{V},\Gamma)$, called the \iw{Chu construction} (or the \iw{Chu category}) over $\ag{\mcg{V}}{\Gamma}$ with the following data:
		\begin{itemize}
			\item Objects: triples $\msf{A}=\agg{A}{r}{X}$, called \iw{Chu objects} or \iw{Chu spaces}, where $A,X$ are $\mcg{V}$-objects, and $r:A\boxtimes X\longrightarrow\Gamma$ is a morphism of $\mcg{V}$.
			
			\item Morphisms: pairs $f=\ag{f^+}{f^-}:\agg{A}{r}{X}\longrightarrow\agg{B}{s}{Y}$, called \iw{Chu morphisms} or \iw{Chu transforms}, where $f^+:A\longrightarrow B$ and $f^-:Y\longrightarrow X$ are $\mcg{V}$-morphisms making Diagram \ref{adj-con-V}~commute. This is called the \iw{adjointness condition} for $f^+,f^-$.
			\begin{figure}[H]
				\begin{center}
					\begin{tikzpicture}[commutative diagrams/every diagram]
					
					\node (N1) at (2,2cm) {$A\boxtimes X$};
					
					\node (N2) at (-2,2cm) {$A\boxtimes Y$};
					
					\node (N3) at (-2,-1cm) {$B\boxtimes Y$};
					
					\node (N4) at (2,-1cm) {$\Gamma$};
					
					\path[commutative diagrams/.cd, every arrow, every label]
					(N2) edge node{$1_A\boxtimes f^-$} (N1)
					(N2) edge node[swap]{$f^+\boxtimes 1_Y$} (N3)
					(N3) edge node[swap]{$s$} (N4)
					(N1) edge node{$r$} (N4);
					\end{tikzpicture}
				\end{center}
				\caption{} \label{adj-con-V}
			\end{figure}
			
			\item For objects $\msf{A}=\agg{A}{r}{X},\msf{B}=\agg{B}{s}{Y},\msf{C}=\agg{C}{t}{Z}$ and morphisms \mbox{$f:\msf{A}\longrightarrow\msf{B},g:\msf{B}\longrightarrow\msf{C}$}, the composition of $f$ and $g$ is defined as \[g\circ f\eqd\ag{g^+\circ f^+}{f^-\circ g^-},\] which satisfies the adjointness condition and, hence, is a Chu transform. To see this, observe that in Diagram \ref{adj-for-composition}, all the triangles and inner quadrilaterals commute and consequently, the outer quadrilateral commutes. 
			\begin{figure}
				\begin{center}
					\begin{tikzpicture}[->,>=angle 90,commutative diagrams/every diagram]
					
					\node (N1) at (5,0cm) {$\Gamma$};
					
					\node (N2) at (5,2cm) {$A\boxtimes X$};
					
					\node (N3) at (2,2cm) {$A\boxtimes Y$};
					
					\node (N4) at (2,0cm) {$B\boxtimes Y$};
					
					\node (N5) at (2,-2cm) {$B\boxtimes Z$};
					
					\node (N6) at (5,-2cm) {$C\boxtimes Z$};
					
					\node (N7) at (-3,0cm) {$A\boxtimes Z$};
					
					\path[commutative diagrams/.cd, every arrow, every label]
					(N2) edge node[swap]{$r$} (N1)
					(N3) edge node[swap]{$1_A\boxtimes f^-$} (N2)
					(N3) edge node[swap]{$f^+\boxtimes 1_Y$} (N4)
					(N5) edge node{$1_B\boxtimes g^-$} (N4)
					(N5) edge node{$g^+\boxtimes 1_Z$} (N6)
					(N4) edge node{$s$} (N1)
					(N6) edge node{$t$} (N1)
					(N7) edge node[swap,near start]{$1_A\boxtimes g^-$} (N3)
					(N7) edge node[near start]{$f^+\boxtimes 1_Z$} (N5);
					
					\draw[bend left=50] (N7) to node [above] {$1_A\boxtimes f^-g^-$} (N2);
					
					\draw[bend right=50] (N7) to  node [below] {$g^+f^+\boxtimes 1_Z$} (N6);
					\end{tikzpicture}
				\end{center}
				\caption{} \label{adj-for-composition}
			\end{figure}
			
			\item Identity morphisms for every $\msf{A}=\agg{A}{r}{X}$ are defined the obvious way: \[1_\msf{A}\eqd\ag{1_A}{1_X}.\]
			
			\item Finally, one can readily verify that the composition so defined is associative (up to natural isomorphism).
		\end{itemize}
	\end{dfn}
	
	Now we show the $*$-autonomous structure of the Chu construction:
	
	\begin{prp}
		\label{Chu-star-auto}
		Let $\agg{\mcg{V}}{\boxtimes}{I}$ be a \cs~with pullbacks, and let $\Gamma$ be an object of $\mcg{V}$.Then the Chu construction $\cat{Chu}(\mcg{V},\Gamma)$ is $*$-autonomous. In particular, it is self-dual.
	\end{prp}
	
	\begin{mypr}
		\label{Chu-star-auto-proof}
		There is an evident self-duality \[(-)^*:\cat{Chu}\op{(\mcg{V},\Gamma)}\longrightarrow\cat{Chu}(\mcg{V},\Gamma)\] which takes an object $\msf{A}=\agg{A}{r}{X}$ to \[\msf{A}^*\eqd\agg{X}{\breve{r}}{A},\] where $\breve{r}\eqd(X\boxtimes A\stackrel{s_{XA}}{\cong}A\boxtimes X\xrightarrow{r}\Gamma)$ is the ``transpose'' of $r$. On morphisms, it takes a pair $f=\ag{f^+}{f^-}$ to $f^*\eqd\ag{f^-}{f^+}$. Now it immediately follows that\\ 
		\begin{center}
			\mbox{
				\begin{tikzpicture}[commutative diagrams/every diagram]
				
				\node (N1) at (-2.3,0cm) {$\cat{Chu}\op{(\mcg{V},\Gamma)}$};
				
				\node (N2) at (2.2,0cm) {$\cat{Chu}(\mcg{V},\Gamma)$};
				
				\node (N3) at (-1.2,0.1cm) {};
				
				\node (N4) at (1.2,0.1cm) {};
				
				\node (N5) at (-1.2,-0.1cm) {};
				
				\node (N6) at (1.2,-0.1cm) {};
				
				\path[commutative diagrams/.cd, every arrow, every label]
				(N3) edge node{$(-)^*$} (N4)
				(N6) edge node{$\op{((-)^*)}$} (N5);
				
				\end{tikzpicture}
			}
		\end{center}
		is an isomorphism of categories.\\
		\par Next, we define a tensor and an internal hom on $\cat{Chu}(\mcg{V},\Gamma)$ and show that these make the Chu construction into a \cs. Then, using the self-duality functor $(-)^*$ defined above, we deduce the $*$-autonomous structure of $\cat{Chu}(\mcg{V},\Gamma)$.\\
		\par Define the Chu object \[\msf{I}\eqd\agg{I}{\lambda}{\Gamma},\] where $\lambda:I\boxtimes\Gamma\cong\Gamma$ is an instance of the left unitor isomorphism. Also, denote the dual of $\msf{I}$ by \[\perp\eqd\msf{I}^*=\agg{\Gamma}{\rho}{I},\] where $\rho:\Gamma\boxtimes I\cong\Gamma$ is the transpose of $\lambda$.
		\par For Chu objects $\msf{A}=\agg{A}{r}{X},\msf{B}=\agg{B}{s}{Y}$, define the $\mcg{V}$-object $H$ as the following pullback:\\
		
		\centering
		\mbox{
			\begin{tikzpicture}[commutative diagrams/every diagram]
			
			\node (N1) at (2,2cm) {$B^A$};
			
			\node (N2) at (-2,2cm) {$H$};
			
			\node (N3) at (-2,-1cm) {$X^Y$};
			
			\node (N4) at (2,-1cm) {$\Gamma^{A\boxtimes Y}$};
			
			\path[commutative diagrams/.cd, every arrow, every label]
			(N2) edge node{} (N1)
			(N2) edge node{} (N3)
			(N3) edge node{$\dot{r}$} (N4)
			(N1) edge node[swap]{$\dot{s}$} (N4);
			\end{tikzpicture}
		}
		\flushleft
		where the exponentials are used to denote the internal homs in $\mcg{V}$, \(\dot{r}\) is the result of currying $r$ to \(X\rightarrow\Gamma^A\) and then exponentiating by $Y$, and similarly for \(\dot{s}\). There is a map \[m:H\boxtimes(A\boxtimes Y)\longrightarrow\Gamma\] obtained by decurrying either leg of the above pullback; so, define \[K\eqd A\boxtimes Y.\] Now define the internal hom $\msf{A}\multimap\msf{B}\eqd\agg{H}{m}{K}$.\\
		\par To form the tensor product $\msf{A}\otimes\msf{B}$, we use the formula \[\msf{A}\otimes\msf{B}\cong(\msf{A}\multimap \msf{B}^*)^*\] as was stated and proved in Proposition \ref{tensor-vs-bullet}. Therefore we define \[\msf{A}\otimes\msf{B}\eqd\agg{A\boxtimes B}{t}{Y^A\times_{\Gamma^{A\boxtimes B}} X^B},\] where the third component is a pullback, and the pairing $t$ is obvious. There are three things to be checked:
		\begin{enumerate}
			\item the presence of a canonical isomorphism \[\cat{Chu}(\mcg{V},\Gamma)(\msf{A}\otimes\msf{B},\msf{C})\cong\cat{Chu}(\mcg{V},\Gamma)(\msf{A},\msf{B}\multimap\msf{C})\] for every $\msf{A,B,C}$;
			\item verifying that $\msf{I}$ is indeed the tensor unit: \[\msf{A\otimes C\cong I\otimes A\cong A}\;\]
			\item verifying that $\perp$ is indeed the dualizing object: \[\msf{A^*\cong A\multimap\perp}.\]
		\end{enumerate}
	\end{mypr}
	
	\begin{prp}
		\label{coreflective-reflective}
		The category $\mcg{V}$ is a coreflective subcategory of $\cat{Chu}(\mcg{V},\Gamma)$, while $\op{\mcg{V}}$ is a reflective subcategory of $\cat{Chu}(\mcg{V},\Gamma)$.
	\end{prp}
	
	\begin{mypr}
		There is a strong monoidal functor \[i:\mcg{V}\longrightarrow\cat{Chu}(\mcg{V},\Gamma)\] taking every $\mcg{V}$-object $C$ to \(\agg{C}{r}{\Gamma^C}\), where $r=(C\boxtimes\Gamma^C\cong\Gamma^C\boxtimes C\xrightarrow{\tu{ev}_{C\Gamma}}\Gamma$. (Note that this does NOT take $\Gamma$ to the dualizing object $\perp$ in $\cat{Chu}(\mcg{V},\Gamma)$, unless of course the canonical morphism \(I\rightarrow\Gamma^\Gamma\) is an isomorphism.) This embedding admits a right adjoint \[p:\cat{Chu}(\mcg{V},\Gamma)\longrightarrow\mcg{V}\] given by the obvious projection, which is also strong monoidal. The unit of the adjunction is an isomorphism, hence $\mcg{V}$ is a coreflective (full) subcategory of $\cat{Chu}(\mcg{V},\Gamma)$.
		\par On the other hand, $\cat{Chu}(\mcg{V},\Gamma)$ is self-dual, hence $\op{\mcg{V}}$ also embeds as a full subcategory of $\cat{Chu}(\mcg{V},\Gamma)$, this time reflectively.
	\end{mypr}
	
	Another fundamental property of the Chu construction is expressed the next proposition:
	
	\begin{prp}
		\label{V-bicomp}
		If $\mcg{V}$ is bicomplete, then so is $\cat{Chu}(\mcg{V},\Gamma)$. The formula for the colimits is the obvious one: \[\Clm_j \agg{A_j}{r_j:A_j\boxtimes X_j\longrightarrow\Gamma}{X_j}=\agg{\Clm_j A_j}{~r}{~\Lm_j X_j}\] where $r$ is the decurrying of \[\Lm_j (X_j\rightarrow \Gamma^{A_j})\cong(\Lm_j X_j\rightarrow\Gamma^{\Clm_j A_j}),\] and the formula for limits is obtained by dualizing the formula for colimits in $\cat{Chu}(\mcg{V},\Gamma)$.
	\end{prp}
	
	\begin{mypr}
		Given bicomplete $\mcg{V}$, it suffices to show that the above formula indeed gives the colimits in $\cat{Chu}(\mcg{V},\Gamma)$. By self-duality of $\cat{Chu}(\mcg{V},\Gamma)$, then, a dual formula for limits follows automatically.\\
		\par Let $D:J\longrightarrow\cat{Chu}(\mcg{V},\Gamma)$ be an arbitrary diagram in $\cat{Chu}(\mcg{V},\Gamma)$, where $J$ is a small category, and let \(\agg{A_j}{r_j}{X_j}\xrightarrow{f_{jk}}\agg{A_k}{r_k}{X_k}\) be a family of arrows in $D(J)$ for indices $j,k$. We have two families of $\mcg{V}$-morphisms \(\left\langle f^+_{jk}:A_j\rightarrow A_k\right\rangle \) and \(\left\langle f^-_{jk}:X_k\rightarrow X_j\right\rangle\). Take the colimit $\ag{C}{c_j:A_j\rightarrow C}$ of the former and the limit $\ag{L}{l_j:L\rightarrow X_j}$ of the latter. Since $C$ is the colimit of $A_j$, \(\ag{\Gamma^C}{\Gamma^{c_j}}\) is the limit of $\Gamma^{A_j}$. For every $j$, the $\mcg{V}$-arrow \(r_j:A_j\boxtimes X_j\rightarrow\Gamma\) has a transpose \(\tilde{r}_j:X_j\rightarrow\Gamma^{A_j}\). Thus, $\ag{L}{\tilde{r}_j l_j}$ is a cone over $\Gamma^{A_j}$ and hence, there exists a unique arrow $s:L\rightarrow\Gamma^C$ such that the following square commutes:
		
		\begin{center}
			\mbox{
				\begin{tikzpicture}[commutative diagrams/every diagram]
				
				\node (N1) at (2,2cm) {$\Gamma^C$};
				
				\node (N2) at (-2,2cm) {$L$};
				
				\node (N3) at (-2,-1cm) {$X_j$};
				
				\node (N4) at (2,-1cm) {$\Gamma^{A_j}$};
				
				\path[commutative diagrams/.cd, every arrow, every label]
				(N2) edge[dashed] node{$s$} (N1)
				(N2) edge node{$l_j$} (N3)
				(N3) edge node{$\tilde{r}_j$} (N4)
				(N1) edge node[swap]{$\Gamma^{c_j}$} (N4);
				\end{tikzpicture}
			}
		\end{center}
		
		\noindent Define the $\mcg{V}$-morphism $r:C\boxtimes L\rightarrow \Gamma$ as the adjoint transpose of $s$, and form the Chu object \(\msf{C}\eqd\agg{C}{r}{L}\). Transpose the above square to obtain the following commutative square:
		
		\begin{center}
			\mbox{
				\begin{tikzpicture}[commutative diagrams/every diagram]
				
				\node (N1) at (2,2cm) {$C\boxtimes L$};
				
				\node (N2) at (-2,2cm) {$A_j\boxtimes L$};
				
				\node (N3) at (-2,-1cm) {$A_j\boxtimes X_j$};
				
				\node (N4) at (2,-1cm) {$\Gamma$};
				
				\path[commutative diagrams/.cd, every arrow, every label]
				(N2) edge node{$c_j\boxtimes 1$} (N1)
				(N2) edge node{$1\boxtimes l_j$} (N3)
				(N3) edge node{$r_j$} (N4)
				(N1) edge node[swap]{$r$} (N4);
				\end{tikzpicture}
			}
		\end{center}
		
		\noindent which proves the adjointness condition for the pair $\ag{c_j}{l_j}$. Whence, we have a Chu morphism \(\ag{c_j}{l_j}:\agg{A_j}{r_j}{X_j}\longrightarrow\msf{C}\) for all $j$. On the other hand, for every $j,k$ the following hold: \[c_j=c_k~f^+_{jk}~~~~~\text{and}~~~~~l_j=f^-_{jk}~l_k.\] Therefore, \[\ag{c_j}{l_j}=\ag{c_k~f^+_{jk}}{f^-_{jk}~l_k}=\ag{c_k}{l_k}\circ\ag{f^+_{jk}}{f^-_{jk}}=\ag{c_k}{l_k}\circ f_{jk}\] and \(\ag{\msf{C}}{\ag{c_j}{l_j}}\) is a cocone under $D(J)$ in $\cat{Chu}(\mcg{V},\Gamma)$. It is straightforward to check the universality of this cocone. We have thus found the desired colimit diagram for $D(J)$.
	\end{mypr}
	\vspace{5mm}
	\section{The Chu construction on the category of sets}
	\label{Chu-Set}
	The category $\cat{Set}$ of sets and functions is well-known in mathematics. Regarding its properties, one may deduce that $\cat{Set}$ is a very good candidate for serving as the base category $\mcg{V}$ in the Chu construction. Thus, we introduce:
	
	\begin{ntn}
		\label{Chu_Gamma}
		Putting $\mcg{V}=\cat{Set}$ in Definition \ref{Chu(V,Gamma)} and considering a nonempty set $\Gamma$, we arrive at a category \[\cat{Chu}(\cat{Set},\Gamma),\] which we denote by $\cat{Chu}_\Gamma$. Whenever $\Gamma$ is clear from the context, we will omit it and denote our category by $\chu$.
	\end{ntn}
	
	Therefore, \chu~is a $*$-autonomous category with the following data:
	
	\begin{itemize}
		\item objects: Chu spaces $\msf{A}=\agg{A}{r}{X}$ with $A,X$ sets, and \(r:A\times X\rightarrow \Gamma\) a map;
		
		\item morphisms: Chu transforms \(f=\ag{f^+}{f^-}\) with $f^+,f^-$ maps satisfying the adjointness condition.
	\end{itemize}
	
	\noindent It is useful to mention that for Chu spaces $\msf{A}=\agg{A}{r}{X},\msf{B}=\agg{B}{s}{Y}$ and a Chu transform $f:\msf{A}\longrightarrow\msf{B}$, the (diagrammatic) adjointness condition\\
	
	\begin{center}
		\mbox{
			\begin{tikzpicture}[commutative diagrams/every diagram]
			
			\node (N1) at (2,2cm) {$A\times X$};
			
			\node (N2) at (-2,2cm) {$A\times Y$};
			
			\node (N3) at (-2,-1cm) {$B\times Y$};
			
			\node (N4) at (2,-1cm) {$\Gamma$};
			
			\path[commutative diagrams/.cd, every arrow, every label]
			(N2) edge node{$1_A\times f^-$} (N1)
			(N2) edge node[swap]{$f^+\times 1_Y$} (N3)
			(N3) edge node{$s$} (N4)
			(N1) edge node[swap]{$r$} (N4);
			\end{tikzpicture}
		}
	\end{center}
	
	can be restated in terms of memberships as: \[\forall a\in A,\forall y\in Y:~ r\left(a,f^-(y) \right)=s\left(f^+(a),y \right),\] which is sometimes more convenient to use.\\
	
	\par \indent From Propositions \ref{coreflective-reflective} and \ref{V-bicomp} we deduce:
	
	\begin{cor}
		\label{set-coref-ref}
		For every nonempty set $\Gamma$, the category $\cat{Set}$ is a coreflective (full) subcategory of \(\chu_\Gamma\), while $\op{\cat{Set}}$ is a reflective subcategory of \(\chu_\Gamma\).
	\end{cor}
	
	\begin{cor}
		\label{set-lim-colim}
		For every nonempty set $\Gamma$, the category \(\chu_\Gamma\) is bicomplete. The formula for colimits in \(\chu_\Gamma\) is \[\Clm_j \agg{A_j}{r_j}{X_j}=\agg{\coprod_j A_j}{r}{\prod_j X_j},\] where $\coprod,\prod$ denote the coproduct and product in $\cat{Set}$, respectively, and $r$ is the decurrying of \[\Lm_j (x_j\rightarrow\Gamma^{A_j})\cong(\prod_j X_j\rightarrow \Gamma^{\coprod_j A_j}).\] A dual formula exists for the case of limits, also.
	\end{cor}
	
	\begin{mypr}
		$\cat{Set}$ is bicomplete, and $\chu_\Gamma$ becomes bicomplete by Proposition \ref{V-bicomp}.
	\end{mypr}
	~\\
	\par There are interesting applications of Chu spaces in the literature. For example, Abramsky \cite{Abram-Toy} views the Hilbert space of quantum mechanics as a Chu space. Also, Papadopoulos and Syropoulos \cite{Papa-Fuzzy} study fuzzy sets and fuzzy relational structures as Chu spaces.
	\subsection{Extensional, separable, and biextensional Chu spaces}
	\label{eChu-sChu-bChu}
	Fix some nonempty set $\Gamma$. In this subsection we discuss some additional properties of $\chu=\chu_\Gamma$.
	
	\begin{dfn}
		\label{carrier-cocarr}
		For a Chu space $\msf{A}=\agg{A}{r}{X}$, the set $A$ is called the \iw{carrier} of $\msf{A}$ while the set $X$ is called the \iw{cocarrier} of $\msf{A}$. Also, the map \(r:A\times X\rightarrow\Gamma\) is called the \iw{matrix} of $\msf{A}$. For arbitrary \(a\in A,x\in X\) we have the functions \[r(a,-):X\longrightarrow\Gamma\] and \[r(-,x):A\longrightarrow\Gamma\] with obvious formulas, which are called a \iw{row} and a \iw{column} (of $r$), respectively.\\
		\par In other words, $r$ is a matrix whose rows are indexed by elements of $A$, whose columns are indexed by elements of $X$, and whose entries are $\Gamma$-elements. It is interesting to note that this way, the dual Chu space of $\msf{A}$ is $\msf{A}^*=\agg{X}{\breve{r}}{A}$, where the matrix \(\breve{r}:X\times A\rightarrow\Gamma\) is just the matrix transpose of $r$ defined by the formula $\breve{r}(x,a)\eqd r(a,x)$. Especially, note that $\msf{A}^{**}=\msf{A}$, that is, the dualization functor on \chu~is an \textit{involution}. Hence the name ``transpose'' is justified.
	\end{dfn}
	
	\begin{dfn}
		\label{ext-sep-bi}
		A Chu space $\msf{A}=\agg{A}{r}{X}$ is called \iw{extensional} if for every $x,y\in X$, \[r(-,x)=r(-,y)~~~~~~~implies~~~~~~~x=y\] (i.e., no repeated columns in $r$). Dually, $\msf{A}$ is called \iw{separable} if for every $a,b\in A$, \[r(a,-)=r(b,-)~~~~~~~implies~~~~~~~a=b\] (i.e., no repeated rows). Also, $\msf{A}$ is called \iw{biextensional} if it is both extensional and separable.
	\end{dfn}
	
	\begin{ntn}
		We write $\cat{eChu}_\Gamma,\cat{sChu}_\Gamma,\cat{bChu}_\Gamma$ for the full subcategories of $\chu_\Gamma$ determined by the extensional, separable, and biextensional Chu spaces, respectively.
	\end{ntn}
	
	\begin{dfn}
		\label{biext-collapse}
		For a given Chu space $\msf{A}=\agg{A}{r}{X}$, we define the equivalence relations $\sim_A$ on $A$ and $\sim_X$ on $X$ as follows:
		\begin{align*}
		a\sim_A b~~~ &\Longleftrightarrow ~~~\forall x\in X~~ r(a,x)=r(b,x)\\
		x\sim_X y~~~ &\Longleftrightarrow ~~~\forall a\in A~~ r(a,x)=r(a,y)
		\end{align*}
		These are evidently equivalence relations. $\msf{A}$ is exactly when $\sim_X$ is the identity. Likewise, $\msf{A}$ is separable exactly when $\sim_A$ is the identity. Finally, $\msf{A}$ is biextensional exactly when both relations are the identities. There are Chu transforms 
		\begin{align*}
		\ag{q_A}{1_X}:\agg{A}{r}{X}&\longrightarrow\agg{\sfr{A}{\sim_A}}{\prm{r}}{X}\\
		\prm{r}([a],x)&=r(a,x)
		\end{align*}
		and
		\begin{align*}
		\ag{1_A}{q_X}:\agg{A}{r}{X}&\longrightarrow\agg{A}{\pprm{r}}{\sfr{X}{\sim_X}}\\
		\pprm{r}(a,[x])&=r(a,x)
		\end{align*}
		and
		\begin{align*}
		\ag{q_A}{q_X}:\agg{A}{r}{X}&\longrightarrow\agg{\sfr{A}{\sim_A}}{\ppprm{r}}{\sfr{X}{\sim_X}}\\
		\ppprm{r}([a],[x])&=r(a,x)
		\end{align*}
		where $q_A,q_X$ are the quotient maps for $\sim_A,\sim_X$, respectively. The Chu space \[\agg{\sfr{A}{\sim_A}}{\ppprm{r}}{\sfr{X}{\sim_X}}\] is called the \iw{biextensional collapse} of $\msf{A}$. We denote this Chu space by the lowercase letter $\msf{a}$.
	\end{dfn}
	
	\begin{prp}
		\label{biext-mor}
		Every Chu transform $f=\ag{f^+}{f^-}:\msf{A}\longrightarrow\msf{B}$ between Chu spaces $\msf{A}=\agg{A}{r}{X},\msf{B}=\agg{B}{s}{Y}$ induces a Chu transform \[\tilde{f}=\ag{\tilde{f}^+}{\tilde{f}^-}:\msf{a\longrightarrow b}\] in a canonical manner.
	\end{prp}
	
	\begin{mypr}
		What we have to show is that there is a pair of functions \(\tilde{f}^+:\sfr{A}{\sim_A}\rightarrow\sfr{B}{\sim_B}\) and \(\tilde{f}^-:\sfr{Y}{\sim_Y}\rightarrow\sfr{X}{\sim_X}\) that satisfy the adjointness condition. This pair of functions can be found via the following.\\
		\par For $a\in A$, define \(\tilde{f}^+([a])\eqd[f^+(a)]\); also, for $y\in Y$ define \(\tilde{f}^-([y])\eqd[f^-(y)]\). We claim that these two definitions are well-defined. We prove our claim for the first one, and the proof for the second one will be similar.\\
		\par Assume $[a]=[\prm{a}]$ for some $a,\prm{a}\in A$; then we have \[\forall x\in X~ r(a,x)=r(\prm{a},x).\] Now let $y\in Y$ be arbitrary. Thus, $f^-(y)\in X$ and \[r(a,f^-(y))=r(\prm{a},f^-(y));\] but from the adjointness condition for the Chu morphism $f$, \[s(f^+(a),y)=r(a,f^-(y))=r(\prm{a},f^-(y))=s(f^+(\prm{a}),y),\] and since $y$ was arbitrary, we have 
		\[[f^+(a)]=[f^+(\prm{a})]\] or \[\tilde{f}^+([a])=\tilde{f}^+([\prm{a}]),\] and we are done.\\
		\par It remains to show the adjointness of $\tilde{f}^+,\tilde{f}^-$:
		\begin{align*}
		\forall[a]\in\sfr{A}{\sim_A}~ \forall[y]\in\sfr{Y}{\sim_Y}:~~ \ppprm{r}\left([a],\tilde{f}^-([y])\right)&=\ppprm{r}\left([a],[f^-(y)] \right)\\
		&=r(a,f^-(y))\\
		&=s(f^+(a),y)\\
		&=\ppprm{s}\left([f^+(a)],[y]\right)\\
		&=\ppprm{s}\left(\tilde{f}^+([a]),[y] \right),
		\end{align*}
		and the proof is complete.
	\end{mypr}
	
	\begin{cor}
		\label{small-chu}
		Biextensional collapses and biextensional morphisms form a category.
	\end{cor}
	
	This category is denoted by \cat{chu} (the ``small Chu''), and the functor indeuced by Proposition \ref{biext-mor} is denoted by $\mathrm{biext}:\cat{Chu}\longrightarrow\cat{chu}$ (the \iw{biextensional collapse functor}).
	
	\par Next, we have \cite{Bif}:
	
	\begin{prp}
		\label{ext-sep-f+-}
		Suppose $f,g:\msf{A\rightrightarrows B}$ are Chu transforms. Then
		\begin{enumerate}
			\item if $\msf{A}$ is extensional, then \(f^+=g^+\) implies \(f^-=g^-\);
			\item if $\msf{B}$ is separable, then \(f^-=g^-\) implies \(f^+=g^+\);
			\item if $\msf{A,B}$ are biextensional, then \(f^+=g^+\) if and only if \(f^-=g^-\).
		\end{enumerate}
	\end{prp}
	Thus, the forward and backward components in a morphism determine each other uniquely in the category \cat{bChu}. [Bif]\\
	
	\begin{mypr}
		Let us write $\msf{A}=\agg{A}{r}{X}$ and $\msf{B}=\agg{B}{s}{Y}$. First we show (1). Suppose $\msf{A}$ is extensional and $f^+=g^+$. Then, for all $y\in Y$ and all $a\in A$ we have 
		\[r(a,f^-(y))=s(f^+(a),y)=s(g^+(a),y)=r(a,g^-(y)).\] Hence $f^-(y)=g^-(y)$ by extensionality of $\msf{A}$. Now (2) follows by dualization, and (1) and (2) together imply (3).
	\end{mypr}
	
	In the next step we take a look at monics and epics \cite{Bif}:
	
	\begin{prp}
		\label{Chu-monic-epic}
		We have:
		\begin{enumerate}
			\item A morphism $f:\msf{A\longrightarrow B}$ in \chu~is monic if and only if $f^+$ is injective and $f^-$ is surjective.
			\item A morphism $f:\msf{A\longrightarrow B}$ in \cat{eChu} is monic if and only if $f^+$ is injective.
			\item A morphism $f:\msf{A\longrightarrow B}$ in \cat{bChu} is monic if and only if $f^+$ is injective.
			\item Suppose $f:\agg{A}{r}{X}\longrightarrow\agg{B}{s}{Y}$ is a morphism in \chu~and $\agg{B}{s}{Y}$ is extensional. If $f^+$ is surjective, then $f^-$ is injective. 
			\item A morphism $f:\msf{A\longrightarrow B}$ in \chu~is epic if and only if $f^+$ is surjective and $f^-$ is injective.
			\item A morphism $f:\msf{A\longrightarrow B}$ in \cat{sChu} is epic if and only if $f^+$ is surjective.
			\item A morphism $f:\msf{A\longrightarrow B}$ in \cat{bChu} is epic if and only if $f^+$ is surjective.
			\item Suppose $f:\agg{A}{r}{X}\longrightarrow\agg{B}{s}{Y}$ is a morphism in \chu~and $\agg{A}{r}{X}$ is separable. If $f^-$ is injective, then $f^+$ is surjective. 
		\end{enumerate}
	\end{prp}
	
	\begin{mypr}
		(1)  The ``If'' part is obvious. We check the ``Only If'' part. Suppose $f:\msf{A\longrightarrow B}$ is such that for any pair of morphisms $g_i:\msf{C\longrightarrow A},~i=1,2$, if \(fg_1=fg_2\), then $g_1=g_2$. We show that $f^+$ is injective and $f^-$ is surjective. Assume \(\msf{A}=\agg{A}{r}{X},\msf{B}=\agg{B}{s}{Y}\).\\
		
		\par Firstly, we show that $f^-$ is surjective. Suppose otherwise. Choose two sets $X_1,X_2$ of the same cardinality as the cardinality of $X\bs f^-(Y)$, in a way that $f^-(Y),X_1,X_2$ be pairwise disjoint. Put \(Z\eqd f^-(Y)\cup X_1\cup X_2\). For $i=1,2$, choose a bijection \(k_i:X\bs f^-(Y)\rightarrow X_i\) and let $\hat{k}_i\eqd 1_{f^-(Y)}\cup k_i:X\rightarrow Z$ be the joint extension of the identity $1_{f^-(Y)}$ and of $k_i$. For \(a\in A,z\in Z\), let 
		\[t(a,z)\eqd\left\lbrace \begin{array}{lr}
		r(a,z)&\text{if~}z\in f^-(Y)\\
		r(a,k_1^{-1}(z))&\text{if~}z\in X_1\\
		r(a,k_2^{-1}(z))&\text{if~}z\in X_2
		\end{array}\right. \]
		Let $\msf{C}\eqd\agg{A}{t}{Z}$. Then \(\langle 1_A,\hat{k}_i\rangle,~i=1,2\) are Chu morphisms from $\msf{C}$ to $\msf{A}$. Indeed, \[r(1_A(a),z)=r(a,z)=t(a,\hat{k}_i(z))\] if \(z\in f^-(Y)\), and \[r(1_A(a),z)=r(a,z)=t(a,k_i(z))=t(a,\hat{k}_i(z))\] also, if $z\in X\bs f^-(Y)$. Moreover, \(\langle 1_A,\hat{k}_i\rangle\) yield the same composition with $f$ because $\hat{k}_1,\hat{k}_2$ behave the same on the image set \(f^-(Y)\). Now $f$ being monic implies that $k_1=k_2$, a contradiction since \(X\bs f^-(Y)\ne\varnothing\). Hence $f^-(Y)=X$. Note that for this construction to work, $\msf{C}$ \underline{cannot} be required to be extensional.\\
		\par The proof for the injectivity of $f^+$ is the same as the ``Only If'' part for item (2), given next.\\
		\par (2) and (3)  \textit{``If''.} Let $f:\msf{A}\longrightarrow\msf{B}$ be a Chu transform and $f^+$ be injective. Consider two arrows \(g_i:\msf{C}\longrightarrow\msf{A},~i=1,2\), which yield the same compositions with $f$. Then $g_1^+=g_2^+$ since $f^+$ is injective. By Proposition \ref{ext-sep-f+-} we have $g_1^-=g_2^-$. Whence, $f$ is monic.\\
		\par \textit{``Only If''}. Let $f:\msf{A}\longrightarrow\msf{B}$ be monic and write \(\msf{A}=\agg{A}{r}{X},\msf{B}=\agg{B}{s}{Y}\). Assume $a_1,a_2\in A$ are such that \(f^+(a_1)=f^+(a_2)=b\). Construct a Chu space $\msf{C}\eqd\agg{C}{t}{Z}$ as follows. Let $C=\{c\}$ be a singleton, and $Z=\Gamma$. Also, let \(t(c,\sigma)=\sigma\) for every $\sigma\in\Gamma$. Clearly, $\msf{C}$ is biextensional.\\
		\par Now define \(g_i:\msf{C\longrightarrow A},~i=1,2\) as follows. Let $g_i^+(c)\eqd a_i$. For $x\in X$, put\\ \(g_i^-(x)\eqd r(a_i,x)\in \Gamma=Z\). Then, \(g_1,g_2:\msf{C\dlong A}\) are Chu transforms. Clearly, $f^+g^+_1=f^+g^+_2$. We claim that $g_1^-f^-=g_2^-f^-$, also. Indeed, if $y\in Y$ and $x=f^-(y)$, then 
		\[g_1^-(x)=r(a_1,x)=s(b,y)=r(a_2,x)=g_2^-(x).\]
		Hence, \(g_1,g_2:\msf{C\dlong A}\) yield the same compositions with $f$. Since $f$ is monic, it follows that \(g_1^+=g_2^+\). Thus, $a_1=a_2$, and $f^+$ is injective.\\
		\par Let $y,\prm{y}\in Y$ with \(f^-(y)=f^-(\prm{y})=x\), say. Choose any $b\in B$ and then $a\in A$ with \(f^+(a)=b\). Then, \(s(b,y)=r(a,x)=s(b,\prm{y})\). Then, $y=\prm{y}$ by extensionality.\\
		\par Parts (5)--(8) are duals to parts (1)--(4), respectively.
	\end{mypr}
	
	~\\
	\par We also have the following remarkable property:
	
	\begin{prp}
		\label{chu-balanced}
		The category \chu~is balanced.
	\end{prp}
	
	\begin{mypr}
		Let \(\msf{A}=\agg{A}{r}{X},~\msf{B}=\agg{B}{s}{Y}\), and let $f:\msf{A\longrightarrow B}$ be a monic epic Chu transform. By Proposition \ref{Chu-monic-epic} parts (1), (5) we find that $f^+$ must be an injective surjective map, i.e., a bijective map. Similarly, $f^-$ must be bijective.\\
		\par Consider the pair of functions $g\eqd\ag{(f^+)^{-1}}{(f^-)^{-1}}$. We claim that $g$ is a Chu transform. If either of $B$ or $X$ is empty then our claim holds trivially. Thus, assume that $B\ne\varnothing,X\ne\varnothing$. Let $b\in B,x\in X$. Since both $f^+$ and $f^-$ are bijective, there exist unique elements $a\in A,y\in Y$ with 
		\[b=f^+(a)~~~ \text{and} ~~~x=f^-(y).\]
		From the adjointness condition for $f$ we have
		\[r(a,f^-(y))=s(f^+(a),y),\] or
		\[r((f^+)^{-1}\left( f^+(a)\right) ,f^-(y))=s(f^+(a),(f^-)^{-1}\left( f^-(y)\right).\]
		Since $b,x$ were arbitrary, we find
		\[\forall b\in B,~\forall x\in X~~~ r((f^+)^{-1}(b) ,x)=s(b,(f^-)^{-1}(x).\]
		Therefore, $g=\ag{(f^+)^{-1}}{(f^-)^{-1}}:\msf{B\longrightarrow A}$ is a Chu transform satisfying
		\[gf=1_\msf{A},~fg=1_\msf{B}.\]
	\end{mypr}
	
	\subsection{Linear logic in $\textbf{Chu}$}
	\label{LL-Chu}
	Barr proposed the Chu construction $\chu(\mcg{V},\Gamma)$ as a source of constructive models of linear logic\index{linear logic} \cite{Barr-LL,Pratt-chu-sp}. The case $\mcg{V}=\cat{Set}$ is particularly important for its combination of simplicity and generality. The latter case was first treated by Lafont and Streicher \cite{Pratt-chu-sp}. Here, we list a number of $\chu_\Gamma$ endofunctors that yield the main operations used in linear logic. We have already defined most of them in the previous sections. Below, $\msf{A}=\agg{A}{r}{X},\msf{B}=\agg{B}{s}{Y},...$ are Chu spaces and $f,g,...$ are Chu transforms.
	
	\subsubsection*{a. Multiplicative connectives:}
	These include the \textit{duality}, \textit{tensor}, \textit{linear implication}, and \textit{par} operations, together with their special objects, namely the tensor unit and perp.\\
	
	\par \noindent \textbf{a1. Duality.}\index{duality connective} The duality functor \[(-)^*:\op{\chu}\longrightarrow\chu,\] which is sometimes referred to as the \textit{perp functor} and denoted by $(-)^\perp$, gives the dualization operation \(\msf{A}^*,\) which is ``involutive'': $\msf{A^{**}=A}$.\\
	
	\par \noindent \textbf{a2. Tensor.}\index{tensor connective} From Proposition \ref{Chu-star-auto}, the tensor product on \chu~is found to be
	\[(-)\otimes(-):\chu\times\chu\longrightarrow\chu\]
	\[\msf{A\otimes B}=\agg{A\times B}{t}{Y^A\times_{\Gamma^{A\times B}}X^B}\]
	where the cocarrier is the pullback depicted in Diagram \ref{tens-pullb},
	\begin{figure}[h]
		\begin{center}
			\begin{tikzpicture}[commutative diagrams/every diagram]
			
			\node (N1) at (2,2cm) {$X^B$};
			
			\node (N2) at (-2,2cm) {$Y^A\times_{\Gamma^{A\times B}}X^B$};
			
			\node (N3) at (-2,-1cm) {$Y^A$};
			
			\node (N4) at (2,-1cm) {$\Gamma^{A\times B}$};
			
			\path[commutative diagrams/.cd, every arrow, every label]
			(N2) edge node{} (N1)
			(N2) edge node{} (N3)
			(N3) edge node{$\check{s}^A$} (N4)
			(N1) edge node[swap]{$\check{r}^B$} (N4);
			\end{tikzpicture}
		\end{center}
		\caption{} \label{tens-pullb}
	\end{figure}
	and the matrix $t$ is defined the obvious way:
	\begin{align*}
	t:(A\times B)\times(Y^A\times_{\Gamma^{A\times B}}X^B)&\longrightarrow\Gamma\\
	\ag{\ag{a}{b}~}{~\ag{\varphi}{\psi}}~~~~~~~~~~~~&\mapsto ~~s(b,\varphi(a))=r(a,\psi(b))
	\end{align*}
	On the other hand, its effect on morphisms $f:\msf{A\longrightarrow\prm{A}},g:\msf{B\longrightarrow\prm{B}}$ can be described as: \[f\otimes g:\msf{A\otimes B\longrightarrow\prm{A}\otimes\prm{B}},\] with
	\begin{align*}
	(f\otimes g)^+:A\times B&\longrightarrow\prm{A}\times\prm{B}\\
	(f\otimes g)^+&=f^+\times g^+,
	\end{align*}
	and
	\begin{align*}
	(f\otimes g)^-:{\prm{Y}}^{\prm{A}} \times_{\Gamma^{\prm{A}\times\prm{B}}} {\prm{X}}^{\prm{B}} &\longrightarrow Y^A \times_{\Gamma^{A\times B}} X^B\\
	\ag{\prm{h}}{\prm{k}}&\mapsto \ag{g^-\prm{h}f^+}{f^-\prm{k}g^+}.
	\end{align*}
	
	\rule{0pt}{1pt}
	\par \noindent \textbf{a3. Tensor unit.}\index{tensor unit} The tensor unit is \[\msf{I}=\agg{\{\bullet\}}{\lambda}{\Gamma}\] with $\{\bullet\}$ the singleton, and $\lambda(\bullet,c)=c$ for every $c\in \Gamma$.\\
	\par \noindent \textbf{a4. Linear implication.}\index{linear implication} This is exactly the internal hom bifunctor
	\[(-)\multimap(-):\op{\chu}\times\chu\longrightarrow\chu.\]
	Using the equation $\msf{A\multimap B=(A\otimes B^*)^*}$ (equation since the dualization on \chu~is actually an involution), we find
	\[\msf{A\multimap B}=\agg{B^A \times_{\Gamma^{A\times Y}} X^Y}{u}{A\times Y},\]
	with
	\begin{align*}
	u:(B^A \times_{\Gamma^{A\times Y}} X^Y) \times (A\times Y) &\longrightarrow \Gamma\\
	\langle \ag{\varphi}{\psi}~,~\ag{a}{y}\rangle~~ &\mapsto s(\varphi(a),y)=r(a,\psi(y)).
	\end{align*}
	On the other hand, for $f:\msf{A\longrightarrow \prm{A},g:B\longrightarrow\prm{B}}$ we have:
	\[f\multimap g:\msf{\prm{A}\multimap B\longrightarrow A\multimap\prm{B}},\] with
	\begin{align*}
	(f\multimap g)^+: B^{\prm{A}} \times_{\Gamma^{\prm{A}\times Y}} {\prm{X}}^Y &\longrightarrow {\prm{B}}^A \times_{\Gamma^{A\times\prm{Y}}} X^{\prm{Y}}\\
	\ag{h}{k} &\mapsto \ag{g^+hf^+}{f^-kg^-},
	\end{align*}
	and
	\[(f\multimap g)^-:\prm{A}\times\prm{Y}\longrightarrow A\times Y\]
	\[(f\multimap g)^-=f^+\times g^-.\]
	
	\rule{0pt}{1pt}
	\par \noindent \textbf{a5. Dualizing object.}\index{dualizing object} The dualizing object is \[\perp=\msf{I^*}=\agg{\Gamma}{\rho}{\{\bullet\}}\]
	with $\rho(c,\bullet)=c$ for all $c\in \Gamma$.\\
	\par \noindent \textbf{a6. Par.}\index{par connective} The par bifunctor $\odot$ is the De Morgan dual to $\otimes$:
	\[(-)\odot(-):\chu\times\chu\longrightarrow\chu.\]
	On objects,
	\begin{align*}
	\msf{A\odot B}&=\msf{(A^*\otimes B^*)^*}\\
	&=\agg{B^X\times_{\Gamma^{X\times Y}}A^Y}{v}{X\times Y}
	\end{align*}
	where 
	\begin{align*}
	v:(B^X\times_{\Gamma^{X\times Y}}A^Y) \times (X\times Y) &\longrightarrow \Gamma\\
	\langle \ag{\varphi}{\psi}~,~\ag{x}{y}\rangle~~ &\mapsto r(\psi(y),x)=s(\varphi(x),y).
	\end{align*}
	On morphisms $f:\msf{A\longrightarrow \prm{A}},g:\msf{B\longrightarrow\prm{B}}$,
	\[f\odot g:\msf{A\odot B\longrightarrow \prm{A}\odot\prm{B}},\]
	with
	\begin{align*}
	(f\odot g)^+:B^X\times_{\Gamma^{X\times Y}}A^Y &\longrightarrow {\prm{B}}^{\prm{X}}\times_{\Gamma^{\prm{X}\times \prm{Y}}}{\prm{A}}^{\prm{Y}}\\
	\ag{h}{k}&\mapsto \ag{g^+hf^-}{f^+kg^-},
	\end{align*}
	and
	\[(f\odot g)^-:\prm{X}\times \prm{Y}\longrightarrow X\times Y\]
	\[(f\odot g)^-=f^-\times g^-.\]
	
	\subsubsection*{b. Additive connectives}
	These include the \textit{plus} and \textit{with} connectives, together with their respective units.\\
	\par \noindent \textbf{b1. Plus.}\index{plus connective} Plus is the binary coproduct bifunctor
	\[(-)\oplus(-):\chu\times\chu\longrightarrow\chu\]
	\[(-)\oplus(-)\eqd\Clm(-,-).\]
	On objects:
	\[\msf{A\oplus B}=\agg{A\sqcup B}{p}{X\times Y},\]
	where $A\sqcup B$ is the coproduct (disjoint union) of the sets $A,B$; also
	\[p:(A\sqcup B)\times(X\times Y)\longrightarrow\Gamma\]
	\[p(m,\ag{x}{y})=\left\lbrace \begin{array}{lr}
	r(a,x)&\text{if~~}m\in A,\\
	s(b,y)&\text{if~~}m\in B.\\
	\end{array}\right. \]
	On morphisms $f:\msf{A\longrightarrow\prm{A}},g:\msf{B\longrightarrow\prm{B}}$,
	\[f\oplus g:\msf{A\oplus B\longrightarrow \prm{A}\oplus \prm{B}}\]
	\begin{align*}
	(f\oplus g)^+&=f^+\sqcup g^+,\\
	(f\oplus g)^-&=f^-\times g^-.
	\end{align*}
	
	\rule{0pt}{1pt}
	\par \noindent \textbf{b2. Plus unit.}\index{plus unit} Defined as:
	\[0\eqd\agg{\varnothing}{!}{\{\bullet\}},\]
	where $!$ is the unique empty matrix. $0$ has the property that \[0\oplus\msf{A}\cong\msf{A}\oplus 0\cong\msf{A}.\]
	
	\rule{0pt}{1pt}
	\par \noindent \textbf{b3. With.}\index{with connective} This is the binary product bifunctor
	\[(-)\&(-):\chu\times\chu\longrightarrow\chu\]
	\[(-)\&(-)\eqd\Lm(-,-),\]
	which is the De Morgan dual to plus:
	\[(-)\&(-)=((-)^*\oplus(-)^*)^*,~~~~~(-)\oplus(-)=((-)^*\&(-)^*)^*.\]
	On objects we have
	\[\msf{A\&B}=\agg{A\times B}{w}{X\sqcup Y}\] so that
	\[w:(A\times B)\times(X\sqcup Y)\longrightarrow\Gamma\]
	\[w(\ag{a}{b},n)=\left\lbrace \begin{array}{lr}
	r(a,n)&\text{if~~}n\in X,\\
	s(b,n)&\text{if~~}n\in Y.\\
	\end{array}\right. \]
	On morphisms $f:\msf{A\longrightarrow\prm{A}},g:\msf{B\longrightarrow\prm{B}}$,
	\[f\&g:\msf{A\&B\longrightarrow\prm{A}\&\prm{B}}\]
	\begin{align*}
	(f\&g)^+&=f^+\times g^+,\\
	(f\&g)^-&=f^-\sqcup g^-.
	\end{align*}
	
	\rule{0pt}{1pt}
	\par \noindent \textbf{b4. With unit.}\index{with unit} Defined as
	\[\top\eqd 0^*=\agg{\{\bullet\}}{!}{\varnothing},\]
	where $!$ is again the unique empty matrix. We have
	\[\msf{\top\&A\cong A\&\top\cong A.}\]
	
	\subsubsection*{c. Other connectives}
	Although \chu~is not cartesian closed, one can still define some sort of ``exponential'' on cartesian closed retractions of this category. The exponential serves syntactically to ``loosen up'' formulas so that they can be ``duplicated'' or ``deleted''. For more information on this subject, see \cite{Pratt-chu-sp}.
	
	
	\subsection{Realizations}
	\label{realize}
	The notion of realization we intend here is the strong version defined by Pultr and Trnkov\'{a} \cite{Pratt-chu-sp}.
	\begin{dfn}
		\label{normal-Chu}
		A Chu space $\agg{A}{r}{X}$ is called \textbf{normal}\index{normal Chu space} whenever $X\subseteq\Gamma^A$ and \(r(a,x)=x(a)\) for every $a\in A,x\in X$. For normal Chu spaces, the matrix $r$ is clearly understood from the context and therefore, the normal Chu space $\agg{A}{r}{X}$ is abbreviated as $\ag{A}{X}$. A \iw{normal Chu transform} is a Chu transform between two normal Chu spaces. Normal spaces and normal transforms form a full subcategory \(\cat{nChu}_\Gamma\) of $\chu_\Gamma$.
	\end{dfn}
	
	\begin{dfn}
		\label{representation}
		A functor $F:\mcg{C\longrightarrow D}$ is a \iw{representation} of objects $C\in\mcg{C}$ by objects \(F(C)\in\mcg{D}\) when $F$ is a full embedding (see Definition \ref{full-faithful-full-embed} and Proposition \ref{iffs-for-full-embed}).
	\end{dfn}
	
	\begin{dfn}
		\label{realization}
		A representation $F$ is a \iw{realization} when in addition, \[U_\mcg{D}F=U_\mcg{C},\] where \(U_\mcg{C}:\mcg{C}\longrightarrow\cat{Set}\) and \(U_\mcg{D}:\mcg{D}\longrightarrow\cat{Set}\) are the corresponding underlying-set functors.\\
		\par The above realization in called a \iw{Chu realization} when $\mcg{D}=\chu$, and is called a \iw{normal Chu realization} if $\mcg{D}=\cat{nChu}$ (where \cat{nChu} abbreviates $\cat{nChu}_\Gamma$).
	\end{dfn}
	
	\begin{rem}
		\label{every-realiz-conc}
		It is clear that every realization is a \textit{concrete functor} (see Chapter 1).
	\end{rem}
	
	Pratt \cite{Pratt-chu-sp,Pratt-rlz} discusses many normal Chu realizations including the realizations of \textit{sets}, \textit{pointed sets}, \textit{preorders}, \textit{topological spaces}, \textit{semilattices}, \textit{distributive lattices}, \textit{Boolean algebras}, etc. Because of their importance in illuminating the power of the Chu construction, we sketch two examples here, namely the cases of sets and topological spaces.\\
	\par Using the conventions $0=\varnothing,n=\{0,1,...,n-1\}$ we have:
	
	\begin{thm}
		\label{Set-nChu_2}
		$\cat{Set}$ is normally realized in $\cat{nChu}_2$.
	\end{thm}
	
	This realization is in fact trivial: any set $A$ is realized as the (normal) Chu space $\msf{A}=\agg{A}{!}{\varnothing}$, and any function $f:A\rightarrow B$ is realized as the (normal) Chu transform \(\msf{f}=\ag{f}{!}:\msf{A\longrightarrow B}\), where ``$!$'' is the unique map $!:\varnothing\rightarrow\varnothing$, doing the double duty of the empty matrix (in $\msf{A,B}$) as well as the empty map (in $\msf{f}$). Note that the composition \(\cat{Set}\rightarrow\cat{nChu}_2\hookrightarrow\chu_2\) yields the (strong monoidal) functor $i$ introduced in Proposition \ref{coreflective-reflective}.
	
	\begin{thm}
		\label{Top-enChu_2}
		A topological space can be viewed as an extensional normal Chu space in \(\cat{nChu}_2\), whose columns are closed under arbitrary union and finite intersection. The normal Chu transforms between topological spaces are exactly the continuous functions between the corresponding topological spaces.
	\end{thm}
	
	Theorem \ref{Top-enChu_2} tells us that normal Chu spaces and transforms (in $\cat{nChu}_2$) can in fact be regarded as some generalization of topological spaces and continuous functions, a generalization in which the conditions of open sets being closed under arbitrary unions and finite intersections are removed.\\
	\par Also:
	
	\begin{thm}
		\label{GCH}
		Assuming the Generalized Continuum Hypothesis, $\chu_\Gamma$ is realizable in \(\chu_{\prm{\Gamma}}\) if and only if $\left| \Gamma\right|\leq\left|\prm{\Gamma} \right|$.
	\end{thm}
	
	Another notable realization in Pratt's work is the normal Chu realization of relational structures:
	
	\begin{thm}
		\label{n-ary-relational}
		Any full subcategory $\mcg{C}$ of the category of n-ary relational structures and their homomorphisms is realized in \(\chu_{2^n}\).
	\end{thm}
	
	Consequently, for example, we have the next corollary:
	
	\begin{cor}
		\label{grps-and-top-grps-in-chu}
		\begin{enumerate}
			\item The category $\cat{Grp}$ of groups and group homomorphisms is realized in $\chu_{2^3}=\chu_8$.
			
			\item The category $\cat{TopGrp}$ of topological groups and continuous group homomorphisms is realized in $\chu_{2^4}=\chu_{16}$.
		\end{enumerate}
	\end{cor}
	
	Also we have:
	
	\begin{thm}
		\label{vect}
		The following functor realizes the category $\cat{Vect}_\Bbbk$ of vector spaces over a field $\Bbbk$ in the category \(\chu_{\left| \Bbbk\right| }\), where $\left| \Bbbk\right|$ denotes the underlying set of $\Bbbk$:
		\begin{align*}
		F:\cat{Vect}_\Bbbk&\longrightarrow\chu_{\left| \Bbbk\right| }\\
		\left( V\xrightarrow{T}W\right) &\mapsto\left( \agg{V}{e_V}{V^*}\xrightarrow{\left\langle T,T^*\right\rangle }\agg{W}{e_W}{W^*}\right) ,
		\end{align*}
		in which $T^*:W^*\rightarrow V^*$ is the dual map to $T$ in $\cat{Vect}_\Bbbk$, and the map $e_V$ (the evaluation) is defined as \(e_V(v_1,v_2^*)\eqd v_2^*(v_1)\in\Bbbk~\) for $v_1\in V,v_2^*\in V^*$.
	\end{thm}
	
	Finally, we mention a remarkable result Pratt has given in [...]:
	\begin{thm}
		\label{every-small-concrete-cat}
		Every small concrete category $\mcg{C}$ embeds fully in \(\chu_{\left| \mcg{C}\right| }\), where $\left| \mcg{C}\right|$ is the number of arrows of $\mcg{C}$.
	\end{thm}

	\subsection{Chu endofunctors}
	\label{endofunc}
	Because of their fundamental roles in universal dialgebra (Chapter 3) and later in the $\bb{DLC}$ construction (Chapter 4), Chu endofunctors deserve special attention. This subsection is meant to provide for a minimum necessary treatment of the issue.\\
	\par First of all, following Proposition \ref{coreflective-reflective}, we introduce a notation that will play a key role in the current thesis:
	
	\begin{ntn}
		\label{F+,F-}
		Let \(\mcg{V}\) be a \cs, and let $F:\cat{Chu}(\mcg{V},\Gamma)\longrightarrow\cat{Chu}(\mcg{V},\Gamma)$ be an endofunctor on $\cat{Chu}(\mcg{V},\Gamma)$. Also, let \(p_1:\cat{Chu}(\mcg{V},\Gamma)\longrightarrow\mcg{V}\) and \(p_2:\cat{Chu}(\mcg{V},\Gamma)\longrightarrow\op{\mcg{V}}\) be the projection functors deduced in the proof of Proposition \ref{coreflective-reflective}. We use the following notations in the sequel:
		\begin{align*}
		F^+&\eqd p_1\circ F:\cat{Chu}(\mcg{V},\Gamma)\longrightarrow\mcg{V},\\
		F^-&\eqd p_2\circ F:\cat{Chu}(\mcg{V},\Gamma)\longrightarrow\op{\mcg{V}}.
		\end{align*}
		In particular, when $\mcg{V}=\cat{Set}$, we have the functors $F^+:\chu_\Gamma\longrightarrow\cat{Set}$ and\\ $F^-:\chu_\Gamma\longrightarrow\op{\cat{Set}}$ for any given Chu endofunctor $F$.
	\end{ntn}
	
	Secondly, we introduce an important notion in the next definition:
	
	\begin{dfn}
		\label{uplifting}
		Let $F:\cat{Set}\longrightarrow \cat{Set}$ be an endofunctor on \cat{Set}. Let \(\Gamma\) be a nonempty set. An \iw{uplifting} of $F$ is a Chu endofunctor $\hat{F}:\chu_\Gamma\longrightarrow\chu_\Gamma$ that makes Diagram \ref{upl} commute. In the diagram, the functor $p_1$ is the same as in Notation \ref{F+,F-}.
		\begin{figure}[H]
			\begin{center}
				\begin{tikzpicture}[commutative diagrams/every diagram]
				
				\node (N1) at (2,2cm) {$\chu_\Gamma$};
				
				\node (N2) at (-2,2cm) {$\chu_\Gamma$};
				
				\node (N3) at (-2,-1cm) {$\cat{Set}$};
				
				\node (N4) at (2,-1cm) {$\cat{Set}$};
				
				\path[commutative diagrams/.cd, every arrow, every label]
				(N2) edge node{$\hat{F}$} (N1)
				(N2) edge node{$p_1$} (N3)
				(N3) edge node{$F$} (N4)
				(N1) edge node{$p_1$} (N4);
				\end{tikzpicture}
			\end{center}
			\caption{} \label{upl}
		\end{figure}
	\end{dfn}
	
	In other words, for an uplifting $\hat{F}$ we have \(\hat{F}^+\eqd p_1\circ\hat{F}=F\circ p_1\). The following theorem demonstrates the usefulness of the notion of uplifting.
	
	\begin{thm}
		\label{upliftable}
		Assume a nonempty set $\Gamma$. Then, every endofunctor \(F:\cat{Set}\longrightarrow\cat{Set}\) can be uplifted to a Chu endofunctor \[\hat{F}:\chu_\Gamma\longrightarrow\chu_\Gamma.\] Moreover, this uplifting can be done at least in two different ways.
	\end{thm}
	
	\begin{mypr}
		Fix a nonempty $\Gamma$. Consider the endofunctor\\ \(F:\cat{Set}\longrightarrow\cat{Set}\). We construct two distinct upliftings $\hat{F}_1,\hat{F}_2$ of $F$.\\
		\par \noindent \textbf{The first construction.} Define 
		
		\begin{align*}
		\hat{F}_1:\chu_\Gamma&\longrightarrow\chu_\Gamma\\
		\agg{A}{r}{X}&\mapsto\agg{F(A)}{!}{\varnothing}\\
		\left( \msf{A}\xrightarrow{\ag{f^+}{f^-}}\msf{B}\right) &\mapsto \left( F(\msf{A})\xrightarrow{\ag{F(f^+)}{!}}F(\msf{B})\right) 
		\end{align*}
		
		$\hat{F}_1$ is a well-defined functor because:
		\begin{align*}
		\forall \agg{A}{r}{X}\in\chu,~ \hat{F}_1(1_{\agg{A}{r}{X}})&=\hat{F}_1(\ag{1_A}{1_X})\\
		&=\ag{F(1_A)}{!}\\
		&=\ag{1_{F(A)}}{1_\varnothing}\\
		&=1_{\agg{F(A)}{!}{\varnothing}};
		\end{align*}
		also, for all consecutive arrows $\cdot\xrightarrow{f}\cdot\xrightarrow{g}\cdot$ in \chu,
		\begin{align*}
		\hat{F}_1(g\circ f)&=\hat{F}_1(\ag{g^+f^+}{f^-g^-})\\
		&=\ag{F(g^+f^+)}{!}\\
		&=\ag{F(g^+)\circ F(f^+)}{!\circ !}\\
		&=\ag{F(g^+)}{!}\circ\ag{F(f^+)}{!}\\
		&=\hat{F}_1(g)\circ\hat{F}_1(f).
		\end{align*}
		
		Now, for any Chu space $\msf{A}=\agg{A}{r}{X}$, \[\hat{F}_1^+(\msf{A})=p_1\circ\hat{F}_1(\msf{A})=p_1(\agg{F(A)}{!}{\varnothing})=F(A)=F\circ p_1(\msf{A});\]
		and for any Chu transform $f$,
		\[\hat{F}_1^+(f)=p_1\circ\hat{F}_1(f)=p_1(\ag{F(f^+)}{!})=F(f^+)=F\circ p_1(f).\]
		Therefore, $\hat{F}_1$ is an uplifting of $F$.
		~\\
		\par \noindent \textbf{The second construction.} Since $\Gamma$ is nonempty, there is some element $c\in\Gamma$. Let $\{\bullet\}$ denote the singleton. Define
		
		\begin{align*}
		\hat{F}_2:\chu_\Gamma&\longrightarrow\chu_\Gamma\\
		\agg{A}{r}{X}&\mapsto\agg{F(A)}{k_c}{\{\bullet\}},~~\text{s.t.~~}k_c(a,\bullet)=c,~\forall a\in A\\
		\left( \msf{A}\xrightarrow{\ag{f^+}{f^-}}\msf{B}\right) &\mapsto \left( F(\msf{A})\xrightarrow{\ag{F(f^+)}{1_{\{\bullet\}}}}F(\msf{B})\right) 
		\end{align*}
		
		Again, $\hat{F}_2$ is a well-defined functor: firstly, for every Chu transform $f$, the pair\\ $F(f)=\ag{F(f^+)}{1_{\{\bullet\}}}$ obviously satisfies the adjointness condition, and so, is a valid Chu morphism; also,
		\begin{align*}
		\forall \agg{A}{r}{X}\in\chu,~ \hat{F}_2(1_{\agg{A}{r}{X}})&=\hat{F}_2(\ag{1_A}{1_X})\\
		&=\ag{F(1_A)}{1_{\{\bullet\}}}\\
		&=\ag{1_{F(A)}}{1_{\{\bullet\}}}\\
		&=1_{\agg{F(A)}{k_c}{\{\bullet\}}};
		\end{align*}
		on the other hand, for all consecutive arrows $\cdot\xrightarrow{f}\cdot\xrightarrow{g}\cdot$ in \chu,
		\begin{align*}
		\hat{F}_2(g\circ f)&=\hat{F}_2(\ag{g^+f^+}{f^-g^-})\\
		&=\ag{F(g^+f^+)}{1_{\{\bullet\}}}\\
		&=\ag{F(g^+)\circ F(f^+)}{1_{\{\bullet\}}\circ 1_{\{\bullet\}}}\\
		&=\ag{F(g^+)}{1_{\{\bullet\}}}\circ\ag{F(f^+)}{1_{\{\bullet\}}}\\
		&=\hat{F}_2(g)\circ\hat{F}_2(f).
		\end{align*}
		
		Now, for any Chu space $\msf{A}=\agg{A}{r}{X}$, \[\hat{F}_2^+(\msf{A})=p_1\circ\hat{F}_2(\msf{A})=p_1(\agg{F(A)}{k_c}{\{\bullet\}})=F(A)=F\circ p_1(\msf{A});\]
		and for any Chu transform $f$,
		\[\hat{F}_2^+(f)=p_1\circ\hat{F}_2(f)=p_1(\ag{F(f^+)}{1_{\{\bullet\}}})=F(f^+)=F\circ p_1(f).\]
		Therefore, $\hat{F}_2$ is another uplifting of $F$.\\
		
		Finally, it is clear that $\hat{F}_1\ne\hat{F}_2$, and we are done.
	\end{mypr}
	
	\rule{0pt}{1pt}
	\begin{rem}
		\label{example-eqv-noneqv}
		In general, a set endofunctor $F:\cat{Set}\longrightarrow\cat{Set}$ may have many upliftings, quite different in structure and behavior. As a remarkable example, consider the endofunctor \((-)\times A:\cat{Set}\longrightarrow\cat{Set}\) for some fixed set $A$ with at least two distinct elements. This functor has at least the following upliftings:
		\begin{itemize}
			\item either of the constructions given in the proof of Theorem \ref{upliftable};
			\item the functor $(-)\otimes\agg{A}{r}{X}$ for some nontrivial $X,r$; and
			\item the functor $(-)\&\agg{A}{r}{X}$, again, for some nontrivial $X,r$.
		\end{itemize}
	\end{rem}
	
	\rule{0pt}{5mm}
	In addition to the notion of uplifting, it is good to point to the following notion:
	
	\begin{dfn}
		\label{bi-uplifting}
		Let $F:\cat{Set}\longrightarrow \cat{Set}$ and $G:\cat{Set}\longrightarrow \cat{Set}$ be endofunctors on \cat{Set}. Let \(\Gamma\) be a nonempty set. A \iw{bi-uplifting} of $F,G$ is a Chu endofunctor $H:\chu_\Gamma\longrightarrow\chu_\Gamma$ that makes both squares in Diagram \ref{bi-upl} commute. In the diagram, the functors $p_1,p_2$ are the same as in Notation \ref{F+,F-}.
		\begin{figure}[H]
			\begin{center}
				\begin{tikzpicture}[commutative diagrams/every diagram]
				
				\node (N1) at (0,-2cm) {$\chu_\Gamma$};
				
				\node (N2) at (0,1cm) {$\chu_\Gamma$};
				
				\node (N3) at (-3,1cm) {$\cat{Set}$};
				
				\node (N4) at (-3,-2cm) {$\cat{Set}$};
				
				\node (M3) at (3,1cm) {$\op{\cat{Set}}$};
				
				\node (M4) at (3,-2cm) {$\op{\cat{Set}}$};
				
				\path[commutative diagrams/.cd, every arrow, every label]
				(N2) edge node{$H$} (N1)
				(N2) edge node[swap]{$p_1$} (N3)
				(N3) edge node{$F$} (N4)
				(N1) edge node[swap]{$p_1$} (N4)
				(N2) edge node{$p_2$} (M3)
				(M3) edge node{$\op{G}$} (M4)
				(N1) edge node{$p_2$} (M4);			
				\end{tikzpicture}
			\end{center}
			\caption{} \label{bi-upl}
		\end{figure}
	\end{dfn}
	
	Every bi-uplifting $H$ of set functors $F,G$ is clearly an uplifting of $F$, also. On the other hand, as an example of a non-trivial bi-uplifting when $\Gamma=2$ and $F=G$ with $F$ a given set functor, see Definition 13 of \cite{Nabl} (termed ``$\msf{F}$-lifting'' therein). In general, characterization of bi-upliftings for arbitrary $\Gamma$ and arbitrary pairs of set functors $F,G$ may be an interesting subject of study by itself.\\
	
	\par In the next chapter, we will introduce the basics of universal dialgebra, with special attention to the theory of universal dialgebra on the Chu construction.

	\chapter{Universal Dialgebra}
	In Section \ref{funda-motivs} of Chapter 1 we discussed the historical origins of universal dialgebra. Now we proceed towards its formalism. Our approach here slightly generalizes that of \cite{Vout}; that is, we study universal dialgbera not on $ \cat{Set} $ but on a given (bicomplete) category $ \mcg{C} $. We will study isomorphisms and bisimulations in $ \dlg{\mcg{C}}{F}{G} $ (Section \ref{iso-bisim}); and we will take a look at limits and colimits in the latter category (Section \ref{lims-and-colims}).	After that, we will pay attention to the problem of dualization (Section \ref{the-prob-of-dualiz}) which we mentioned in Chapter 1. Finally, we will specifically focus on universal dialgebra on $ \chu $ (Section \ref{unidialg-on-chu}).
	
	\section{Basic notions}
	Many of the familiar structures in mathematics can be expressed in terms of ``$ F $-algebras'' and/or ``$ G $-coalgebras'', in which $ F,G $ are set functors. As a few examples:
	\begin{itemize}
		\item a \textit{semigroup} $ \ag{A}{*:A^2\longrightarrow A} $ may be viewed as an $ F $-algebra of the form $ \ag{A}{*:FA\longrightarrow A} $ in which $ F=(-)^2:\cat{Set}\longrightarrow\cat{Set} $;
		
		\item a \textit{partial function} $ \ag{B}{p:B\longrightarrow B\sqcup \{\bullet\}} $ may be viewed as a $ G $-coalgebra of the form $ \ag{B}{p:B\longrightarrow GB} $ where $ G=(-)\sqcup\{\bullet\}:\cat{Set}\longrightarrow\cat{Set} $.
	\end{itemize}
	On the other hand, the above examples can also be viewed as ``$ \ag{F}{G} $-dialgebras'' such that $ F,G $ are again set functors (assuming $ G=1_{\cat{Set}} $ for the first example and $ F=1_{\cat{Set}} $ for the second). As a non-trivial example of $ \ag{F}{G} $-dialgebras, one may take a \textit{field}
	\[\langle\Bbbk,+,0,-(-),\times,1,(-)^{-1}\rangle.\]
	This is an $ \ag{F}{G} $-dialgebra for some non-identity set functors $ F,G $.
	\par In addition to the above, one may think of generalized $ \ag{F}{G} $-dialgebras for endofunctors
	\[F,G:\mcg{C}\dlong\mcg{C},\]
	where $ \mcg{C} $ is an arbitrary category. In the following, we study such formalisms.
	
	\begin{dfn}
		\label{F,G-dialg}
		Let $\mcg{C}$ be a category. Let $F,G:\mcg{C}\dlong\mcg{C}$ be two endofunctors on $\mcg{C}$. An $\langle F,G\rangle$-\iw{dialgebra}\index{$\langle F,G\rangle$-dialgebra} \(\eb{A}\) is a pair \(\ag{A}{\alpha}\) consisting of a $\mcg{C}$-object $A$, together with a $\mcg{C}$-morphism \(\alpha:FA\longrightarrow GA\). Whenever the functors $F,G$ are clear from the context, we may simply refer to $\eb{A}$ as a ``dialgebra''.
	\end{dfn}
	
	\begin{dfn}
		\label{dialg-homo}
		Let $\eb{A}=\ag{A}{\alpha},\eb{B}=\ag{B}{\beta}$ be two $\ag{F}{G}$-dialgebras. An\\ $\ag{F}{G}$-\iw{dialgebra homomorphism} \(\eb{h:A\longrightarrow B}\) is a $\mcg{C}$-arrow \(h:A\longrightarrow B\) making Diagram \ref{dialg-homo-diag} commute.
		\begin{figure}[H]
			\begin{center}
				\begin{tikzpicture}[commutative diagrams/every diagram]
				
				\node (N1) at (3,2cm) {$GA$};
				
				\node (N2) at (-1,2cm) {$FA$};
				
				\node (N3) at (-1,-1cm) {$FB$};
				
				\node (N4) at (3,-1cm) {$GB$};
				
				\node (N5) at (-3,2cm) {$A$};
				
				\node (N6) at (-3,-1cm) {$B$};
				
				\path[commutative diagrams/.cd, every arrow, every label]
				(N2) edge node{$\alpha$} (N1)
				(N2) edge node{$Fh$} (N3)
				(N3) edge node{$\beta$} (N4)
				(N1) edge node{$Gh$} (N4)
				(N5) edge node{$h$} (N6);
				\end{tikzpicture}
			\end{center}
			\caption{} \label{dialg-homo-diag}
		\end{figure}
	\end{dfn}
	
	\begin{rem}
		\label{stem-yield}
		The arrow $h:A\longrightarrow B$ is called the \iw{stem} of the arrow \(\eb{h:A\longrightarrow B}\), or equivalently, one may say that the homomorphism $\eb{h}$ \textbf{stems from} $h$ (via $F,G$). Yet another equivalent statement is that $h$ \textbf{yields} a homomorphism $\eb{h}$ (via $F,G$). We use boldface characters to represent the corresponding homomorphisms of the stems.
	\end{rem}
	
	Identity homomorphisms are $\ag{F}{G}$-dialgebra homomorphisms and, given two dialgebra homomorphisms \(\eb{f:A\longrightarrow B}\) and \(\eb{g:B\longrightarrow C}\), the composition $g\circ f:A\longrightarrow C$ yields a dialgebra homomorphism \(\eb{g\circ f:A\longrightarrow C}\) (Diagram \ref{dialg-homo-compos}). This composition is clearly associative, also. Therefore:
	
	\begin{figure}[H]
		\begin{center}
			\begin{tikzpicture}[commutative diagrams/every diagram]
			
			\node (N1) at (4,1.5cm) {$FC$};
			
			\node (N2) at (0,1.5cm) {$FB$};
			
			\node (N3) at (-4,1.5cm) {$FA$};
			
			\node (N4) at (-4,-1cm) {$GA$};
			
			\node (N5) at (0,-1cm) {$GB$};
			
			\node (N6) at (4,-1cm) {$GC$};
			
			\node (N7) at (4,3cm) {$C$};
			
			\node (N8) at (0,3cm) {$B$};
			
			\node (N9) at (-4,3cm) {$A$};
			
			\path[commutative diagrams/.cd, every arrow, every label]
			(N2) edge node{$Fg$} (N1)
			(N3) edge node{$Ff$} (N2)
			(N3) edge node{$\alpha$} (N4)
			(N4) edge node{$Gf$} (N5)
			(N5) edge node{$Gg$} (N6)
			(N1) edge node{$\gamma$} (N6)
			(N2) edge node{$\beta$} (N5)
			(N8) edge node{$g$} (N7)
			(N9) edge node{$f$} (N8);
			\end{tikzpicture}
		\end{center}
		\caption{} \label{dialg-homo-compos}
	\end{figure}
	
	\begin{prp}
		\label{C^F_G-prp}
		Given endofunctors \(F,G:\mcg{C\dlong C}\) on some category $\mcg{C}$, $\ag{F}{G}$-dialgebras together with $\ag{F}{G}$-dialgebra homomorphisms form a category.
	\end{prp}
	
	\begin{ntn}
		\label{C^F_G-base}
		This category is denoted by $\mcg{C}^F_G$. We may refer to $\mcg{C}$ as the \iw{base category} for $\mcg{C}^F_G$. To simplify notation for special cases, let $\mcg{C}^F\eqd \mcg{C}^F_{1_\mcg{C}}$ and $\mcg{C}_G\eqd \mcg{C}_G^{1_\mcg{C}}$.
	\end{ntn}
	
	\begin{rem}
		\label{special-cases-alg-coalg}
		$\mcg{C}^F$ is essentially the same thing as the category of $F$-algebras over $\mcg{C}$ while $\mcg{C}_G$ is just the category of $G$-coalgebras over $\mcg{C}$. This way, universal algebra and universal coalgebra may be viewed as two special cases of universal dialgebra.
	\end{rem}
	
	\par The reader may refer to \cite{Blok,ModAr,Palm,Poll,Vout} for various examples of categories of dialgebras. For algebras and coalgebras, references \cite{Burr,Dene,Grat,Hugh,Jaco,Rutt} are suggested. Also, Abramsky \cite{Abram-Coalg} gives an interesting example of connections between coalgebras and the Chu construction. Chang and Keisler \cite{Chang-Mod} point to the connection between universal algebra and model theory. Pavlovi\'{c} and Pratt \cite{Conti-final} introduce the continuum $ \bb{R} $ of real numbers as a final coalgebra (i.e. as a terminal object in the category of coalgebras). Gumm \cite{Gumm-FunCo} studies functors for coalgebras. Moss \cite{Moss} studies coalgebraic logic. Palmigiano \cite{Palm} introduces abstract logics as dialgebras. Martins et al \cite{ModAr} suggest dialgebras as appropriate models for computational processes which are combinations of ``algebraic construction'' and ``coalgebraic observation''. Finally, Rodrigues \cite{Rod} develops a dialgebraic theory over the category $ \cat{Rel} $ (the category of sets and binary relations).
	
	\section{Isomorphisms and bisimulations}
	\label{iso-bisim}
	In this section, $\mcg{C}$ is a fixed category and $F,G$ are two fixed $\mcg{C}$-endofunctors.
	
	\begin{dfn}
		\label{dialg-iso}
		A dialgebra homomorphism \(\eb{h:A\longrightarrow B}\) is called a \iw{dialgebra isomorphism} if it is an isomorphism arrow in the category $\dlg{\mcg{C}}{F}{G}$.
	\end{dfn}
	
	The following proposition shows that the existence of an inverse homomorphism is equivalent to the existence of an inverse for the stem.
	
	\begin{prp}
		\label{h-iso-conversely}
		\begin{enumerate}
			\item If $\eb{h}$ is an isomorphism in $\dlg{\mcg{C}}{F}{G}$, then its stem $h:A\longrightarrow B$ is an isomorphism in $\mcg{C}$.
			
			\item Conversely, if $h$ has an inverse $k:B\longrightarrow A$, then $k$ yields a dialgebra homomorphism $\eb{k:B\longrightarrow A}$ such that $\eb{k}=\eb{h}^{-1}$. That is, $\eb{h}$ is a dialgebra isomorphism.\\
		\end{enumerate}
		
		\par Therefore, a dialgebra homomorphism is a dialgebra isomorphism if and only if its stem is an isomorphism arrow in the corresponding base category.
	\end{prp}
	
	\begin{mypr}
		(1) is obvious. We only prove (2):\\
		\par Consider
		\begin{center}
			\mbox{
				\begin{tikzpicture}[commutative diagrams/every diagram]
				
				\node (N1) at (3,2cm) {$GA$};
				\node (V1) at (2.9,1.8cm) {};
				\node (V2) at (3.1,1.8cm) {};
				
				\node (N2) at (-1,2cm) {$FA$};
				\node (V3) at (-1.1,1.8cm) {};
				\node (V4) at (-0.9,1.8cm) {};
				
				\node (N3) at (-1,-1cm) {$FB$};
				\node (V5) at (-1.1,-0.8cm) {};
				\node (V6) at (-0.9,-0.8cm) {};
				
				\node (N4) at (3,-1cm) {$GB$};
				\node (V7) at (2.9,-0.8cm) {};
				\node (V8) at (3.1,-0.8cm) {};
				
				\node (N5) at (-3,2cm) {$A$};
				\node (V9) at (-3.1,1.8cm) {};
				\node (V10) at (-2.9,1.8cm) {};
				
				\node (N6) at (-3,-1cm) {$B$};
				\node (V11) at (-3.1,-0.8cm) {};
				\node (V12) at (-2.9,-0.8cm) {};
				
				\path[commutative diagrams/.cd, every arrow, every label]
				(N2) edge node{$\alpha$} (N1)
				(N3) edge node{$\beta$} (N4)
				(V9) edge node[swap]{$h$} (V11)
				(V12) edge node[swap]{$k$} (V10)
				(V3) edge node[swap]{$Fh$} (V5)
				(V6) edge node[swap]{$Fk$} (V4)
				(V1) edge node[swap]{$Gh$} (V7)
				(V8) edge node[swap]{$Gk$} (V2);
				\end{tikzpicture}
			}
		\end{center}
		We have $hk=1_B, kh=1_A$; thus
		\begin{align*}
		G(k)\beta&=G(k)\beta F(h)F(k)\\
		&=G(k)G(h)\alpha F(k)\\
		&=\alpha F(k).
		\end{align*}
		Therefore, $\eb{k:B\longrightarrow A}$ is a homomorphism, and clearly, $\eb{k}=\eb{h}^{-1}$.
	\end{mypr}
	
	\begin{dfn}
		\label{dialg-mono}
		A homomorphism $\eb{h:A\longrightarrow B}$ is called a \textbf{(dialgebra) monomorphism}\index{dialgebra monomorphism} if it is a monic arrow in $\dlg{\mcg{C}}{F}{G}$, i.e., if for every dialgebra $\eb{C}=\ag{C}{\gamma}$ and all homomorphisms \(\eb{f,g:C\dlong A}\), \[\eb{hf=hg}~~~~\text{implies}~~~~\eb{f=g}.\]
		Dually, $\eb{h}$ is said to be an(a) (dialgebra) epimorphism whenever $\eb{h}$ is an epic arrow in $\dlg{\mcg{C}}{F}{G}$.
	\end{dfn}
	
	\begin{prp}
		\label{stem-monic-epic}
		\begin{enumerate}
			\item Every dialgebra homomorphism stemming from a monic arrow is monic.
			\item Dually, every dialgebra homomorphism stemming from an epic arrow is epic.
		\end{enumerate}
	\end{prp}
	
	\begin{mypr}
		(1)  For a homomorphism $\eb{h}$ with monic stem $h$, the equality $\eb{hf=hg}$ for some homomorphisms $\eb{f,g}$ implies $hf=hg$, hence $f=g$, and consequently $\eb{f=g}$.
	\end{mypr}
	
	\begin{rem}
		Contrary to the case of Proposition \ref{h-iso-conversely}, the converse statements for Proposition \ref{stem-monic-epic} are false: there are counterexamples monomorphisms (epimorphisms) $\eb{h}$ for which the stem $h$ is \textit{not} a monic (epic) arrow in the base category $\mcg{C}$. For instance, take the category $\cat{Ring}=\cat{Set}^\mathfrak{R}$ of rings, where $\mathfrak{R}$ is the correponding endofunctor. Here we have an epimorphism $\eb{i}:\ag{\bb{Z}}{\fk{R}\bb{Z}\rightarrow\bb{Z}}\longrightarrow\ag{\bb{Q}}{\fk{R}\bb{Q}\rightarrow\bb{Q}}$; however, its stem is the inclusion map $i:\bb{Z\hookrightarrow Q}$ which is obviously not surjective (hence not epic) in \cat{Set}.
	\end{rem}
	
	\begin{dfn}
		\label{dialg-section-retrac}
		A homomorphism $\eb{h:A\longrightarrow B}$ is called a \textbf{(dialgebra) section}\index{dialgebra section} provided that it is a section arrow in $\dlg{\mcg{C}}{F}{G}$, i.e., if there exists a homomorphism $\eb{k:B\longrightarrow A}$ with $\eb{kh}=1_\eb{A}$.
		\par Dually, $\eb{h}$ is a \textbf{(dialgebra) retraction}\index{dialgebra retraction} if it is a retraction arrow in $\dlg{\mcg{C}}{F}{G}$, that is, if there exists $\eb{k}$ with $\eb{hk}=1_\eb{B}$.
	\end{dfn}
	
	\begin{prp}
		\label{section-retrac-stem}
		\begin{enumerate}
			\item If $\eb{h}$ is a section (retraction) in $\dlg{\mcg{C}}{F}{G}$, then its stem $h:A\longrightarrow B$ is also a section (retraction) in $\mcg{C}$.
			\item If $k:B\longrightarrow A$ is a $\mcg{C}$-arrow that makes $h$ a section (retraction), and if a homomorphism $\eb{k:B\longrightarrow A}$ stems from $k$, then the homomorphism $\eb{h}$ is a section (retraction) as well.
		\end{enumerate}
	\end{prp}
	
	\begin{dfn}
		\label{subdialg-quotient}
		A dialgebra $\eb{A}=\ag{A}{\alpha}$ is called a \iw{subdialgebra} of a dialgebra $\eb{B}=\ag{B}{\beta}$ if there exists a homomorphism $\eb{m:A\longrightarrow B}$ that stems from a monic $m:A\rightarrowtail B$. By Proposition \ref{stem-monic-epic} (1) we know that such $\eb{m}$ must be a dialgebra monomorphism.
		\par Dually, $\eb{B}$ is called a \iw{quotient dialgebra} of $\eb{A}$ provided that there exists a homomorphism $\eb{q:A\longrightarrow B}$ such that its stem $q:A\twoheadrightarrow B$ is epic. Again, from Proposition \ref{stem-monic-epic} (2) we know that such $\eb{q}$ is necessarily a dialgebra epimorphism.
	\end{dfn}
	
	\rule{0pt}{0.5cm}
	Now we turn to the fundamental concept of ``bisimulation'' for dialgebras. In \cite{Vout}, the concept was defined and investigated for the case $\mcg{C}=\cat{Set}$. Here, we generalize the concept for categories with finite products, possibly non-cartesian-closed cases. Before proceeding we need a categorical definition of the concept of ``reltion''. The following material is a generalization of what can be found in \cite{Borc2}.
	\begin{dfn}
		\label{cat-relation}
		Let $\mcg{C}$ be a category. By a \iw{relation} between objects $A,B\in\mcg{C}$ we mean an object \(R\in\mcg{C}\) together with a pair of arrows \[A\xleftarrow{~\pi_1~}R\xrightarrow{~\pi_2~}B,\]
		which form a \iw{monomorphic pair of arrows}, that is, given an object $X$ and parallel arrows $f,g:X\dlong R$, 
		\[f=g ~~~~~\text{iff}~~~~~ \pi_if=\pi_ig, ~i=1,2.\]
		For every object $X\in\mcg{C}$ we write
		\[R_X\eqd\{\ag{\pi_1x}{\pi_2x}\mid x\in\mcg{C}(X,R)\}\subseteq\mcg{C}(X,A)\times\mcg{C}(X,B)\] 
		which is a relation in the usual sense (i.e. set-theoretic), and call it the \iw{corresponding relation} between $\mcg{C}(X,A)$ and $\mcg{C}(X,B)$. Also, by a \textbf{relation on} $A$ we mean a relation between $A$ and itself.
	\end{dfn}
	
	Recall that an ordinary (set-thoretic) relation $R$ between sets $A,B$ was defined as a \textit{subset} $R\subseteq A\times B$ of the (cartesian) product of $A$ and $B$. In category-theoretic terms, this might be expressed as an injective arrow $R\rightarrowtail A\times B$ in the category \cat{Set}. The following proposition generalizes this idea to categories with finite products:
	
	\begin{prp}
		\label{relation-as-monic}
		Let $\mcg{C}$ be a category with finite products, and let $A,B\in \mcg{C}$. Then, every relation $\agg{R}{\pi_1}{\pi_2}$ between $A,B$ can be viewed as a monic arrow \[\varphi:R\rightarrowtail A\times B\] and vice versa.
	\end{prp}
	
	\begin{mypr}
		($\Longrightarrow$)  Let \(A\xleftarrow{~p_1~}A\times B\xrightarrow{~p_2~}B\) be the product of $A$ and $B$. Then by universality of the product, there exists a unique arrow $\varphi:R\longrightarrow A\times B$ such that \[\pi_i=p_i\varphi,~i=1,2.\]
		To see $\varphi$ being monic, assume parallel arrows $f,g:X\dlong R$ with \(\varphi f=\varphi g\). Then
		\[p_i\varphi f=p_i\varphi g,\] or
		\[\pi_i f=\pi_i g\] for $i=1,2$. But since $\pi_i$ are a monomorphic pair of arrows, we have $f=g$.\\
		\par ($\Longleftarrow$)  Conversely, assume that the product $\agg{A\times B}{p_1}{p_2}$ and a mono $\varphi:R\rightarrowtail A\times B$ are given. Define \(\pi_i\eqd p_i\varphi\) for each $i$. For any parallel pair $f,g:X\dlong R$, the equalities \(\pi_i f=\pi_i g,~\forall i~\) imply
		\[p_i\varphi f=p_i\varphi g,~ \forall i.\]
		Since $\agg{A\times B}{p_1}{p_2}$ is itself a relation between $A,B$, (because $p_1,p_2$ form a monomorphic family), we must have \[\varphi f=\varphi g.\]
		By hypothesis, $\varphi$ is monic; thus $f=g$.
	\end{mypr}
	
	~\\
	\par Given a relation $\agg{R}{\pi_1}{\pi_2}$ on an object $A$ of some category $\mcg{C}$, it is now possible to require the classical properties on the various relations $R_X$ on the sets $\mcg{C}(X,R)$.
	
	\begin{dfn}
		\label{equiv-R_X}
		By an \iw{equivalence relation} on an object $A$ of a category $\mcg{C}$ we mean a relation $\agg{R}{\pi_1}{\pi_2}$ on $A$ such that for every object $X\in \mcg{C}$, the corresponding relation $R_X$ on the set $\mcg{C}(X,A)$ is an equivalence relation. More generally, the relation $R$ is \textbf{reflexive} (resp. \textbf{symmetric, transitive, antisymmetric, etc.}) when each corresponding realtion $R_X$ is so.
	\end{dfn}
	
	In \cat{Set}, given an equivalence relation $R\subseteq A\times A$, one can take the \textit{quotient} of $A$ by $R$, arriving at a diagram \[R\raisebox{-1pt}{$\def\arraystretch{0.5}\begin{array}{c} \xrightarrow{~\pi_1~}\\ \xrightarrow[~\pi_2~]{} \end{array}$}A\xrightarrow{~q~}A/R.\]
	The coequalizer of $\pi_i$ is the quotient of $A$ by the equivalence relation generated by the pairs \(\ag{\pi_1(x)}{\pi_2(x)},x\in R\), which is $q$. On the other hand, \(q(a)=q(\prm{a})\) iff \(\ag{a}{\prm{a}}\in R\), which indicates that $\ag{\pi_1}{\pi_2}$ is the kernel pair of $q$.
	
	\begin{dfn}
		\label{eqv-effective}
		An equivalence relation $\agg{R}{\pi_1}{\pi_2}$ on an object $A$ of a category $\mcg{C}$ is \iw{effective} when the coequalizer $q$ of $\ag{\pi_1}{\pi_2}$ exists, and, $\ag{\pi_1}{\pi_2}$ is the kernel pair of $q$.
	\end{dfn}
	
	For more on categorical equivalence relations, see \cite{Borc2}.\\
	\par Now we return to dialgebras and bisimulations:
	
	\begin{dfn}
		\label{F,G-bisimulation}
		Let $\eb{A}=\ag{A}{\alpha},\eb{B}=\ag{B}{\beta}$ be two dialgebras in $\dlg{\mcg{C}}{F}{G}$. An $\ag{F}{G}$-\iw{bisimulation}\index{$\ag{F}{G}$-bisimulation} between $\eb{A}$ and $\eb{B}$ is a dialgebra $\eb{R}=\ag{R}{\rho:FR\rightarrow GR}$, where \(\agg{R}{\pi_1}{\pi_2}\) is a relation between $A,B$, such that both \(\boldsymbol{\uppi}_1:\eb{R\longrightarrow A}\) and \(\boldsymbol{\uppi}_2:\eb{R\longrightarrow B}\) are dialgebra homomorphisms. In other words, Diagram \ref{F,G-bisimulation-diag} commutes.\\
		\begin{figure}[H]
			\begin{center}
				\begin{tikzpicture}[commutative diagrams/every diagram]
				
				\node (N1) at (3,2cm) {$FB$};
				
				\node (N2) at (0,2cm) {$FR$};
				
				\node (N3) at (0,-1cm) {$GR$};
				
				\node (N4) at (3,-1cm) {$GB$};
				
				\node (N5) at (-3,2cm) {$FA$};
				
				\node (N6) at (-3,-1cm) {$GA$};
				
				\path[commutative diagrams/.cd, every arrow, every label]
				(N2) edge node{$F\pi_2$} (N1)
				(N2) edge node{$\rho$} (N3)
				(N3) edge node{$G\pi_2$} (N4)
				(N1) edge node{$\beta$} (N4)
				(N5) edge node{$\alpha$} (N6)
				(N2) edge node[swap]{$F\pi_1$} (N5)
				(N3) edge node[swap]{$G\pi_1$} (N6);
				\end{tikzpicture}
			\end{center}
			\caption{} \label{F,G-bisimulation-diag}
		\end{figure}
		\par We omit the functors $F,G$ and use the term \textbf{bisimulation} whenever $F$ and $G$ are clear from the context. Also, an $\ag{F}{G}$-\textbf{bisimulation} on $\eb{A}$ is an $\ag{F}{G}$-bisimulation between $\eb{A}$ and itself. Finally, an $\ag{F}{G}$-\iw{bisimulation equivalnce} on $\eb{A}$ is a bisimulation on $\eb{A}$ that is also an equivalence relation on $\eb{A}$.
	\end{dfn}
	
	Next, another important concept is introduced. Following McLarty \cite{Lart}:
	
	\begin{dfn}
		\label{graph-arr}
		Suppose that $\mcg{C}$ has finite products, and take any arrow $f:A\longrightarrow B$. The \iw{graph arrow} of  $f$, denoted by $\Gr_f$, is the arrow that makes Diagram \ref{graph-arr-diag} commute.\\
		\begin{figure}[H]
			\begin{center}
				\begin{tikzpicture}[commutative diagrams/every diagram]
				
				\node (N2) at (0,1.5cm) {$A$};
				
				\node (N3) at (0,-1cm) {$A\times B$};
				
				\node (N4) at (3,-1cm) {$B$};
				
				\node (N6) at (-3,-1cm) {$A$};
				
				\path[commutative diagrams/.cd, every arrow, every label]
				(N2) edge node{$\Gr_f$} (N3)
				(N2) edge node[swap]{$1_A$} (N6)
				(N2) edge node{$f$} (N4)
				(N3) edge node[swap]{$p_1$} (N6)
				(N3) edge node{$p_2$} (N4);
				\end{tikzpicture}
			\end{center}
			\caption{} \label{graph-arr-diag}
		\end{figure}
	\end{dfn}
	
	\begin{prp}
		\label{Gr-eq}
		$\Gr_f$ as defined above is an equalizer for $fp_1$ and $p_2$.
	\end{prp}
	
	\begin{mypr}
		We have \[fp_1\Gr_f=f1_a=f=p_2\Gr_f.\] 
		Also, for any other $\gamma:C\longrightarrow A\times B$ with $fp_1\gamma=p_2\gamma$, consider the composition \[p_1\gamma:C\longrightarrow A.\]
		We have \[\left\lbrace \begin{array}{l}
		p_1(\Gr_fp_1\gamma)=(p_1\Gr_f)p_1\gamma=p_1\gamma,\\ 
		p_2(\Gr_fp_1\gamma)=fp_1\gamma=p_2\gamma, \end{array}\right.\]
		and then by universality of the product $A\times B$,
		\[\Gr_fp_1\gamma=\gamma.\]
		Finally, suppose that there is some other arrow $\theta:C\longrightarrow A$ with $\Gr_f\theta=\gamma$. Then \[\theta=p_1\Gr_f\theta=p_1\gamma,\] so $p_1\gamma$ is the unique arrow with the above property. This completes the proof.
	\end{mypr}
	
	\begin{cor}
		\label{Gr-monic}
		The graph arrow $\Gr_f$ defined above is always monic.
	\end{cor}
	
	Consequently, by Proposition \ref{relation-as-monic}:
	
	\begin{cor}
		\label{Gr-relation}
		For any $f$, the graph arrow $\Gr_f:A\rightarrowtail A\times B$ is a relation $\agg{A}{1_A}{f}$ between $A,B$.
	\end{cor}
	
	In addition, we have:
	
	\begin{prp}
		\label{diagonal-Gr}
		For any $A$ with a product with itself, the diagonal \(\triangle_A:A\rightarrowtail A\times A\) equals the graph arrow of $1_A$: \[\triangle_A=\Gr_{1_A}.\]
	\end{prp}
	
	\rule{0pt}{0.5cm}
	Now, equipped with all the above, we state and prove the following theorem, which generalizes Theorem 9 in \cite{Vout}.
	
	\begin{thm}
		\label{homo-iff-Gr-bisim}
		Assume that $\mcg{C}$ has finite products, and let $\eb{A}=\ag{A}{\alpha},\eb{B}=\ag{B}{\beta}$ be two $\ag{F}{G}$-dialgebras in \(\dlg{\mcg{C}}{F}{G}\). Then an arrow $f:A\longrightarrow B$ yields a dialgebra homomorphism if and only if its graph arrow $\Gr_f$ induces an $\ag{F}{G}$-bisimulation.
	\end{thm}
	
	\begin{mypr}
		By Corollary \ref{Gr-relation} we know that for any arrow $f:A\longrightarrow B$, the graph arrow $\Gr_f$ is equivalent to a relation \(\agg{A}{1_A}{f}\) between $A,B$. Since $\eb{1}_A$ is always a dialgebra homomorphism, it obviously follows that $f$ yields a homomorphism if and only if \(\agg{A}{1_A}{f}\) yields a bisimulation between $\eb{A,B}$.
	\end{mypr}
	
	\begin{cor}
		\label{diagonal-bisim-eqv}
		Let $\eb{A}=\ag{A}{\alpha}$ be a dialgebra. The diagonal $\triangle_A$ of $A$ yields a bisimulation equivalence on $\eb{A}$.
	\end{cor}
	
	\begin{mypr}
		Use Proposition \ref{diagonal-Gr} together with Theorem \ref{homo-iff-Gr-bisim}.
	\end{mypr}
	\vspace{5mm}
	\section{Limits and colimits}
	\label{lims-and-colims}
	Now let $\mcg{C}$ be a bicomplete category. In this section, a study of limits and colimits in the category $\dlg{\mcg{C}}{F}{G}$ is undertaken. The ultimate goal is to prove that $\dlg{\mcg{C}}{F}{G}$ has all limits that are preserved by G, and all colimits that are preserved by $F$.
	\par Recall that in the category theory literature, especially in the study of limits and colimits, a functor $d:\mcg{D\longrightarrow C}$ is sometimes called a \iw{diagram}, with the category $ \mcg{D} $ being called the \iw{scheme} of $ d $ \cite{Adam}. Now we sate another closely-related definition \cite{Vout}.
	
	\begin{dfn}
		\label{type-lim-colim}
		By the \iw{type} of a limit or colimit, we mean the isomorphism class, in the category \cat{Grph} of graphs and graph homomorphisms, of the base graph of the limit or colimit, respectively. A functor $H:\mcg{C\longrightarrow C}$ is said to \iw{preserve} a type $D$ of limit (colimit), if for every diagram $d:D\longrightarrow \mcg{C}$ with limit (colimit) \(\ag{L}{\{f_v \mid v\in V(D)\}}\),
		\[\ag{H(L)}{\{H(f_v)\mid v\in V(D)\}}\] is the limit (colimit) of the diagram $H\circ d:D\longrightarrow\mcg{C}$.
	\end{dfn}
	
	Let $U:\dlg{\mcg{C}}{F}{G}\longrightarrow\mcg{C}$ be the forgetful functor from the category of $\ag{F}{G}$-dialgebras to $\mcg{C}$. For every dialgebra \(\eb{A}=\ag{A}{\alpha}\), $U(\eb{A})=A$ and, for every dialgebra homomorphism $\eb{h}$, \(U(\eb{h})=h\). For this functor we consider the following definition.
	
	\begin{dfn}
		\label{create-type-lim}
		The functor $U$ is said to \iw{create a type} $D$ \textbf{limit (colimit)} if for every diagram $d:D\longrightarrow \dlg{\mcg{C}}{F}{G}$, its limit (colimit) is constructed by first taking the limit (colimit) of \(U\circ d:D\longrightarrow\mcg{C}\) in $\mcg{C}$, and then supplying it in a unique way with an $\ag{F}{G}$-dialgebra structure.
	\end{dfn}
	
	\begin{thm}
		\label{forget-crea-pres-G-F}
		Assume that $\mcg{C}$ is bicomplete. Then:
		\begin{enumerate}
			\item The forgetful functor $U:\dlg{\mcg{C}}{F}{G}\longrightarrow\mcg{C}$ creates and preserves all types of limits that the functor $G:\mcg{C}\longrightarrow\mcg{C}$ preserves.
			
			\item The forgetful functor $U:\dlg{\mcg{C}}{F}{G}\longrightarrow\mcg{C}$ creates and preserves all types of colimits that the functor $F:\mcg{C}\longrightarrow\mcg{C}$ preserves.
		\end{enumerate}
	\end{thm}
	
	\begin{mypr}
		We only prove part (1). Part (2) is automatically proved by dualization.\\
		\par Let $D=\ag{V(D)}{E(D)}$ be a graph and $d:D\longrightarrow\dlg{\mcg{C}}{F}{G}$ a diagram in \(\dlg{\mcg{C}}{F}{G}\). Suppose that $G$ preserves limits of type $D$. Consider the diagram $U\circ d:D\longrightarrow\mcg{C}$. Since $\mcg{C}$ is complete, $Ud$ has a limit \(\ag{L}{\{l_v:L\rightarrow Ud(v)\mid v\in V(D)\}}\) in $\mcg{C}$, such that for all $e:v_1\longrightarrow v_2$ in $E(D)$, the following diagram commutes:
		\begin{center}
			\mbox{
				\begin{tikzpicture}[commutative diagrams/every diagram]
				
				\node (N1) at (0,2cm) {$L$};
				
				\node (N2) at (-2,0cm) {$Ud(v_1)$};
				
				\node (N3) at (2,0cm) {$Ud(v_2)$};
				
				\path[commutative diagrams/.cd, every arrow, every label]
				(N1) edge node[swap]{$l_{v_1}$} (N2)
				(N1) edge node{$l_{v_2}$} (N3)
				(N2) edge node{$Ud(e)$} (N3);
				\end{tikzpicture}
			}
		\end{center}
		and such that, for all other cones \(\ag{M}{\{f_v:M\rightarrow Ud(v)\mid v\in V(D)\}}\) with analogous commutativity conditions as above, there exists a unique morphism $f:M\longrightarrow L$ such that the following triangle commutes for all $v\in V(D)$:
		\begin{center}
			\mbox{
				\begin{tikzpicture}[commutative diagrams/every diagram]
				
				\node (N1) at (0,-2cm) {$Ud(v)$};
				
				\node (N2) at (-2,0cm) {$M$};
				
				\node (N3) at (2,0cm) {$L$};
				
				\path[commutative diagrams/.cd, every arrow, every label]
				(N2) edge node[swap]{$f_v$} (N1)
				(N3) edge node{$l_v$} (N1)
				(N2) edge node{$f$} (N3);
				\end{tikzpicture}
			}
		\end{center}
		\rule{0pt}{0.5cm}
		Now consider the diagram
		\begin{center}
			\mbox{
				\begin{tikzpicture}[commutative diagrams/every diagram]
				
				\node (N1) at (0,2cm) {$FL$};
				
				\node (N2) at (-2.5,0cm) {$FUd(v_1)$};
				
				\node (N3) at (2.5,0cm) {$FUd(v_2)$};
				
				\node (N4) at (0,-1cm) {$GL$};
				
				\node (N5) at (-2.5,-3cm) {$GUd(v_1)$};
				
				\node (N6) at (2.5,-3cm) {$GUd(v_2)$};
				
				\path[commutative diagrams/.cd, every arrow, every label]
				(N1) edge node[swap]{$F(l_{v_1})$} (N2)
				(N1) edge node{$F(l_{v_2})$} (N3)
				(N2) edge node{$F(Ud(e))$} (N3)
				(N4) edge node[swap]{$F(l_{v_1})$} (N5)
				(N4) edge node{$F(l_{v_2})$} (N6)
				(N5) edge node{$F(Ud(e))$} (N6)
				(N2) edge node[swap]{$\delta_{d(v_1)}$} (N5)
				(N3) edge node{$\delta_{d(v_2)}$} (N6);
				\end{tikzpicture}
			}
		\end{center}
		\rule{0pt}{0.5cm}
		Since $G$ preserves limits of type $D$, the cone \(\ag{G(L)}{\{G(l_v)\mid v\in V(D)\}}\) is a limiting cone in $C$. Therefore, since
		\begin{align*}
		GUd(e)\delta_{d(v_1)}F(l_{v_1})&=\delta_{v_2}FUd(e)F(l_{v_1})\\
		&=\delta_{v_2}F(l_{v_2}),
		\end{align*}
		there exists a unique map $\lambda:FL\longrightarrow GL$, such that the following diagram commutes in $\mcg{C}$:
		\begin{center}
			\mbox{
				\begin{tikzpicture}[commutative diagrams/every diagram]
				
				\node (N1) at (0,2cm) {$FL$};
				
				\node (N2) at (-2.5,0cm) {$FUd(v_1)$};
				\node (R2) at (-1,0.2cm) {};
				
				\node (N3) at (2.5,0cm) {$FUd(v_2)$};
				\node (R3) at (1,0.2cm) {};
				
				\node (N4) at (0,-1cm) {$GL$};
				
				\node (N5) at (-2.5,-3cm) {$GUd(v_1)$};
				
				\node (N6) at (2.5,-3cm) {$GUd(v_2)$};
				
				\path[commutative diagrams/.cd, every arrow, every label]
				(N1) edge node[swap]{$F(l_{v_1})$} (N2)
				(N1) edge node{$F(l_{v_2})$} (N3)
				(N4) edge node[swap]{$F(l_{v_1})$} (N5)
				(N4) edge node{$F(l_{v_2})$} (N6)
				(N2) edge node[swap]{$\delta_{d(v_1)}$} (N5)
				(N3) edge node{$\delta_{d(v_2)}$} (N6)
				(N1) edge[dashed] node[near start]{$\lambda$} (N4)
				(R2) edge[-,line width=20pt,draw=white] node{} (R3)
				(N2) edge node{$F(Ud(e))$} (N3)
				(N5) edge node{$F(Ud(e))$} (N6);
				\end{tikzpicture}
			}
		\end{center}
		\rule{0pt}{0.5cm}
		Clearly, $\eb{L}\eqd\ag{L}{\lambda}$ is a dialgebra and, for all $v \in V(D)$, \(\eb{l}_v:\eb{L}\longrightarrow d(v)\) is a dialgebra homomorphism.\\
		\par To show that $\ag{\eb{L}}{\{\eb{l}_v\mid v\in V(D)\}}$ is a limiting cone in $\dlg{\mcg{C}}{F}{G}$, consider any other cone $\ag{\eb{M}=\ag{M}{\mu}}{\{\eb{f}_v\mid v\in V(D)\}}$ in $\dlg{\mcg{C}}{F}{G}$. Since \(\ag{L}{\{l_v\mid v\in V(D)\}}\) is a limiting cone in $\mcg{C}$, there exists a unique arrow \(f:UM\longrightarrow L\) in $\mcg{C}$ such that for all $v \in V(D)$ \[Uf_v=l_vf.\]
		It suffices to show that $f$ yields a dialgebra homomorphism $\eb{f:M\longrightarrow L}$. For all $v$ we have:
		\begin{center}
			\mbox{
				\begin{tikzpicture}[commutative diagrams/every diagram]
				
				\node (N1) at (0,2cm) {$FL$};
				
				\node (N2) at (-2.5,2cm) {$FUM$};
				
				\node (N3) at (-2.5,0cm) {$GUM$};
				
				\node (N4) at (0,0cm) {$GL$};
				
				\node (N5) at (2.5,0cm) {$GUd(v)$};
				
				\path[commutative diagrams/.cd, every arrow, every label]
				(N2) edge node{$Ff$} (N1)
				(N3) edge node{$Gf$} (N4)
				(N2) edge node[swap]{$\mu$} (N3)
				(N4) edge node{$Gl_v$} (N5)
				(N1) edge node[swap]{$\lambda$} (N4);
				\end{tikzpicture}
			}
		\end{center}
		\rule{0pt}{0.5cm}
		\begin{align*}
		G(l_v)\lambda F(f)&=\delta_{d(v)}F(l_v)F(f)\\
		&=\delta_{d(v)}F(l_vf)\\
		&=\delta_{d(v)}F(U(f_v))\\
		&=G(U(f_v))\mu\\
		&=G(l_vf)\mu\\
		&=G(l_v)G(f)\mu;
		\end{align*}
		whence, by the universal arrow property of $\ag{GL}{\{Gl_v\mid v\in V(D)\}}$, it now follows that $\lambda F(f)=G(f)\mu$.
	\end{mypr}
	~\\
	\par Following the lines of the proof of Theorem \ref{forget-crea-pres-G-F}, the next theorem may also be proved.
	
	\begin{thm}
		\label{pres-weak-lims-colims}
		Assume that $\mcg{C}$ is bicomplete. Let $D$ be a graph, and \(d:D\longrightarrow\dlg{\mcg{C}}{F}{G}\) be a diagram in $\dlg{\mcg{C}}{F}{G}$.
		\begin{enumerate}
			\item If the functor $G$ preserves weak limits of type $D$, then the limit \(\ag{L}{\{l_v\mid v\in V(D)\}}\) of $Ud:D\longrightarrow \mcg{C}$ in $\mcg{C}$ may be endowed with a dialgebra structure $\eb{L}=\ag{L}{\lambda}$, such that Diagram \ref{pres-weak-lims-colims-d1} commutes in $\dlg{\mcg{C}}{F}{G}$, for all $e:v_1\longrightarrow v_2$ in \(E(D)\).
			\begin{figure}[h]
				\begin{center}
					\begin{tikzpicture}[commutative diagrams/every diagram]
					
					\node (N1) at (0,2cm) {$\eb{L}$};
					
					\node (N2) at (-2,0cm) {$d(v_1)$};
					
					\node (N3) at (2,0cm) {$d(v_2)$};
					
					\path[commutative diagrams/.cd, every arrow, every label]
					(N1) edge node[swap]{$\eb{l}_{v_1}$} (N2)
					(N1) edge node{$\eb{l}_{v_2}$} (N3)
					(N2) edge node{$d(e)$} (N3);
					\end{tikzpicture}
				\end{center}
				\caption{} \label{pres-weak-lims-colims-d1}
			\end{figure}
			
			\item If the functor $F$ preserves weak colimits of type $D$, then the colimit \(\ag{C}{\{c_v\mid v\in V(D)\}}\) of $Ud:D\longrightarrow \mcg{C}$ in $\mcg{C}$ may be endowed with a dialgebra structure $\eb{C}=\ag{C}{\kappa}$, such that Diagram \ref{pres-weak-lims-colims-d2} commutes in $\dlg{\mcg{C}}{F}{G}$, for all $e:v_1\longrightarrow v_2$ in \(E(D)\).
			\begin{figure}[h]
				\begin{center}
					\begin{tikzpicture}[commutative diagrams/every diagram]
					
					\node (N1) at (0,-2cm) {$\eb{C}$};
					
					\node (N2) at (-2,0cm) {$d(v_1)$};
					
					\node (N3) at (2,0cm) {$d(v_2)$};
					
					\path[commutative diagrams/.cd, every arrow, every label]
					(N2) edge node[swap]{$\eb{c}_{v_1}$} (N1)
					(N3) edge node{$\eb{c}_{v_2}$} (N1)
					(N2) edge node{$d(e)$} (N3);
					\end{tikzpicture}
				\end{center}
				\caption{} \label{pres-weak-lims-colims-d2}
			\end{figure}
		\end{enumerate}
	\end{thm}
	
	The following theorem exhibits some of the fruitful interaction of the notion of bisimulation with the weak preservation of limits by the functor $G$.
	
	\begin{thm}
		\label{G-pres-weak-pullbacks}
		Assume that $\mcg{C}$ is bicomplete. Let $\eb{A}=\ag{A}{\alpha},\eb{B}=\ag{B}{\beta},\eb{C}=\ag{C}{\gamma}$ be $\ag{F}{G}$-dialgebras, and \(\eb{f:A\longrightarrow C},\eb{g:B\longrightarrow C}\) be dialgebra homomorphisms in $\dlg{\mcg{C}}{F}{G}$. If $G$ preserves weak pullbacks, then the pullback of $f:A\longrightarrow C,g:B\longrightarrow C$ in $\mcg{C}$ yields a bisimulation from $\eb{A}$ to $\eb{B}$ in $\dlg{\mcg{C}}{F}{G}$.
	\end{thm}
	
	\begin{mypr}
		Let $\eb{f:A\longrightarrow C}$ and \(\eb{g:B\longrightarrow C}\) be homomorphisms. Consider the pullback diagram in $\mcg{C}$:
		\begin{center}
			\mbox{
				\begin{tikzpicture}[commutative diagrams/every diagram]
				
				\node (N1) at (1.5,2cm) {$B$};
				
				\node (N2) at (-1.5,2cm) {$A\times_C B$};
				
				\node (N3) at (-1.5,0cm) {$A$};
				
				\node (N4) at (1.5,0cm) {$C$};
				
				\path[commutative diagrams/.cd, every arrow, every label]
				(N2) edge node{$\pi_2$} (N1)
				(N3) edge node{$f$} (N4)
				(N2) edge node[swap]{$\pi_1$} (N3)
				(N1) edge node[swap]{$g$} (N4);
				\end{tikzpicture}
			}
		\end{center}
		\rule{0pt}{0.5cm}
		Then, since $G$ preserves weak pullbacks, by Theorem \ref{pres-weak-lims-colims} (1), there exists an arrow $\pi:F(A\times_C B)\longrightarrow G(A\times_C B)$ such that both \(\boldsymbol{\uppi}_1:\ag{A\times_C B}{\pi}\longrightarrow\eb{A}\) and \(\boldsymbol{\uppi}_2:\ag{A\times_C B}{\pi}\longrightarrow\eb{B}\) are homomorphisms.
	\end{mypr}
	\vspace{5mm}
	\section{The problem of dualization}
	\label{the-prob-of-dualiz}
	As we discussed earlier (see \ref{funda-motivs}), there is a problem concerning the dualization of dialgebras. We explain this more here. Firstly, we give a statement of the problem:\\
	
	\par \noindent \textbf{Statement of the problem.} \textit{Let $ F,G $ be endofunctors on $ \mcg{C} $. Can the results in $ \dlg{\mcg{C}}{F}{G} $ be translated into corresponding ones in $ \dlg{\mcg{C}}{G}{F} $?}\\
	
	\par To answer the above question, we pay attention to the fact that there is a formal relationship between the category of $ F $-coalgebras and the category of $ \op{F} $-algebras \cite{Hugh}:
	
	\begin{prp}
		\label{formal-rel-alg-coalg}
		Let $ \mcg{C} $ be any category. Then the category $ \mcg{C}_F $ of $ F $-coalgebras arises formally as the category $ \op{((\op{\mcg{C}})^{\op{F}})} $.
	\end{prp}
	
	This fact can be generalized to arbitrary dialgebras:
	
	\begin{prp}
		\label{formal-dlz-dialg}
		Let $ \mcg{C} $ be any category. Then the category $ \dlg{\mcg{C}}{G}{F} $ of $ \ag{G}{F} $-dialgebras arises formally as the category $ \op{(\dlg{(\op{\mcg{C}})}{\op{F}}{\op{G}})} $.
	\end{prp}
	
	Therefore, whenever $ \mcg{C} $ is non-self-dual (e.g. when $ \mcg{C}=\cat{Set} $), our problem seems to have no simple answer (also, see \cite{Gumm-UnCoLo}). However, self-dual bases provide us with an advantage:
	
	\begin{prp}
		\label{slfd-advant}
		Let $ \mcg{C} $ be self-dual (see Definition \ref{selfdual}), and let $ F,G $ be endofunctors on $ \mcg{C} $. Then
		\[\dlg{\mcg{C}}{G}{F}\cong\op{(\dlg{\mcg{C}}{F}{G})}.\]
	\end{prp}
	
	\vspace{2mm}	
	In other words, we have a (contravariant) isomorphism of categories of dialgebras whenever the base category is isomorphic to its dual.
	
	\section{Universal dialgebra on the Chu construction}
	\label{unidialg-on-chu}
	Eventually, we fix $\mcg{C}=\chu=\chu_\Gamma$ for some fixed nonempty set $\Gamma$, and study some specific properties of the category $\dlg{\chu}{F}{G}$ for some fixed endofunctors $F,G$. One pleasant thing about \chu~is that it is bicomplete (and also, balanced). Therefore, all the statements proven in the previous sections apply to it. Below we give more detailed statements regarding the above.\\
	\par First of all, let us take a look at the detailed structure of a relation between two Chu spaces. For Chu spaces $\msf{A}=\agg{A}{r}{X},\msf{B}=\agg{B}{s}{Y}$, the categorical product $\msf{A\&B}$ has the form \(\agg{A\times B}{w}{X\sqcup Y}\), as was discussed before. Now, the following corollary to Proposition \ref{relation-as-monic} characterizes the relations between Chu spaces.
	
	\begin{cor}
		\label{Chu-relation-as-monic}
		Let $\msf{A}=\agg{A}{r}{X},\msf{B}=\agg{B}{s}{Y}$ be Chu spaces. Then every relation \(\agg{\msf{R}}{\pi_1}{\pi_2}\) between $\msf{A,B}$ can be viewed as a monic arrow \(\varphi:\msf{R\rightarrowtail A\&B}\) and vice versa. For $\msf{R}=\agg{C}{t}{Z}$ and monic \(\varphi:\msf{R\rightarrowtail A\&B}\), the function $\varphi^+:C\rightarrow A\times B$ is injective while $\varphi^-:X\sqcup Y\rightarrow Z$ is surjective.
	\end{cor}
	
	\begin{mypr}
		Use Propositions \ref{relation-as-monic} and \ref{Chu-monic-epic}.
	\end{mypr}
	
	~\\
	\par Next, the graph arrow of Chu transforms:
	
	\begin{prp}
		\label{Gr-Chu-transforms}
		Let $\msf{A}=\agg{A}{r}{X},\msf{B}=\agg{B}{s}{Y}$ be Chu spaces. For every Chu transform $f:\msf{A\longrightarrow B}$ we have
		\[\Gr_f=\ag{\Gr_{f^+}}{\left[1_X,f^- \right]},\] where $\Gr_f$ and $\Gr_{f^+}$ denote the graph arrows of $f$ and \(f^+\) in the categories \chu~and \cat{Set}, respectively, and $\left[1_X,f^- \right]$ denotes the unique arrow that makes the following diagram commute. In this diagram, the top arrows are the obvious inclusions:
		\begin{center}
			\mbox
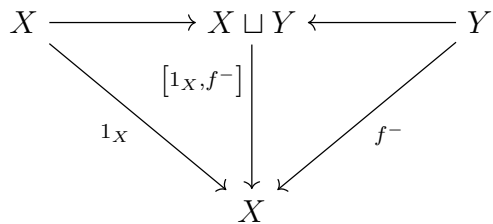
{
				\begin{tikzpicture}[commutative diagrams/every diagram]
				
				\node (N1) at (0,0cm) {$X\sqcup Y$};
				
				\node (N2) at (-3,0cm) {$X$};
				
				\node (N3) at (0,-2.5cm) {$X$};
				
				\node (N4) at (3,0cm) {$Y$};
				
				\path[commutative diagrams/.cd, every arrow, every label]
				(N2) edge node{} (N1)
				(N1) edge node[swap,near start]{$\left[1_X,f^- \right]$} (N3)
				(N2) edge node[swap]{$1_X$} (N3)
				(N4) edge node{} (N1)
				(N4) edge node{$f^-$} (N3);
				\end{tikzpicture}
			}
		\end{center}
	\end{prp}
	
	\begin{mypr}
		This easily follows from the properties of the graph arrow in Diagram \ref{graph-arr-diag}.
	\end{mypr}
	
	~\\
	\par Now we turn to limits and colimits. Following Theorem \ref{forget-crea-pres-G-F}:
	
	\begin{cor}
		\label{Chu-forget-crea-pres-G-F}
		\begin{enumerate}
			\item The forgetful functor $U:\dlg{\chu}{F}{G}\longrightarrow\chu$ creates and preserves all types of limits that the functor $G$ preserves.
			\item The forgetful functor $U:\dlg{\chu}{F}{G}\longrightarrow\chu$ creates and preserves all types of colimits that the functor $F$ preserves.
		\end{enumerate}
	\end{cor}
	
	In the same manner, one can deduce corollaries from Theorems \ref{pres-weak-lims-colims} and \ref{G-pres-weak-pullbacks} when $\mcg{C}=\chu$.\\
	
	\par Finally, we state the analog of Proposition \ref{slfd-advant}:
	
	\begin{prp}
		\label{slfd-advant-Chu}
		We have \[\dlg{\chu}{G}{F}\cong\op{(\dlg{\chu}{F}{G})}.\]
	\end{prp}

	\chapter{The Double Category of Paired Dialgebras on Chu}
	As we mentioned in Chapter 1, the theory of universal dialgebra was first developed over the category \cat{Set} \cite{Vout}. There are a number of possible stages of expanding the theory in various ways. In the first stage of expansion, the base category \cat{Set} might well be generalized to some other bicomplete category. What we did in the previous chapter was taking a few steps just in that direction. More precisely, we observed how the category \chu~can serve as a substitute for \cat{Set}.\\
	\par In the second stage, one may try to generalize the concept of ``universal dialgebra'' itself to an even broader framework. Here, equipped with all the eesential tools having been developed in the previous chapters, we are to enter that second stage, which constitues the main part of the current thesis.
	\par In this chapter we introduce double categories firstly (Section \ref{intro-double-cats}). Then, we will introduce the main formalism of the present thesis, i.e., the $ \dlc_{\Gamma,\Sigma} $ construction (Section \ref{the-main-forma}), and we will show that the construction is well-defined. Lastly, we will investigate some of the basic properties of the formalism (Section \ref{some-elem-props}).
	
	\section{Double categories}
	\label{intro-double-cats}
	Our ``second-stage'' generalization of universal dialgebra utilizes the language of \textit{double categories}. Double categories were introduced by C. Ehresmann and have been further developed by several people since then (see \cite{Fiore-Thoma} and the references therein). The notion of double category might be seen as a ``two-dimensional'' generalization of the very notion of category. That is to say, whereas in ordinary cateogries we used to deal with objects (``zero-dimensional points'', diagrammatically) and morphisms (``one-dimensional arrows''), in double categories we have objects (the so-called ``0-cells''), \textit{horizontal morphisms} (``horizontal 1-cells''), \textit{vertical morphisms} (``vertical 1-cells''), and 2-cells. This way, then, ordianry categories can be seen as special cases of double categories when some of the double-categorical structure is trivialized. It is worth noting that there also exist other two-dimensional generalizations of ordinary categories, namely \textit{2-categories} and \textit{bicategories}. We intend to mention the fact that double categories--together with their ``weak'' versions--include 2-categories and bicategories as special cases (see Remark \ref{double-to-2-cat}).\\
	
	\par It is important to know that development of a theory of double categories can be done in various ways. The usual method is to introduce them as \textit{internal categories} either in \cat{CAT} or in \cat{Cat}. This is exactly the approach we take in the present work. However, there is another approach in the literature that develops the theory based on the theory of \textit{indexed categories} (see \cite{GranPar-Lims}); last but not least, some authors believe that double categories are worth an independent treatment, directly founded on basic category theory \cite{GranPar-Lims}.\\
	
	\par In Section \ref{internal-cats-section} we introduced internal categories. Now, we proceed as the following. The material in this section is borrowed from \cite{Doub-Q-Goid,Fiore-Adj,Fiore-Model,Fiore-Thoma,GranPar-Lims,Nief-Expo,Nief-Span,Schul,nLab-doub,Catsters}.
	
	\begin{dfn}
		\label{double-cat-def}
		\begin{enumerate}
			\item A \iw{(large) double category} is an internal category $\bb{D}=\langle \mcg{A,B,s,t,i,\bowtie}\rangle$ in \cat{CAT}.
			
			\item Likewise, a \iw{small double category} is an internal category in \cat{Cat}.
			
			\item A double category is said to be \textbf{properly large} if it is not small.
		\end{enumerate}
	\end{dfn}
	
	\rule{0pt}{1mm}
	
	\begin{rem}
		\label{explain-double-cats}
		The following explanations may help the reader obtain a better understanding of the notion of double category. A double category \[\bb{D}=\langle \mcg{A,B,s,t,i,\bowtie}\rangle\] consists of an ``object of objects'' $\mcg{A}$, an ``object of arrows'' $\mcg{B}$, a ``source'' $\mcg{s}$, a ``target'' $\mcg{t}$, an ``identity'' $\mcg{i}$, and a ``composition'' $\bowtie$. Then, as $\bb{D}$ is internal to \cat{CAT}, it follows that $\mcg{A}$ and $\mcg{B}$ are themselves categories, and $\mcg{s,t,i,\bowtie}$ are themselves functors.\\
		
		\par Therefore, $\mcg{B}$ itself has a class of objects $\uo{\mcg{B}}$ and a class of morphisms $\um{\mcg{B}}$, together with source and target maps
		\[\mcg{b_0,b_1}:\um{\mcg{B}}\dlong\uo{\mcg{B}},\]
		and similarly for $\mcg{A}$:
		\[\mcg{a_0,a_1}:\um{\mcg{A}}\dlong\uo{\mcg{A}}.\]
		Also, each of the functors $\mcg{s}=\ag{\mcg{s_0}}{\mcg{s_1}}$ and $\mcg{t}=\ag{\mcg{t_0}}{\mcg{t_1}}$ consists of a map on objects and a map on morphisms. Diagram \ref{A-B-obj-mor-s-t} depicts all of these maps together.
		\begin{figure}[H]
			\begin{center}
				\begin{tikzpicture}[commutative diagrams/every diagram]
				
				\node (N1) at (1.8,1.5cm) {$\uo{\mcg{B}}$};
				\node (V1) at (1.4,1.6cm) {};
				\node (V2) at (1.4,1.4cm) {};
				\node (V9) at (1.7,1.3cm) {};
				\node (V10) at (1.9,1.3cm) {};
				
				\node (N2) at (-1.8,1.5cm) {$\um{\mcg{B}}$};
				\node (V3) at (-1.4,1.6cm) {};
				\node (V4) at (-1.4,1.4cm) {};
				\node (V11) at (-1.7,1.3cm) {};
				\node (V12) at (-1.9,1.3cm) {};
				
				\node (N3) at (-1.8,-1.5cm) {$\um{\mcg{A}}$};
				\node (V5) at (-1.4,-1.6cm) {};
				\node (V6) at (-1.4,-1.4cm) {};
				\node (V13) at (-1.7,-1.3cm) {};
				\node (V14) at (-1.9,-1.3cm) {};
				
				\node (N4) at (1.8,-1.5cm) {$\uo{\mcg{A}}$};
				\node (V7) at (1.4,-1.6cm) {};
				\node (V8) at (1.4,-1.4cm) {};
				\node (V15) at (1.7,-1.3cm) {};
				\node (V16) at (1.9,-1.3cm) {};
				
				\path[commutative diagrams/.cd, every arrow, every label]
				(V3) edge node{$\mcg{b_0}$} (V1)
				(V4) edge node[swap]{$\mcg{b_1}$} (V2)
				(V5) edge node[swap]{$\mcg{a_1}$} (V7)
				(V6) edge node{$\mcg{a_0}$} (V8)
				(V12) edge node[swap]{$\mcg{s_1}$} (V14)
				(V11) edge node{$\mcg{t_1}$} (V13)
				(V9) edge node[swap]{$\mcg{s_0}$} (V15)
				(V10) edge node{$\mcg{t_0}$} (V16);
				\end{tikzpicture}
			\end{center}
			\caption{} \label{A-B-obj-mor-s-t}
		\end{figure}
		In this diagram, the following commutativity conditions hold:
		\begin{align}
		\mcg{a}_k\mcg{s}_1&=\mcg{s}_0\mcg{b}_k,~~k=0,1;\tag{$*$}\\
		\mcg{a}_k\mcg{t}_1&=\mcg{t}_0\mcg{b}_k,~~k=0,1.\tag{$**$}
		\end{align}
		
		\par On the other hand, each of $\mcg{i},\bowtie$ is itself composed of a function on objects as well as a function on arrows. That is, we have functions
		\[\mcg{i}_0:\uo{\mcg{A}}\longrightarrow\uo{\mcg{B}},\]
		\[\mcg{i}_1:\um{\mcg{A}}\longrightarrow\um{\mcg{B}},\]
		and also
		\[\bowtie_0:\uo{(\pb{\mcg{B}}{\mcg{A}}{\mcg{B}})}\longrightarrow\uo{\mcg{B}},\]
		\[\bowtie_1:\um{(\pb{\mcg{B}}{\mcg{A}}{\mcg{B}})}\longrightarrow\um{\mcg{B}}.\]
		Let the compositions in $\mcg{A,B}$ be denoted by $\circ_\mcg{A},\circ_\mcg{B}$, respectively. Then, by functoriality of $\mcg{s,t}$ we have for any $\bm{c,d}\in\um{\mcg{B}}$,
		\[\um{\mcg{s}}(\bm{d}\circ_\mcg{B}\bm{c})=\um{\mcg{s}}(\bm{d})\circ_\mcg{A}\um{\mcg{s}}(\bm{c}),\]
		\[\um{\mcg{t}}(\bm{d}\circ_\mcg{B}\bm{c})=\um{\mcg{t}}(\bm{d})\circ_\mcg{A}\um{\mcg{t}}(\bm{c}).\]
		From the functoriality of $\mcg{a_0,a_1,b_0,b_1}$ we find
		\[\mcg{b}_k\mcg{i_1}=\mcg{i_0}\mcg{a}_k,~~k=0,1.\]
		Also, the functoriality of $\bowtie$ amounts to
		\[\mcg{b}_k\bowtie_1(\bm{d,c})=\bt_0(\mcg{b}_k(\bm{d}),\mcg{b}_k(\bm{c}),~~k=0,1\]
		for $\bm{c,d}\in\um{\mcg{B}}$, and
		\[\bt_1(\bm{d}\circ_\mcg{B}\bm{c},\bm{f}\circ_\mcg{B}\bm{e})=\bt_1(\bm{d,f})\circ_\mcg{B}\bt_1(\bm{c,e}),\]
		for $\bm{c,d,e,f}\in\um{\mcg{B}}$, whenever all the above compositions are defined.
	\end{rem}
	
	Now we analyze the subject from a different viewpoint. What the above formalism actually gives us is four classes of ``things''. We interpret these ``things'' in the following manner:
	\begin{itemize}
		\item[\checkmark] we regard $\uo{\mcg{A}}$ as the class of \iw{0-cells};
		
		\item[\checkmark] we regard $\um{\mcg{A}}$ as the class of \iw{horizontal 1-cells};
		
		\item[\checkmark] we regard $\uo{\mcg{B}}$ as the class of \iw{vertical 1-cells}; and
		
		\item[\checkmark] we regard $\um{\mcg{B}}$ as the class of \iw{2-cells}.
	\end{itemize}
	
	Let us look at what a 2-cell in $\um{\mcg{B}}$ looks like. It has a source and target in $\um{\mcg{A}}$, so it has a ``source horizontal 1-cell'' and a ``target horizontal 1-cell''; but it also has a source and target in $\uo{\mcg{B}}$, thus it has a ``source vertical 1-cell'' and a ``target vertical 1-cell''. The commutativity conditions ($*$) and ($**$) of Diagram \ref{A-B-obj-mor-s-t} tell us that \textit{the corners of the adjacent 1-cells match up}. The result looks like the \iw{cell structure} sketched in Diagram \ref{cell-struc}. Notice that the horizontal and vertical arrows in that diagram are depicted by different shapes of arrowheads, in order to distinguish between the horizontal and vertical 1-cells.\\
	\begin{figure}[h]
		\begin{center}
			\begin{tikzpicture}[commutative diagrams/every diagram]
			
			\node (N1) at (2,2cm)[circle,draw,fill] {};
			\node (V1) at (1.7,2cm) {};
			\node (V2) at (2,1.7cm) {};
			
			\node (N2) at (-2,2cm)[circle,draw,fill] {};
			\node (V3) at (-1.7,2cm) {};
			\node (V4) at (-2,1.7cm) {};
			
			\node (N3) at (-2,-2cm)[circle,draw,fill] {};
			\node (V5) at (-1.7,-2cm) {};
			\node (V6) at (-2,-1.7cm) {};
			
			\node (N4) at (2,-2cm)[circle,draw,fill] {};
			\node (V7) at (1.7,-2cm) {};
			\node (V8) at (2,-1.7cm) {};
			
			\draw[line width=5pt, -to] (V3) to node {} (V1);
			\draw[line width=5pt, -left to] (V4) to node {} (V6);
			\draw[line width=5pt, -to] (V5) to node {} (V7);
			\draw[line width=5pt, -left to] (V2) to node {} (V8);
			
			\node (B1) at (5,0cm) {\fbox{\begin{tabular}{c}
					target vertical 1-cell\\
					$\in\uo{\mcg{B}}$
					\end{tabular}}};
			\node (X1) at (2,0cm) {};
			
			\node (B2) at (0,4.5cm) {\fbox{\begin{tabular}{c}
					source horizontal 1-cell\\
					$\in\um{\mcg{A}}$
					\end{tabular}}};
			\node (X2) at (0,2cm) {};
			
			\node (B3) at (-5,0cm) {\fbox{\begin{tabular}{c}
					source vertical 1-cell\\
					$\in\uo{\mcg{B}}$
					\end{tabular}}};
			\node (X3) at (-2,0cm) {};
			
			\node (B4) at (0,-4.5cm) {\fbox{\begin{tabular}{c}
					target horizontal 1-cell\\
					$\in\um{\mcg{A}}$
					\end{tabular}}};
			\node (X4) at (0,-2cm) {};
			
			\node (B5) at (4,2.5cm) {\fbox{\begin{tabular}{c}
					0-cell\\
					$\in\uo{\mcg{A}}$
					\end{tabular}}};
			\node (X5) at (2,2cm) {};
			
			\node (B6) at (-4,2.5cm) {\fbox{\begin{tabular}{c}
					0-cell\\
					$\in\uo{\mcg{A}}$
					\end{tabular}}};
			\node (X6) at (-2,2cm) {};
			
			\node (B7) at (-4,-2.5cm) {\fbox{\begin{tabular}{c}
					0-cell\\
					$\in\uo{\mcg{A}}$
					\end{tabular}}};
			\node (X7) at (-2,-2cm) {};
			
			\node (B8) at (4,-2.5cm) {\fbox{\begin{tabular}{c}
					0-cell\\
					$\in\uo{\mcg{A}}$
					\end{tabular}}};
			\node (X8) at (2,-2cm) {};
			
			\node (W1) at (1.7,0cm) {};
			\node (W2) at (-1.7,0cm) {};
			
			\draw[ -implies, double distance=30pt,line width=5pt] (W2) to node {} (W1);
			
			\node (B9) at (-0.5,0cm) {\mbox{\begin{tabular}{c}
					2-cell\\
					$\in\um{\mcg{B}}$
					\end{tabular}}};
			
			\draw[dashed, line width=1pt] (B1) to node {} (X1);
			\draw[dashed, line width=1pt] (B2) to node {} (X2);
			\draw[dashed, line width=1pt] (B3) to node {} (X3);
			\draw[dashed, line width=1pt] (B4) to node {} (X4);
			\draw[dashed, line width=1pt] (B5) to node {} (X5);
			\draw[dashed, line width=1pt] (B6) to node {} (X6);
			\draw[dashed, line width=1pt] (B7) to node {} (X7);
			\draw[dashed, line width=1pt] (B8) to node {} (X8);
			\end{tikzpicture}
		\end{center}
		\caption{The cell structure} \label{cell-struc}
	\end{figure}
	
	\par Now let us see what kinds of compositions we have here. There are two things that are going on here. First of all, the fact that $\mcg{B}$ and $\mcg{A}$ are themselves categories means that each of them is equipped with a composition, namely $\circ_\mcg{B}$ and $\circ_\mcg{A}$. The composition $\circ_\mcg{B}$ composes the adjacent 2-cells along the vertical 1-cells, while the composition $\circ_\mcg{A}$ composes the adjacent horizontal 1-cells along the 0-cells. The commutativity conditions ($*$) and ($**$) of Diagram \ref{A-B-obj-mor-s-t} ensure that these two compositions are compatible with each other. The result is the \textit{horizontal composition} sketched in Diagram \ref{hor-compo}.\\
	
	\begin{figure}[H]
		\begin{center}
			\begin{tikzpicture}[commutative diagrams/every diagram]
			\node (N1) at (2,1cm) {$\bullet$};
			\node (N2) at (0,1cm) {$\bullet$};
			\node (N3) at (-2,1cm) {$\bullet$};
			\node (N4) at (2,-1cm) {$\bullet$};
			\node (N5) at (0,-1cm) {$\bullet$};
			\node (N6) at (-2,-1cm) {$\bullet$};
			
			\draw[line width=1pt, -to] (N2) to node[auto] {$g$} (N1);
			
			\draw[line width=1pt, -left to] (N1) to node {} (N4);
			
			\draw[line width=1pt, -to] (N3) to node[auto] {$f$} (N2);
			
			\draw[line width=1pt, -left to] (N2) to node {} (N5);
			
			\draw[line width=1pt, -to] (N5) to node[auto,swap] {$k$} (N4);
			
			\draw[line width=1pt, -left to] (N3) to node {} (N6);
			
			\draw[line width=1pt, -to] (N6) to node[auto,swap] {$h$} (N5);
			
			\node (V1) at (-1.8,0cm) {};
			\node (V2) at (-0.2,0cm) {};
			\node (V3) at (0.2,0cm) {};
			\node (V4) at (1.8,0cm) {};
			
			\draw[ -implies, double distance=3pt, line width=1pt] (V1) to node[auto] {$\bm{c}$} (V2);
			
			\draw[ -implies, double distance=3pt, line width=1pt] (V3) to node[auto] {$\bm{d}$} (V4);
			
			\node (E1) at (3,0cm) {};
			\node (E2) at (4,0cm) {};
			
			\draw[-, double distance=4pt, line width=3pt] (E1) to node {} (E2);
			
			\node (P1) at (5,1cm) {$\bullet$};
			\node (P2) at (5,-1cm) {$\bullet$};
			\node (P3) at (8,-1cm) {$\bullet$};
			\node (P4) at (8,1cm) {$\bullet$};
			
			\draw[line width=1pt, -to] (P1) to node[auto] {$g\circ_\mcg{A}f$} (P4);
			
			\draw[line width=1pt, -left to] (P1) to node[auto] {} (P2);
			
			\draw[line width=1pt, -to] (P2) to node[auto,swap] {$k\circ_\mcg{A}h$} (P3);
			
			\draw[line width=1pt, -left to] (P4) to node[auto] {} (P3);
			
			\node (W1) at (5.2,0cm) {};
			\node (W2) at (7.8,0cm) {};
			
			\draw[ -implies, double distance=3pt, line width=1pt] (W1) to node[auto] {$\bm{d}\circ_\mcg{B}\bm{c}$} (W2);
			\end{tikzpicture}
		\end{center}	
		\caption{} \label{hor-compo}
	\end{figure}
	
	\par On the other hand, we have the internal composition bifunctor $\bowtie=\ag{\bowtie_0}{\bowtie_1}$. This internal composition provides a way for \textit{vertical composition} of stacked cell structures, as sketched in Diagram \ref{ver-compo}.\\
	
	\begin{figure}[H]
		\begin{center}
			\begin{tikzpicture}[commutative diagrams/every diagram]
			\node (N1) at (1,2cm) {$\bullet$};
			\node (N2) at (-1,2cm) {$\bullet$};
			\node (N3) at (1,0cm) {$\bullet$};
			\node (N4) at (-1,0cm) {$\bullet$};
			\node (N5) at (1,-2cm) {$\bullet$};
			\node (N6) at (-1,-2cm) {$\bullet$};
			
			\draw[line width=1pt, -to] (N2) to node {} (N1);
			
			\draw[line width=1pt, -left to] (N2) to node[auto,swap] {$u$} (N4);
			
			\draw[line width=1pt, -to] (N4) to node {} (N3);
			
			\draw[line width=1pt, -left to] (N4) to node[auto,swap] {$v$} (N6);
			
			\draw[line width=1pt, -to] (N6) to node {} (N5);
			
			\draw[line width=1pt, -left to] (N1) to node[auto] {$\prm{u}$} (N3);
			
			\draw[line width=1pt, -left to] (N3) to node[auto] {$\prm{v}$} (N5);
			
			\node (V1) at (-0.8,1cm) {};
			\node (V2) at (0.8,1cm) {};
			\node (V3) at (-0.8,-1cm) {};
			\node (V4) at (0.8,-1cm) {};
			
			\draw[ -implies, double distance=3pt, line width=1pt] (V1) to node[auto] {$\bm{c}$} (V2);
			
			\draw[ -implies, double distance=3pt, line width=1pt] (V3) to node[auto] {$\bm{d}$} (V4);
			
			\node (E1) at (2.5,0cm) {};
			\node (E2) at (3.5,0cm) {};
			
			\draw[-, double distance=4pt, line width=3pt] (E1) to node {} (E2);
			
			\node (P1) at (5,2cm) {$\bullet$};
			\node (P2) at (5,-2cm) {$\bullet$};
			\node (P3) at (7.3,-2cm) {$\bullet$};
			\node (P4) at (7.3,2cm) {$\bullet$};
			
			\draw[line width=1pt, -to] (P1) to node {} (P4);
			
			\draw[line width=1pt, -left to] (P1) to node[auto,swap, near start] {$u\bowtie_0 v$} (P2);
			
			\draw[line width=1pt, -to] (P2) to node {} (P3);
			
			\draw[line width=1pt, -left to] (P4) to node[auto] {$\prm{u}\bowtie_0 \prm{v}$} (P3);
			
			\node (W1) at (5.1,0cm) {};
			\node (W2) at (7.2,0cm) {};
			
			\draw[ -implies, double distance=15pt, line width=1pt] (W1) to node {} (W2);
			
			\node (X) at (5.9,0cm) {$\bm{c}\bowtie_1\bm{d}$};
			
			\end{tikzpicture}
		\end{center}
		\caption{} \label{ver-compo}	
	\end{figure}
	
	\begin{rem}
		\label{diff-ntns-compo-pasting}
		Notice that we have adopted different notational conventions for horizontal and vertical compositions; namely, we have used the ``pasting order'' for $\bowtie$ (so that $\bm{c\bowtie d}$ abbreviates $\bowtie(\bm{c,d})$), while we have used the usual ``algebraic order'' for $\circ_\mcg{A},\circ_\mcg{B}$.
	\end{rem}
	
	Now, these horizontal and vertical compositions must satisfy the previously-mentioned equations (Remark \ref{explain-double-cats}). For example, for 2-cells $\bm{c,d,e,f}$ in Diagram \ref{compat-cond} we have
	\[(\bm{d}\circ_\mcg{B}\bm{c})\bowtie(\bm{f}\circ_\mcg{B}\bm{e})=(\bm{c\bowtie e})\circ_\mcg{B}(d\bowtie f).\]
	
	\begin{figure}[H]
		\begin{center}
			\begin{tikzpicture}[commutative diagrams/every diagram]
			\node (N1) at (2,1cm) {$\bullet$};
			\node (N2) at (0,1cm) {$\bullet$};
			\node (N3) at (-2,1cm) {$\bullet$};
			\node (N4) at (2,-1cm) {$\bullet$};
			\node (N5) at (0,-1cm) {$\bullet$};
			\node (N6) at (-2,-1cm) {$\bullet$};
			\node (N7) at (2,-3cm) {$\bullet$};
			\node (N8) at (0,-3cm) {$\bullet$};
			\node (N9) at (-2,-3cm) {$\bullet$};
			
			\draw[line width=1pt, -to] (N2) to node {} (N1);
			\draw[line width=1pt, -to] (N3) to node {} (N2);
			\draw[line width=1pt, -to] (N5) to node {} (N4);
			\draw[line width=1pt, -to] (N6) to node {} (N5);
			\draw[line width=1pt, -to] (N8) to node {} (N7);
			\draw[line width=1pt, -to] (N9) to node {} (N8);
			
			\draw[line width=1pt, -left to] (N1) to node {} (N4);
			\draw[line width=1pt, -left to] (N2) to node {} (N5);
			\draw[line width=1pt, -left to] (N3) to node {} (N6);
			\draw[line width=1pt, -left to] (N4) to node {} (N7);
			\draw[line width=1pt, -left to] (N5) to node {} (N8);
			\draw[line width=1pt, -left to] (N6) to node {} (N9);
			
			\node (V1) at (-1.8,0cm) {};
			\node (V2) at (-0.2,0cm) {};
			\node (V3) at (0.2,0cm) {};
			\node (V4) at (1.8,0cm) {};
			\node (V5) at (-1.8,-2cm) {};
			\node (V6) at (-0.2,-2cm) {};
			\node (V7) at (0.2,-2cm) {};
			\node (V8) at (1.8,-2cm) {};
			
			\draw[ -implies, double distance=3pt, line width=1pt] (V1) to node[auto] {$\bm{c}$} (V2);
			
			\draw[ -implies, double distance=3pt, line width=1pt] (V3) to node[auto] {$\bm{d}$} (V4);
			
			\draw[ -implies, double distance=3pt, line width=1pt] (V5) to node[auto] {$\bm{e}$} (V6);
			
			\draw[ -implies, double distance=3pt, line width=1pt] (V7) to node[auto] {$\bm{f}$} (V8);
			\end{tikzpicture}
		\end{center}
		\caption{} \label{compat-cond}
	\end{figure}
	
	Next, we pay attention to further structures that ``emerge'' from our formalism of double categories.
	
	\begin{prp}
		\label{A-obj-B-obj-cat}
		The $\mcg{A}$-objects (i.e. the 0-cells) together with the $\mcg{B}$-objects (i.e. the vertical 1-cells) form a category.
	\end{prp}
	
	\begin{mypr}
		For every $\mcg{B}$-object $u$, seen as an arrow between two $\mcg{A}$-objects, the domain and codomain are given by $\uo{\mcg{s}}(u)$ and $\uo{\mcg{s}}(u)$, respectively. Also, consecutive $\mcg{B}$-objects $u,v$ compose as $u\bowtie_0 v$. Identity morphisms are given by \(\mcg{i_0}:\uo{\mcg{A}}\longrightarrow\uo{\mcg{B}}.\) Finally, the composition $\bowtie_0$ is associative.
	\end{mypr}
	
	By similar reasoning one can prove:
	
	\begin{prp}
		\label{A-mor-B-mor-cat}
		The class of $\mcg{A}$-morphisms (i.e. the horizontal 1-cells) together with the class of $\mcg{B}$-morphisms (i.e. the 2-cells) form a category, for which the composition is given by $\bowtie_1$.
	\end{prp}
	
	\begin{rem}
		\label{double-to-2-cat}
		It is worth noting here that if all the vertical 1-cells happen to be identities (i.e., if $\mcg{A}$ is a discrete category), then the double cateogry reduces to what is called a \textit{2-category} in the literature. Also, if all the horizontal 1-cells are identities, then the result will be another 2-category. Thus, every double category yields two different 2-categories as special cases. For more on 2-cateogries and their relationships to double categories and other closely-related concepts (such as their respective ``weak'' versions), the reader is referred to \cite{HDC,BasBi,nLab-2-cat,nLab-bicat,nLab-hicat}.
	\end{rem}
	
	Next, we give another important definition.
	
	\begin{dfn}
		\label{transpose-def}
		The \iw{transpose} of a double category $\bb{D}=\langle \mcg{A,B,s,t,i,\bowtie}\rangle$ is the double category
		\[\bb{D^\T}\eqd\langle \prm{\mcg{A}},\prm{\mcg{B}},\prm{\mcg{s}},\prm{\mcg{t}},\prm{\mcg{i}},\prm{\bowtie}\rangle,\]
		where $\prm{\mcg{A}}$ is the category introduced in Proposition \ref{A-obj-B-obj-cat}, and $\prm{\mcg{B}}$ is the category introduced in Proposition \ref{A-mor-B-mor-cat}; the source functor $\prm{\mcg{s}}$ is given as the pair $\ag{\mcg{a_0}}{b_0}$ while $\prm{\mcg{t}}$ is defined as $\langle \mcg{a_1,b_1}\rangle$; the functor $\prm{\mcg{i}}$ is made of the identity functions of $\mcg{A,B}$, and finally, the composition bifunctor $\prm{\bowtie}$ comes from the compositions of the categories $\mcg{A,B}$.
	\end{dfn}
	
	In other words, $\bb{D^\T}$ is the double category which consists of the same 0-cells, 1-cells, and 2-cells as of $\bb{D}$, but in which the roles of horizontal and vertical categories are interchanged.
	
	\begin{rem}
		\label{T-of-T}
		It follows immediately that by transposing the transposed category $\bb{D^\T}$ we return to the original cateogry $\bb{D}$:
		\[(\bb{D^\T})^\T=\bb{D}.\]
	\end{rem}
	
	\begin{dfn}
		\label{boundary}
		Let $\bb{D}$ be a double category, and let Diagram \ref{boundary-diag} be a cell structure in $\bb{D}$.
		\begin{figure}[H]
			\begin{center}
				\begin{tikzpicture}[commutative diagrams/every diagram]
				\node (N1) at (2,1cm) {$N$};
				\node (N2) at (0,1cm) {$M$};
				\node (N3) at (0,-1cm) {$\prm{M}$};
				\node (N4) at (2,-1cm) {$\prm{N}$};
				\draw[line width=1pt, -to] (N2) to node[auto] {$f$} (N1);
				\draw[line width=1pt, -to] (N3) to node[auto,swap] {$g$} (N4);
				\draw[line width=1pt, -left to] (N1) to node[auto] {$v$} (N4);
				\draw[line width=1pt, -left to] (N2) to node[auto,swap] {$u$} (N3);
				\node (V1) at (0.2,0cm) {};
				\node (V2) at (1.8,0cm) {};
				\draw[ -implies, double distance=3pt, line width=1pt] (V1) to node[auto] {$\bm{c}$} (V2);
				\end{tikzpicture}
			\end{center}
			\caption{} \label{boundary-diag}
		\end{figure}
		The \iw{boundary} of the 2-cell $\bm{c}$ is defined as the 4-tuple
		\[\left(\begin{matrix}
		&f& \\ u& &v\\ &g& 
		\end{matrix} \right),\]
		and the situation is denoted by 
		\[\left(\begin{matrix}
		&f& \\ u&\bm{c}&v\\ &g& 
		\end{matrix} \right).\]
	\end{dfn}
	
	\rule{0pt}{0.1cm}
	\begin{dfn}
		\label{flat-def}
		A double category $\bb{D}$ is said to be \iw{flat} if its 2-cells are completely determined by their boundaries. In such case, then, for a cell structure like Diagram \ref{boundary-diag} we write:
		\[\bm{c}=\left(\begin{matrix}
		&f& \\ u& &v\\ &g& 
		\end{matrix} \right).\]
	\end{dfn}
	
	\rule{0pt}{1mm}
	\begin{dfn}
		\label{hor-dual-def}
		Let $\bb{D}=\langle \mcg{A,B,s,t,i,\bowtie}\rangle$ be a double category. The \iw{horizontal dual} of $\bb{D}$, denoted by $\bb{D^\leftrightarrow}$, is defined as
		\[\bb{D^\leftrightarrow}\eqd\langle \op{\mcg{A}},\op{\mcg{B}},\mcg{s^\leftrightarrow,t^\leftrightarrow,i^\leftrightarrow,\bowtie^\leftrightarrow}\rangle,\]
		for which we have
		\[\op{\uo{\mcg{A}}}=\uo{\mcg{A}},~~~~\op{\uo{\mcg{B}}}=\uo{\mcg{B}};\]
		\rule{0pt}{1mm}
		\[\forall f:(f\in\um{\mcg{A}}~~\text{iff}~~\op{f}\in\op{\um{\mcg{A}}}),\]
		\[\forall \bm{c}:~(\bm{c}\in\um{\mcg{B}}~~\text{iff}~~\op{\bm{c}}\in\op{\um{\mcg{B}}});\]
		\rule{0pt}{1mm}
		\[\mcg{s^\leftrightarrow_0}\eqd\mcg{s_0},\]
		\[\forall\op{\bm{c}}\in\op{\um{\mcg{B}}}:~~\mcg{s^\leftrightarrow_1}(\op{\bm{c}})\eqd\op{(\mcg{s_1}(\bm{c}))};\]
		\rule{0pt}{1mm}
		\[\mcg{t^\leftrightarrow_0}\eqd\mcg{t_0},\]
		\[\forall\op{\bm{c}}\in\op{\um{\mcg{B}}}:~~\mcg{t^\leftrightarrow_1}(\op{\bm{c}})\eqd\op{(\mcg{t_1}(\bm{c}))};\]
		\rule{0pt}{1mm}
		\[\mcg{i^\leftrightarrow_0}\eqd\mcg{i_0},\]
		\[\forall\op{f}\in\op{\um{\mcg{A}}}:~~\mcg{i^\leftrightarrow_1}(\op{f})\eqd\op{(\mcg{i_1}(f))};\]
		\rule{0pt}{1mm}
		\[\mcg{\bowtie^\leftrightarrow_0}\eqd\mcg{\bowtie_0},\]
		\[\forall\op{\bm{c}},\op{\bm{d}}\in\op{\um{\mcg{B}}}:~~\op{\bm{c}}\mcg{\bowtie^\leftrightarrow_1}\op{\bm{d}}\eqd\op{(\bm{c}\mcg{\bowtie_1}\bm{d})}.\]
	\end{dfn}
	
	\begin{dfn}
		\label{ver-dual-def}
		Let $\bb{D}=\langle\mcg{A,B,s,t,i,\bowtie}\rangle$ be a double category. The \iw{vertical dual} of $\bb{D}$, denoted by $\bb{D^\updownarrow}$, is defined as
		\[\bb{D^\updownarrow}\eqd((\bb{D^\T})^\leftrightarrow)^\T.\]
	\end{dfn}
	
	\begin{rem}
		\label{ver-dual-data}
		Therefore, the vertical dual $\bb{D^\updownarrow}$ will have the following data:
		\[\bb{D^\updownarrow}=\langle\mcg{A,B^\updownarrow,s^\updownarrow,t^\updownarrow,i^\updownarrow,\bowtie^\updownarrow}\rangle,\]
		where $\mcg{B^\updownarrow}$ is the category consisting of a class of objects $\uo{\mcg{B^\updownarrow}}$ and a class of morphisms $\um{\mcg{B^\updownarrow}}$ such that
		\[\uo{\mcg{B^\updownarrow}}\eqd\op{\um{(\prm{\mcg{A}})}}\]
		where $\prm{\mcg{A}}$ was given in Definition \ref{transpose-def}, and
		\[\forall v:~(v^\updownarrow\in\uo{\mcg{B^\updownarrow}}~~\text{iff}~~v\in\uo{\mcg{B}})\]
		\[\forall v^\updownarrow\in\uo{\mcg{B^\updownarrow}}:~\mcg{s^\updownarrow}(v^\updownarrow)=\mcg{t}(v),\]
		\[\forall v^\updownarrow\in\uo{\mcg{B^\updownarrow}}:~\mcg{t^\updownarrow}(v^\updownarrow)=\mcg{s}(v),\]
		\[\forall X\in\uo{\mcg{A}}:~\mcg{i^\updownarrow}(X)=(\mcg{i}(X))^\updownarrow,\]
		\[\forall\ag{v^\updownarrow}{w^\updownarrow}\in\pb{\uo{\mcg{B^\updownarrow}}}{\uo{\mcg{A}}}{\uo{\mcg{B^\updownarrow}}}:~v^\updownarrow\bowtie^\updownarrow_0w^\updownarrow=(w\bowtie_0 v)^\updownarrow.\]
		Also,
		\[\um{\mcg{B}}^\updownarrow=\op{\um{(\prm{\mcg{B}})}},\]
		where $\prm{\mcg{B}}$ was given in Definition \ref{transpose-def}, and
		\[\forall\bm{c}:~((v^\updownarrow\stackrel{\bm{c}^\updownarrow}{\Longrightarrow}w^\updownarrow)\in\um{\mcg{B^\updownarrow}}~~\text{iff}~~(v\stackrel{\bm{c}}{\Longrightarrow}w)\in\um{\mcg{B}}),\]
		\[\forall\bm{c^\updownarrow}\in\um{\mcg{B^\updownarrow}}:~\mcg{s^\updownarrow}(\bm{c^\updownarrow})=\mcg{t}(\bm{c}),\]
		\[\forall\bm{c^\updownarrow}\in\um{\mcg{B^\updownarrow}}:~\mcg{t^\updownarrow}(\bm{c^\updownarrow})=\mcg{s}(\bm{c}),\]
		\[\forall f\in\um{\mcg{A}}:~\mcg{i^\updownarrow}(f)=(\mcg{i}(f))^\updownarrow;\]
		further, for every composable arrows $\bm{c^\updownarrow},\bm{d^\updownarrow}$ in $\mcg{B^\updownarrow}$ there is a horizontal composition as
		\[\bm{d^\uda}\circ_{\mcg{B^\uda}}\bm{c^\uda}=(\bm{d}\circ_\mcg{B}\bm{c})^\uda;\]
		finally, for every $\langle\bm{c^\uda,d^\uda}\rangle\in\pb{\um{\mcg{B^\uda}}}{\um{\mcg{A}}}{\um{\mcg{B^\uda}}}$ there is a vertical composition given by 
		\[\bm{c}^\uda\bowtie^\uda\bm{d}^\uda=(\bm{d\bowtie c})^\uda.\]
	\end{rem}
	
	\begin{dfn}
		\label{cen-dual-def}
		Let $\bb{D}=\langle\mcg{A,B,s,t,i,\bowtie}\rangle$ be a double category. The \iw{central dual} of $\bb{D}$, denoted by $\bb{D^\boxplus}$, is defined as
		\[\bb{D}^\boxplus\eqd(\bb{D}^\lra)^\uda.\]
	\end{dfn}
	
	\begin{rem}
		\label{D8-group}
		Let $\tu{D}_8$ denote the dihedral group of order 8, corresponding to the \textit{symmetry group of a square} \cite{Arms,Rob}; that is
		\[\tu{D}_8\eqd\langle a,b\mid a^4=b^2=e,~bab=a^{-1}\rangle,\]
		where $a,b$ are two generators and $e$ is the identity element of the group. As stated in \cite{GranPar-Lims}, for any double category $\bb{D}$, the group $\tu{D}_8$ acts on $\bb{D}$. In particular, we have
		\[\bb{D^\boxplus}=(\bb{D}^\lra)^\uda=(\bb{D}^\uda)^\lra;\]
		that is, horizontal dualization commutes with vertical dualization. Other useful relations are
		\[(\bb{D^\lra})^\T=\bb{(D^\T)}^\uda;~~(\bb{D^\T})^\lra=\bb{(D^\uda)}^\T;\]
		\[(\bb{D^\lra})^\lra=(\bb{D^\uda})^\uda=(\bb{D^\boxplus})^\boxplus=(\bb{D^\T})^\T=\bb{D}.\]
	\end{rem}
	
	\rule{0pt}{1mm}
	\begin{rem}
		\label{better-view-cell-struc}
		It is sometimes better to view the notion of double category together with its transpose, its horizontal, vertical, and central duals in the following schematic way. This viewpoint has the advantage that it respects, and clearly manifests, the two-dimensional nature of the cell structures in double categories.\\
		\par Let $\bb{D}=\langle\mcg{A,B,s,t,i,\bowtie}\rangle$ be a double category, and $\bb{D^\T}=\langle\prm{\mcg{A}},\prm{\mcg{B}},\prm{\mcg{s}},\prm{\mcg{t}},\prm{\mcg{i}},\prm{\bowtie}\rangle$ be its transpose. We know that
		\[\um{\mcg{B}}=\um{\prm{\mcg{B}}},\]
		that is, every 2-cell $\bm{c}\in\um{\mcg{B}}$ is at the same time a 2-cell in $\um{\prm{\mcg{B}}}$ and vice versa. Therefore, we can depict $\bm{c}$ as an inclined double-arrow in the cell structure of Diagram \ref{inclined-1}.
		\begin{figure}[h]
			\begin{center}
				\begin{tikzpicture}[commutative diagrams/every diagram]
				\node (N1) at (1,1cm) {$\bullet$};
				\node (N2) at (-1,1cm) {$\bullet$};
				\node (N3) at (-1,-1cm) {$\bullet$};
				\node (N4) at (1,-1cm) {$\bullet$};
				\draw[line width=1pt, -to] (N2) to node[auto] {} (N1);
				\draw[line width=1pt, -to] (N3) to node[auto,swap] {} (N4);
				\draw[line width=1pt, -left to] (N1) to node[auto] {} (N4);
				\draw[line width=1pt, -left to] (N2) to node[auto,swap] {} (N3);
				\node (V1) at (-0.7,0.7cm) {};
				\node (V2) at (0.7,-0.7cm) {};
				\draw[ -implies, double distance=9pt, line width=1pt] (V1) to node {$\bm{c}$} (V2);
				\end{tikzpicture}
			\end{center}
			\caption{} \label{inclined-1}
		\end{figure}
		In this diagram, the inclined double-arrow represents the 2-cell $\bm{c}$ \textit{both as a member of} $\um{\mcg{B}}$ \textit{within} $\bb{D}$ (i.e., as an arrow from the leftmost vertical arrow to the rightmost vertical arrow) \textit{and as a member of} $\prm{\um{\mcg{B}}}$ within $\bb{D^\T}$ (i.e., as an arrow from the uppermost horizontal arrow to the lowermost horizontal arrow). This way, the 2-cell $\bm{c}$ is truly regarded as a ``two-dimensional arrow'' capable of doing two different duties in two different dimensions.\\
		\par Consequently, the cell structures in $\bb{D}$ together with those in its horizontal, vertical, and central duals can easily be depicted as the following (Diagram \ref{inclined-2}).
		\begin{figure}[h]
			\begin{center}
				\begin{tikzpicture}[commutative diagrams/every diagram]
				\node (N1) at (3,3cm) {$\bullet$};
				\node (N2) at (1,3cm) {$\bullet$};
				\node (N3) at (1,1cm) {$\bullet$};
				\node (N4) at (3,1cm) {$\bullet$};
				\draw[line width=1pt, -to] (N1) to node {} (N2);
				\draw[line width=1pt, -to] (N4) to node {} (N3);
				\draw[line width=1pt, -left to] (N1) to node {} (N4);
				\draw[line width=1pt, -left to] (N2) to node {} (N3);
				\node (V1) at (2.7,2.7cm) {};
				\node (V2) at (1.3,1.3cm) {};
				\draw[ -implies, double distance=9pt, line width=1pt] (V1) to node {} (V2);
				
				\node (M1) at (-3,3cm) {$\bullet$};
				\node (M2) at (-1,3cm) {$\bullet$};
				\node (M3) at (-1,1cm) {$\bullet$};
				\node (M4) at (-3,1cm) {$\bullet$};
				\draw[line width=1pt, -to] (M1) to node {} (M2);
				\draw[line width=1pt, -to] (M4) to node {} (M3);
				\draw[line width=1pt, -left to] (M1) to node {} (M4);
				\draw[line width=1pt, -left to] (M2) to node {} (M3);
				\node (U1) at (-2.7,2.7cm) {};
				\node (U2) at (-1.3,1.3cm) {};
				\draw[ -implies, double distance=9pt, line width=1pt] (U1) to node {} (U2);
				
				\node (L1) at (-3,-3cm) {$\bullet$};
				\node (L2) at (-1,-3cm) {$\bullet$};
				\node (L3) at (-1,-1cm) {$\bullet$};
				\node (L4) at (-3,-1cm) {$\bullet$};
				\draw[line width=1pt, -to] (L1) to node[auto] {} (L2);
				\draw[line width=1pt, -to] (L4) to node[auto,swap] {} (L3);
				\draw[line width=1pt, -left to] (L1) to node[auto] {} (L4);
				\draw[line width=1pt, -left to] (L2) to node[auto,swap] {} (L3);
				\node (T1) at (-2.7,-2.7cm) {};
				\node (T2) at (-1.3,-1.3cm) {};
				\draw[ -implies, double distance=9pt, line width=1pt] (T1) to node {} (T2);
				
				\node (O1) at (3,-3cm) {$\bullet$};
				\node (O2) at (1,-3cm) {$\bullet$};
				\node (O3) at (1,-1cm) {$\bullet$};
				\node (O4) at (3,-1cm) {$\bullet$};
				\draw[line width=1pt, -to] (O1) to node[auto] {} (O2);
				\draw[line width=1pt, -to] (O4) to node[auto,swap] {} (O3);
				\draw[line width=1pt, -left to] (O1) to node[auto] {} (O4);
				\draw[line width=1pt, -left to] (O2) to node[auto,swap] {} (O3);
				\node (W1) at (2.7,-2.7cm) {};
				\node (W2) at (1.3,-1.3cm) {};
				\draw[ -implies, double distance=9pt, line width=1pt] (W1) to node {} (W2);
				
				\node (C1) at (-2,0.3cm) {$\bb{D}$};
				\node (C2) at (2,0.3cm) {$\bb{D^\lra}$};
				\node (C3) at (-2,-3.7cm) {$\bb{D^\uda}$};
				\node (C4) at (2,-3.7cm) {$\bb{D^\boxplus}$};
				\end{tikzpicture}
			\end{center}
			\caption{} \label{inclined-2}
		\end{figure}
		(There are also four other cell structures, with the vertical and horizontal 1-cells interchanged, corresponding to the transposes of each of the above cell structures.)
	\end{rem}
	
	Next, we take a look at horizontal isomorphisms.
	
	\begin{dfn}
		\label{hor-iso}
		Let the following be a cell structure in a double category\\ $\bb{D}=\langle\mcg{A,B,s,t,i,\bowtie}\rangle$:
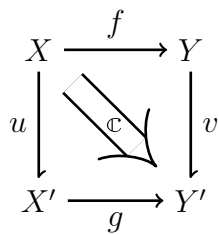
		\begin{center}
			\mbox{
				\begin{tikzpicture}[commutative diagrams/every diagram]
				\node (N1) at (1,1cm) {$Y$};
				\node (N2) at (-1,1cm) {$X$};
				\node (N3) at (-1,-1cm) {$\prm{X}$};
				\node (N4) at (1,-1cm) {$\prm{Y}$};
				\draw[line width=1pt, -to] (N2) to node[auto] {$f$} (N1);
				\draw[line width=1pt, -to] (N3) to node[auto,swap] {$g$} (N4);
				\draw[line width=1pt, -left to] (N1) to node[auto] {$v$} (N4);
				\draw[line width=1pt, -left to] (N2) to node[auto,swap] {$u$} (N3);
				\node (V1) at (-0.7,0.7cm) {};
				\node (V2) at (0.7,-0.7cm) {};
				\draw[ -implies, double distance=9pt, line width=1pt] (V1) to node {$\bm{c}$} (V2);
				\end{tikzpicture}
			}
		\end{center}
		Then, the 2-cell $\bm{c}$ is called a \iw{horizontal isomorphism from} $u$ \textbf{to} $v$ if 
		\[\bm{c}:u\cong v\]
		is an isomorphism in $\mcg{B}$, with inverse $\bm{c}^{-1}$. Note that it follows from this definition and from the functoriality of $\mcg{s,t}$ that:
		\begin{enumerate}
			\item both of the horizontal 1-cells $f,g$ must be isomorphisms in $\mcg{A}$, with inverses $f^{-1},g^{-1}$, respectively; and
			
			\item the boundary of $\bm{c}^{-1}$ shall be $\left(\begin{smallmatrix}
			&f^{-1}& \\u~& &~v\\ &g^{-1}&
			\end{smallmatrix} \right)$.
		\end{enumerate}
		In such case, then, the 2-cell $\left( \begin{smallmatrix}
		&f^{-1}& \\u~&\bm{c}^{-1} &~v\\ &g^{-1}&
		\end{smallmatrix}\right) $ is called the \iw{horizontal inverse} to the 2-cell $\left( \begin{smallmatrix}
		&f& \\u~&\bm{c} &~v\\ &g&
		\end{smallmatrix}\right)$.
	\end{dfn}
	
	\begin{rem}
		\label{simp-hor-iso-flat}
		For \textit{flat} double categories, the above situation simplifies to:
		\[\bm{c}=\begin{pmatrix}
		&f& \\u& &v\\ &g& 
		\end{pmatrix};~~~~\bm{c}^{-1}=\begin{pmatrix}
		&f^{-1}& \\u& &v\\ &g^{-1}& 
		\end{pmatrix}.\]
	\end{rem}
	
	\rule{0pt}{1mm}
	\begin{dfn}
		\label{ver-cen-iso}
		\begin{enumerate}
			\item For the double category $\bb{D}$ as above, the 2-cell $\bm{c}$ is called a \iw{vertical isomorphism} \textbf{from} $f$ \textbf{to} $g$ if $\bm{c}$ is a horizontal isomorphism in \(\bb{D^\T}\).
			
			\item For the double category $\bb{D}$ as above, the 2-cell $\bm{c}$ is called a a \iw{central isomorphism} \textbf{from} $\ag{u}{f}$ \textbf{to} $\ag{v}{g}$ if $\bm{c}$ is both a horizontal isomorphism from $u$ to $v$ and a vertical isomorphism from $f$ to $g$.
		\end{enumerate}
	\end{dfn}
	
	\begin{rem}
		\label{other-isos}
		There are also other kinds of isomorphisms and equivalences definable in double categories. See for example \cite{GranPar-Lims}.
	\end{rem}
	
	\rule{0pt}{5mm}
	Now we introduce another fundamental concept:
	
	\begin{dfn}
		\label{hor-functor}
		Let $\bb{D}=\langle\mcg{A,B,s,t,i,\bowtie}\rangle$ and $\bb{E}=\langle\mcg{X,Y,p,q,j,\circledast}\rangle$ be two double categories. A \iw{double functor} $\mcg{F}:\bb{D\longrightarrow E}$ consists of the following data:
		\[\mcg{F}\eqd\ag{\mcg{F^0}}{\mcg{F^1}},\]
		where $\mcg{F^0:A\longrightarrow X}$ and $\mcg{F^1:B\longrightarrow Y}$ are functors that make Diagram \ref{hor-functor-diag} commute.
		\begin{figure}[h]
			\begin{center}
				\begin{tikzpicture}[commutative diagrams/every diagram]
				
				\node (N1) at (0.8,0.3cm) {$\mcg{Y}$};
				
				\node (N2) at (-1.3,-0.3cm) {$\mcg{B}$};
				
				\node (N3) at (4.7,-0.3cm) {$\mcg{Y}$};
				
				\node (N4) at (-5.2,0.3cm) {$\mcg{B}$};
				
				\node (N5) at (3,2cm) {$\pb{\mcg{Y}}{\mcg{X}}{\mcg{Y}}$};
				
				\node (N6) at (-3.5,-2cm) {$\mcg{A}$};
				
				\node (N7) at (-3,2cm) {$\pb{\mcg{B}}{\mcg{A}}{\mcg{B}}$};
				
				\node (N8) at (2.5,-2cm) {$\mcg{X}$};
				
				\node (N9) at (3,5cm) {$\mcg{Y}$};
				
				\node (N10) at (-3.5,-5cm) {$\mcg{B}$};
				
				\node (N11) at (-3,5cm) {$\mcg{B}$};
				
				\node (N12) at (2.5,-5cm) {$\mcg{Y}$};
				
				\path[commutative diagrams/.cd, every arrow, every label]
				(N7) edge node{$\bowtie$} (N11)
				(N5) edge node[swap]{$\circledast$} (N9)
				(N7) edge node[swap]{$\pi_1$} (N4)
				(N7) edge node{$\pi_2$} (N2)
				(N5) edge node[swap]{$\prm{\pi}_1$} (N1)
				(N5) edge node{$\prm{\pi}_2$} (N3)
				(N4) edge node[swap]{$\mcg{t}$} (N6)
				(N2) edge node{$\mcg{s}$} (N6)
				(N1) edge node[swap]{$\mcg{q}$} (N8)
				(N3) edge node{$\mcg{p}$} (N8)
				(N6) edge node[swap]{$\mcg{i}$} (N10)
				(N8) edge node{$\mcg{j}$} (N12)
				(N4) edge[-,line width=6pt,draw=white] node{} (N1)
				(N2) edge[-,line width=6pt,draw=white] node{} (N3)
				(N4) edge[line width=1pt] node[near end]{$\mcg{F^1}$} (N1)
				(N2) edge[line width=1pt] node[near start,swap]{$\mcg{F^1}$} (N3)
				(N11) edge[line width=1pt] node{$\mcg{F^1}$} (N9)
				(N10) edge[line width=1pt] node{$\mcg{F^1}$} (N12)
				(N6) edge[line width=1pt] node{$\mcg{F^0}$} (N8)
				(N7) edge[line width=1pt] node{$(\mcg{F^1}\pi_1,\mcg{F^1}\pi_2)$} (N5);
				\end{tikzpicture}
			\end{center}
			\caption{} \label{hor-functor-diag}
		\end{figure}
		In this diagram, $\agg{\pb{\mcg{B}}{\mcg{A}}{\mcg{B}}}{\pi_1}{\pi_2}$ is the pullback of $\mcg{t,s}$, while $\agg{\pb{\mcg{Y}}{\mcg{X}}{\mcg{Y}}}{\prm{\pi}_1}{\prm{\pi}_2}$ is the pullback of $\mcg{q,p}$; also, \[(\mcg{F_1}\pi_1,\mcg{F_1}\pi_2):\pb{\mcg{B}}{\mcg{A}}{\mcg{B}}\longrightarrow\pb{\mcg{Y}}{\mcg{X}}{\mcg{Y}}\] is the unique arrow according to the universality of $\agg{\pb{\mcg{Y}}{\mcg{X}}{\mcg{Y}}}{\prm{\pi}_1}{\prm{\pi}_2}$.
	\end{dfn}
	
	\begin{rem}
		\label{horfun-also-diagram}
		A double functor $\mcg{F}:\bb{D\longrightarrow E}$ is sometimes called a \iw{double diagram} in $\bb{E}$, especially when working with limits and colimits. The double diagram $\mcg{F}$ is said to be \textbf{small} whenever $\bb{D}$ is a small double category.
	\end{rem}
	
	\begin{dfn}
		\label{identity-fun-hor}
		\begin{enumerate}
			\item Let $\bb{D}=\langle\mcg{A,B,s,t,i,\bowtie}\rangle$ be a double category. The double functor 
			\[\tu{Id}_\bb{D}\eqd\ag{1_\mcg{A}}{1_\mcg{B}}\]
			is called the \iw{identity double functor} on $\bb{D}$.
			
			\item We define a \iw{composition} for double functors; that is, for consecutive $\mcg{F,G}$:
			\[\mcg{G\circ F}\eqd\ag{\mcg{G^0F^0}}{\mcg{G^1F^1}}.\]
			It is immediate that $\mcg{G\circ F}$ satisfies the commutativity conditions of Diagram \ref{hor-functor-diag}, and, hence, is well-defined.
			
			\item A double functor $\mcg{F}:\bb{D\longrightarrow E}$ is called an \iw{isomorphism double functor} if there exists another double functor $\mcg{G}:\bb{E\longrightarrow D}$ such that 
			\[\mcg{F\circ G}=\tu{Id}_\bb{E},~~~~\mcg{G\circ F}=\tu{Id}_\bb{D}.\]
			In such case, we may denote the isomorphism situation by $\mcg{F}:\bb{D\cong E}$.
			
			\item A \iw{constant double functor} $\mcg{K}=\ag{\mcg{K^0}}{\mcg{K^1}}:\bb{D\longrightarrow E}$ is a double functor in which both $\mcg{K^0,K^1}$ are constant functors.
		\end{enumerate}
	\end{dfn}
	
	\begin{prp}
		\label{hor-iso-iff-both-iso}
		A double functor $\mcg{F}=\ag{\mcg{F^0}}{\mcg{F^1}}$ is an isomorphism if and only if $\mcg{F^0,F^1}$ are isomorphisms as ordinary functors.
	\end{prp}
	
	\begin{mypr}
		The ``Only if'' side is obvious. We prove the ``If'' side. Assume that $\bb{D,E}$ are double categories and $\mcg{F}=\ag{\mcg{F^0}}{\mcg{F^1}}:\bb{D\longrightarrow E}$ is a double functor from $\bb{D}$ to $\bb{E}$ such that both $\mcg{F^0,F^1}$ are (ordinary) isomorphism functors. Thus, $\mcg{F^0,F^1}$ have corresponding inverses, namely $(\mcg{F^0})^{-1},(\mcg{F^1})^{-1}$. It suffices to show that the pair 
		\[\mcg{F}^{-1}\eqd\ag{(\mcg{F^0})^{-1}}{(\mcg{F^1})^{-1}}\]
		actually forms a double functor from $\bb{E}$ to $\bb{D}$. Referring to Diagram \ref{hor-functor-diag}, we do the following observation: by replacing all $\mcg{F^0,F^1}$ with $(\mcg{F^0})^{-1},(\mcg{F^1})^{-1}$, the unique arrow $(\mcg{F_1}\pi_1,\mcg{F_1}\pi_2):\pb{\mcg{B}}{\mcg{A}}{\mcg{B}}\longrightarrow\pb{\mcg{Y}}{\mcg{X}}{\mcg{Y}}$ finds an inverse 
		\[\left( (\mcg{F^1})^{-1}\prm{\pi}_1,(\mcg{F^1})^{-1}\prm{\pi}_2\right) :\pb{\mcg{Y}}{\mcg{X}}{\mcg{Y}}\longrightarrow\pb{\mcg{B}}{\mcg{A}}{\mcg{B}}.\]
		Consequently, the whole diagram again commutes, and therefore, \(\mcg{F}^{-1}:\bb{E\longrightarrow D}\) is a well-defined double functor.
	\end{mypr}
	
	\begin{dfn}
		\label{transpose-of-functor}
		Let $\mcg{F}:\bb{D\longrightarrow E}$ be a double functor between double categories $\bb{D}=\langle\mcg{A,B,s,t,i,\bowtie}\rangle$ and $\bb{E}=\langle\mcg{X,Y,p,q,j,\circledast}\rangle$. The \iw{transpose} of $\mcg{F}$ is the double functor $\mcg{F}^\T:\bb{D^\T\longrightarrow E^\T}$ defined canonically. That is,
		\[\mcg{F^\T}\eqd\ag{(\mcg{F^\T})^0}{(\mcg{F^\T})^1}:\langle\prm{\mcg{A}},\prm{\mcg{B}},\prm{\mcg{s}},\prm{\mcg{t}},\prm{\mcg{i}},\prm{\bowtie}\rangle\longrightarrow\langle\prm{\mcg{X}},\prm{\mcg{Y}},\prm{\mcg{p}},\prm{\mcg{q}},\prm{\mcg{j}},\prm{\circledast}\rangle,\]
		with functors $(\mcg{F^\T})^0:\prm{\mcg{A}}\longrightarrow\prm{\mcg{X}}$ and $(\mcg{F^\T})^1:\prm{\mcg{B}}\longrightarrow\prm{\mcg{Y}}$ defined as the following:
		\begin{align*}
		\forall C\in\uo{\prm{\mcg{A}}}:~~(\mcg{F^\T})^0(C)&\eqd\mcg{F^0}(C);\\
		\forall u\in\um{\prm{\mcg{A}}}:~~(\mcg{F^\T})^0(u)&\eqd\mcg{F^1}(u);\\
		\forall f\in\uo{\prm{\mcg{B}}}:~~(\mcg{F^\T})^1(f)&\eqd\mcg{F^0}(f);\\
		\forall \bm{c}\in\um{\prm{\mcg{B}}}:~~(\mcg{F^\T})^1(\bm{c})&\eqd\mcg{F^1}(\bm{c}).
		\end{align*}
	\end{dfn}
	
	It is easy to check the commutativity conditions of \ref{hor-functor-diag} for $\mcg{F^\T}$. For example, the equation $(\mcg{F^\T})^0\prm{\mcg{s}}=\prm{\mcg{p}}(\mcg{F^\T})^1$ follows from the definition of $\prm{\mcg{s}},\prm{\mcg{p}}$ and from the functoriality of $\mcg{F^0,F^1}$, and so on.\\
	
	Now we introduce a concept from group theory \cite{Arms,Rob}.
	
	\begin{dfn}
		\label{Klein-def}
		The \iw{Klein four-group} is the group defined as
		\[\tu{V}_4\eqd\langle a,b \mid a^2=b^2=[a,b]=e\rangle,\]
		where $a,b$ are two generators, and $[a,b]\eqd aba^{-1}b^{-1}$ is the \textit{commutator} of $a,b$. The underlying set of $\tu{V}_4$ has four elements (hence the name ``four-group''):
		\[\tu{V}_4=\{e,a,b,ab\}.\]
	\end{dfn}
	
	$\tu{V}_4$ is \textit{abelian}, a \textit{subgroup} of the dihedral group $\tu{D}_8$, and the \textit{smallest non-cyclic group}. Also, geometrically, $\tu{V}_4$ is the \textit{symmetry group} of a rhombus and of a rectangle which are \textit{not} squares, the four elements of the symmetry group being identity, vertical reflection, horizontal reflection, and 180\textdegree~ rotation (Diagram \ref{rhombu-rectang}).
	\begin{figure}[H]
		\begin{center}
			\begin{tikzpicture}
			
			\node (C1) at (3.5,0cm) {$\bullet$};
			\node (C2) at (-3.5,0cm) {$\bullet$};
			
			\draw (5.5,1) -- (1.5,1) -- (1.5,-1) -- (5.5,-1) -- (5.5,1);
			
			\draw (-2.5,0) -- (-3.5,2) -- (-4.5,0) -- (-3.5,-2) -- (-2.5,0);
			
			\draw[dashed] (3.5,2) -- (3.5,-2);
			\draw[dashed] (0.5,0) -- (6.5,0);
			
			\draw[dashed] (-3.5,3) -- (-3.5,-3);
			\draw[dashed] (-1.5,0) -- (-5.5,0);			
			\end{tikzpicture}
		\end{center}
		\caption{} \label{rhombu-rectang}
	\end{figure}
	
	We need the group $\tu{V}_4$ for the next definition:
	
	\begin{dfn}
		\label{Klein-inv}
		A double category $\bb{D}$ is said to be \iw{Klein-invariant} or $\tu{V}_4$-\textbf{invariant}\index{$\tu{V}_4$-invariant} if there exist the following isomorphism double functors:
		\begin{align*}
		(-)^\#:\bb{D^\lra}&\cong \bb{D},\\
		(-)^\star:\bb{D^\uda}&\cong \bb{D}.
		\end{align*}
		$(-)^\#$ and $(-)^\star$ are called the \iw{horizontal self-duality} and the \iw{vertical self-duality} of $\bb{D}$, respectively.
	\end{dfn}
	
	\begin{rem}
		\label{cen-self-dua}
		It follows immediately that for every Klein-invariant double category $\bb{D}$ we have
		\[\bb{D\cong D^\lra\cong D^\uda\cong D^\boxplus},\]
		where the composite isomorphism
		\[(-)^\circledcirc\eqd(-)^\#\circ(-)^\star=(-)^\star\circ(-)^\#:\bb{D^\boxplus\cong D}\]
		is called the \iw{central self-duality} of $\bb{D}$.
	\end{rem}
	
	There is another useful definition (which is not used in the present work but may be useful for future development):
	
	\begin{dfn}
		\label{self-tr-def}
		\begin{enumerate}
			\item A double category $\bb{D}$ is said to be \iw{self-transpose} if \[\bb{D^\T\cong D}.\]
			
			\item $\bb{D}$ is said to be \iw{strictly self-transpose} if \[\bb{D^\T}=\bb{D}.\]
		\end{enumerate}
	\end{dfn}
	
	Next, we turn to transformations between double functors \cite{GranPar-Lims}.
	
	\begin{dfn}
		\label{hor-NT}
		Let $\mcg{F,G}:\bb{D\dlong E}$ be double functors from $\bb{D}=\langle\mcg{A,B,s,t,i,\bowtie}\rangle$ to $\bb{E}=\langle\mcg{X,Y,p,q,j,\circledast}\rangle$. A \iw{horizontal (natural) transformation} $\fk{h}:\mcg{F}\Rrightarrow\mcg{G}$ assigns
		\begin{enumerate}
			\item a family of $\mcg{X}$-morphisms (horizontal 1-cells) $\fk{h}_C:\mcg{F^0}(C)\longrightarrow \mcg{G^0}(C)$ for each $\mcg{A}$-object (0-cell) $C$,
			
			\item a family of $\mcg{Y}$-morphisms (2-cells) $\left(\begin{smallmatrix}
			&\fk{h}_C& \\ \mcg{F^1}(u)&\fk{h}_u&\mcg{G^1}(u)\\ &\fk{h}_\prm{C}& 
			\end{smallmatrix} \right) $ for each $\mcg{B}$-object (vertical 1-cell) $u:C\rightharpoonup\prm{C}$,
		\end{enumerate}
		so that the following preservation and naturality conditions hold:
		\begin{itemize}
			\item[\textbf{(HT1)}~~] $\left(\begin{smallmatrix}
			&\fk{h}_C& \\ \mcg{F^1}(\mcg{i}(C))&\fk{h}_u&\mcg{G^1}(\mcg{i}(C))\\ &\fk{h}_\prm{C}& 
			\end{smallmatrix} \right) = \left(\begin{smallmatrix}
			&\fk{h}_C& \\ \mcg{j}(\mcg{F^0}(C))&\fk{h}_u&\mcg{j}(\mcg{G^0}(C))\\ &\fk{h}_\prm{C}& 
			\end{smallmatrix} \right) $ for every $C\in\mcg{A}$;
			
			\item[\textbf{(HT2)}~~] $\fk{h}_{u\bowtie v}=\fk{h}_u\circledast\fk{h}_v$ for all $u,v\in\uo{\mcg{B}}$;
			
			\item[\textbf{(HT3)}~~] $\fk{h}_v\circ\mcg{F^1}(\bm{c})=\mcg{G^1}(\bm{c})\circ\fk{h}_u$ for every $\left(\begin{smallmatrix}
			&f& \\u&\bm{c}&v\\ &g& 
			\end{smallmatrix} \right) $ in $\bb{D}$ (Diagram \ref{hor-NT-diag}).
			\begin{figure}[h]
				\begin{center}
					\begin{tikzpicture}[commutative diagrams/every diagram]
					\node (N1) at (3,1cm) {$\bullet$};
					\node (N2) at (0,1cm) {$\bullet$};
					\node (N3) at (-3,1cm) {$\bullet$};
					\node (N4) at (3,-1cm) {$\bullet$};
					\node (N5) at (0,-1cm) {$\bullet$};
					\node (N6) at (-3,-1cm) {$\bullet$};
					
					\draw[line width=1pt, -to] (N2) to node[auto] {} (N1);
					
					\draw[line width=1pt, -left to] (N1) to node[auto] {$\mcg{G^1}(v)$} (N4);
					
					\draw[line width=1pt, -to] (N3) to node[auto] {} (N2);
					
					\draw[line width=1pt, -left to] (N2) to node {$\mcg{F^1}(v)$} (N5);
					
					\draw[line width=1pt, -to] (N5) to node[auto,swap] {} (N4);
					
					\draw[line width=1pt, -left to] (N3) to node[auto,swap] {$\mcg{F^1}(u)$} (N6);
					
					\draw[line width=1pt, -to] (N6) to node[auto,swap] {} (N5);
					
					\node (V1) at (1.5,0cm) {$\fk{h}_v$};
					\node (V2) at (-1.5,0cm) {$\mcg{F^1}(\bm{c})$};

					\node (E1) at (0,-1.5cm) {};
					\node (E2) at (0,-2.5cm) {};
					
					\draw[-, double distance=4pt, line width=3pt] (E1) to node {} (E2);
					
					\node (P1) at (3,-3cm) {$\bullet$};
					\node (P2) at (0,-3cm) {$\bullet$};
					\node (P3) at (-3,-3cm) {$\bullet$};
					\node (P4) at (3,-5cm) {$\bullet$};
					\node (P5) at (0,-5cm) {$\bullet$};
					\node (P6) at (-3,-5cm) {$\bullet$};
					
					\draw[line width=1pt, -to] (P2) to node[auto] {} (P1);
					
					\draw[line width=1pt, -left to] (P1) to node[auto] {$\mcg{G^1}(v)$} (P4);
					
					\draw[line width=1pt, -to] (P3) to node[auto] {} (P2);
					
					\draw[line width=1pt, -left to] (P2) to node {$\mcg{G^1}(u)$} (P5);
					
					\draw[line width=1pt, -to] (P5) to node[auto,swap] {} (P4);
					
					\draw[line width=1pt, -left to] (P3) to node[auto,swap] {$\mcg{F^1}(u)$} (P6);
					
					\draw[line width=1pt, -to] (P6) to node[auto,swap] {} (P5);
					
					\node (W1) at (1.5,-4cm) {$\mcg{G^1}(\bm{c})$};
					\node (W2) at (-1.5,-4cm) {$\fk{h}_u$};
					\end{tikzpicture}
				\end{center}
				\caption{} \label{hor-NT-diag}
			\end{figure}
		\end{itemize}
	\end{dfn}
	
	\begin{rem}
		\label{produce-NT}
		It follows that the arrows $\fk{h}_C:\mcg{F^0}C\longrightarrow\mcg{G^0}C$ give a natural transformation $\mcg{F^0}\xrightarrow{~\bullet~}\mcg{G^0}$.
	\end{rem}
	
	\begin{rem}
		\label{flat-HTx}
		If $\bb{E}$ is flat, then \textbf{(HT1)} and \textbf{(HT2)} are trivially satisfied while \textbf{(HT3)} reduces to ordinary naturality. A horizontal transformation $\fk{h}:\mcg{F\Rrightarrow G}$ reduces thus to a natural transformation $\langle\fk{h}_C:\mcg{F^0}C\longrightarrow\mcg{G^0}C\rangle_C$ such that, for every vertical 1-cell $u:C\rightharpoonup\prm{C}$ in $\bb{D}$, the boundary 
		\[\begin{pmatrix}
		&\fk{h}_C& \\ \mcg{F^1}u& &\mcg{G^1}u\\ &\fk{h}_\prm{C}&
		\end{pmatrix}\]
		determines a unique 2-cell.
	\end{rem}
	
	Horizontal transformations compose canonically:
	
	\begin{dfn}
		\label{horNT-compo}
		Given $\fk{h}=\ag{\langle \fk{h}_C\rangle_C}{\langle \fk{h}_u\rangle_u}\mcg{F}\Rrightarrow\mcg{G}$ and $\fk{k}=\ag{\langle \fk{k}_C\rangle_C}{\langle \fk{k}_u\rangle_u}\mcg{G}\Rrightarrow\mcg{H}$, their \iw{(horizontal) composition} is defined as
		\[\fk{k\circ h}\eqd\ag{\langle\fk{k}_C\circ\fk{h}_C \rangle_C}{\langle \fk{k}_u\circ\fk{h}_u\rangle_u}:\mcg{F}\Rrightarrow\mcg{H}.\]
	\end{dfn}
	
	\rule{0pt}{1mm}
	Likewise, one can define \textit{identity horizontal (natural) transformations}, \textit{horizontal natural isomorphisms}, etc.\\
	
	\par Next, we take a look at (horizontal) limits and colimits.
	
	\begin{dfn}
		\label{hor-source-sink}
		Let $\mcg{F}:\bb{D\longrightarrow E}$ be a double diagram in $\bb{E}$ (see Remark \ref{horfun-also-diagram}). Also, let $\tu{K}:\bb{D\longrightarrow E}$ be a constant double functor.
		\begin{enumerate}
			\item A \iw{horizontal natural source} for $\mcg{F}$ is a horizontal transformation $\fk{h}:\tu{K}\Rrightarrow\mcg{F}$.
			
			\item A \iw{horizontal natural sink} for $\mcg{F}$ is a horizontal transformation $\fk{h}:\mcg{F}\Rrightarrow\tu{K}$.
			
			\item A \iw{horizontal limit} for $\mcg{F}$ is a horizontal natural source $\tu{L}\stackrel{~\fk{l}~}{\Rrightarrow}\mcg{F}$ which is \textit{universal} in the sense that for any other horizontal natural source $\tu{K}\stackrel{~\fk{k}~}{\Rrightarrow}\mcg{F}$, there exists a \textit{unique} horizontal transformation $\fk{h}:\tu{K}\Rrightarrow\tu{L}$ such that \[\fk{k=l\circ h}.\]
			
			\item A \iw{horizontal colimit} for $\mcg{F}$ is a horizontal natural sink $\mcg{F}\stackrel{~\fk{c}~}{\Rrightarrow}\tu{C}$ which is \textit{universal} in the sense that for any other horizontal natural sink $\mcg{F}\stackrel{~\fk{k}~}{\Rrightarrow}\tu{K}$, there exists a \textit{unique} horizontal transformation $\fk{h}:\tu{C}\Rrightarrow\tu{K}$ such that \[\fk{k=h\circ c}.\]
		\end{enumerate}
	\end{dfn}
	
	This way, we can have (binary) horizontal products and coproducts, horizontal pullbacks and pushouts, etc.
	
	\rule{0pt}{5mm}
	
	\section{The main formalism}
	\label{the-main-forma}
	First of all, we introduce an auxiliary category, which will be used in the sequel.
	
	\begin{dfn}
		\label{fluid-def}
		Consider the category $\chu=\chu_\Gamma$ for some fixed nonempty $\Gamma$. By the \iw{fluid category} on \chu, denoted by $\cat{Fluid}_\chu$, we mean the category which consists of the following data:
		\begin{itemize}
			\item Objects: all dialgebras of the form $\langle\msf{A},F\msf{A}\xrightarrow{~\alpha~}G\msf{A}\rangle$, where $\msf{A}\in\chu$ is a Chu space, $F$ and $G$ are Chu endofunctors, and $\alpha$ is a Chu transform. More precisely, we take
			\[\uo{\left( \cat{Fluid}_\chu\right) }~\eqd~\dju_{F,G}\uo{\left( \dlg{\chu}{F}{G}\right) },\]
			where the right hand side is the disjoint union of classes.
			
			\item Morphisms: if $\ag{\msf{A}}{\alpha}, \ag{\msf{B}}{\beta}$ are two objects with $\ag{\msf{A}}{\alpha}\in\dlg{\chu}{F^1}{G^1}$ and\\ $\ag{\msf{B}}{\beta}\in\dlg{\chu}{F^2}{G^2}$ for endofunctors $F^1,F^2,G^1,G^2$, then 
			\begin{itemize}
				\item for $\msf{A\ne B}$, there are no morphisms between $\ag{\msf{A}}{\alpha}, \ag{\msf{B}}{\beta}$;
				
				\item but for $\msf{A=B}$, the morphisms are pairs of the form\\ $\ag{\varphi}{\psi}:\ag{\msf{A}}{\alpha}\longrightarrow \ag{\msf{B}}{\beta}$, where 
				\[\varphi:F^1\msf{A}\longrightarrow F^2\msf{A},~\psi:G^1\msf{A}\longrightarrow G^2\msf{A}\]
				are Chu transforms that make Diagram \ref{fluid-def-diag} commute.
				\begin{figure}[H]
					\begin{center}
						\begin{tikzpicture}[commutative diagrams/every diagram]
						
						\node (N1) at (3,2cm) {$G^1\msf{A}$};
						
						\node (N2) at (-1,2cm) {$F^1\msf{A}$};
						
						\node (N3) at (-1,-1cm) {$F^2\msf{A}$};
						
						\node (N4) at (3,-1cm) {$G^2\msf{A}$};
						
						\node (N5) at (-3,2cm) {$\msf{A}$};
						
						\node (N6) at (-3,-1cm) {$\msf{A}$};
						
						\path[commutative diagrams/.cd, every arrow, every label]
						(N2) edge node{$\alpha$} (N1)
						(N2) edge node[swap]{$\varphi$} (N3)
						(N3) edge node{$\beta$} (N4)
						(N1) edge node{$\psi$} (N4)
						(N5) edge node[swap]{$1_\msf{A}$} (N6);
						\end{tikzpicture}
					\end{center}
					\caption{} \label{fluid-def-diag}
				\end{figure}
				The pair $\ag{\varphi}{\psi}$ is called a \iw{local transformation} from $\langle F^1,G^1\rangle$ to $\langle F^2,G^2\rangle$.
			\end{itemize}
			
			\item Composition: the pairs $\ag{\varphi}{\psi}:\ag{\msf{A}}{\alpha^1}\longrightarrow\ag{\msf{A}}{\alpha^2}$ and $\ag{\prm{\varphi}}{\prm{\psi}}:\ag{\msf{A}}{\alpha^2}\longrightarrow\ag{\msf{A}}{\alpha^3}$ compose componentwise:
			\[\ag{\prm{\varphi}}{\prm{\psi}}\circ\ag{\varphi}{\psi}\eqd\ag{\prm{\varphi}\varphi}{\prm{\psi}\psi}.\]
			
			\item Identities: pairs $\ag{1_{F\msf{A}}}{1_{G\msf{A}}}$ for each object $\ag{\msf{A}}{F\msf{A}\xrightarrow{~\alpha~}G\msf{A}}$.
		\end{itemize}
	\end{dfn}
	
	\rule{0pt}{2mm}
	
	The category $\cat{Fluid}_\chu$ is intended to serve as one of the resources for our main formalism:
	
	\begin{dfn}
		\label{DLC-def}
		Consider the category $\chu=\chu_\Gamma$ for some fixed nonempty set $\Gamma$. Also, let $\Sigma\in\cat{Set}$ be a fixed set (possibly empty). The (properly large) double category $\dlc_{\Gamma,\Sigma}$ is defined by the following data:
		\[\dlc_{\Gamma,\Sigma}\eqd\langle\mcg{A,B,s,t,i,\bowtie}\rangle,\]
		such that categories $\mcg{A,B}$ and functors $\mcg{s,t,i,\bowtie}$ are defined as below:\\
		
		\par \noindent \textbf{a. The category} $\mcg{A}$\textbf{.} The category $\mcg{A}$ (the ``object of objects'') is defined as the following:\\
		
		\par \noindent \textbf{a.1. Objects of} $\mcg{A}$\textbf{.} The $\mcg{A}$-objects (the 0-cells) are pairs
		\[\eb{P}\eqd\eb{\ag{D}{S}}\]
		such that\\
		\textbf{a.1.1.} $\eb{D}=\left[ \begin{smallmatrix}
		\eb{A}\\ \eb{B}
		\end{smallmatrix}\right] $ is a pair called a \iw{dialgebra column}, where $\eb{A}=\ag{\msf{A}}{\alpha}$ is an\\ $\ag{F}{G}$-dialgebra in $\dlg{\chu}{F}{G}$ for endofunctors $F,G$, while $\eb{B}=\ag{\msf{B}}{\beta}$ is a $\ag{K}{H}$-dialgebra in $\dlg{\chu}{K}{H}$ for endofunctors $K,H$;\\
		\textbf{a.1.2.} $\eb{S}=\ag{S^+}{S^-}$ is a pair of functions
		\begin{align*}
		S^+:\Sigma&\longrightarrow F^+\msf{A}\sqcup K^+\msf{B},\\
		S^-:\Sigma&\longrightarrow G^-\msf{A}\sqcup H^-\msf{B};
		\end{align*}
		$\eb{S}$ is called a \iw{super-matrix}. Furthermore, the \iw{functorial profile} of $\eb{P}$ is defined as
		\[\boldsymbol{\Phi}_\eb{P}\eqd\begin{bmatrix}
		F&G\\ K&H
		\end{bmatrix}.\]
		
		\rule{0pt}{1mm}
		\par \noindent \textbf{a.2. Morphisms of} $\mcg{A}$\textbf{.} For 0-cells $\eb{P}_1=\ag{\eb{D}_1}{\eb{S}_1}$ and $\eb{P}_2=\ag{\eb{D}_2}{\eb{S}_2}$, if $\boldsymbol{\Phi}_{\eb{P}_1}\ne\boldsymbol{\Phi}_{\eb{P}_2}$, then there are no $\mcg{A}$-morphisms between them; but if $\boldsymbol{\Phi}_{\eb{P}_1}=\boldsymbol{\Phi}_{\eb{P}_2}=\left[\begin{smallmatrix}
		F&G\\ K&H
		\end{smallmatrix} \right] $ with
		\[\eb{D}_j=\begin{bmatrix}
		\eb{A}_j\\ \eb{B}_j
		\end{bmatrix},~~\eb{A}_j\in\dlg{\chu}{F}{G},~~ \eb{B}_j\in\dlg{\chu}{K}{H},~~ j=1,2,\] 
		the $\mcg{A}$-morphisms (the horizontal 1-cells) are pairs of the form
		\[\bm{h}:\eb{P}_1\longrightarrow\eb{P}_2,\]
		where $\bm{h}=\left[ \begin{smallmatrix}
		\eb{f}\\\eb{g}
		\end{smallmatrix}\right] $ is a pair such that $\eb{f}:\eb{A}_1\longrightarrow\eb{A}_2$ is an $\ag{F}{G}$-dialgebra homomorphism in $\dlg{\chu}{F}{G}$ while $\eb{g}:\eb{B}_2\longrightarrow\eb{B}_1$ is a $\ag{K}{H}$-dialgebra homomorphism in $\dlg{\chu}{K}{H}$, and Diagram \ref{PHS-NHS} commutes, in which the left and right diamonds are called the \textbf{Positive} and \textbf{Negative Horizontal Super-adjointness} conditions, or briefly \textbf{PHS}\index{Positive Horizontal Super-adjointness Condition (PHS)} and \textbf{NHS}\index{Negative Horizontal Super-adjointness Condition (NHS)}, respectively (in this diagram, also, $f,g$ are the stems of $\eb{f,g}$, respectively).\\
		\begin{figure}[h]
			\begin{center}
				\begin{tikzpicture}[commutative diagrams/every diagram]
				
				\node (N1) at (0,0cm) {$\Sigma$};
				
				\node (N2) at (-2.4,2cm) {$F^+\msf{A}_1\sqcup K^+\msf{B}_1$};
				
				\node (N3) at (-4.8,0cm) {$F^+\msf{A}_2\sqcup K^+\msf{B}_1$};
				
				\node (N4) at (-2.4,-2cm) {$F^+\msf{A}_2\sqcup K^+\msf{B}_2$};
				
				\node (N5) at (2.4,2cm) {$G^-\msf{A}_1\sqcup H^-\msf{B}_1$};
				
				\node (N6) at (4.8,0cm) {$G^-\msf{A}_1\sqcup H^-\msf{B}_2$};
				
				\node (N7) at (2.4,-2cm) {$G^-\msf{A}_2\sqcup H^-\msf{B}_2$};
				
				\node (N8) at (-2.4,0cm) {\fbox{\textbf{PHS}}};
				
				\node (N9) at (2.4,0cm) {\fbox{\textbf{NHS}}};
				
				\path[commutative diagrams/.cd, every arrow, every label]
				(N1) edge node[swap]{$S^+_1$} (N2)
				(N2) edge node[swap]{$F^+f\sqcup 1$} (N3)
				(N1) edge node{$S^+_2$} (N4)
				(N4) edge node{$1\sqcup K^+g$} (N3)
				(N1) edge node{$S^-_1$} (N5)
				(N5) edge node{$1\sqcup H^-g$} (N6)
				(N1) edge node[swap]{$S_2^-$} (N7)
				(N7) edge node[swap]{$G^-f\sqcup 1$} (N6);
				\end{tikzpicture}
			\end{center}
			\caption{} \label{PHS-NHS}
		\end{figure}
		
		\par \noindent \textbf{a.3. Composition in} $\mcg{A}$\textbf{.} For horizontal 1-cells $\bm{h,k}$ with $\bm{h}=\left[\begin{smallmatrix}
		\eb{f}_1\\\eb{g}_1
		\end{smallmatrix} \right], \bm{k}=\left[\begin{smallmatrix}
		\eb{f}_2\\\eb{g}_2
		\end{smallmatrix} \right]$, their composition is given by
		\[\bm{k\circ h}\eqd\begin{bmatrix}
		\eb{f}_2\eb{f}_1\\\eb{g}_1\eb{g}_2
		\end{bmatrix} ,\]
		with composition in each row being done in the corresponding category. The above composition satisfies the positive and negative super-adjointness conditions and, hence, is well-defined (see Step 1 of the proof of Proposition \ref{DLC-welldef} below).\\
		
		\par \noindent \textbf{a.4. Identities in} $\mcg{A}$\textbf{.} For a 0-cell $\eb{P}=\ag{\eb{D}}{\eb{S}}$,
		\[1_\eb{P}\eqd\begin{bmatrix}
		1_\eb{A}\\1_\eb{B}
		\end{bmatrix},\]
		where $1_\eb{A},1_\eb{B}$ are the identity dialgebra homomorphisms on $\eb{A},\eb{B}$, respectively.\\
		
		\rule{0pt}{1mm}
		\par \noindent \textbf{b. The category} $\mcg{B}$\textbf{.} The category $\mcg{B}$ (the ``object of arrows'') is defined as the following. Before proceeding we introduce an auxiliary category:\\
		
		\par \noindent \textbf{b.1. The auxiliary category} $\prm{\mcg{A}}$\textbf{.} Consider the category $\prm{\mcg{A}}$ consisting of the following data:
		\begin{itemize}
			\item $\prm{\mcg{A}}$-objects: the same as objects of $\mcg{A}$:
			\[\uo{\prm{\mcg{A}}}\eqd\uo{\mcg{A}}.\]
			
			\item $\prm{\mcg{A}}$-morphisms: for 0-cells $\eb{P}^1=\ag{\eb{D}^1}{\eb{S}^1}$ and $\eb{P}^2=\ag{\eb{D}^2}{\eb{S}^2}$ with\\ \(\eb{D}^i=\left[\begin{smallmatrix}
			\eb{A}^i\\\eb{B}^i
			\end{smallmatrix} \right], \eb{A}^i=\ag{\msf{A}}{\alpha^i}\in\dlg{\chu}{F^i}{G^i},\eb{B}^i=\ag{\msf{B}}{\beta^i}\in\dlg{\chu}{K^i}{H^i},i=1,2 \), the morphisms are pairs of the form 
			\[\bm{v}\eqd\begin{bmatrix}
			\ag{\mu}{\nu}\\\ag{\theta}{\zeta}
			\end{bmatrix}:\eb{P}^1\rightharpoonup\eb{P}^2,\]
			where $\ag{\mu}{\nu}:\ag{\msf{A}}{\alpha^1}\longrightarrow\ag{\msf{A}}{\alpha^2}$ and \(\ag{\theta}{\zeta}:\ag{\msf{B}}{\beta^1}\longrightarrow\ag{\msf{B}}{\beta^2}\) are arrows in $\cat{Fluid}_\chu$, such that Diagram \ref{PVS-NVS} commutes (the left and right triangles are called the \textbf{Positive} and \textbf{Negative Vertical Super-adjointness} conditions, or briefly \textbf{PVS}\index{Positive Vertical Super-adjointness Condition (PVS)} and \textbf{NVS}\index{Negative Vertical Super-adjointness Condition (NVS)}, respectively).
			\begin{figure}[h]
				\begin{center}
					\begin{tikzpicture}[commutative diagrams/every diagram]
					
					\node (N1) at (0,0cm) {$\Sigma$};
					
					\node (N2) at (-4,2cm) {${F^1}^+\msf{A}\sqcup {K^1}^+\msf{B}$};
					
					\node (N3) at (-4,-2cm) {${F^2}^+\msf{A}\sqcup {K^2}^+\msf{B}$};
					
					\node (N4) at (4,-2cm) {${G^2}^-\msf{A}\sqcup {H^2}^-\msf{B}$};
					
					\node (N5) at (4,2cm) {${G^1}^-\msf{A}\sqcup {H^1}^-\msf{B}$};
					
					\node (N6) at (-2.4,0cm) {\fbox{\textbf{PVS}}};
					
					\node (N7) at (2.4,0cm) {\fbox{\textbf{NVS}}};
					
					\path[commutative diagrams/.cd, every arrow, every label]
					(N1) edge node[swap]{${S^1}^+$} (N2)
					(N2) edge node[swap]{$\mu^+\sqcup \theta^+$} (N3)
					(N1) edge node{${S^2}^+$} (N3)
					(N1) edge node{${S^1}^-$} (N5)
					(N4) edge node[swap]{$\nu^-\sqcup\zeta^-$} (N5)
					(N1) edge node[swap]{${S^2}^-$} (N4);
					\end{tikzpicture}
				\end{center}
				\caption{} \label{PVS-NVS}
			\end{figure}
			
			\item $\prm{\mcg{A}}$-composition: $\prm{\mcg{A}}$-morphisms compose componentwise: for $\bm{v}=\left[\begin{smallmatrix}
			\mcg{l}\\\mcg{m}
			\end{smallmatrix} \right]:\eb{P}^1\rightharpoonup\eb{P}^2$ and $\bm{w}=\left[\begin{smallmatrix}
			\prm{\mcg{l}}\\\prm{\mcg{m}}
			\end{smallmatrix} \right]:\eb{P}^2\rightharpoonup\eb{P}^3$,
			\[\bm{w\circ v}\eqd\begin{bmatrix}
			\prm{\mcg{l}}\mcg{l}\\\prm{\mcg{m}}\mcg{m}
			\end{bmatrix}:\eb{P}^1\rightharpoonup\eb{P}^3.\]
			
			\item $\prm{\mcg{A}}$-identities: The identities are defined componentwise; for an object $\eb{P=\ag{D}{S}}$,
			\[\prm{1}_\eb{P}\eqd\begin{bmatrix}
			\ag{1_{F\msf{A}}}{1_{G\msf{A}}}\\\ag{1_{K\msf{B}}}{1_{H\msf{B}}}
			\end{bmatrix};\]
			that is, $\prm{1}_\eb{P}$ is defined as a pair of identity $\cat{Fluid}_\chu$-arrows (the notation $\prm{1}_\eb{P}$ is used in order to avoid confusion with the identity $1_\eb{P}$ defined in part (a.4) in the above).
		\end{itemize}
		Also, the above composition is associative because of associativity of $\cat{Fluid}_\chu$-arrows. Therefore, $\prm{\mcg{A}}$ is indeed a category.\\
		
		\par \noindent \textbf{b.2. Objects of} $\mcg{B}$\textbf{.} The $\mcg{B}$-objects (the vertical 1-cells) are exactly the $\prm{\mcg{A}}$-morphisms defined above:
		\[\uo{\mcg{B}}\eqd\um{\prm{\mcg{A}}}.\]
		
		\par \noindent \textbf{b.3. Morphisms of} $\mcg{B}$\textbf{.} For $\mcg{A}$-objects
		\[\eb{P}^i_j=\ag{\eb{D}^i_j}{\eb{S}^i_j},~\eb{D}^i_j=\begin{bmatrix}
		\eb{A}^i_j\\\eb{B}^i_j
		\end{bmatrix},~\eb{A}^i_j=\ag{\msf{A}_j}{\alpha^i_j}\in\dlg{\chu}{F^i}{G^i},\]
		\[\eb{B}^i_j=\ag{\msf{B}_j}{\beta^i_j}\in\dlg{\chu}{K^i}{H^i},~i,j\in\{1,2\},\]
		for $\mcg{A}$-morphisms
		\[\bm{h}^i=\begin{bmatrix}
		\eb{f}^i\\\eb{g}^i
		\end{bmatrix}:\eb{P}^i_1\longrightarrow\eb{P}^i_2,~i\in\{1,2\},\]
		and for $\mcg{B}$-objects
		\[\bm{v}_j=\begin{bmatrix}
		\ag{\mu_j}{\nu_j}\\\ag{\theta_j}{\zeta_j}
		\end{bmatrix}:\eb{P}^1_j\rightharpoonup\eb{P}^2_j,~j\in\{1,2\},\]
		the $\mcg{B}$-morphisms (the 2-cells), also called the \textbf{cubicles}\index{cubicle}, are of the form
		\[\bm{c}\eqd\begin{pmatrix}
		&\bm{h}^1& \\\bm{v}_1& &\bm{v}_2\\ &\bm{h}^2& 
		\end{pmatrix}:\bm{v}_1\Longrightarrow\bm{v}_2,\]
		provided that the following four parallelograms commute (Diagram \ref{four-parall}):
		\begin{figure}[H]
			\begin{center}
				\begin{tikzpicture}[commutative diagrams/every diagram]
				
				\node (F11) at (-10,3.5cm) {$F^1\msf{A}_1$};
				
				\node (G11) at (-7,3.5cm) {$G^1\msf{A}_1$};
				
				\node (F12) at (-8.5,2cm) {$F^1\msf{A}_2$};
				
				\node (G12) at (-5.5,2cm) {$G^1\msf{A}_2$};
				
				\node (F21) at (-10,1cm) {$F^2\msf{A}_1$};
				
				\node (G21) at (-7,1cm) {$G^2\msf{A}_1$};
				
				\node (F22) at (-8.5,-0.5cm) {$F^2\msf{A}_2$};
				
				\node (G22) at (-5.5,-0.5cm) {$G^2\msf{A}_2$};
				
				\node (K11) at (-2,3.5cm) {$K^1\msf{B}_1$};
				
				\node (H11) at (1,3.5cm) {$H^1\msf{B}_1$};
				
				\node (K12) at (-0.5,2cm) {$K^1\msf{B}_2$};
				
				\node (H12) at (2.5,2cm) {$H^1\msf{B}_2$};
				
				\node (K21) at (-2,1cm) {$K^2\msf{B}_1$};
				
				\node (H21) at (1,1cm) {$H^2\msf{B}_1$};
				
				\node (K22) at (-0.5,-0.5cm) {$K^2\msf{B}_2$};
				
				\node (H22) at (2.5,-0.5cm) {$H^2\msf{B}_2$};
				
				\path[commutative diagrams/.cd, every arrow, every label]
				(F11) edge node[swap]{$\mu_1$} (F21)
				(G11) edge node[swap]{$\nu_1$} (G21)
				(F12) edge node{$\mu_2$} (F22)
				(G12) edge node{$\nu_2$} (G22)
				(K11) edge node[swap]{$\theta_1$} (K21)
				(H11) edge node[swap]{$\zeta_1$} (H21)
				(K12) edge node{$\theta_2$} (K22)
				(H12) edge node{$\zeta_2$} (H22)
				(F11) edge node{$F^1f^1$} (F12)
				(G11) edge node{$G^1f^1$} (G12)
				(F21) edge node[swap]{$F^2f^2$} (F22)
				(G21) edge node[swap]{$G^2f^2$} (G22)
				(K12) edge node[swap]{$K^1g^1$} (K11)
				(H12) edge node[swap]{$H^1g^1$} (H11)
				(K22) edge node{$K^2g^2$} (K21)
				(H22) edge node{$H^2g^2$} (H21);
				\end{tikzpicture}
			\end{center}
			\caption{} \label{four-parall}
		\end{figure}
		Notice the fact that the above conditions imply the commutativity of both of the cubes in the following diagram (Diagram \ref{cube-cube}).
		\begin{figure}[H]
			\begin{center}
				\begin{tikzpicture}[commutative diagrams/every diagram]
				
				\node (F11) at (-10,3.5cm) {$F^1\msf{A}_1$};
				
				\node (G11) at (-7,3.5cm) {$G^1\msf{A}_1$};
				
				\node (F12) at (-8.5,2cm) {$F^1\msf{A}_2$};
				
				\node (G12) at (-5.5,2cm) {$G^1\msf{B}_2$};
				
				\node (F21) at (-10,1cm) {$F^2\msf{A}_1$};
				
				\node (G21) at (-7,1cm) {$G^2\msf{A}_1$};
				
				\node (F22) at (-8.5,-0.5cm) {$F^2\msf{A}_2$};
				
				\node (G22) at (-5.5,-0.5cm) {$G^2\msf{A}_2$};
				
				\node (K11) at (-2,3.5cm) {$K^1\msf{B}_1$};
				
				\node (H11) at (1,3.5cm) {$H^1\msf{B}_1$};
				
				\node (K12) at (-0.5,2cm) {$K^1\msf{B}_2$};
				
				\node (H12) at (2.5,2cm) {$H^1\msf{B}_2$};
				
				\node (K21) at (-2,1cm) {$K^2\msf{B}_1$};
				
				\node (H21) at (1,1cm) {$H^2\msf{B}_1$};
				
				\node (K22) at (-0.5,-0.5cm) {$K^2\msf{B}_2$};
				
				\node (H22) at (2.5,-0.5cm) {$H^2\msf{B}_2$};
				
				\path[commutative diagrams/.cd, every arrow, every label]
				(F11) edge node[swap]{$\mu_1$} (F21)
				(G11) edge node[swap,near start]{$\nu_1$} (G21)
				(F21) edge node[swap,near start]{$\alpha^2_1$} (G21)
				(G12) edge node{$\nu_2$} (G22)
				(K11) edge node[swap]{$\theta_1$} (K21)
				(H11) edge node[swap,near start]{$\zeta_1$} (H21)
				(K21) edge node[swap,near start]{$\beta^2_1$} (H21)
				(H12) edge node{$\zeta_2$} (H22)
				(F11) edge node{$F^1f^1$} (F12)
				(G11) edge node{$G^1f^1$} (G12)
				(F21) edge node[swap]{$F^2f^2$} (F22)
				(G21) edge node[swap]{$G^2f^2$} (G22)
				(K12) edge node[swap]{$K^1g^1$} (K11)
				(H12) edge node[swap]{$H^1g^1$} (H11)
				(K22) edge node{$K^2g^2$} (K21)
				(H22) edge node{$H^2g^2$} (H21)
				(F11) edge node{$\alpha^1_1$} (G11)
				(F12) edge[-,line width=6pt,draw=white] node{} (G12)
				(F12) edge node[swap,near end]{$\alpha^1_2$} (G12)
				(F22) edge node[swap]{$\alpha^2_2$} (G22)
				(F12) edge[-,line width=6pt,draw=white] node{} (F22)
				(F12) edge node[near end]{$\mu_2$} (F22)
				(K11) edge node{$\beta^1_1$} (H11)
				(K12) edge[-,line width=6pt,draw=white] node{} (H12)
				(K12) edge node[swap,near end]{$\beta^1_2$} (H12)
				(K22) edge node[swap]{$\beta^2_2$} (H22)
				(K12) edge[-,line width=6pt,draw=white] node{} (K22)
				(K12) edge node[near end]{$\theta_2$} (K22);
				\end{tikzpicture}
			\end{center}
			\caption{} \label{cube-cube}
		\end{figure}
		Also, notice that the above definition means that we are defining a \textit{flat} double category since every 2-cell is completely determined by its boundary (see Definitions \ref{boundary} and \ref{flat-def}).\\
		
		\par \noindent \textbf{b.4. Composition in} $\mcg{B}$\textbf{.} For cubicles 
		\[\begin{pmatrix}
		&\bm{h}^1& \\\bm{v}_1& &\bm{v}_2\\ &\bm{h}^2& 
		\end{pmatrix}:\bm{v}_1\Longrightarrow\bm{v}_2~~~~\text{and}~~~~\begin{pmatrix}
		&\bm{k}^1& \\\bm{v}_2& &\bm{v}_3\\ &\bm{k}^2& 
		\end{pmatrix}:\bm{v}_2\Longrightarrow\bm{v}_3\] 
		the composition is defined as
		\[\begin{pmatrix}
		&\bm{k}^1& \\\bm{v}_2& &\bm{v}_3\\ &\bm{k}^2& 
		\end{pmatrix}\circ\begin{pmatrix}
		&\bm{h}^1& \\\bm{v}_1& &\bm{v}_2\\ &\bm{h}^2& 
		\end{pmatrix}\eqd\begin{pmatrix}
		&\bm{k}^1\bm{h}^1& \\\bm{v}_1& &\bm{v}_3\\ &\bm{k}^2\bm{h}^2& 
		\end{pmatrix}:\bm{v}_1\Longrightarrow\bm{v}_3.\]
		
		\par \noindent \textbf{b.5. Identities in} $\mcg{B}$\textbf{.} For every vertical 1-cell $\bm{v}:\ag{\eb{D}^1}{\eb{S}^1}\rightharpoonup\ag{\eb{D}^2}{\eb{S}^2}$, the identity 2-cell is defiend as
		\[\bm{1_v}\eqd\begin{pmatrix}
		&1_{\eb{D}^1}& \\\bm{v}& &\bm{v}\\ &1_{\eb{D}^2}& 
		\end{pmatrix}.\]
		
		\rule{0pt}{1mm}
		\par \noindent \textbf{c. The internal source and target functors.} The functors $\mcg{s,t:B\dlong A}$ are defined as the following:\\
		
		\par \noindent \textbf{c.1. Internal source and target on} $\mcg{B}$\textbf{-objects.} For a $\mcg{B}$-object $\bm{v}:\eb{P}^1\rightharpoonup\eb{P}^2$,
		\[\mcg{s}(\bm{v})\eqd\eb{P}^1;~~~~\mcg{t}(\bm{v})\eqd\eb{P}^2.\]
		
		\par \noindent \textbf{c.1. Internal source and target on} $\mcg{B}$\textbf{-morphisms.} For a $\mcg{B}$-morphism\\ $\bm{c}=\left( \begin{smallmatrix}
		&\bm{h}^1& \\\bm{v}_1& &\bm{v}_2\\ &\bm{h}^2& 
		\end{smallmatrix}\right):\bm{v}_1\Longrightarrow\bm{v}_2 $ with $\bm{v}_j:\eb{P}^1_j\rightharpoonup\eb{P}^2_j,~j=1,2$ and $\bm{h}^i:\eb{P}^i_1\longrightarrow\eb{P}^i_2,~i=1,2$,
		\[\mcg{s}(\bm{c})\eqd\bm{h}^1;~~~~\mcg{t}(\bm{c})\eqd\bm{h}^2.\]
		
		\par \noindent \textbf{d. The internal identity functor.} The functor $\mcg{i:A\longrightarrow B}$ is defined as the following:\\
		
		\par \noindent \textbf{d.1. Internal identity on} $\mcg{A}$\textbf{-objects.} For $\eb{P=\ag{D}{S}}$,
		\[\mcg{i}(P)\eqd\prm{1}_\eb{P},\]
		where $\prm{1}_\eb{P}$ was defined in part (b.1).\\
		
		\par \noindent \textbf{d.2. Internal identity on} $\mcg{A}$\textbf{-morphisms.} For $\bm{h}:\eb{P}_1\longrightarrow\eb{P}_2$,
		\[\mcg{i}(\bm{h})\eqd\begin{pmatrix}
		&\bm{h}& \\\prm{1}_{\eb{P}_1}& &\prm{1}_{\eb{P}_2}\\ &\bm{h}& 
		\end{pmatrix}:\prm{1}_{\eb{P}_1}\Longrightarrow\prm{1}_{\eb{P}_2}.\]
		
		\par \noindent \textbf{e. The internal composition bifunctor.} For the pullback of Diagram \ref{pullb-s-t-use},
		\begin{figure}[h]
			\begin{center}
				\begin{tikzpicture}[commutative diagrams/every diagram]
				
				\node (N1) at (2,1cm) {$\mcg{B}$};
				
				\node (N2) at (-1,1cm) {$\mcg{B\times_A B}$};
				
				\node (N3) at (-1,-1cm) {$\mcg{B}$};
				
				\node (N4) at (2,-1cm) {$\mcg{A}$};
				
				\path[commutative diagrams/.cd, every arrow, every label]
				(N2) edge node{$\pi_2$} (N1)
				(N1) edge node{$\mcg{s}$} (N4)
				(N3) edge node{$\mcg{t}$} (N4)
				(N2) edge node[swap]{$\pi_1$} (N3);
				\end{tikzpicture}
			\end{center}
			\caption{} \label{pullb-s-t-use}
		\end{figure}
		the composition bifunctor
		\[\bt:\mcg{\pb{B}{A}{B}\longrightarrow B}\]
		is defined as the following:\\
		
		\par \noindent \textbf{e.1. Internal composition on objects.} Considering a pair
		\[\bm{\ag{v}{w}}\in\ob(\mcg{\pb{B}{A}{B}}),~\bm{v}=\begin{bmatrix}
		\mcg{l}\\\mcg{m}
		\end{bmatrix}:\eb{P}^1\rightharpoonup\eb{P}^2,~\bm{w}=\begin{bmatrix}
		\prm{\mcg{l}}\\\prm{\mcg{m}}
		\end{bmatrix}:\eb{P}^2\rightharpoonup\eb{P}^3,\]
		and abbreviating $\bt(\bm{v,w})$ as $\bm{v\bt w}$, we define
		\[\bm{v\bt w}\eqd\bm{w}\circ_{\prm{\mcg{A}}}\bm{v}=\begin{bmatrix}
		\prm{\mcg{l}}\mcg{l}\\\prm{\mcg{m}}\mcg{m}
		\end{bmatrix}:\eb{P}^1\rightharpoonup\eb{P}^3,\]
		where $\circ_{\prm{\mcg{A}}}$ is composition in the category $\prm{\mcg{A}}$ (see parts (b.1), (b.2)).\\
		
		\par \noindent \textbf{e.1. Internal composition on morphisms.} Considering a pair
		\[\bm{\ag{c}{d}}\in\mor(\mcg{\pb{B}{A}{B}})\]
		with
		\[\bm{c}=\begin{pmatrix}
		&\bm{h}^1& \\\bm{v}_1& &\bm{v}_2\\ &\bm{h}^2& 
		\end{pmatrix}:\bm{v}_1\Longrightarrow\bm{v}_2,~\bm{v}_j:\eb{P}^1_j\rightharpoonup\eb{P}^2_j,~j=1,2,\]
		\[\bm{d}=\begin{pmatrix}
		&\bm{h}^2& \\\bm{w}_1& &\bm{w}_2\\ &\bm{h}^3& 
		\end{pmatrix}:\bm{w}_1\Longrightarrow\bm{w}_2,~\bm{w}_j:\eb{P}^2_j\rightharpoonup\eb{P}^3_j,~j=1,2,\]
		we define
		\[\bm{c\bt d}\eqd\begin{pmatrix}
		&\bm{h}^1& \\\bm{v}_1\bt\bm{w}_1& &\bm{v}_2\bt\bm{w}_2\\ &\bm{h}^3& 
		\end{pmatrix}:\bm{v}_1\bt\bm{w}_1\Longrightarrow\bm{v}_2\bt\bm{w}_2.\]
	\end{dfn}
	
	\rule{0pt}{5mm}
	Whenever $\Gamma,\Sigma$ are clear from the context, we denote the double category $\dlc_{\Gamma,\Sigma}$ by $\dlc$. The initialism ``DLC'' stands for ``the properly large flat \textbf{D}ouble category of paired \textbf{D}ialgebras and dialgebra homomorphisms, paired \textbf{L}ocal transformations, and \textbf{C}ubicles on \chu''.\\
	
	\par We need to show that Definition \ref{DLC-def} is well-defined.
	
	\begin{prp}
		\label{DLC-welldef}
		$\dlc$ as defined in Definition \ref{DLC-def} is well-defined.
	\end{prp}
	
	\begin{mypr}
		We prove this through the following steps:\\
		
		\par \noindent \textbf{Step 1.} \textit{Composition in} $\mcg{A}$ \textit{is well-defined.}\\
		
		\par The composition of $\bm{h}=\left[ \begin{smallmatrix}
		\eb{A}_1\xrightarrow{~\eb{f}~}\eb{A}_2\\\eb{B}_2\xrightarrow{~\eb{g}~}\eb{B}_1
		\end{smallmatrix}\right]$ and $\prm{\bm{h}}=\left[ \begin{smallmatrix}
		\eb{A}_2\xrightarrow{~\prm{\eb{f}}~}\eb{A}_3\\\eb{B}_3\xrightarrow{~\prm{\eb{g}}~}\eb{B}_2
		\end{smallmatrix}\right]$ was given by
		\[\prm{\bm{h}}\circ_\mcg{A}\bm{h}=\left[ \begin{smallmatrix}
		\eb{A}_1\xrightarrow{~\prm{\eb{f}}\eb{f}~}\eb{A}_3\\\eb{B}_3\xrightarrow{~\eb{g}\prm{\eb{g}}~}\eb{B}_1
		\end{smallmatrix}\right].\]
		Well-definedness of the above composition is proved by the commutativity of the following diagrams (Diagrams \ref{PHS-in-A-proof}, \ref{NHS-in-A-proof}). The former diagram guarantees the PHS condition to hold for $\prm{\bm{h}}\circ_\mcg{A}\bm{h}$, while the latter verifies NHS.\\
		\begin{figure}[h]
			\begin{center}
				\begin{tikzpicture}[->,>=angle 90,commutative diagrams/every diagram]
				
				\node (N1) at (5,0cm) {$\Sigma$};
				
				\node (N2) at (5,3cm) {$F^+\msf{A}_1\sqcup K^+\msf{B}_1$};
				
				\node (N3) at (0,3cm) {$F^+\msf{A}_2\sqcup K^+\msf{B}_1$};
				
				\node (N4) at (0,0cm) {$F^+\msf{A}_2\sqcup K^+\msf{B}_2$};
				
				\node (N5) at (0,-3cm) {$F^+\msf{A}_3\sqcup K^+\msf{B}_2$};
				
				\node (N6) at (5,-3cm) {$F^+\msf{A}_3\sqcup K^+\msf{B}_3$};
				
				\node (N7) at (-5,0cm) {$F^+\msf{A}_3\sqcup K^+\msf{B}_1$};
				
				\path[commutative diagrams/.cd, every arrow, every label]
				(N1) edge node{$S^+_1$} (N2)
				(N2) edge node{$F^+f\sqcup 1_{K^+\msf{B}_1}$} (N3)
				(N4) edge node[swap]{$1_{F^+\msf{A}_2}\sqcup K^+g$} (N3)
				(N4) edge node{$F^+\prm{f}\sqcup 1_{K^+\msf{B}_2}$} (N5)
				(N6) edge node[swap]{$1_{F^+\msf{A}_3}\sqcup K^+\prm{g}$} (N5)
				(N1) edge node[swap]{$S^+_2$} (N4)
				(N1) edge node[swap]{$S^+_3$} (N6)
				(N3) edge node[near end]{$F^+\prm{f}\sqcup 1_{K^+\msf{B}_1}$} (N7)
				(N5) edge node[swap,near end]{$1_{F^+\msf{A}_3}\sqcup K^+g$} (N7);
				
				\draw[bend right=50] (N2) to node [above] {$F^+(\prm{f}f)\sqcup 1_{K^+\msf{B}_1}$} (N7);
				
				\draw[bend left=50] (N6) to  node [below] {$1_{F^+\msf{A}_3}\sqcup K^+(g\prm{g})$} (N7);
				\end{tikzpicture}
			\end{center}
			\caption{} \label{PHS-in-A-proof}
		\end{figure}
		\begin{figure}[h]
			\begin{center}
				\begin{tikzpicture}[->,>=angle 90,commutative diagrams/every diagram]
				
				\node (N1) at (5,0cm) {$\Sigma$};
				
				\node (N2) at (5,3cm) {$G^-\msf{A}_3\sqcup H^-\msf{B}_3$};
				
				\node (N3) at (0,3cm) {$G^-\msf{A}_2\sqcup H^-\msf{B}_3$};
				
				\node (N4) at (0,0cm) {$G^-\msf{A}_2\sqcup H^-\msf{B}_2$};
				
				\node (N5) at (0,-3cm) {$G^-\msf{A}_1\sqcup H^-\msf{B}_2$};
				
				\node (N6) at (5,-3cm) {$G^-\msf{A}_1\sqcup H^-\msf{B}_1$};
				
				\node (N7) at (-5,0cm) {$G^-\msf{A}_1\sqcup H^-\msf{B}_3$};
				
				\path[commutative diagrams/.cd, every arrow, every label]
				(N1) edge node{$S^-_3$} (N2)
				(N2) edge node{$G^-\prm{f}\sqcup 1_{H^-\msf{B}_3}$} (N3)
				(N4) edge node[swap]{$1_{G^-\msf{A}_2}\sqcup H^-\prm{g}$} (N3)
				(N4) edge node{$G^-f\sqcup 1_{H^-\msf{B}_2}$} (N5)
				(N6) edge node[swap]{$1_{G^-\msf{A}_1}\sqcup H^-g$} (N5)
				(N1) edge node[swap]{$S^-_2$} (N4)
				(N1) edge node[swap]{$S^-_1$} (N6)
				(N3) edge node[near end]{$G^-f\sqcup 1_{H^-\msf{B}_3}$} (N7)
				(N5) edge node[swap,near end]{$1_{G^-\msf{A}_1}\sqcup H^-\prm{g}$} (N7);
				
				\draw[bend right=50] (N2) to node [above] {$G^-(f\prm{f})\sqcup 1_{H^-\msf{B}_3}$} (N7);
				
				\draw[bend left=50] (N6) to  node [below] {$1_{G^-\msf{A}_1}\sqcup H^-(\prm{g}g)$} (N7);
				\end{tikzpicture}
			\end{center}
			\caption{} \label{NHS-in-A-proof}
		\end{figure}
		
		\par \noindent \textbf{Step 2.} \textit{The identity and associativity propoerties hold for} $\mcg{A}$\textit{-morphisms.} Given $\eb{P}_j=\ag{\eb{D}_j}{\eb{S}_j},~\eb{D}_j=\left[\begin{smallmatrix}
		\eb{A}_j\\\eb{B}_j
		\end{smallmatrix} \right],~j=1,2$, it is clear that $1_{\eb{P}_j}=\left[\begin{smallmatrix}
		1_{\eb{A}_j}\\1_{\eb{B}_j}
		\end{smallmatrix} \right],~j=1,2 $ are actually the identity morphisms on $\eb{P}_j$ since
		\begin{align*}
		\bm{h}\circ 1_{\eb{P}_1}&=\begin{bmatrix}
		\eb{f}\circ 1_{\eb{A}_1}\\\eb{g}\circ 1_{\eb{B}_1}
		\end{bmatrix}\\
		&=\begin{bmatrix}
		\eb{f}\\\eb{g}
		\end{bmatrix}\\
		&=\begin{bmatrix}
		1_{\eb{A}_2}\circ\eb{f}\\1_{\eb{B}_2}\circ\eb{g}
		\end{bmatrix}\\
		&=1_{\eb{P}_2}\circ\bm{h}
		\end{align*}
		for every $\bm{h}=\left[ \begin{smallmatrix}
		\eb{f}\\\eb{g}
		\end{smallmatrix}\right]:\eb{P}_1\longrightarrow\eb{P}_2$.\\
		\par With a similar line of reasoning, as $\eb{f},\eb{g}$ are dialgebra homomorphisms for every\\ $\bm{h}=\left[ \begin{smallmatrix}
		\eb{f}\\\eb{g}
		\end{smallmatrix}\right]\in\mcg{A}$, the associativity of $\mcg{A}$-morphisms is verified.\\
		
		\rule{0pt}{1mm}
		\par \noindent \textbf{Step 3.} \textit{Composition in} $\prm{\mcg{A}}$ \textit{is well-defined.}\\
		
		\par The composition of $\bm{v}=\left[\begin{smallmatrix}
		\ag{\mu}{\nu}\\\ag{\theta}{\zeta}
		\end{smallmatrix} \right]:\eb{P}^1\rightharpoonup\eb{P}^2$ and $\bm{w}=\left[\begin{smallmatrix}
		\ag{\prm{\mu}}{\prm{\nu}}\\\ag{\prm{\theta}}{\prm{\zeta}}
		\end{smallmatrix} \right]:\eb{P}^2\rightharpoonup\eb{P}^3$ was given by
		\[\bm{w}\circ_\mcg{\prm{A}}\bm{v}=\begin{bmatrix}
		\ag{\prm{\mu}\mu}{\prm{\nu}\nu}\\ \ag{\prm{\theta}\theta}{\prm{\zeta}\zeta}
		\end{bmatrix}.\]
		Well-definedness of this composition is proved by the commutativity of Diagram \ref{A-prime-compo-well}. In the diagram, the left half guarantees the PVS condition to hold for $\bm{w}\circ_\mcg{\prm{A}}\bm{v}$ while right half verifies NVS.
		\begin{figure}[h]
			\begin{center}
				\begin{tikzpicture}[->,>=angle 90,commutative diagrams/every diagram]
				
				\node (C) at (0,0cm) {$\Sigma$};
				\node (T) at (0,3cm) {$\Sigma$};
				\node (B) at (0,-3cm) {$\Sigma$};
				
				\node (N2) at (-5,3cm) {${F^1}^+\msf{A}\sqcup {K^1}^+\msf{B}$};
				
				\node (N3) at (-3,0cm) {${F^2}^+\msf{A}\sqcup {K^2}^+\msf{B}$};
				
				\node (M1) at (-5,-3cm) {${F^3}^+\msf{A}\sqcup {K^3}^+\msf{B}$};
				
				\node (N4) at (3,0cm) {${G^2}^-\msf{A}\sqcup {H^2}^-\msf{B}$};
				
				\node (N5) at (5,3cm) {${G^1}^-\msf{A}\sqcup {H^1}^-\msf{B}$};
				
				\node (M2) at (5,-3cm) {${G^3}^-\msf{A}\sqcup {H^3}^-\msf{B}$};
				
				\path[commutative diagrams/.cd, every arrow, every label]
				(T) edge node[swap]{${S^1}^+$} (N2)
				(N2) edge node{$\mu^+\sqcup \theta^+$} (N3)
				(C) edge node[swap]{${S^2}^+$} (N3)
				(T) edge node{${S^1}^-$} (N5)
				(N4) edge node{$\nu^-\sqcup\zeta^-$} (N5)
				(C) edge node{${S^2}^-$} (N4)
				(B) edge node{${S^3}^+$} (M1)
				(B) edge node[swap]{${S^3}^-$} (M2)
				(N3) edge node{$(\prm{\mu})^+\sqcup (\prm{\theta})^+$} (M1)
				(M2) edge node{$(\prm{\nu})^-\sqcup(\prm{\zeta})^-$} (N4)
				(N2) edge node[swap]{\rotatebox{-90}{$(\prm{\mu}\mu)^+\sqcup(\prm{\theta}\theta)^+$}} (M1)
				(M2) edge node[swap]{\rotatebox{90}{$(\prm{\nu}\nu)^-\sqcup(\prm{\zeta}\zeta)^-$}} (N5);
				
				\draw[-, double distance=2pt,line width=0.5pt] (T) to node {} (C);				
				\draw[-, double distance=2pt,line width=0.5pt] (C) to node {} (B);
				\end{tikzpicture}
			\end{center}
			\caption{} \label{A-prime-compo-well}
		\end{figure}

		\rule{0pt}{1mm}
		\par \noindent \textbf{Step 4.} $\mcg{B}$\textit{-morphisms} \textit{are well-defined.}\\
		
		\par To show that the cubicles are well-defined, we observe that for the cubicle\\ $\bm{c}=\left(\begin{smallmatrix}
		&\bm{h}^1& \\\bm{v}_1& &\bm{v}_2\\ &\bm{h}^2& 
		\end{smallmatrix} \right)$ defined in part (b.2) of Definition \ref{DLC-def}, the horizontal and vertical super-adjointness conditions, PHS, NHS, PVS, and NVS are compatible with each other. In other words, Diagram \ref{PHS-NHS-PVS-NVS-for-cubicles} commutes.
		\begin{figure}[h]
			\begin{center}
				\begin{tikzpicture}[commutative diagrams/every diagram]
				
				\node (U0) at (0,3cm) {$\Sigma$};
				\node (D0) at (0,-3cm) {$\Sigma$};
				
				\node (U11p) at (-4,4cm) {${F^1}^+\msf{A}_1\sqcup {K^1}^+\msf{B}_1$};	
				\node (D11p) at (-4,-2cm) {${F^2}^+\msf{A}_1\sqcup {K^2}^+\msf{B}_1$};
				
				\node (U12p) at (-2,1cm) {${F^1}^+\msf{A}_2\sqcup {K^1}^+\msf{B}_2$};	
				\node (D12p) at (-2,-5cm) {${F^2}^+\msf{A}_2\sqcup {K^2}^+\msf{B}_2$};
				
				\node (Uxp) at (-6,2.2cm) {${F^1}^+\msf{A}_2\sqcup {K^1}^+\msf{B}_1$};	
				\node (Dxp) at (-6,-3.8cm) {${F^2}^+\msf{A}_2\sqcup {K^2}^+\msf{B}_1$};
				
				\node (U11m) at (4,4cm) {${G^1}^-\msf{A}_1\sqcup {H^1}^-\msf{B}_1$};
				\node (D11m) at (4,-2cm) {${G^2}^-\msf{A}_1\sqcup {H^2}^-\msf{B}_1$};
				
				\node (U12m) at (2,1cm) {${G^1}^-\msf{A}_2\sqcup {H^1}^-\msf{B}_2$};
				\node (D12m) at (2,-5cm) {${G^2}^-\msf{A}_2\sqcup {H^2}^-\msf{B}_2$};
				
				\node (Uxm) at (6,2.2cm) {${G^1}^-\msf{A}_1\sqcup {H^1}^-\msf{B}_2$};
				\node (Dxm) at (6,-3.8cm) {${G^2}^-\msf{A}_1\sqcup {H^2}^-\msf{B}_2$};
				
				\path[commutative diagrams/.cd, every arrow, every label]
				(U0) edge node[swap]{${S^1_1}^+$} (U11p)
				(U0) edge node[swap]{${S^1_2}^+$} (U12p)
				(U0) edge node{${S^1_1}^-$} (U11m)
				(U0) edge node{${S^1_2}^-$} (U12m)
				(D0) edge node[near start]{${S^2_1}^+$} (D11p)
				(D0) edge node[swap]{${S^2_2}^+$} (D12p)
				(D0) edge node[swap,near start]{${S^2_1}^-$} (D11m)
				(D0) edge node{${S^2_2}^-$} (D12m)
				(U11p) edge node[swap]{${F^1}^+f^1\sqcup 1$} (Uxp)
				(U11p) edge node[swap,near end]{$\mu_1^+\sqcup\theta_1^+$} (D11p)
				(U12p) edge[-,line width=6pt,draw=white] node{} (Uxp)
				(U12p) edge node[swap,near start]{$1\sqcup {K^1}^+g^1$} (Uxp)
				(U11m) edge node{$1\sqcup {H^1}^-g^1$} (Uxm)
				(D11m) edge node[swap,near start]{$\nu_1^-\sqcup\zeta_1^-$} (U11m)
				(U12m) edge[-,line width=6pt,draw=white] node{} (Uxm)
				(U12m) edge node[near start]{${G^1}^-f^1\sqcup 1$} (Uxm)
				(D11p) edge node{${F^2}^+f^2\sqcup 1$} (Dxp)
				(D12p) edge node[swap]{$1\sqcup {K^2}^+g^2$} (Dxp)
				(D11m) edge node[swap]{$1\sqcup {H^2}^-g^2$} (Dxm)
				(D12m) edge node[swap]{${G^2}^-f^2\sqcup 1$} (Dxm)
				(U12p) edge[-,line width=6pt,draw=white] node{} (D12p)
				(U12p) edge node[near start]{$\mu_2^+\sqcup\theta_2^+$} (D12p)
				(Uxp) edge node[swap]{$\mu_2^+\sqcup\theta_1^+$} (Dxp)
				(D12m) edge[-,line width=6pt,draw=white] node{} (U12m)
				(D12m) edge node[near end]{$\nu_2^-\sqcup\zeta_2^-$} (U12m)
				(Dxm) edge node[swap]{$\nu_1^-\sqcup\zeta_2^-$} (Uxm);
				
				\draw[-, double distance=2pt,line width=0.5pt] (U0) to node {} (D0);
				\end{tikzpicture}
			\end{center}
			\caption{} \label{PHS-NHS-PVS-NVS-for-cubicles}
		\end{figure}
		In this diagram, the upper-left, upper-right, lower-left, and lower-right diamonds are $\tu{PHS}^1,~\tu{NHS}^1,~\tu{PHS}^2,$ and $\tu{NHS}^2$, respectively (the superscripts refer to the corresponding horizontal 1-cells). Also, the front and rear parallelograms with common segment $\begin{matrix}
		\Sigma\\\parallel\\\Sigma
		\end{matrix}$ are $\tu{PVS}_1,~\tu{NVS}_1,~\tu{PVS}_2$, and $\tu{NVS}_2$, respectively (the subscripts referring to the corresponding vertical 1-cells). Finally, the commutativity of the front and rear lateral parallelograms follows immediately from the commutativity of Diagram \ref{four-parall}. Whence, the whole shape commutes, and consequently, cubicle $\bm{c}$ is a well-defined arrow from $\bm{v}_1$ to $\bm{v}_2$.\\
		
		\rule{0pt}{1mm}
		\par \noindent \textbf{Step 5.} \textit{Axiom \textbf{A1} of Definition \ref{internal-cat-def} holds for} $\dlc$ \textit{as an internal category in} $\cat{CAT}$\textit{.}\\
		
		\par Consider the functors $\mcg{si:A\longrightarrow A}$ and $\mcg{ti:A\longrightarrow A}$. For every $\mcg{A}$-object
		\[\eb{P=\ag{D}{S}},~~\eb{D}=\begin{bmatrix}
		\ag{\msf{A}}{F\msf{A}\xrightarrow{~\alpha~}G\msf{A}}\\\ag{\msf{B}}{K\msf{B}\xrightarrow{~\beta~}H\msf{B}}
		\end{bmatrix}\]
		we have
		\[\mcg{si}(\eb{P})=s\begin{bmatrix}
		\ag{1_{F\msf{A}}}{1_{G\msf{A}}}\\\ag{1_{K\msf{B}}}{1_{H\msf{B}}}
		\end{bmatrix}=\eb{P},\]
		\[\mcg{ti}(\eb{P})=t\begin{bmatrix}
		\ag{1_{F\msf{A}}}{1_{G\msf{A}}}\\\ag{1_{K\msf{B}}}{1_{H\msf{B}}}
		\end{bmatrix}=\eb{P}.\]
		Also, for every $\mcg{A}$-morphism $\bm{h}:\eb{P}_1\longrightarrow\eb{P}_2$,
		\[\mcg{si}(\bm{h})=\mcg{s}\begin{pmatrix}
		&\bm{h}& \\\prm{1}_{\eb{P}_1}& &\prm{1}_{\eb{P}_2}\\ &\bm{h}& 
		\end{pmatrix}=\bm{h},\]
		\[\mcg{ti}(\bm{h})=\mcg{t}\begin{pmatrix}
		&\bm{h}& \\\prm{1}_{\eb{P}_1}& &\prm{1}_{\eb{P}_2}\\ &\bm{h}& 
		\end{pmatrix}=\bm{h}.\]
		Therefore, 
		\[\mcg{si=ti}=1_\mcg{A}.\]
		
		\rule{0pt}{1mm}
		\par \noindent \textbf{Step 6.} \textit{Axiom \textbf{A2} of Definition \ref{internal-cat-def} holds for} $\dlc$ \textit{as an internal category in} $\cat{CAT}$\textit{.}\\
		
		\par For any pair $\bm{\ag{v}{w}}\in\uo{(\bb{\pb{B}{A}{B}})}$ with $\bm{v}:\eb{P}^1\rightharpoonup\eb{P}^2$ and $\bm{w}:\eb{P}^2\rightharpoonup\eb{P}^3$ we have
		\[\mcg{t}\pi_2(\bm{v,w})=\mcg{t}(\bm{w})=\eb{P}^3=\mcg{t}(\bm{v\bt w}),\]
		\[\mcg{s}\pi_1(\bm{v,w})=\mcg{s}(\bm{v})=\eb{P}^1=\mcg{s}(\bm{v\bt w}).\]
		On the other hand, for a pair $\bm{\ag{c}{d}}\in\um{(\bb{\pb{B}{A}{B}})}$ with $\bm{v}_1\xR{c}\bm{v}_2$ and $\bm{v}_2\xR{d}\bm{v}_3$,
		\[\mcg{t}\pi_2(\bm{c,d})=\mcg{t}(\bm{d})=\bm{v}_3=\mcg{t}(\bm{c\bt d}),\]
		\[\mcg{s}\pi_1(\bm{c,d})=\mcg{s}(\bm{c})=\bm{v}_1=\mcg{s}(\bm{c\bt d}).\]
		So, 
		\[\mcg{t}\pi_2=\mcg{t}\bt~~~~\text{and}~~~~ \mcg{s}\pi_1=\mcg{s}\bt,\]
		as desired.\\
		
		\rule{0pt}{1mm}
		\par \noindent \textbf{Step 7.} \textit{Axiom \textbf{A3} of Definition \ref{internal-cat-def} holds for} $\dlc$ \textit{as an internal category in} $\cat{CAT}$\textit{.}\\
		
		\par Let $\bm{v}:\eb{P}^1\rightharpoonup\eb{P}^2$ be a $\mcg{B}$-object with
		\[\eb{P}^i=\ag{\eb{D}^i}{\eb{S}^i},~~\eb{D}^i=\begin{bmatrix}
		\ag{\msf{A}}{F^i\msf{A}\xrightarrow{~\alpha^i~}G^i\msf{A}}\\\ag{\msf{B}}{K^i\msf{B}\xrightarrow{~\beta^i~}H^i\msf{B}}
		\end{bmatrix},~~i=1,2.\]
		Then, $\mcg{s}(\bm{v})=\eb{P}^1$ and $\mcg{t}(\bm{v})=\eb{P}^2$, thus
		\[\mcg{is}(\bm{v})=\prm{1}_{\eb{P}^1}=\begin{bmatrix}
		\ag{1_{F^1\msf{A}}}{1_{G^1\msf{A}}}\\\ag{1_{K^1\msf{B}}}{1_{H^1\msf{B}}}
		\end{bmatrix}:\eb{P}^1\rightharpoonup\eb{P}^1,\]
		\[\mcg{it}(\bm{v})=\prm{1}_{\eb{P}^2}=\begin{bmatrix}
		\ag{1_{F^2\msf{A}}}{1_{G^2\msf{A}}}\\\ag{1_{K^2\msf{B}}}{1_{H^2\msf{B}}}
		\end{bmatrix}:\eb{P}^2\rightharpoonup\eb{P}^2.\]
		Whence,
		\[\begin{bmatrix}
		\ag{1_{F^1\msf{A}}}{1_{G^1\msf{A}}}\\\ag{1_{K^1\msf{B}}}{1_{H^1\msf{B}}}
		\end{bmatrix}\bt\bm{v}=\bm{v},\]
		\[\bm{v}\bt\begin{bmatrix}
		\ag{1_{F^2\msf{A}}}{1_{G^2\msf{A}}}\\\ag{1_{K^2\msf{B}}}{1_{H^2\msf{B}}}
		\end{bmatrix}=\bm{v}.\]
		These imply that we have the following agreement of functors on every $\mcg{B}$-object $\bm{v}$:
		\begin{align}
		\bt\circ(\mcg{is},1_\mcg{B})(\bm{v})&=1_\mcg{B}(\bm{v})=\bt\circ(1_\mcg{B},\mcg{it})(\bm{v}). \tag{$*$}
		\end{align}
		
		\par On the other hand, assume a $\mcg{B}$-morphism $\bm{c}=\left( \begin{smallmatrix}
		&\bm{h}^1& \\\bm{v}_1& &\bm{v}_2\\ &\bm{h}^2& 
		\end{smallmatrix}\right)$ between $\mcg{B}$-objects
		\[\bm{v}_j:\eb{P}^1_j\rightharpoonup\eb{P}^2_j,~~~j=1,2,\]
		together with $\mcg{A}$-morphisms
		\[\bm{h}^i:\eb{P}^i_1\longrightarrow\eb{P}^i_2,~~~i=1,2\]
		such that
		\[\mcg{s}(\bm{c})=\bm{h}^1,~~~~\mcg{t}(\bm{c})=\bm{h}^2.\]
		We have
		\[\mcg{is}(\bm{c})=\mcg{i}(\bm{h}^1)=\begin{pmatrix}
		&\bm{h}^1& \\\prm{1}_{\eb{P}^1_1}& &\prm{1}_{\eb{P}^1_2}\\ &\bm{h}^1& 
		\end{pmatrix},\]
		\[\mcg{it}(\bm{c})=\mcg{i}(\bm{h}^2)=\begin{pmatrix}
		&\bm{h}^2& \\\prm{1}_{\eb{P}^2_1}& &\prm{1}_{\eb{P}^2_2}\\ &\bm{h}^2& 
		\end{pmatrix}.\]
		Therefore,
		\begin{align*}
		\begin{pmatrix}
		&\bm{h}^1& \\\prm{1}_{\eb{P}^1_1}& &\prm{1}_{\eb{P}^1_2}\\ &\bm{h}^1& 
		\end{pmatrix} \bt \bm{c}&=\begin{pmatrix}
		&\bm{h}^1& \\\prm{1}_{\eb{P}^1_1}& &\prm{1}_{\eb{P}^1_2}\\ &\bm{h}^1& 
		\end{pmatrix} \bt \begin{pmatrix}
		&\bm{h}^1& \\\bm{v}_1& &\bm{v}_2\\ &\bm{h}^2& 
		\end{pmatrix}\\
		&=\begin{pmatrix}
		&\bm{h}^1& \\\bm{v}_1& &\bm{v}_2\\ &\bm{h}^2& 
		\end{pmatrix}\\
		&=\bm{c};
		\end{align*}
		
		\begin{align*}
		\bm{c} \bt \begin{pmatrix}
		&\bm{h}^2& \\\prm{1}_{\eb{P}^2_1}& &\prm{1}_{\eb{P}^2_2}\\ &\bm{h}^2& 
		\end{pmatrix}&=\begin{pmatrix}
		&\bm{h}^1& \\\bm{v}_1& &\bm{v}_2\\ &\bm{h}^2& 
		\end{pmatrix} \bt \begin{pmatrix}
		&\bm{h}^2& \\\prm{1}_{\eb{P}^2_1}& &\prm{1}_{\eb{P}^2_2}\\ &\bm{h}^2& 
		\end{pmatrix}\\
		&=\begin{pmatrix}
		&\bm{h}^1& \\\bm{v}_1& &\bm{v}_2\\ &\bm{h}^2& 
		\end{pmatrix}\\
		&=\bm{c}.
		\end{align*}
		These imply that we have the following agreement of functors on every $\mcg{B}$-morphism $\bm{c}$:
		\begin{align}
		\bt\circ(\mcg{is},1_\mcg{B})(\bm{c})&=1_\mcg{B}(\bm{c})=\bt\circ(1_\mcg{B},\mcg{it})(\bm{c}). \tag{$**$}
		\end{align}
		Now, results $(*),(**)$ imply
		\[\bt\circ(\mcg{is},1_\mcg{B})=1_\mcg{B}=\bt\circ(1_\mcg{B},\mcg{it}).\]
		
		\rule{0pt}{1mm}
		\par \noindent \textbf{Step 8.} \textit{Axiom \textbf{A4} of Definition \ref{internal-cat-def} holds for} $\dlc$ \textit{as an internal category in} $\cat{CAT}$\textit{.}\\
		
		\par Finally, we verify the internal associativity. Consider the bifunctors
		\[\mcg{\pb{1_B}{A}{\bt}}:\mcg{\pb{B}{A}{(\pb{B}{A}{B})}}\longrightarrow\mcg{\pb{B}{A}{B}},\]
		\[\mcg{\pb{\bt}{A}{1_B}}:\mcg{\pb{(\pb{B}{A}{B})}{A}{B}}\longrightarrow\mcg{\pb{B}{A}{B}},\]
		where the pullback $\mcg{\pb{B}{A}{(\pb{B}{A}{B})}}$ is that of $\mcg{t,s}\pi_1$, while the pullback $\mcg{\pb{(\pb{B}{A}{B})}{A}{B}}$ is that of $\mcg{t}\pi_2,\mcg{s}$ (see Propositions \ref{slice-prof-iff-pullb} and \ref{pullb-as-bif}). Assume that
		\[\bm{v}=\begin{bmatrix}
		\mcg{l}_1\\\prm{\mcg{l}}_1
		\end{bmatrix},~~\bm{w}=\begin{bmatrix}
		\mcg{m}_1\\\prm{\mcg{m}}_1
		\end{bmatrix},~~\bm{x}=\begin{bmatrix}
		\mcg{n}_1\\\prm{\mcg{n}}_1
		\end{bmatrix}\]
		are $\mcg{B}$-objects such that 
		\[\ag{\bm{x}}{\bm{\ag{w}{v}}}\in\mcg{\pb{B}{A}{(\pb{B}{A}{B})}}.\]
		This means
		\[\mcg{t}(\bm{x})=\mcg{s}\pi_1(\bm{w,v})=\mcg{s}(w),\]
		and
		\[\mcg{t}(\bm{w})=\mcg{s}(\bm{v}).\]
		But these two conditions together are equivalent to stating that
		\[\bm{\ag{\ag{x}{w}}{v}}\in\mcg{\pb{(\pb{B}{A}{B})}{A}{B}};\]
		therefore,
		\begin{align*}
		\bt\circ(\pb{1_\mcg{B}}{\mcg{A}}{\bt})(\bm{x,\ag{v}{w}})&=\bm{x\bt(w\bt v)}\\
		&=\begin{bmatrix}
		\mcg{n}_1\\\prm{\mcg{n}}_1
		\end{bmatrix} \bt \begin{bmatrix}
		\mcg{l}_1\mcg{m}_1\\\prm{\mcg{l}}_1\prm{\mcg{m}}_1
		\end{bmatrix}\\
		&=\begin{bmatrix}
		\mcg{l}_1\mcg{m}_1\mcg{n}_1\\\prm{\mcg{l}}_1\prm{\mcg{m}}_1\prm{\mcg{n}}_1
		\end{bmatrix}\\
		&=\begin{bmatrix}
		\mcg{m}_1\mcg{n}_1\\\prm{\mcg{m}}_1\prm{\mcg{n}}_1
		\end{bmatrix} \bt \begin{bmatrix}
		\mcg{l}_1\\\prm{\mcg{l}}_1
		\end{bmatrix}\\
		&=\bm{(x\bt w)\bt v}\\
		&=\bt \circ (\pb{\bt}{\mcg{A}}{1_\mcg{B}})(\bm{\ag{x}{w},v}).
		\end{align*}
		
		Similarly, for $\mcg{B}$-morphisms $\bm{c,d,e}$ it can shown that
		\[\bm{\ag{c}{\ag{d}{e}}}\in\mcg{\pb{B}{A}{(\pb{B}{A}{B})}}~~~\text{iff}~~~~\bm{\ag{\ag{c}{d}}{e}}\in\mcg{\pb{(\pb{B}{A}{B})}{A}{B}}.\]
		Now let
		\[\bm{c}=\begin{pmatrix}
		&\bm{h}^1& \\\bm{v}& &\bm{w}\\ &\bm{h}^2& 
		\end{pmatrix},~~\bm{d}=\begin{pmatrix}
		&\bm{h}^2& \\\prm{\bm{v}}& &\prm{\bm{w}}\\ &\bm{h}^3& 
		\end{pmatrix},~~\bm{e}=\begin{pmatrix}
		&\bm{h}^3& \\ \bm{v}^{\prime\prime}& &\bm{w}^{\prime\prime}\\ &\bm{h}^4& 
		\end{pmatrix}\]
		so that $\bm{\ag{c}{\ag{d}{e}}}$ is a morphism in $\mcg{\pb{B}{A}{(\pb{B}{A}{B})}}$. Then
		\begin{align*}
		\bt\circ(\pb{1_\mcg{B}}{\mcg{A}}{\bt})(\bm{c,\ag{d}{e}})&=\bm{c\bt(d\bt e)}\\
		&=\begin{pmatrix}
		&\bm{h}^1& \\\bm{v}& &\bm{w}\\ &\bm{h}^2& 
		\end{pmatrix} \bt \begin{pmatrix}
		&\bm{h}^2& \\ \prm{\bm{v}}\bt\bm{v}^{\prime\prime}& &\prm{\bm{w}}\bt\bm{w}^{\prime\prime}\\ &\bm{h}^4& 
		\end{pmatrix}\\
		&=\begin{pmatrix}
		&\bm{h}^1& \\ \bm{v}\bt(\prm{\bm{v}}\bt\bm{v}^{\prime\prime})& &\bm{w}\bt(\prm{\bm{w}}\bt\bm{w}^{\prime\prime})\\ &\bm{h}^4& 
		\end{pmatrix}\\
		&=\begin{pmatrix}
		&\bm{h}^1& \\ (\bm{v}\bt\prm{\bm{v}})\bt\bm{v}^{\prime\prime}& &(\bm{w}\bt\prm{\bm{w}})\bt\bm{w}^{\prime\prime}\\ &\bm{h}^4& 
		\end{pmatrix}\\
		&=\begin{pmatrix}
		&\bm{h}^1& \\ \bm{v}\bt\prm{\bm{v}}& &\bm{w}\bt\prm{\bm{w}}\\ &\bm{h}^3& 
		\end{pmatrix} \bt \begin{pmatrix}
		&\bm{h}^3& \\ \bm{v}^{\prime\prime}& &\bm{w}^{\prime\prime}\\ &\bm{h}^4& 
		\end{pmatrix}\\
		&=\bm{(c\bt d)\bt e}\\
		&=\bt \circ (\pb{\bt}{\mcg{A}}{1_\mcg{B}})(\bm{\ag{c}{d},e}).
		\end{align*}
		From the above results on objects and morphisms we conclude:
		\[\bt\circ(\pb{1_\mcg{B}}{\mcg{A}}{\bt})=\bt \circ (\pb{\bt}{\mcg{A}}{1_\mcg{B}}),\]
		and the proof is complete.
	\end{mypr}
	
	\begin{rem}
		\label{DLC-T}
		At this point, it is interesting to take a look at the transpose double category $\dlc^\T$. It consists of the following data:
		\[\dlc^\T=\langle\mcg{\prm{A},\prm{B},\prm{s},\prm{t},\prm{i},\prm{\bt}}\rangle,\]
		where $\mcg{\prm{A}}$ is exactly the auxiliary category defined in part (b.1) of Definition \ref{DLC-def}; objects of $\mcg{\prm{B}}$ are $\mcg{A}$-morphisms $\bm{h}=\left[\begin{smallmatrix}
		\eb{f}\\ \eb{g}
		\end{smallmatrix} \right] $, morphisms of $\prm{\mcg{B}}$ are the $\mcg{B}$-morphisms (i.e. the cubicles)--this time viewed as arrows between $\mcg{A}$-morphisms--and the functors $\mcg{\prm{s},\prm{t},\prm{i},\prm{\bt}}$ are determined according to Definition \ref{transpose-def}.
	\end{rem}
	\rule{0pt}{5mm}

	\section{Some of the basic properties}
	\label{some-elem-props}
	Now we state and prove some basic properties of $\dlc$. First and foremost, we have the fundamental property of Klein-invariance. We prove this as a corollary of Theorems \ref{DLC-hor-sd} and \ref{DLC-ver-sd}; in these theorems, respectively, we explicitly construct the isomorphism double functors \((-)^\#:\dlc^\lra\cong\dlc\) and \((-)^\star:\dlc^\uda\cong\dlc\). The Klein-invariance of $\dlc$ follows immediately afterwards. Next, we will take a look at (binary) horizontal products and coproducts in $ \dlc $.
	
	\subsection{Klein-invariance}
	
	\begin{thm}
		\label{DLC-hor-sd}
		For every set $\Sigma$ and every nonempty set $\Gamma$, the double category $\dlc_{\Gamma,\Sigma}$ is horizontally self-dual.
	\end{thm}
	
	\begin{mypr}
		Assume an arbitrary set $\Sigma$ and an arbitrary nonempty set $\Gamma$. Consider the double category $\dlc=\dlc_{\Gamma,\Sigma}$. Observe that for $\dlc=\langle\mcg{A,B,s,t,i,\bowtie}\rangle$, the double category $\dlc^\lra$ consists of the following data:
		\[\dlc^\lra=\langle\op{\mcg{A}},\op{\mcg{B}},\mcg{s^\lra,t^\lra,i^\lra,\bt^\lra}\rangle.\]
		Now we define
		\[(-)^\#\eqd\ag{(-)^{\#^0}}{(-)^{\#^1}}:\dlc^\lra\longrightarrow\dlc,\]
		where $(-)^{\#^0},(-)^{\#^1}$ are functors we define below. Firstly, \[(-)^{\#^0}:\op{\mcg{A}}\longrightarrow\mcg{A}\] has the following data:
		\begin{itemize}
			\item on $\op{\mcg{A}}$-objects: for any $\eb{P=\ag{D}{S}}$ with $\eb{D}=\left[\begin{smallmatrix}
			\eb{A}\\\eb{B}
			\end{smallmatrix} \right]$ we define
			\[\eb{P}^{\#^0}\eqd\ag{\begin{bmatrix}
				\eb{B}\\\eb{A}
				\end{bmatrix}}{\eb{S}},\]
			in which we have identified the isomorphic coproducts 
			\[F^+\msf{A}\sqcup K^+\msf{B}\cong K^+\msf{B}\sqcup F^+\msf{A}\]
			and also
			\[G^-\msf{A}\sqcup H^-\msf{B}\cong H^-\msf{B}\sqcup G^-\msf{A},\]
			so that the super-matrix $\eb{S}=\langle\Sigma\xrightarrow{~S^+~}F^+\msf{A}\sqcup K^+\msf{B},~\Sigma\xrightarrow{~S^-~}G^-\msf{A}\sqcup H^-\msf{B}\rangle$ remains the same;
			
			\item on $\op{\mcg{A}}$-morphisms: for any $\op{\bm{h}}:\eb{P}_2\longrightarrow\eb{P}_1$, where $\bm{h}=\left[\begin{smallmatrix}
			\eb{f}\\\eb{g}
			\end{smallmatrix} \right] :\eb{P}_1\longrightarrow\eb{P}_2$ is an $\mcg{A}$-morphism, we define
			\[(\op{\bm{h}})^{\#^0}\eqd\begin{bmatrix}
			\eb{g}\\\eb{f}
			\end{bmatrix}:\eb{P}_2^{\#^0}\longrightarrow\eb{P}_1^{\#^0}.\]
			Note that $(\op{\bm{h}})^{\#^0}$ satisfies the PHS and NHS conditions (in fact, the corresponding commutative diagram  for $(\op{\bm{h}})^{\#^0}$ is a ``vertical flip'' of Diagram \ref{PHS-NHS}), hence it is a well-defined $\mcg{A}$-morphism.
		\end{itemize}
		On the other hand, $(-)^{\#^1}:\op{\mcg{B}}\longrightarrow\mcg{B}$ is such that:
		\begin{itemize}
			\item on $\op{\mcg{B}}$-objects: for any $\bm{v}=\left[\begin{smallmatrix}
			\ag{\mu}{\nu}\\\ag{\theta}{\zeta}
			\end{smallmatrix} \right] :\eb{P}^1\rightharpoonup\eb{P}^2$ we define
			\[\bm{v}^{\#^1}\eqd\begin{bmatrix}
			\ag{\theta}{\zeta}\\\ag{\mu}{\nu}
			\end{bmatrix}:(\eb{P}^1)^{\#^0}\rightharpoonup(\eb{P}^2)^{\#^0}\]
			which satisfies the PVS and NVS conditions (with commutativity diagrams the same as Diagram \ref{PVS-NVS}), and is thus a well-defined $\mcg{B}$-object;
			
			\item on $\op{\mcg{B}}$-morphisms: for any $\op{\bm{c}}:\bm{v}_2\Longrightarrow\bm{v}_1$ where 
			\[\bm{c}=\begin{pmatrix}
			&\bm{h}^1& \\\bm{v}_1& &\bm{v}_2\\ &\bm{h}^2& 
			\end{pmatrix}:\bm{v}_1\Longrightarrow\bm{v}_2\]
			is a $\mcg{B}$-morphism, we define
			\[(\op{\bm{c}})^{\#^1}\eqd\begin{pmatrix}
			&(\op{(\bm{h}^1)})^{\#^0}& \\\bm{v}_2& &\bm{v}_1\\ &(\op{(\bm{h}^2)})^{\#^0}& 
			\end{pmatrix}:\bm{v}_2\Longrightarrow\bm{v}_1,\]
			which clearly satisfies the commutativity conditions of Diagram \ref{four-parall} and is, therefore, well-defined.
		\end{itemize}
		
		\par Now, from the definition of $\mcg{s^\lra}$, the following holds for every $\op{\mcg{B}}$-object\\ $\bm{v}=\left[ \begin{smallmatrix}
		\ag{\mu}{\nu}\\\ag{\theta}{\zeta}
		\end{smallmatrix}\right]:\eb{P}^1\rightharpoonup\eb{P}^2 $:
		\begin{align}
		(\uo{\mcg{s}}^\lra(\bm{v}))^{\#^0}=(\uo{\mcg{s}}(\bm{v}))^{\#^0}=(\eb{P}_1)^{\#^0}=\uo{\mcg{s}}(\bm{v}^{\#^1}). \tag{1}
		\end{align}
		Also, for every $\op{\mcg{B}}$-morphism $\op{\bm{c}}=\op{\left(\begin{smallmatrix}
			&\bm{h}^1& \\\bm{v}_1& &\bm{v}_2\\ &\bm{h}^2& 
			\end{smallmatrix} \right)}:\bm{v}_2\Longrightarrow\bm{v}_1$,
		\begin{align}
		(\um{\mcg{s}}^\lra(\op{\bm{c}}))^{\#^0}=(\op{(\um{\mcg{s}}(\bm{c}))})^{\#^0}=(\op{(\bm{h}^1)})^{\#^0}=\um{\mcg{s}}((\op{\bm{c}})^{\#^1}). \tag{2}
		\end{align}
		Equations (1) and (2) together imply
		\begin{align}
		(\mcg{s^\lra}(-))^{\#^0}=\mcg{s}((-)^{\#^1}). \tag{3}
		\end{align}
		Likewise, one can deduce
		\begin{align}
		(\mcg{t^\lra}(-))^{\#^0}=\mcg{t}((-)^{\#^1}). \tag{4}
		\end{align}
		\par Next, from the definition of $\mcg{i^\lra}$, for every $\op{\mcg{A}}$-object $\eb{P}$ we have
		\begin{align}
		(\uo{\mcg{i}}^\lra(\eb{P}))^{\#^1}=(\uo{\mcg{i}}(\eb{P}))^{\#^1}=(\prm{1}_\eb{P})^{\#^1}=\prm{1}_{\eb{P}^{\#^0}}. \tag{5}
		\end{align}
		Also, for every $\op{\mcg{A}}$-morphism $\op{\bm{h}}:\eb{P}_2\longrightarrow\eb{P}_1$ where $\bm{h}=\left[\begin{smallmatrix}
		\eb{f}\\ \eb{g}
		\end{smallmatrix}\right]:\eb{P}_1\longrightarrow\eb{P}_2$,
		\begin{align}
		(\um{\mcg{i}}^\lra(\op{\bm{h}}))^{\#^1}&=(\op{(\um{\mcg{i}}(\bm{h}))})^{\#^1}=\left(\op{\begin{pmatrix}
			&\bm{h}& \\\prm{1}_{\eb{P}_1}& &\prm{1}_{\eb{P}_2}\\ &\bm{h}& 
			\end{pmatrix}} \right)^{\#^1} \notag\\
		&=\begin{pmatrix}
		&(\op{\bm{h}})^{\#^0}& \\\prm{1}_{\eb{P}_2}& &\prm{1}_{\eb{P}_1}\\ &(\op{\bm{h}})^{\#^0}& 
		\end{pmatrix}=\um{\mcg{i}}((\op{\bm{h}})^{\#^0}). \tag{6}
		\end{align}
		Equations (5) and (6) together imply
		\begin{align}
		(\mcg{i}^\lra(-))^{\#^1}=\mcg{i}((-)^{\#^0}). \tag{7}
		\end{align}
		Now consider Diagram \ref{hor-sd-diag}.
		\begin{figure}[h]
			\begin{center}
				\begin{tikzpicture}[commutative diagrams/every diagram]
				
				\node (N1) at (0.8,0.3cm) {$\mcg{B}$};
				
				\node (N2) at (-1.3,-0.3cm) {$\op{\mcg{B}}$};
				
				\node (N3) at (4.7,-0.3cm) {$\mcg{B}$};
				
				\node (N4) at (-5.2,0.3cm) {$\op{\mcg{B}}$};
				
				\node (N5) at (3,2cm) {$\pb{\mcg{B}}{\mcg{A}}{\mcg{B}}$};
				
				\node (N6) at (-3.5,-2cm) {$\op{\mcg{A}}$};
				
				\node (N7) at (-3,2cm) {$\pb{\op{\mcg{B}}}{\op{\mcg{A}}}{\op{\mcg{B}}}$};
				
				\node (N8) at (2.5,-2cm) {$\mcg{A}$};
				
				\node (N9) at (3,5cm) {$\mcg{B}$};
				
				\node (N10) at (-3.5,-5cm) {$\op{\mcg{B}}$};
				
				\node (N11) at (-3,5cm) {$\op{\mcg{B}}$};
				
				\node (N12) at (2.5,-5cm) {$\mcg{B}$};
				
				\path[commutative diagrams/.cd, every arrow, every label]
				(N7) edge node{$\bowtie^\lra$} (N11)
				(N5) edge node[swap]{$\bt$} (N9)
				(N7) edge node[swap]{$\rho_1$} (N4)
				(N7) edge node{$\rho_2$} (N2)
				(N5) edge node[swap]{$\pi_1$} (N1)
				(N5) edge node{$\pi_2$} (N3)
				(N4) edge node[swap]{$\mcg{t}^\lra$} (N6)
				(N2) edge node{$\mcg{s}^\lra$} (N6)
				(N1) edge node[swap]{$\mcg{t}$} (N8)
				(N3) edge node{$\mcg{s}$} (N8)
				(N6) edge node[swap]{$\mcg{i}^\lra$} (N10)
				(N8) edge node{$\mcg{i}$} (N12)
				(N4) edge[-,line width=6pt,draw=white] node{} (N1)
				(N2) edge[-,line width=6pt,draw=white] node{} (N3)
				(N4) edge node[near end]{$(-)^{\#^1}$} (N1)
				(N2) edge node[near start,swap]{$(-)^{\#^1}$} (N3)
				(N11) edge node{$(-)^{\#^1}$} (N9)
				(N10) edge node{$(-)^{\#^1}$} (N12)
				(N6) edge node{$(-)^{\#^0}$} (N8)
				(N7) edge[dashed] node{$\left( (-)^{\#^1}\circ\rho_1,(-)^{\#^1}\circ\rho_2\right) $} (N5);
				\end{tikzpicture}
			\end{center}
			\caption{} \label{hor-sd-diag}
		\end{figure}
		In this diagram, $\rho_1,\rho_2$ are the projections for\\ $\pb{\op{\mcg{B}}}{\op{\mcg{A}}}{\op{\mcg{B}}}$, making the left hand side diamond commute. Now, in addition, from (3), (4), and (7) we know that, respectively, the lower-right parallelogram, the lower-left parallelogram, and the bottom rectangle all commute. Therefore, we have the following equations:
		\begin{align*}
		\mcg{t}\circ(-)^{\#^1}\circ\rho_1&=(-)^{\#^0}\circ\mcg{t}^\lra\circ\rho_1\\
		&=(-)^{\#^0}\circ\mcg{s}^\lra\circ\rho_2\\
		&=\mcg{s}\circ(-)^{\#^1}\circ\rho_2.
		\end{align*}
		So, by universality of the pullback $\agg{\mcg{\pb{B}{A}{B}}}{\pi_1}{\pi_2}$, there exists a unique arrow
		\[\left( (-)^{\#^1}\circ\rho_1,(-)^{\#^1}\circ\rho_2\right):\pb{\op{\mcg{B}}}{\op{\mcg{A}}}{\op{\mcg{B}}}\longrightarrow\mcg{\pb{B}{A}{B}}\]
		making the upper-left and upper-right parallelograms commute. Now, in order to prove the double functoriality of $(-)^{\#}$, it remains to prove the commutativity of the top rectangle.\\
		\par On objects we have
		\begin{align*}
		(\bm{v}\bt^\lra_0\bm{w})^{\#^1}&=(\bm{v}\bt_0\bm{w})^{\#^1}\\
		&=\bm{v}^{\#^1}\bt_0\bm{w}^{\#^1}\\
		&=\bt_0\circ\left((-)^{\#^1}\circ\rho_1,(-)^{\#^1}\circ\rho_2 \right)(\bm{v,w});
		\end{align*}
		on morphisms,
		\begin{align*}
		(\op{\bm{c}}\bt^\lra_1\op{\bm{d}})^{\#^1}&=(\op{(\bm{c}\bt_1\bm{d})})^{\#^1}\\
		&=(\op{\bm{c}})^{\#^1}\bt_1(\op{\bm{d}})^{\#^1}\\
		&=\bt_1\circ\left((-)^{\#^1}\circ\rho_1,(-)^{\#^1}\circ\rho_2 \right)(\op{\bm{c}},\op{\bm{d}}).
		\end{align*}
		These complete the proof of commutativity of Diagram \ref{hor-sd-diag}. Thus, $(-)^{\#}$ is actually a double functor.\\
		
		\par The final step of the proof is to show that $(-)^{\#}$ is an isomorphism. For this, according to Proposition \ref{hor-iso-iff-both-iso}, it suffices to show that both $(-)^{\#^0}$ and $(-)^{\#^1}$ are isomorphisms (as ordinary functors). More precisely, we show that the \textit{opposite} functors
		\[(-)^{\op{(\#^0)}}:\mcg{A}\longrightarrow\op{\mcg{A}}~~~~\text{and}~~~~(-)^{\op{(\#^1)}}:\mcg{B}\longrightarrow\op{\mcg{B}}\]
		are inverses to $(-)^{\#^0}$ and $(-)^{\#^1}$, respectively.\\
		\par For any $\op{\mcg{A}}$-object $\eb{P=\ag{D}{S}}$ with $\eb{D}=\left[ \begin{smallmatrix}
		\eb{A}\\ \eb{B}
		\end{smallmatrix}\right] $,
		\[(\eb{P}^{\#^0})^{\op{(\#^0)}}=\ag{\begin{bmatrix}
			\eb{B}\\\eb{A}
			\end{bmatrix}}{\eb{S}}^{\op{(\#^0)}}=\ag{\begin{bmatrix}
			\eb{A}\\\eb{B}
			\end{bmatrix}}{\eb{S}}=\eb{P};\]
		thus, on $\op{\mcg{A}}$-objects, $((-)^{\#^0})^{\op{(\#^0)}}=1_{\op{\mcg{A}}}$. Similarly, on $\mcg{A}$-objects, $((-)^{\op{(\#^0)}})^{\#^0}=1_{\mcg{A}}$.\\
		\par On the other hand, for any $\op{\mcg{A}}$-morphism $\op{\bm{h}}:\eb{P}_2\longrightarrow\eb{P}_1$, where $\bm{h}=\left[\begin{smallmatrix}
		\eb{f}\\\eb{g}
		\end{smallmatrix} \right]:\eb{P}_1\longrightarrow\eb{P}_2 $ is an $\mcg{A}$-morphism,
		\[((\op{\bm{h}})^{\#^0})^{\op{(\#^0)}}=\begin{bmatrix}
		\eb{g}\\\eb{f}
		\end{bmatrix}^{\op{(\#^0)}}=\op{\begin{bmatrix}
			\eb{f}\\\eb{g}
			\end{bmatrix}}=\op{\bm{h}};\]
		that is, on $\op{\mcg{A}}$-morphisms, $((-)^{\#^0})^{\op{(\#^0)}}=1_{\op{\mcg{A}}}$. Similarly, on $\mcg{A}$-morphisms,\\ $((-)^{\op{(\#^0)}})^{\#^0}=1_\mcg{A}$. We conclude from all the above that $(-)^{\#^0}$ is an isomorphism.\\
		
		\par For the case of $(-)^{\#^1}$ the reasoning goes as the following: for any $\op{\mcg{B}}$-object $\bm{v}=\left[\begin{smallmatrix}
		\ag{\mu}{\nu}\\ \ag{\theta}{\zeta}
		\end{smallmatrix} \right] $,
		\[(\bm{v}^{\#^1})^{\op{(\#^1)}}=\begin{bmatrix}
		\ag{\theta}{\zeta}\\ \ag{\mu}{\nu}
		\end{bmatrix}^{\op{(\#^1)}}=\begin{bmatrix}
		\ag{\mu}{\nu}\\ \ag{\theta}{\zeta}
		\end{bmatrix}=\bm{v}.\]
		On the other hand, for any $\op{\mcg{B}}$-morphism $\op{\bm{c}}$ where $\bm{c}=\left(\begin{smallmatrix}
		&\bm{h}^1& \\ \bm{v}_1& &\bm{v}_2\\ &\bm{h}^2& 
		\end{smallmatrix} \right) $ is a $\mcg{B}$-morphism,
		\[((\op{\bm{c}})^{\#^1})^{\op{(\#^1)}}=\begin{pmatrix}
		&(\op{(\bm{h}^1)})^{\#^0}& \\ \bm{v}_2& &\bm{v}_1\\ &(\op{(\bm{h}^2)})^{\#^0}& 
		\end{pmatrix}^{\op{(\#^1)}}=\op{\begin{pmatrix}
			&\bm{h}^1& \\ \bm{v}_1& &\bm{v}_2\\ &\bm{h}^2& 
			\end{pmatrix}}=\op{\bm{c}},\]
		so, $((-)^{\#^1})^{\op{(\#^1)}}=1_{\op{\mcg{B}}}$. By a similar argument, $((-)^{\op{(\#^1)}})^{\#^1}=1_{\mcg{B}}$. Consequently, $(-)^{\#^1}$, too, is an isomorphism. This completes the proof of the horizontal self-duality of $(-)^{\#}$.
	\end{mypr}
	
	\rule{0pt}{5mm}
	\par For the vertical self-duality a different line of reasoning is to be taken. First of all, recall that the $*$-autonomous structure of \chu~includes a (covariantly written) self-duality functor
	\[(-)^*:\op{\chu}\longrightarrow\chu\]
	that sends every Chu space $\msf{A}=\agg{A}{r}{X}$ to the Chu space $\msf{A}^*=\agg{X}{\breve{r}}{A}$, and sends every $\op{f}=\op{\ag{f^+}{f^-}}:\msf{B}\longrightarrow\msf{A}$ to the Chu transform $(\op{f})^*=\ag{f^-}{f^+}:\msf{B}^*\longrightarrow\msf{A}^*$. To be able to work with $\ag{F}{G}$-dialgebras, we need a \textit{contravariant} version of the above self-duality. Therefore, we define:
	
	\begin{dfn}
		\label{contra-tilde-def}
		The \textit{contravariant} self-duality functor $\widetilde{(-)}$ is defined as
		\[\widetilde{(-)}\eqd(\rvr{\chu})^*:\chu\longrightarrow\chu,\]
		for which the reverser $\rvr{(-)}$ was introduced in Definition \ref{the-reverser}.
	\end{dfn}
	
	Further, we expand the above definition to dialgebras:
	
	\begin{dfn}
		\label{dialg-tilde}
		Let $F,G$ be two Chu endofunctors, and let $\eb{A}=\ag{\msf{A}}{F\msf{A}\xrightarrow{~\alpha~}G\msf{A}}$ be an $\ag{F}{G}$-dialgebra. Then we define
		\[\widetilde{\eb{A}}\eqd\ag{\msf{A}}{\widetilde{a}:\wt{G}\msf{A}\longrightarrow\wt{F}\msf{A}},\]
		where $\wt{F}\msf{A}$ abbreviates $\wt{(F\msf{A})}$, and so on.
	\end{dfn}
	
	\begin{prp}
		\label{chuFG-iff-chuGtFt}
		Let $\ag{\msf{A}}{\alpha},\ag{\msf{B}}{\beta}\in\dlg{\chu}{F}{G}$. A Chu transform $h:\msf{A\longrightarrow B}$ yields an $\ag{F}{G}$-dialgebra homomorphism $\eb{h:A\longrightarrow B}$ if and only if $h$ yields an $\ag{\wt{G}}{\wt{F}}$-dialgebra homomorphism $\wt{\eb{h}}:\wt{\eb{B}}\longrightarrow\wt{\eb{A}}$. Moreover, $\dlg{\chu}{F}{G}\cong\dlg{\chu}{\wt{G}}{\wt{F}}$.
	\end{prp}
	
	\begin{mypr}
		In Diagram \ref{chuFG-iff-chuGtFt-diag}, the commutativity of the left square is equivalent to the commutativity of the right one.\\
		\begin{figure}[h]
			\begin{center}
				\begin{tikzpicture}[commutative diagrams/every diagram]
				
				\node (N1) at (-2,2cm) {$G\msf{A}$};
				
				\node (N2) at (-5,2cm) {$F\msf{A}$};
				
				\node (N3) at (-5,-1cm) {$F\msf{B}$};
				
				\node (N4) at (-2,-1cm) {$G\msf{B}$};
				
				\node (M1) at (1,2cm) {$\wt{F}\msf{A}$};
				
				\node (M2) at (4,2cm) {$\wt{G}\msf{A}$};
				
				\node (M3) at (4,-1cm) {$\wt{G}\msf{B}$};
				
				\node (M4) at (1,-1cm) {$\wt{F}\msf{B}$};
				
				\path[commutative diagrams/.cd, every arrow, every label]
				(N2) edge node[swap]{$\alpha$} (N1)
				(N2) edge node{$Fh$} (N3)
				(N3) edge node[swap]{$\beta$} (N4)
				(N1) edge node{$Gh$} (N4)
				(M2) edge node[swap]{$\wt{\alpha}$} (M1)
				(M4) edge node{$\wt{F}h$} (M1)
				(M3) edge node[swap]{$\wt{\beta}$} (M4)
				(M3) edge node{$\wt{G}h$} (M2);
				\end{tikzpicture}
			\end{center}
			\caption{} \label{chuFG-iff-chuGtFt-diag}
		\end{figure}
		\par \noindent It is obvious that $\wt{(\wt{(\eb{A})})}=\eb{A}$ for every $\eb{A}$ and $\wt{(\wt{(\eb{h})})}=\eb{h}$ for every $\eb{h}$. So, the map sending every $\eb{A}\in\dlg{\chu}{F}{G}$~ to ~$\wt{\eb{A}}\in\dlg{\chu}{\wt{G}}{\wt{F}}\text{\tiny~}$, together with the map which sends every $\ag{F}{G}$-dialgebra homomorphism $\eb{h}$ to the $\langle\wt{G},\wt{F}\rangle$-dialgebra homomorphism $\wt{\eb{h}}$, constitute an isomorphism functor $\dlg{\chu}{F}{G}\cong\dlg{\chu}{\wt{G}}{\wt{F}}$.
	\end{mypr}
	
	Even more---by some abuse of notation---we expand the definition of $\wt{(-)}$ to the case of the category $\mcg{A}$:
	
	\begin{dfn}
		\label{expansion-to-A}
		For every $\mcg{A}$-object $\eb{P=\ag{D}{S}}$ with $\eb{D}=\left[\begin{smallmatrix}
		\eb{A}\\\eb{B}
		\end{smallmatrix} \right],\eb{S}=\ag{S^+}{S^-} $ we define
		\[\wt{\eb{P}}\eqd\ag{\begin{bmatrix}
			\wt{\eb{A}}\\ \wt{\eb{B}}
			\end{bmatrix}}{\wt{\eb{S}}},\]
		where $\wt{\eb{S}}\eqd\ag{S^-}{S^+}$. Also, for every $\mcg{A}$-morphism $\bm{h}=\left[ \begin{smallmatrix}
		\eb{f}\\ \eb{g}
		\end{smallmatrix}\right]:\eb{P}_1\longrightarrow\eb{P}_2 $, we define
		\[\wt{\bm{h}}\eqd\begin{bmatrix}
		\wt{\eb{f}}\\ \wt{\eb{g}}
		\end{bmatrix}:\wt{\eb{P}_2}\longrightarrow\wt{\eb{P}_1}.\]
	\end{dfn}
	
	It is evident that $\wt{\bm{h}}$ satisfies the PHS and NHS conditions (the corresponding commutative diagram for $\wt{\bm{h}}$ is a ``180\textdegree~ rotation'' of Diagram \ref{PHS-NHS}), thus it is a well-defined $\mcg{A}$-morphism. This way, $\wt{(-)}:\mcg{A\longrightarrow A}$ can be seen as an endofunctor.
	
	\begin{prp}
		\label{endo-tilde-sd-A}
		The endofunctor $\wt{(-)}:\mcg{A\longrightarrow A}$ is a (contravariant self-duality) isomorphism.
	\end{prp}
	
	\begin{mypr}
		By definition, for every $\eb{P}\in\mcg{A}$, $\wt{\left( \wt{\eb{P}}\right) }=\eb{P}$; also, for every $\bm{h}\in\mcg{A}$ we have \(\wt{\left( \wt{\bm{h}}\right) }=\bm{h}\).
		
	\end{mypr}
	
	\begin{prp}
		\label{induces-contra-sd-Fluid}
		The endofunctor $\wt{(-)}:\chu\longrightarrow\chu$ induces a (contravariant self-duality) isomorphism on $\cat{Fluid}_\chu$.
	\end{prp}
	
	\begin{mypr}
		Let an endofunctor $\ol{\ag{-}{-}}:\flu\longrightarrow\flu$ be defined as the following: for every dialgebra $\eb{A}\in\flu$,
		\[\ol{\eb{A}}\eqd\wt{\eb{A}};\]
		on the other hand, for every $\flu$-morphism $\ag{\varphi}{\psi}:\eb{A}^1\longrightarrow\eb{A}^2$ we define
		\[\ol{\langle\varphi,\psi\rangle}\eqd\langle\wt{\psi},\wt{\varphi}\rangle:\wt{\eb{A}^2}\longrightarrow\wt{\eb{A}^1}.\]
		(Notice the reversal of the order of the letters $\varphi,\psi$.) Now it is evident that for every $\langle\varphi,\psi\rangle$,
		\[\ol{(\ol{\langle\varphi,\psi\rangle})}=\langle\varphi,\psi\rangle,\]
		as desired.
	\end{mypr}
	
	Now, equipped with all the above, we are ready to state and prove our next theorem:
	
	\begin{thm}
		\label{DLC-ver-sd}
		For every set $\Sigma$ and every nonempty set $\Gamma$, the double category $\dlc_{\Gamma,\Sigma}$ is vertically self-dual.
	\end{thm}
	
	\begin{mypr}
		Assume that an arbitrary set $\Sigma$ and an arbitrary nonempty set $\Gamma$ are given. Consider the double category $\dlc=\dlc_{\Gamma,\Sigma}$. Observe that the double category $\dlc^\uda$ consists of the following data:
		\[\dlc^\uda=\langle\mcg{A,B^\uda,s^\uda,t^\uda,i^\uda,\bt^\uda}\rangle.\]
		Now we define:
		\[(-)^\star\eqd\ag{(-)^{\star^0}}{(-)^{\star^1}}:\dlc^\uda\longrightarrow\dlc,\]
		where
		\[(-)^{\star^0}\eqd\wt{(-)}:\mcg{A\longrightarrow A},\]
		and the functor $(-)^{\star^1}:\mcg{B^\uda\longrightarrow B}$ is defined as the following:
		\begin{itemize}
			\item for any $\mcg{B^\uda}$-object $\bm{v^\uda}:\eb{P}^2\rightharpoonup\eb{P}^1$ where $\bm{v}=\left[\begin{smallmatrix}
			\ag{\mu}{\nu}\\ \ag{\theta}{\zeta}
			\end{smallmatrix} \right]:\eb{P}^1\rightharpoonup\eb{P}^2$ is a $\mcg{B}$-object (equivalently, an $\prm{\mcg{A}}$-morphism),
			\[(\bm{v^\uda})^{\star^1}\eqd\begin{bmatrix}
			\ol{\ag{\mu}{\nu}}\\ \ol{\ag{\theta}{\zeta}}
			\end{bmatrix}=\begin{bmatrix}
			\ag{\wt{\nu}}{\wt{\mu}}\\ \ag{\wt{\zeta}}{\wt{\theta}}
			\end{bmatrix}:\wt{\eb{P}^2}\rightharpoonup\wt{\eb{P}^1},\]
			which satisfies the PVS and NVS conditions (with commutativity diagram being a ``180\textdegree~ rotation'' of Diagram \ref{PVS-NVS}), and is thus a well-defined $\mcg{B}$-object;
			
			\item for any $\mcg{B^\uda}$-morphism $\bm{v}_1^\uda\xR{c^\uda}\bm{v}_2^\uda$ where $\bm{c}=\left(\begin{smallmatrix}
			&\bm{h}^1& \\\bm{v}_1& &\bm{v}_2\\ &\bm{h}^2& 
			\end{smallmatrix} \right) :\bm{v}_1\Longrightarrow\bm{v}_2$ is a $\mcg{B}$-morphism,
			\[(\bm{c^\uda})^{\star^1}\eqd\begin{pmatrix}
			&\wt{\bm{h}^2}& \\ (\bm{v}_1^\uda)^{\star^1}& &(\bm{v}_2^\uda)^{\star^1}\\ &\wt{\bm{h}^1}& 
			\end{pmatrix}:(\bm{v}_1^\uda)^{\star^1}\Longrightarrow(\bm{v}_2^\uda)^{\star^1},\]
			which satisfies the commutativity conditions of Diagram 15 and is therefore, well-defined.
		\end{itemize}
		Now, from the definition of $\mcg{s^\uda}$, the following holds for every $\mcg{B^\uda}$-object $\bm{v^\uda}:\eb{P}^2\rightharpoonup\eb{P}^1$ where $\bm{v}=\left[\begin{smallmatrix}
		\ag{\mu}{\nu}\\ \ag{\theta}{\zeta}
		\end{smallmatrix} \right] $:
		\begin{align}
		\left(\mcg{s}^\uda_0(\bm{v^\uda})\right)^{\star^0}=\wt{\mcg{t}_0(\bm{v})}=\wt{\eb{P}^2}=\mcg{s}_0\left((\bm{v^\uda})^{\star^1} \right); \tag{1}
		\end{align}
		also, for every $\mcg{B^\uda}$-morphism $\bm{v}_1^\uda\xR{c^\uda}\bm{v}_2^\uda$ where $\bm{c}=\left(\begin{smallmatrix}
		&\bm{h}^1& \\\bm{v}_1& &\bm{v}_2\\ &\bm{h}^2& 
		\end{smallmatrix} \right)$,
		\begin{align}
		\left(\mcg{s}^\uda_1(\bm{c^\uda})\right)^{\star^0}=\wt{\mcg{t}_1(\bm{c})}=\wt{\bm{h}^2}=\mcg{s}_1\left((\bm{c^\uda})^{\star^1} \right). \tag{2}
		\end{align}
		Equations (1) and (2) together imply
		\begin{align}
		\left( \mcg{s^\uda}(-)\right)^{\star^0}=\mcg{s}\left((-)^{\star^1} \right). \tag{3}
		\end{align}
		Likewise, one can deduce
		\begin{align}
		\left( \mcg{t^\uda}(-)\right)^{\star^0}=\mcg{t}\left((-)^{\star^1} \right). \tag{4}
		\end{align}
		Next, from the definition of $\mcg{i^\uda}$, for every $\mcg{A}$-object $\eb{P}$ we have
		\begin{align}
		\left(\mcg{i}^\uda_0 (\eb{P}) \right)^{\star^1}=\left(\left( \mcg{i}_0 (\eb{P})\right)^\uda \right)^{\star^1}= \left(\left( \prm{1}_\eb{P}\right)^\uda \right)^{\star^1}=\prm{1}_{\eb{P}^{\star^0}}; \tag{5}
		\end{align}
		also, for every $\mcg{A}$-morphism $\bm{h}=\left[\begin{smallmatrix}
		\eb{f}\\ \eb{g}
		\end{smallmatrix} \right] :\eb{P}_1\longrightarrow\eb{P}_2$,
		\begin{align}
		\left(\mcg{i}^\uda_0 (\bm{h}) \right)^{\star^1}&=\left(\left( \mcg{i}_0 (\bm{h})\right)^\uda \right)^{\star^1} \notag \\
		&=\left(\begin{pmatrix}
		&\bm{h}& \\ \prm{1}_{\eb{P}_1}& &\prm{1}_{\eb{P}_2}\\ &\bm{h} &
		\end{pmatrix}^\uda \right)^{\star^1} \notag \\
		&=\begin{pmatrix}
		&\wt{\bm{h}}& \\ \prm{1}_{\eb{P}_1^{\star^0}}& &\prm{1}_{\eb{P}_2^{\star^0}}\\ &\wt{\bm{h}}& 
		\end{pmatrix} \notag \\
		&=\mcg{i_1}\left( \bm{h}^{\star^0}\right). \tag{6}
		\end{align}
		Equations (5) and (6) together imply
		\begin{align}
		\left(\mcg{i^\uda}(-) \right)^{\star^1} = \mcg{i}\left((-)^{\star^0} \right).  \tag{7}
		\end{align}
		Now we consider Diagram \ref{ver-sd-diag}.
		\begin{figure}[h]
			\begin{center}
				\begin{tikzpicture}[commutative diagrams/every diagram]
				
				\node (N1) at (0.8,0.3cm) {$\mcg{B}$};
				
				\node (N2) at (-1.3,-0.3cm) {$\mcg{B}^\uda$};
				
				\node (N3) at (4.7,-0.3cm) {$\mcg{B}$};
				
				\node (N4) at (-5.2,0.3cm) {$\mcg{B}^\uda$};
				
				\node (N5) at (3,2cm) {$\pb{\mcg{B}}{\mcg{A}}{\mcg{B}}$};
				
				\node (N6) at (-3.5,-2cm) {$\mcg{A}$};
				
				\node (N7) at (-3,2cm) {$\pb{\mcg{B}^\uda}{\mcg{A}}{\mcg{B}^\uda}$};
				
				\node (N8) at (2.5,-2cm) {$\mcg{A}$};
				
				\node (N9) at (3,5cm) {$\mcg{B}$};
				
				\node (N10) at (-3.5,-5cm) {$\mcg{B}^\uda$};
				
				\node (N11) at (-3,5cm) {$\mcg{B}^\uda$};
				
				\node (N12) at (2.5,-5cm) {$\mcg{B}$};
				
				\path[commutative diagrams/.cd, every arrow, every label]
				(N7) edge node{$\bowtie^\uda$} (N11)
				(N5) edge node[swap]{$\bt$} (N9)
				(N7) edge node[swap]{$\sigma_1$} (N4)
				(N7) edge node{$\sigma_2$} (N2)
				(N5) edge node[swap]{$\pi_1$} (N1)
				(N5) edge node{$\pi_2$} (N3)
				(N4) edge node[swap]{$\mcg{t}^\uda$} (N6)
				(N2) edge node{$\mcg{s}^\uda$} (N6)
				(N1) edge node[swap]{$\mcg{t}$} (N8)
				(N3) edge node{$\mcg{s}$} (N8)
				(N6) edge node[swap]{$\mcg{i}^\uda$} (N10)
				(N8) edge node{$\mcg{i}$} (N12)
				(N4) edge[-,line width=6pt,draw=white] node{} (N1)
				(N2) edge[-,line width=6pt,draw=white] node{} (N3)
				(N4) edge node[near end]{$(-)^{\star^1}$} (N1)
				(N2) edge node[near start,swap]{$(-)^{\star^1}$} (N3)
				(N11) edge node{$(-)^{\star^1}$} (N9)
				(N10) edge node{$(-)^{\star^1}$} (N12)
				(N6) edge node{$(-)^{\star^0}$} (N8)
				(N7) edge[dashed] node{$\left( (-)^{\star^1}\circ\sigma_1,(-)^{\star^1}\circ\sigma_2\right) $} (N5);
				\end{tikzpicture}
			\end{center}
			\caption{} \label{ver-sd-diag}
		\end{figure}
		In this diagram, $\sigma_1,\sigma_2$ are the projections for $\mcg{\pb{B^\uda}{A}{B^\uda}}$, making the left hand side diamond commute. Now, in addition, from (3), (4), and (7) we know that, respectively, the lower-right parallelogram, the lower-left parallelogram, and the bottom rectangle all commute. Therefore, we have the following equations:
		\begin{align*}
		\mcg{t}\circ(-)^{\star^1}\circ\sigma_1&=(-)^{\star^0}\circ\mcg{t^\uda}\circ\sigma_1\\
		&=(-)^{\star^0}\circ\mcg{s^\uda}\circ\sigma_2\\
		&=\mcg{s}\circ(-)^{\star^1}\circ\sigma_2.\\
		\end{align*}
		So, by universality of the pullback $\agg{\mcg{\pb{B}{A}{B}}}{\pi_1}{\pi_2}$, there exists a unique arrow
		\[\left( (-)^{\star^1}\circ\sigma_1,(-)^{\star^1}\circ\sigma_2\right):\mcg{\pb{B^\uda}{A}{B^\uda}}\longrightarrow\mcg{\pb{B}{A}{B}}\]
		making the upper-left and the upper-right parallelograms commute. Now, in order to prove the double functoriality of $(-)^\star$, it remains to show the commutativity of the top rectangle.\\
		\par On objects we have
		\begin{align*}
		\left(\bm{v}^\uda\bt^\uda_0 \bm{w}^\uda \right)^{\star^1}&= \left((\bm{w}\bt_0 \bm{v})^\uda \right)^{\star^1}\\
		&=\left( \bm{v}^\uda\right)^{\star^1} \bt_0 \left( \bm{w}^\uda\right)^{\star^1}\\
		&=\bt_0\circ \left( (-)^{\star^1}\circ\sigma_1,(-)^{\star^1}\circ\sigma_2\right)\left(\bm{v^\uda},\bm{w^\uda} \right); 
		\end{align*}
		for morphisms,
		\begin{align*}
		\left(\bm{c}^\uda\bt^\uda_1 \bm{d}^\uda \right)^{\star^1}&= \left((\bm{d}\bt_1 \bm{c})^\uda \right)^{\star^1}\\
		&=\left( \bm{c}^\uda\right)^{\star^1} \bt_1 \left( \bm{d}^\uda\right)^{\star^1}\\
		&=\bt_1\circ \left( (-)^{\star^1}\circ\sigma_1,(-)^{\star^1}\circ\sigma_2\right)\left(\bm{c^\uda},\bm{d^\uda} \right). 
		\end{align*}
		These complete the proof of commutativity of Diagram \ref{ver-sd-diag}, and $(-)^\star$ is actually a double functor (from $\dlc^\uda$ to $\dlc$).\\
		
		\par The final step of the proof is to show that $(-)^\star$ is an isomorphism. For this, according to Proposition \ref{hor-iso-iff-both-iso}, it suffices to show that both $(-)^{\star^0}$ and $(-)^{\star^1}$ are isomorphisms as ordinary functors. By Proposition \ref{endo-tilde-sd-A} we already know that the functor $(-)^{\star^0}=\wt{(-)}:\mcg{A\longrightarrow A}$ is inverse to itself. On the other hand, we show that the functor
		\[(-)^{{\star^1}^\uda}:\mcg{B\longrightarrow B^\uda},\]
		which is defined below, is inverse to $(-)^{\star^1}$.\\
		
		\par We define $(-)^{{\star^1}^\uda}$ as the functor which sends every $\mcg{B}$-object\\ $\bm{v}=\left[\begin{smallmatrix}
		\ag{\mu}{\nu}\\ \ag{\theta}{\zeta}
		\end{smallmatrix} \right]:\eb{P}^1\rightharpoonup\eb{P}^2 $ to the $\mcg{B^\uda}$-object
		\[\bm{v}^{{\star^1}^\uda}\eqd\left(\left( \bm{v^\uda}\right)^{\star^1}  \right)^\uda =\begin{bmatrix}
		\ol{\ag{\mu}{\nu}} \\ \ol{\ag{\theta}{\zeta}}
		\end{bmatrix}^\uda = \begin{bmatrix}
		\ag{\wt{\nu}}{\wt{\mu}} \\ \ag{\wt{\zeta}}{\wt{\theta}}
		\end{bmatrix}^\uda : \wt{\eb{P}^1}\rightharpoonup\wt{\eb{P}^2}, \]
		and which sends every $\mcg{B}$-morphism $\bm{c}=\left(\begin{smallmatrix}
		&\bm{h}^1& \\ \bm{v}_1& &\bm{v}_2\\ &\bm{h}^2& 
		\end{smallmatrix} \right):\bm{v}_1\Longrightarrow\bm{v}_2 $ to the $\mcg{B}^\uda$-morphism
		\[\bm{c}^{{\star^1}^\uda}\eqd\left(\left(\bm{c^\uda} \right)^{\star^1} \right)^\uda=\begin{pmatrix}
		&\wt{\bm{h}^2}& \\ \left(\bm{v}^\uda_1 \right)^{\star^1} & &\left(\bm{v}^\uda_2 \right)^{\star^1}\\ &\wt{\bm{h}^1}& 
		\end{pmatrix}^\uda:\left(\bm{v}^\uda_1 \right)^{\star^1}\Longrightarrow\left(\bm{v}^\uda_2 \right)^{\star^1}. \]
		Now, for any $\mcg{B}$-object $\bm{v}=\left[\begin{smallmatrix}
		\ag{\mu}{\nu}\\ \ag{\theta}{\zeta}
		\end{smallmatrix} \right]$,
		\begin{align*}
		\left(\bm{v}^{{\star^1}^\uda}\right)^{\star^1}&=\left(\left(\left(\bm{v^\uda} \right)^{\star^1} \right)^\uda \right)^{\star^1} \\
		&=\left(\begin{bmatrix}
		\ol{\ag{\mu}{\nu}}\\ \ol{\ag{\theta}{\zeta}}
		\end{bmatrix}^\uda \right)^{\star^1}\\
		&=\begin{bmatrix}
		\ol{(\ol{\ag{\mu}{\nu}})}\\ \ol{(\ol{\ag{\theta}{\zeta}})}
		\end{bmatrix}\\
		&=\begin{bmatrix}
		\ag{\mu}{\nu}\\ \ag{\theta}{\zeta}
		\end{bmatrix}\\
		&=\bm{v}; \tag{8}
		\end{align*}
		also, using (8), we have for any $\mcg{B^\uda}$-object $\bm{v}^\uda=\left[ \begin{smallmatrix}
		\ag{\mu}{\nu}\\ \ag{\theta}{\zeta}
		\end{smallmatrix}\right]^\uda $,
		\begin{align}
		\left( \left( \bm{v^\uda}\right)^{\star^1} \right)^{{\star^1}^\uda}=\left(\left(\left(\left(\bm{v^\uda} \right)^{\star^1} \right)^\uda \right)^{\star^1} \right)^\uda=\bm{v^\uda}. \tag{9}
		\end{align}
		On the other hand, for any $\mcg{B}$-morphism $\bm{c}=\left(\begin{smallmatrix}
		&\bm{h}^1& \\ \bm{v}_1& &\bm{v}_2\\ &\bm{h}^2& 
		\end{smallmatrix} \right):\bm{v}_1\Longrightarrow\bm{v}_2 $,
		\begin{align*}
		\left( \bm{c}^{{\star^1}^\uda}\right)^{\star^1}&=\left(\left(\left(\bm{c^\uda} \right)^{\star^1} \right)^\uda \right)^{\star^1}\\
		&=\left(\begin{pmatrix}
		&\wt{\bm{h}^2}& \\ \left( \bm{v}_1^\uda\right)^{\star^1}& &\left( \bm{v}_2^\uda\right)^{\star^1}\\ &\wt{\bm{h}^1}& 
		\end{pmatrix}^\uda \right)^{\star^1}\\
		&=\begin{pmatrix}
		&\wt{(\wt{\bm{h}^1})}& \\ \left(\left(\left(\bm{v}^\uda_1 \right)^{\star^1} \right)^\uda \right)^{\star^1}& &\left(\left(\left(\bm{v}^\uda_2 \right)^{\star^1} \right)^\uda \right)^{\star^1}\\ &\wt{(\wt{\bm{h}^2})}& 
		\end{pmatrix}\\
		&=\begin{pmatrix}
		&\bm{h}^1& \\ \bm{v}_1& &\bm{v}_2\\ &\bm{h}^2& 
		\end{pmatrix}\\
		&=\bm{c}; \tag{10}
		\end{align*}
		also, using (10), we have for any $\mcg{B^\uda}$-morphism $\bm{c^\uda}=\left(\begin{smallmatrix}
		&\bm{h}^1& \\ \bm{v}_1& &\bm{v}_2\\ &\bm{h}^2& 
		\end{smallmatrix} \right)^\uda:\bm{v}_1^\uda\Longrightarrow\bm{v}_2^\uda $,
		\begin{align}
		\left( \left( \bm{c^\uda}\right)^{\star^1} \right)^{{\star^1}^\uda}=\left(\left(\left(\left(\bm{c^\uda} \right)^{\star^1} \right)^\uda \right)^{\star^1} \right)^\uda=\bm{c^\uda}. \tag{11}
		\end{align}
		Now, Equations (8), (10) together imply
		\[\left((-)^{{\star^1}^\uda} \right)^{\star^1}=1_\mcg{B}. \]
		Likewise, Equations (9), (11) together imply
		\[\left( (-)^{\star^1} \right)^{{\star^1}^\uda}=1_{\mcg{B^\uda}}. \]
		Consequently, $(-)^{\star^1}$, too, is an isomorphism. This completes the proof of the vertical self-duality of $(-)^\star$.
	\end{mypr}
	
	~\\
	\par Theorems \ref{DLC-hor-sd} and \ref{DLC-ver-sd} together yield one of the fundamental features of $\dlc$:
	
	\begin{cor}
		\label{DLC-Klein-inv}
		For every set $\Sigma$ and every nonempty set $\Gamma$, the double category $\dlc_{\Gamma,\Sigma}$ is Klein-invariant.
	\end{cor}
	
	Corollary \ref{DLC-Klein-inv} is so important that it deserves special attention. What it tells is that we have the following four isomorphic double categories:
	\[\dlc\cong\dlc^\lra\cong\dlc^\uda\cong\dlc^\boxplus.\]
	This means that whenever one proves a result in any of the above four double categories, that result automatically yields three other analogous results \textit{in the same double category} by horizontal, vertical, and central self-dualizations, respectively. This situation is comparable to (and is a generalization of) the self-duality of \chu.
	
	\subsection{Binary horizontal products and coproducts}
	In this subsection we have both negative and positive results. We start with negative results.
	
	\begin{prp}
		\label{no-bin-prod-coprod}
		Assume a nonempty set $\Gamma$ and a (possibly empty) set $\Sigma$. Then:
		\begin{enumerate}
			\item there exist 0-cells $\eb{P}_1,\eb{P}_2$ in $\dlc_{\Gamma,\Sigma}$ such that for any two vertical 1-cells\\ $\bm{v}_1:\eb{P}_1\rightharpoonup\prm{\eb{P}}_1$ and $\bm{v}_2:\eb{P}_2\rightharpoonup\prm{\eb{P}}_2$  there exist no horizontal binary product and no horizontal binary coproduct;
			
			\item there exist 0-cells $\eb{P}_1,\eb{P}_2$ in $\dlc_{\Gamma,\Sigma}$ such that for any two vertical 1-cells\\ $\bm{v}_1:\prm{\eb{P}}_1\rightharpoonup\eb{P}_1$ and $\bm{v}_2:\prm{\eb{P}}_2\rightharpoonup\eb{P}_2$  there exist no horizontal binary product and no horizontal binary coproduct.
		\end{enumerate}		
	\end{prp}
	
	\begin{mypr}
		(1)~ Consider any two 0-cells $\eb{P}_1,\eb{P}_2$ with different functorial profiles
		\[\boldsymbol{\Phi}_{\eb{P}_1}\ne\boldsymbol{\Phi}_{\eb{P}_2},\]
		and consider vertical 1-cells $\bm{v}_1:\eb{P}_1\rightharpoonup\prm{\eb{P}}_1,~\bm{v}_2:\eb{P}_2\rightharpoonup\prm{\eb{P}}_2$ for some $\prm{\eb{P}}_1,\prm{\eb{P}}_2$. Also, assume that
		\[\bm{w}:\eb{Q}\rightharpoonup\prm{\eb{Q}}\]
		is the horizontal product of $\bm{v}_1,\bm{v}_2$. This implies that the 0-cell $\eb{Q}$ must be the (ordinary) product of $\eb{P}_1,\eb{P}_2$. On the other hand, we know that at least one of the following holds:
		\[\boldsymbol{\Phi}_\eb{Q}\ne\boldsymbol{\Phi}_{\eb{P}_1}~~~~\text{or}~~~~\boldsymbol{\Phi}_\eb{Q}\ne\boldsymbol{\Phi}_{\eb{P}_2}.\]
		Assume, for example, $\boldsymbol{\Phi}_\eb{Q}\ne\boldsymbol{\Phi}_{\eb{P}_1}$. But this means that there can be no horizontal 1-cells from $\eb{Q}$ to $\eb{P}_1$; thus, we have no projection from $\eb{Q}$ to $\eb{P}_1$, and consequently, no projection 2-cell from $\bm{w}$ to $\bm{v}_1$. Therefore, $\eb{Q}$ cannot be a product of $\eb{P}_1,\eb{P}_2$, and $\bm{w}$ cannot be a horizontal product of $\bm{v}_1,\bm{v}_2$. Contradiction.\\
		\par A similar argument works for the coproduct of $\bm{v}_1,\bm{v}_2$, as well.\\
		
		\par (2)~ For this part, just apply the vertical self-duality of $\dlc$ on part (1).
	\end{mypr}
	
	~\\
	\par In addition, when $\Sigma\ne\varnothing$, the situation is even worse:	
	
	\begin{prp}
		\label{no-bin-prod-coprod-for-nonempty}
		Assume nonempty sets $\Gamma,\Sigma$. Then, in the double category $\dlc_{\Gamma,\Sigma}$:
		\begin{enumerate}
			\item there exist 0-cells $\eb{P}_1,\eb{P}_2$ with $\boldsymbol{\Phi}_{\eb{P}_1}=\boldsymbol{\Phi}_{\eb{P}_2}$, such that for any pair of 0-cells $\prm{\eb{P}}_1,\prm{\eb{P}}_2$ with $\boldsymbol{\Phi}_{\prm{\eb{P}}_1}=\boldsymbol{\Phi}_{\prm{\eb{P}}_2}$ and vertical 1-cells $\bm{v}_1:\eb{P}_1\rightharpoonup\prm{\eb{P}}_1,\bm{v}_2:\eb{P}_2\rightharpoonup\prm{\eb{P}}_2$ no horizontal binary product exists;
			
			\item there exist 0-cells $\eb{Q}_1,\eb{Q}_2$ with $\boldsymbol{\Phi}_{\eb{Q}_1}=\boldsymbol{\Phi}_{\eb{Q}_2}$, such that for any pair of 0-cells $\prm{\eb{Q}}_1,\prm{\eb{Q}}_2$ with $\boldsymbol{\Phi}_{\prm{\eb{Q}}_1}=\boldsymbol{\Phi}_{\prm{\eb{Q}}_2}$ and vertical 1-cells $\bm{w}_1:\eb{Q}_1\rightharpoonup\prm{\eb{Q}}_1,\bm{w}_2:\eb{Q}_2\rightharpoonup\prm{\eb{Q}}_2$ no horizontal binary coproduct exists;
			
			\item there exist 0-cells $\eb{P}_1,\eb{P}_2$ with $\boldsymbol{\Phi}_{\eb{P}_1}=\boldsymbol{\Phi}_{\eb{P}_2}$, such that for any pair of 0-cells $\prm{\eb{P}}_1,\prm{\eb{P}}_2$ with $\boldsymbol{\Phi}_{\prm{\eb{P}}_1}=\boldsymbol{\Phi}_{\prm{\eb{P}}_2}$ and vertical 1-cells $\bm{v}_1:\prm{\eb{P}}_1\rightharpoonup\eb{P}_1,\bm{v}_2:\prm{\eb{P}}_2\rightharpoonup\eb{P}_2$ no horizontal binary product exists;
			
			\item there exist 0-cells $\eb{Q}_1,\eb{Q}_2$ with $\boldsymbol{\Phi}_{\eb{Q}_1}=\boldsymbol{\Phi}_{\eb{Q}_2}$, such that for any pair of 0-cells $\prm{\eb{Q}}_1,\prm{\eb{Q}}_2$ with $\boldsymbol{\Phi}_{\prm{\eb{Q}}_1}=\boldsymbol{\Phi}_{\prm{\eb{Q}}_2}$ and vertical 1-cells $\bm{w}_1:\prm{\eb{Q}}_1\rightharpoonup\eb{Q}_1,\bm{w}_2:\prm{\eb{Q}}_2\rightharpoonup\eb{Q}_2$ no horizontal binary coproduct exists.
		\end{enumerate}
	\end{prp}
	
	\rule{0pt}{1mm}
	
	\begin{mypr}
		We only prove part (1); parts (2), (3), and (4) are automatically proved by horizontal, vertical, and central self-dualities, respectively.\\
		\par Assume nonempty sets $\Gamma,\Sigma$, and consider the double category
		\[\dlc_{\Gamma,\Sigma}=\langle\mcg{A,B,s,t,i,\bt}\rangle.\]
		Consider $\chu_\Gamma$-objects $\msf{A}_i=\agg{A_i}{r_i}{X_i},~\msf{B}_i=\agg{B_i}{s_i}{Y_i},~i=1,2$ with nonempty carriers $A_i,B_i$. Also, consider four $\chu_\Gamma$-endofunctors $F,G,H,K$ such that for every Chu space $\msf{C}=\agg{C}{t}{Z}$ with a nonempty carrier $C$,
		\[F^+\msf{C}\ne\varnothing,~~G^+\msf{C}\ne\varnothing,~~H^+\msf{C}\ne\varnothing,~~K^+\msf{C}\ne\varnothing\]
		(for example, take $F=G=H=K=1_{\chu_\Gamma}$). Assume that
		\[\ag{\msf{A}_i}{F\msf{A}_i\xrightarrow{~\alpha_i~} G\msf{A}_i},~~\ag{\msf{B}_i}{K\msf{B}_i\xrightarrow{~\beta_i~} H\msf{B}_i},~~i=1,2\]
		are dialgebras. 
		\par Fix $m\in\Sigma$. Define 0-cell $\eb{P}_1$ as
		\[\eb{P}_1\eqd\ag{\begin{bmatrix}
			\ag{\msf{A}_1}{F\msf{A}_1\xrightarrow{~\alpha_1~} G\msf{A}_1}\\ \ag{\msf{B}_1}{K\msf{B}_1\xrightarrow{~\beta_1~} H\msf{B}_1}
			\end{bmatrix}}{\ag{S_1^+}{S_1^-}}\]
		such that
		\begin{align*}
		S_1^+:\Sigma&\longrightarrow F^+\msf{A}_1\sqcup K^+\msf{B}_1\\
		S_1^+(m)&\in F^+\msf{A}_1.
		\end{align*}
		Also, define 0-cell $\eb{P}_2$ as
		\[\eb{P}_2\eqd\ag{\begin{bmatrix}
			\ag{\msf{A}_2}{F\msf{A}_2\xrightarrow{~\alpha_2~} G\msf{A}_2}\\ \ag{\msf{B}_2}{K\msf{B}_2\xrightarrow{~\beta_2~} H\msf{B}_2}
			\end{bmatrix}}{\ag{S_2^+}{S_2^-}}\]
		such that
		\begin{align*}
		S_2^+:\Sigma&\longrightarrow F^+\msf{A}_2\sqcup K^+\msf{B}_2\\
		S_2^+(m)&\in K^+\msf{B}_2.
		\end{align*}
		Now, consider any pair of 0-cells $\prm{\eb{P}}_1,\prm{\eb{P}}_1$ with $\boldsymbol{\Phi}_{\prm{\eb{P}}_1}=\boldsymbol{\Phi}_{\prm{\eb{P}}_2}$, together with vertical 1-cells
		\[\bm{v}_1:\eb{P}_1\rightharpoonup\prm{\eb{P}}_1,~~\bm{v}_2:\eb{P}_2\rightharpoonup\prm{\eb{P}}_2.\]
		We claim that there are no horizontal products of $\bm{v}_1,\bm{v}_2$.\\
		\par We prove the claim by contradiction. Assume there exists some vertical 1-cell\\ $\bm{u}:\eb{P}\rightharpoonup\prm{\eb{P}}$ together with, say, cubicles
		\[\bm{v}_1\stackrel{~~\bm{c}_1~~}{\Longleftarrow}\bm{u}\stackrel{~~\bm{c}_2~~}{\Longrightarrow}\bm{v}_2,\]
		\[\bm{c}_i=\begin{pmatrix}
		&\left[ \begin{smallmatrix}
		\eb{f}_i\\\eb{g}_i
		\end{smallmatrix}\right] & \\ \bm{u}& &\bm{v}_i\\ &\left[ \begin{smallmatrix}
		\prm{\eb{f}}_i\\\prm{\eb{g}}_i
		\end{smallmatrix}\right]& 
		\end{pmatrix},~~i=1,2\]
		that possess the universal property of horizontal product. Then,
		\begin{align}
		\eb{P}_1\xleftarrow{~\left[ \begin{smallmatrix}
			\eb{f}_1\\\eb{g}_1
			\end{smallmatrix}\right]~}\eb{P}\xrightarrow{~\left[ \begin{smallmatrix}
			\eb{f}_2\\\eb{g}_2
			\end{smallmatrix}\right]~}\eb{P}_2 \tag{$*$}
		\end{align}
		must be a (product) diagram in $\mcg{A}$. Let $\eb{P}=\ag{\left[\begin{smallmatrix}
			\ag{\msf{C}}{\gamma}\\ \ag{\msf{D}}{\delta}
			\end{smallmatrix} \right] }{\ag{T^+}{T^-}}$. Then, the above diagram implies that for the horizontal 1-cells $\left[ \begin{smallmatrix}
		\eb{f}_1\\\eb{g}_1
		\end{smallmatrix}\right]$ and $\left[ \begin{smallmatrix}
		\eb{f}_2\\\eb{g}_2
		\end{smallmatrix}\right]$ we must have the following positive horizontal super-adjointness conditions $\tu{PHS}_1$ and $\tu{PHS}_2$, respectively:
		\begin{center}
			\mbox{
				\begin{tikzpicture}[commutative diagrams/every diagram]
				
				\node (N1) at (0,0cm) {$\Sigma$};
				
				\node (N2) at (0,3cm) {$F^+\msf{C}\sqcup K^+\msf{D}$};
				
				\node (N3) at (-4,3cm) {$F^+\msf{A}_1\sqcup K^+\msf{D}$};
				
				\node (N4) at (-4,0cm) {$F^+\msf{A}_1\sqcup K^+\msf{B}_1$};
				
				\node (M3) at (4,3cm) {$F^+\msf{A}_2\sqcup K^+\msf{D}$};
				
				\node (M4) at (4,0cm) {$F^+\msf{A}_2\sqcup K^+\msf{B}_2$};
				
				\node (N8) at (-2,1.5cm) {\fbox{$\tu{PHS}_1$}};
				
				\node (M8) at (2,1.5cm) {\fbox{$\tu{PHS}_2$}};	
				
				\path[commutative diagrams/.cd, every arrow, every label]
				(N1) edge node[swap]{$T^+$} (N2)
				(N2) edge node[swap]{$F^+f_1\sqcup 1$} (N3)
				(N1) edge node[swap]{$S^+_1$} (N4)
				(N4) edge node{$1\sqcup K^+g_1$} (N3)
				(N2) edge node{$F^+f_2\sqcup 1$} (M3)
				(N1) edge node{$S^+_2$} (M4)
				(M4) edge node[swap]{$1\sqcup K^+g_2$} (M3);
				\end{tikzpicture}
			}
		\end{center}
		Since $S^+_1(m)\in F^+\msf{A}_1$, we obtain from $\tu{PHS}_1$ that
		\[(1\sqcup K^+g_1)\circ S^+_1(m)\in F^+\msf{A}_1,\]
		or
		\begin{align}
		(F^+f_1\sqcup 1)\circ T^+(m)\in F^+\msf{A}_1. \tag{1}
		\end{align}
		Now if $T^+(m)\in K^+\msf{D}$, it follows that $(F^+f_1\sqcup 1)\circ T^+(m)\in K^+\msf{D}$ which contradicts (1) because the coproduct $F^+\msf{A}_1\sqcup K^+\msf{D}$ is the \textit{disjoint union} of the sets $F^+\msf{A}_1$ and $K^+\msf{D}$. Therefore, we must have
		\begin{align}
		T^+(m)\in F^+\msf{C}. \tag{2}
		\end{align}
		On the other hand, considering $\tu{PHS}_2$, we find from $S_2^+(m)\in K^+\msf{B}_2$ that
		\[(1\sqcup K^+g_2)\circ S^+_2(m)\in K^+\msf{D},\]
		or
		\begin{align}
		(F^+f_2\sqcup 1)\circ T^+(m)\in K^+\msf{D}. \tag{3}
		\end{align}
		But Equation (2) says that $T^+(m)\in F^+\msf{C}$, thus
		\begin{align}
		(F^+f_2\sqcup 1)\circ T^+(m)\in F^+\msf{A}_2, \tag{4}
		\end{align}
		which clearly contradicts (3). Consequently, there can be no such product diagram as $(*)$, and we are done.
	\end{mypr}
	
	~\\
	Finally, we have a positive result for the case $\Sigma=\varnothing$. 
	
	\begin{prp}
		\label{positive-res-for-empty}
		Assume a nonempty set $\Gamma$, and let $\Sigma=\varnothing$. Assume that $\eb{P}_1,\eb{P}_2,\prm{\eb{P}}_1,\prm{\eb{P}}_2$ are 0-cells in the double category $\dlc_{\Gamma,\varnothing}$ such that
		\begin{align*}
		\boldsymbol{\Phi}_{\eb{P}_1}=\boldsymbol{\Phi}_{\eb{P}_2}&=\begin{bmatrix}
		F&G\\K&H
		\end{bmatrix},\\
		\boldsymbol{\Phi}_{\prm{\eb{P}}_1}=\boldsymbol{\Phi}_{\prm{\eb{P}}_2}&=\begin{bmatrix}
		\prm{F}&\prm{G}\\\prm{K}&\prm{H}
		\end{bmatrix}.
		\end{align*}
		Then, for every pair of vertical 1-cells $\bm{v}_1:\eb{P}_1\rightharpoonup\prm{\eb{P}}_1,\bm{v}_2:\eb{P}_2\rightharpoonup\prm{\eb{P}}_2$ we have the following:
		\begin{enumerate}
			\item If the endofunctors $F,G,\prm{F},\prm{G}$ preserve binary products and the endofunctors $H,K,\prm{H},\prm{K}$ preserve binary coproducts, then there exists a horizontal product for $\bm{v}_1,\bm{v}_2$; that is, there exists a vertical 1-cell $\bm{w}:\eb{Q}\rightharpoonup\prm{\eb{Q}}$ for some $\eb{Q},\prm{\eb{Q}}$, together with cubicles
			\[\bm{v}_1\stackrel{~~\bm{p}_1~~}{\Longleftarrow}\bm{w}\stackrel{~~\bm{p}_2~~}{\Longrightarrow}\bm{v}_2\]
			with the following universal property: for any other vertical 1-cell $\bm{u}$ and any diagram 
			\[\bm{v}_1\stackrel{~~\bm{d}_1~~}{\Longleftarrow}\bm{u}\stackrel{~~\bm{d}_2~~}{\Longrightarrow}\bm{v}_2,\]
			there exists a unique cubicle $\bm{d}:\bm{u}\Longrightarrow\bm{w}$ such that
			\[\bm{d}_k=\bm{p}_k\circ\bm{d},~~k=1,2.\]
			
			\item If the endofunctors $F,G,\prm{F},\prm{G}$ preserve binary coproducts and the endofunctors $H,K,\prm{H},\prm{K}$ preserve binary products, then there exists a horizontal coproduct for $\bm{v}_1,\bm{v}_2$; that is, there exists a vertical 1-cell $\bm{w}:\eb{Q}\rightharpoonup\prm{\eb{Q}}$ for some $\eb{Q},\prm{\eb{Q}}$, together with cubicles
			\[\bm{v}_1\stackrel{~~\bm{c}_1~~}{\Longrightarrow}\bm{w}\stackrel{~~\bm{c}_2~~}{\Longleftarrow}\bm{v}_2\]
			with the following universal property: for any other vertical 1-cell $\bm{u}$ and any diagram 
			\[\bm{v}_1\stackrel{~~\bm{d}_1~~}{\Longrightarrow}\bm{u}\stackrel{~~\bm{d}_2~~}{\Longleftarrow}\bm{v}_2,\]
			there exists a unique cubicle $\bm{d}:\bm{w}\Longrightarrow\bm{u}$ such that
			\[\bm{d}_k=\bm{d}\circ\bm{c}_k,~~k=1,2.\]
		\end{enumerate}
	\end{prp}
	
	\begin{mypr}
		We only prove part (1); part (2) is automatically proved by horizontal self-duality.\\
		
		\par Note that since $\Sigma=\varnothing$, the positive and negative horizontal super-adjointness conditions are trivially satisfied for all horizontal 1-cells, and the positive and negative vertical super-adjointness conditions are trivially satisfied for all vertical 1-cells. Thus, whenever we define either a horizontal 1-cell or a vertical 1-cell, there will be no need to worry about any of the conditions PHS, NHS, PVS, or NVS.\\
		
		\par Assume that $\eb{P}_k=\ag{\left[\begin{smallmatrix}
			\alpha_k:F\msf{A}_k\longrightarrow G\msf{A}_k\\
			\beta_k:K\msf{B}_k\longrightarrow H\msf{B}_k
			\end{smallmatrix} \right]}{\eb{S}_k} $ and $\prm{\eb{P}}_k=\ag{\left[\begin{smallmatrix}
			\prm{\alpha}_k:\prm{F}\msf{A}_k\longrightarrow \prm{G}\msf{A}_k\\
			\prm{\beta}_k:\prm{K}\msf{B}_k\longrightarrow \prm{H}\msf{B}_k
			\end{smallmatrix} \right]}{\prm{\eb{S}}_k} $ for\\ $k=1,2$, where each of $\eb{S}_k,\prm{\eb{S}}_k$ is a pair of empty functions. Since $F,G,\prm{F},\prm{G}$ preserve binary products, using Corollary \ref{Chu-forget-crea-pres-G-F}, we find that there exist product dialgebras
		\[\alpha=\alpha_1\&\alpha_2:F\msf{A}_1\& F\msf{A}_2\longrightarrow G\msf{A}_1\& G\msf{A}_2~~~~\text{in}~~\dlg{\chu}{F}{G}\]
		and
		\[\prm{\alpha}=\prm{\alpha}_1\&\prm{\alpha}_2:\prm{F}\msf{A}_1\& \prm{F}\msf{A}_2\longrightarrow \prm{G}\msf{A}_1\& \prm{G}\msf{A}_2~~~~\text{in}~~\dlg{\chu}{\prm{F}}{\prm{G}}.\]
		These products are equipped with projection dialgebra homomorphisms
		\[\alpha_1\xleftarrow{~\boldsymbol{\uppi}_1~}\alpha\xrightarrow{~\boldsymbol{\uppi}_2~}\alpha_2~~~~\text{and}~~~~\prm{\alpha}_1\xleftarrow{~\prm{\boldsymbol{\uppi}}_1~}\prm{\alpha}\xrightarrow{~\prm{\boldsymbol{\uppi}}_2~}\prm{\alpha}_2.\]
		On the other hand, since $H,K,\prm{H},\prm{K}$ preserve binary coproducts, again using Corollary \ref{Chu-forget-crea-pres-G-F}, we obtain coproduct dialgebras
		\[\beta=\beta_1\oplus\beta_2:K\msf{B}_1\oplus K\msf{B}_2\longrightarrow H\msf{B}_1\oplus H\msf{B}_2~~~~\text{in}~~\dlg{\chu}{K}{H}\]
		and
		\[\prm{\beta}=\prm{\beta}_1\oplus\prm{\beta}_2:\prm{K}\msf{B}_1\oplus \prm{K}\msf{B}_2\longrightarrow \prm{H}\msf{B}_1\oplus \prm{H}\msf{B}_2~~~~\text{in}~~\dlg{\chu}{\prm{K}}{\prm{H}}.\]
		These too, are equipped with some injection dialgebra homomorphisms
		\[\beta_1\xrightarrow{~\textbf{i}_1~}\beta\xleftarrow{~\textbf{i}_2~}\beta_2~~~~\text{and}~~~~\prm{\beta}_1\xrightarrow{~\prm{\textbf{i}}_1~}\prm{\beta}\xleftarrow{~\prm{\textbf{i}}_2~}\prm{\beta}_2.\]
		Let $\eb{T}=\ag{!^+}{!^-}$ and $\prm{\eb{T}}=\langle{\prm{!}}^+,{\prm{!}}^-\rangle$ where
		\begin{align*}
		!^+:\varnothing&\longrightarrow F^+(\msf{A}_1\& \msf{A}_2)\sqcup K^+(\msf{B}_1\oplus \msf{B}_2),\\
		!^-:\varnothing&\longrightarrow G^-(\msf{A}_1\& \msf{A}_2)\sqcup H^-(\msf{B}_1\oplus \msf{B}_2),\\
		{\prm{!}}^+:\varnothing&\longrightarrow {\prm{F}}^+(\msf{A}_1\& \msf{A}_2)\sqcup {\prm{K}}^+(\msf{B}_1\oplus \msf{B}_2),\\
		{\prm{!}}^-:\varnothing&\longrightarrow {\prm{G}}^-(\msf{A}_1\& \msf{A}_2)\sqcup {\prm{H}}^-(\msf{B}_1\oplus \msf{B}_2)
		\end{align*}
		are all empty functions. Define
		\begin{align*}
		\eb{Q}&\eqd\ag{\begin{bmatrix}
			\alpha\\ \beta
			\end{bmatrix}}{\eb{T}};\\
		\prm{\eb{Q}}&\eqd\ag{\begin{bmatrix}
			\prm{\alpha}\\ \prm{\beta}
			\end{bmatrix}}{\prm{\eb{T}}};\\
		\bm{q}_k&\eqd\begin{bmatrix}
		\boldsymbol{\uppi}_k\\ \textbf{i}_k
		\end{bmatrix}:\eb{Q}\longrightarrow\eb{P}_k,~~k=1,2;\\
		\prm{\bm{q}}_k&\eqd\begin{bmatrix}
		\prm{\boldsymbol{\uppi}}_k\\ \prm{\textbf{i}}_k
		\end{bmatrix}:\prm{\eb{Q}}\longrightarrow\prm{\eb{P}}_k,~~k=1,2.
		\end{align*}
		Now, for any pair of vertical 1-cells $\bm{v}_1:\eb{P}_1\rightharpoonup\prm{\eb{P}}_1,\bm{v}_2:\eb{P}_2\rightharpoonup\prm{\eb{P}}_2$ with
		\[\bm{v}_k=\begin{bmatrix}
		\ag{\mu_k}{\nu_k}\\ \ag{\theta_k}{\zeta_k}
		\end{bmatrix},~~k=1,2,\]
		define
		\[\bm{w}\eqd\begin{bmatrix}
		\ag{\mu_1\&\mu_2}{\nu_1\&\nu_2}\\ \ag{\theta_1\oplus\theta_2}{\zeta_1\oplus\zeta_2}
		\end{bmatrix}:\eb{Q}\rightharpoonup\prm{\eb{Q}}.\]
		Also, define the cubicles
		\[\bm{p}_k\eqd\begin{pmatrix}
		&\bm{q}_k& \\
		\bm{w}& &\bm{v}_k\\
		&\prm{\bm{q}}_k& 
		\end{pmatrix}:\bm{w}\Longrightarrow\bm{v}_k,~~k=1,2,\]
		We claim that the triple $\agg{\bm{w}}{\bm{p}_1}{\bm{p}_2}$ is a horizontal product for $\bm{v}_1,\bm{v}_2$.\\
		\par To see this, let
		\[\bm{v}_1\stackrel{~~\bm{d}_1~~}{\Longleftarrow}\bm{u}\stackrel{~~\bm{d}_2~~}{\Longrightarrow}\bm{v}_2\]
		be a double diagram in $\dlc_{\Gamma,\varnothing}$, such that $\bm{u}:\eb{R}\rightharpoonup\prm{\eb{R}}$ is a vertical 1-cell between 0-cells
		\[\eb{R}=\ag{\begin{bmatrix}
			\gamma\\ \delta
			\end{bmatrix}}{\eb{U}},~~~~\prm{\eb{R}}=\ag{\begin{bmatrix}
			\prm{\gamma}\\ \prm{\delta}
			\end{bmatrix}}{\prm{\eb{U}}},\]
		for some dialgebras $\gamma,\delta,\prm{\gamma},\prm{\delta}$, and pairs of empty functions $\eb{U},\prm{\eb{U}}$, and such that
		\[\bm{d}_k=\begin{pmatrix}
		&\bm{h}_k& \\
		\bm{u}& &\bm{v}_k\\
		&\prm{\bm{h}}_k& 
		\end{pmatrix},~~k=1,2,\]
		with
		\begin{align*}			
		\bm{h}_k&=\begin{bmatrix}
		\gamma\xrightarrow{~\eb{f}_k~}\alpha_k\\ \delta\xleftarrow{~\eb{g}_k~}\beta_k
		\end{bmatrix},~~k=1,2;\\
		\prm{\bm{h}}_k&=\begin{bmatrix}
		\prm{\gamma}\xrightarrow{~\prm{\eb{f}}_k~}\prm{\alpha}_k\\ \prm{\delta}\xleftarrow{~\prm{\eb{g}}_k~}\prm{\beta}_k
		\end{bmatrix},~~k=1,2.
		\end{align*}
		Then, by universalities of dialgebras $\alpha,\beta,\prm{\alpha},\prm{\beta}$, there exist unique dialgebra homomorphisms
		\begin{align*}
		\boldsymbol{\upvarphi}&:\gamma\longrightarrow\alpha,\\
		\boldsymbol{\uppsi}&:\beta\longrightarrow\delta,\\
		\prm{\boldsymbol{\upvarphi}}&:\prm{\gamma}\longrightarrow\prm{\alpha},\\
		\prm{\boldsymbol{\uppsi}}&:\prm{\beta}\longrightarrow\prm{\delta}
		\end{align*}
		satisfying the following equations for $k=1,2$:
		\begin{align*}
		\eb{f}_k&=\boldsymbol{\uppi}_k\boldsymbol{\upvarphi},\\
		\eb{g}_k&=\boldsymbol{\uppsi}\eb{i}_k,\\
		\prm{\eb{f}}_k&=\prm{\boldsymbol{\uppi}}_k\prm{\boldsymbol{\upvarphi}},\\
		\prm{\eb{g}}_k&=\prm{\boldsymbol{\uppsi}}\prm{\eb{i}}_k.
		\end{align*}
		Therefore, we have a unique cubicle
		\[\bm{d}\eqd\begin{pmatrix}
		&\begin{bmatrix}
		\boldsymbol{\upvarphi}\\ \boldsymbol{\uppsi}
		\end{bmatrix}& \\
		\bm{u}& &\bm{w}\\
		&\begin{bmatrix}
		\prm{\boldsymbol{\upvarphi}}\\ \prm{\boldsymbol{\uppsi}}
		\end{bmatrix}& 
		\end{pmatrix}:\bm{u}\Longrightarrow\bm{w}\]
		which satisfies
		\[\bm{d}_k=\bm{p}_k\circ\bm{d},~~k=1,2.\]
	\end{mypr}
	
	~\\
	\par In the next chapter, we will outline some of the themes for future development of the formalism presented in the current chapter.

	\chapter{Conclusion and Suggestions for Future Work}
	In the final chapter of this work, we sum up what has been presented so far, and give some guidelines for future theoretical developments as well as for future practical implementations of the new formalism.
	
	\section{Conclusion}
	In the previous chapters, we took the initial steps of developing what may be called ``$ \dlc $ theory'' or ``cubicle theory''. First of all, we introduced the preliminary concepts and definitions, including the essentials of basic category theory and internal categories in Chapter 1. Then, in Chapter 2, we did the same for monoidal categories and $*$-autonomous categories. We then proceeded towards the Chu construction, we focused on $\chu=\chu(\cat{Set},\Gamma)$, and we gave a number of its important properties. We observed the bicompleteness of the Chu construction in general and of $\chu$ in particular; we dealt with extensional, separable, and biextensional Chu spaces, and we proved that \chu~is balanced.
	\par Next, we paid attention to linear logic in Chu. We gave exact formulations for multiplicative as well as additive connectives of linear logic in Chu. Multiplicative connectives consist of \textit{dualization}, \textit{tensor}, \textit{linear implication}, and \textit{par} operations, together with the \textit{tensor unit} and the \textit{dualizing object}. On the other hand, additive connectives consist of \textit{plus} and \textit{with} operations (which are exactly the categorical binary coproduct and binary product, respectively), together with the \textit{plus unit} and the \textit{with unit}.
	\par The issue of \textit{realization} was our next subject of discussion. There, we quoted some remarkable results indicating the vast potentials of the Chu construction to realize concrete mathematics: $\cat{Set}$ and $\cat{Top}$ are realized in $\chu_2$, $\cat{Grp}$ is realized in $\chu_8$, $\cat{TopGrp}$ is realized in $\chu_{16}$, $n$-ary relational structures in general are realized in $\chu_{2^n}$, $\cat{Vect}_\Bbbk$ is realized in $\chu_{\left| \Bbbk\right| }$, and so on.
	\par In the final subsection of Chapter 2, we focused on Chu endofunctors, establishing a general method for uplifting arbitrary set functors to Chu functors. Also, we mentioned bi-uplifting as a generalization of the notion of uplifting.
	\par In Chapter 3, we turned to universal dialgebra. We gave the basic notions of the theory for a general (bicomplete) base category $\mcg{C}$, and we stated and proved some of the fundamental properties of the categories of the form $\dlg{\mcg{C}}{F}{G}$.
	\par In Chapter 4, we firstly introduced double categories as categories internal to $ \cat{CAT} $. Next, we introduced the main formalism of the present work, namely the $ \dlc $ construction; we proved its well-definedness, and we studied a number of its basic properties. The most important property was Klein-invariance, which we likened to a two-dimensional generalization of the self-duality of the Chu construction. Also, we observed that although horizontal products and coproducts do not exist in $ \dlc $ in many situations, they do exist at least when certain conditions are met (concerning the functorial profiles and the supermatrices).
	\par With all the above, we finished the present work. Now it is time to point to some opportunities for future work.
	
	\section{Future development}
	For this, the following suggestions are given.
	
	\subsubsection{Horizontal and vertical limits/colimits in general}
	A very special kind of limits/colimits was studied in this work, namely binary horizontal products and coproducts. The study may well be continued to cover arbitrary horizontal as well as vertical limits/colimits in $ \dlc $.
	
	\subsubsection{Other kinds of double-categorical constructions}
	Other double-categorical concepts are to be developed in $ \dlc $, also. These include:
	\begin{itemize}
		\item various subcategories of $ \dlc $;
		
		\item ``exponentials'' (see \cite{GranPar-Lims});
		
		\item ``comma double categories ([ibid]);
		
		\item ``double isomorphisms'' and ``sesqui-isomorphisms'' ([ibid]);
		
		\item ``double limits/colimits'' ([ibid]);
		
		\item various ``forgetful'' and ``free functors'';
		
		\item ``double adjunctions'' and ``free monads'' (\cite{Fiore-Adj});
		
		\item ``Kan extensions'' (\cite{CWM,GranPar-Kan,Koud}) and ``Kan lifts'' (\cite{nLab-lift}) in $ \dlc $;
		
		\item etc.
	\end{itemize}
	
	Certainly, there are other theoretical as well as practical aspects that may be suggested for future development of cubicle theory; however, the author prefers to postpone mentioning them to other occasions. 

	\bibliographystyle{plain}
	\bibliography{refs-en}
	\printindex

\end{document}